\documentclass{amsart}
\usepackage{amssymb,latexsym,graphics}
\usepackage[mathscr]{eucal}
\usepackage{enumerate}
\newtheorem{thm}{Theorem}[section]
\newtheorem{cor}[thm]{Corollary}
\newtheorem{lem}[thm]{Lemma}
\newtheorem{prop}[thm]{Proposition}

\theoremstyle{remark}
\newtheorem{rem}[thm]{Remark}

\newtheorem{defn}[thm]{Definition}
\theoremstyle{definition}

\newcommand{\step}[2]{\noindent\textbf{Step #1:}(\emph{#2})}

\newcommand{\tphi}{\widetilde{\phi}}
\newcommand{\llp}[2]{|\!|\!| \, #1 \, |\!|\!|_{#2}}
\newcommand{\R}{\mathbb{R}}
\newcommand{\Z}{\mathbb{Z}}

\newcommand{\lp}[2]{\Vert \, #1 \, \Vert_{#2}}

\newcommand{\td}{\widetilde}

\newcommand{\dint}{ {\int\!\!\int} }
\newcommand{\snabla}{ {\slash\!\!\!\nabla} }
\newcommand{\Anabla}{ { {}^A\!\nabla } }
\newcommand{\uX}{ {\underline{X}} }
\newcommand{\uS}{ {\underline{S}} }
\newcommand{\uE}{ {\underline{E}} }
\newcommand{\tc}{\tilde{c}}
\newcommand{\lsf}{\lesssim_F}
\newcommand{\lse}{\lesssim_E}

\newcommand{\Ut}{U^\dagger}

\newcommand{\W}{\mathcal{W}}

\newcommand{\ret}{\vspace{.2cm}}

\begin{document}

\title[Energy Dispersed Wave Maps]
{Energy dispersed large data wave maps in $2+1$ dimensions}
\author{Jacob Sterbenz}
\address{Department of Mathematics,
University of California, San Diego, CA 92093-0112}
\email{jsterben@math.ucsd.edu}
\thanks{The first author was supported in part by the NSF grant DMS-0701087}
\author{Daniel Tataru}
\address{Department of Mathematics, University of California, Berkeley, CA 94720-3840}
\email{tataru@math.berkeley.edu}
\thanks{The second author was
  supported in part by the NSF grant DMS-0801261.}
\date{}
\dedicatory{}
\commby{}

%%% ----------------------------------------------------------------------

\begin{abstract}
  In this article we consider large data Wave-Maps from $\R^{2+1}$
  into a compact Riemannian manifold $(\mathcal{M},g)$, and we prove that
  regularity and dispersive bounds persist as long as a certain type
  of bulk (non-dispersive) concentration is absent. This is a
  companion to our concurrent article \cite{ST_1}, which together with
  the present work establishes a full regularity theory for large data
  Wave-Maps.
\end{abstract}

%%% ----------------------------------------------------------------------
\maketitle
%%% ----------------------------------------------------------------------
%%%%%%%%%%%%%%%%%%%%%%%%%%%%%%%%%%%%%%%%%%%%%%%%%%%%%%%%%%%%%%%%%%%%%%%%%%
%%%-----------------------------------------------------------------------

%-------------------------------------------------------------------------
%-------------------------------------------------------------------------
%-------------------------------------------------------------------------
%-------------------------------------------------------------------------
%-------------------------------------------------------------------------
%%%%%%%%%%%%%%%%%%%%%%%%%%%%%%%%%%%%%%%%%%%%%%%%%%%%%%%%%%%%%%%%%%%%%%%%%%
%-------------------------------------------------------------------------

\tableofcontents

\section{Introduction}
In this article we consider finite energy large data Wave-Maps from
$\R^{2+1}$ into a compact Riemannian manifold $(\mathcal{M},g)$.
Our main result asserts that regularity and dispersive
bounds persist as long as a certain type of bulk
concentration is absent. The results proved here are used
in the companion article \cite{ST_1} to establish a full
regularity theory for large data Wave-Maps.\ret

The set-up we consider is the same as the one in \cite{Tataru_WM1},
using the so-called extrinsic formulation of the Wave-Maps equation.
Precisely, we consider the target manifold $(\mathcal{M},g)$ as an
isometrically embedded submanifold of $\R^N$. Then we can view the
$\mathcal{M}$ valued functions as $\R^N$ valued functions whose range
is contained in $\mathcal{M}$.  Such an embedding always exists by
Nash's theorem~\cite{MR17:782b} (see also Gromov~\cite{MR43:1212} and
G\"unther~\cite{MR93b:53049}).  In this context the Wave-Maps equation
can be expressed in a form which involves the second fundamental form
$\mathcal{S}$ of $\mathcal{M}$, viewed as a symmetric bilinear form:
\begin{align}
    \mathcal{S}:\, T\mathcal{M} \times T\mathcal{M}  \ &\to  \ N\mathcal{M}\ ,
    &\langle S(X,Y),N\rangle  \ &=  \ \langle \partial_X
    N, Y\rangle \ , \notag
\end{align}
The Cauchy problem for the wave maps equation
has the form:
\begin{subequations}\label{main_eq}
\begin{align}
    \Box \phi^a  \ &=  \ - \mathcal{S}^{a}_{bc}(\phi) \partial^\alpha \phi^b
        \partial_\alpha\phi^c \ , &\phi  &\in  \R^{N} \ , \label{main_eq1}\\
    \phi(0,x)\ &= \ \phi_0(x) \ , \ \ \partial_t\phi(0,x) \ = \ \dot{\phi}_0(x) \ , \label{main_eq2}
\end{align}
\end{subequations}
where the initial data $(\phi_0,\dot{\phi}_0)$ is chosen to obey the
constraint:
\begin{align}
    \phi_0(x) \ &\in\  \mathcal{M} \ ,
    &\dot{\phi}_0(x) \ &\in \ T_{\phi_0(x)}\mathcal{M} \ ,
    &x &\in \R^2 \ . \notag
\end{align}
In the sequel, it will be convenient for us to use the notation
$\phi[t] = (\phi(t), \partial_t \phi(t))$.  The system of equations
\eqref{main_eq} admits a conserved quantity, namely the Dirichlet
energy:
\begin{equation}
    E[\phi(t)] \ := \ \int_{\R^2} |\partial_t \phi(t)|^2 + |\nabla_x
    \phi(t)|^2 dx \ := \ \lp{\phi[t]}{\dot{H}^1 \times L^2}^2 \ =
    \ E \ . \label{energy}
\end{equation}
 Finite energy solutions
for \eqref{main_eq} correspond to initial data in the energy space,
namely $\phi[t] \in \dot{H}^{1} \times L^2$. We call a
Wave-Map ``classical'' on a bounded time interval
$(t_0,t_1)\times\R^2$ if $\nabla_{x,t} \phi(t)$ belongs to
the Schwartz class for  all $t\in  (t_0,t_1)$.\ret

The Wave-Maps equation is also invariant with respect to the change
of scale $\phi(t,x) \to \phi(\lambda t, \lambda x)$ for any positive
$\lambda \in \R$. In $(2+1)$ dimensions, it is easy to see that the
energy $E[\phi]$ is dimensionless with respect to this scale
transformation. For this reason, the problem we consider is called
{\em energy critical}.\ret

For the evolution \eqref{main_eq}, a
local well-posedness theory in Sobolev spaces ${H^s} \times {H}^{s+1}$
for $s$ above scaling, $s > 1$, was established some time ago.
See \cite{MR94h:35137} and \cite{MR99c:35163}, and references therein.
 The small data Cauchy-problem in the scale invariant Sobolev space is, by now,
also well understood. Following work of the second author
\cite{Tataru_WM2} for initial data in a scale invariant Besov space,
Tao was the first to consider the wave map equation with small energy
data.  In the case when the target manifold is a sphere, Tao
\cite{Tao_WM} proved global regularity and scattering for small energy
solutions.  This result was extended to the case of arbitrary compact
target manifolds by the second author in \cite{Tataru_WM1}.  Finite
energy solutions were also introduced in \cite{Tataru_WM1} as unique
strong limits of classical solutions, and the continuous dependence of
the solutions with respect to the initial data was established. The
case when the target is the hyperbolic plane was handled by Krieger
\cite{MR2094472}. There is also an extensive literature devoted to the
more tractable higher dimensional case; we refer the reader to
\cite{MR99j:58209}, \cite{MR2001m:35200}, \cite{MR2002h:58051},
\cite{MR1990880}, \cite{MR2016196}, and \cite{MR1890048} for more information.

To measure the dispersive properties of solutions $\phi$ to the
Wave-Maps equation, we shall use a variant of the standard dispersive
norm $S$ from \cite{Tataru_WM1}. This was originally defined in
\cite{Tao_WM} by modifying a construction in \cite{Tataru_WM2}.  $S$
is used together with its companion space $N$ which has the linear
property (precise definitions will be given shortly):
\begin{equation}
    \lp{\phi}{S[I]} \ \lesssim \ \lp{\phi}{L^\infty_t(L^\infty_x)[I]} + \lp{\phi[0]}{\dot H^1 \times L^2}
    + \lp{\Box\phi}{N[I]} \ . \notag
\end{equation}
The main result in \cite{Tataru_WM1} asserts that global regularity and
scattering hold for the small energy critical problem:\ret

\begin{thm}\label{Tat_Thm}
The wave maps equation \eqref{main_eq} is globally well-posed for
small initial data  $ \phi[0] \in \dot H^{1} \times L^2$ in the
following sense:\ret

\begin{enumerate}[(i)]
\item (Classical Solutions) If the initial data $\phi[0]$ is constant
  outside of a compact set and $C^\infty$, then there is a global
  classical solution $\phi$ with this data.\ret
%  which for $s \geq 1$
%   satisfies the global bounds:
% \begin{equation}
% \|\phi[t]\|_{L^\infty(\dot H^s \times \dot H^{s-1})} \lesssim
% \|\phi[0]\|_{\dot H^s \times\dot  H^{s-1}}
% \end{equation}
\item (Finite Energy Solutions) For each small initial data set in
  $\phi[0] \in \dot H^{1} \times L^2$ there is a global solution $\phi
  \in S$, obtained as the unique $S$ limit of classical solutions, so
  that:
    \begin{equation}
        \lp{\phi}{S} \ \lesssim \ \lp{\phi[0]}{\dot{H}^{1} \times L^2}
    \end{equation}\label{small_data}
\item (Continuous dependence) The solution map $\phi[0] \to \phi$  from
a small ball in   $\dot H^{1} \times L^2$ to $S$ is continuous.

\end{enumerate}
\end{thm}\ret

We remark that due to the finite speed of propagation one can also
state a local version of the above result, where the small energy
initial data is taken in a ball, and the solution is defined in the
corresponding uniqueness cone. This allows one to define large data
finite energy solutions:

\begin{defn}
  Let $I$ be a time interval. We say that $\phi$ is a finite energy
  wave map in $I$ if $\phi[\cdot] \in C(I;\dot H^{1} \times L^2)$
and, for each $(t_0,x_0) \in I$ and $r > 0$ so that
  $E[\phi(t_0)](B(x_0,r))$ is small enough, the solution $\phi$
  coincides with the one given by Theorem~\ref{Tat_Thm} in the
  uniqueness cone $I \cap \{ |x-x_0| + |t-t_0| \leqslant r\}$.
\end{defn}

In this work we consider a far more subtle case, which is a
conditional version of the large data problem. It is first important
to  observe that for general targets the above theorem cannot be
extended to arbitrarily large $C^\infty$ initial data, and that this
failure can be attributed to several different mechanisms. For
instance any harmonic map $\phi_0: \R^2 \to \mathcal{M}$ yields a time
independent wave-map which does not decay in time therefore it does
not belong to $S$. More interesting is that for certain non-convex
targets, for example when we take $\mathcal{M}=\mathbb S^2$, finite
time blow-up of smooth solutions is possible (see \cite{KST}, \cite{RS}).
In this latter case, the blow-up occurs along a family of rescaled
harmonic maps.
% Our primary goal  will be to see how far one can push the regularity theory for
% \eqref{main_eq} using \emph{only} two things:
% \begin{itemize}
%     \item Peturbative estimates applied to \eqref{main_eq}.
%     \item The conservation of energy \eqref{energy}.
% \end{itemize}
% In particular, we will make no use here of special monotonicity formulas (e.g. Morawetz type
% estimates).
To avoid such Harmonic-Map based solutions, as well as other
possible concentration scenarios,
in this article we prove a  conditional regularity theorem:\ret

\begin{thm}[Energy Dispersed Regularity Theorem]\label{main_thm}
There exist two functions $1\ll F(E)$ and $0 < \epsilon(E)\ll 1$ of
the energy \eqref{energy} such that the following statement is true.
If $\phi$ is a finite energy solution to \eqref{main_eq} on the open interval
$(t_1,t_2)$ with energy $E$ and:
\begin{equation}
     \sup_k \lp{P_k\phi}{L^\infty_{t,x}[(t_1,t_2)\times \R^2]} \ \leqslant \  \epsilon(E)
     \label{energy_disp}
\end{equation}
then one also has:
\begin{equation}
    \lp{\phi}{S(t_1,t_2)} \ \leqslant \ F(E) \ . \label{S_est}
\end{equation}
Finally, such a solution $\phi(t)$ extends in a regular way to a neighborhood
of the interval $I=[t_1,t_2]$.
\end{thm} \ret

\begin{rem}
In  Section \ref{s:main_proof}, Theorem \ref{main_thm_FE}, we
shall state a slightly stronger version of this result which uses the language of
frequency envelopes from \cite{Tao_WM}. In particular, we will show
the energy dispersion bound \eqref{energy_disp} implies that a
certain range of subcritical Sobolev norms may only grow by a
universal energy dependent factor. Put
another way, one may interpret this restatement of Theorem
\ref{main_thm} as saying that  \emph{in the energy dispersed
scenario, the Wave-Maps equation becomes subcritical in the sense
that there is a quasi-conserved norm of \textbf{higher} regularity
than the physical energy}. This information, coupled with the
standard regularity theory for Wave-Maps (e.g. see \cite{Tataru_WM1}) provides
us with the continuation property.
\end{rem} \ret

\begin{rem}
  The result in this article is stated and proved in space dimension
  $d=2$. However, given its perturbative nature, one would expect to
  have a similar result in higher dimension $d \geqslant 3$ as well. That
  is indeed the case. There are two reasons why we have decided to stay with
  $d=2$ here. One is to fix the notations. The second, and the more
  important reason, is to avoid lengthening the paper with an additional
  argument in Section~\ref{s:main_proof}, which is the only place in
  the article where the conservation of energy is used. In higher
dimensions, this aspect would have to be replaced by an almost conservation
of energy, with errors controlled by the energy dispersion parameter $\epsilon$.
\end{rem} \ret

\begin{rem}
The proof of  Theorem \ref{main_thm} allows us to obtain  explicit formulas
for $F(E)$ and $\epsilon(E)$. Precisely, in the conclusion of the
proof of Corollary \ref{cor:maint} below, we show that these parameters
may be chosen of the form:
\begin{equation}
    F(E)  \ =  \ e^{C e^{E^M}} \ ,
    \qquad \epsilon(E)  \ =  \ e^{-C e^{E^M}} \ , \notag
\end{equation}
with $C$ and $M$ sufficiently large.
\end{rem}\ret

As a consequence of the frequency envelope version of this result
in  Theorem \ref{main_thm_FE} we can also state a weaker non-conditional
version of the above result:

\begin{cor} \label{maint_nonpert}
There exists two functions $1\ll F(E)$ and $0 < \epsilon(E)\ll 1$ of
the energy \eqref{energy} such that for each initial data $\phi[0]$
satisfying:
\begin{equation}
     \sup_k \lp{P_k\phi[0]}{\dot H^1 \times L^2} \ \leqslant \  \epsilon(E)
     \label{energy_dispa}
\end{equation}
there exists a unique global finite energy solution $\phi \in S$,
satisfying:
\begin{equation}
    \lp{\phi}{S} \ \leqslant \ F(E) \ . \label{S_esta}
\end{equation}
which depends continuously on the initial data. If in addition the initial
data is smooth, then the solution is also smooth.
\end{cor}\ret

Our main interest in Theorem \ref{main_thm} is to combine it  with
the results of our concurrent work \cite{ST_1}, which together implies a
full regularity theory for Wave-Maps. In this context, one may  view
Theorem \ref{main_thm} as providing  a  ``compactness continuation''
principle, which roughly states that there is the following
dichotomy for  classical Wave-Maps defined on the open time interval
$(t_0,t_1)\times \R^2$:
\begin{enumerate}
    \item The solution $\phi$ continues to a neighborhood  of the closed
    time interval $[t_0,t_1]$ as a classical Wave-Map.
    \item The solution $\phi$ exhibits a compactness property on a sequence of rescaled times.
\end{enumerate}
In particular, the second case may used with the energy estimates
from \cite{ST_1} to conclude that a portion of any singular Wave-Map
must become stationary, and via compactness must therefore rescale
to a Harmonic-Map of non-trivial energy. This was known as the
\emph{bubbling conjecture} (see the introduction of \cite{ST_1} for more
background).\ret

Finally, we would like to remark that results similar in spirit to
the ones of this paper and \cite{ST_1} have been recently announced.
In the case where $\mathcal{M}=\mathbb{H}^n$, the hyperbolic spaces,
globally regularity and scattering follows from the program of Tao
\cite{Tao_LWM1}, \cite{Tao_LWM2}, \cite{Tao_LWM3}, \cite{Tao_LWM4}, \cite{Tao_LWM5}
and \cite{Tao_LWM6}. In the case where the target
$\mathcal{M}$ is a negatively curved Riemann surface, Krieger and Schlag \cite{KSc}
provide global regularity and scattering via a modification of the
Kenig-Merle method \cite{KM},
which uses as a key component suitably defined Bahouri-Gerard \cite{BG}
type decompositions. \ret

\textbf{Acknowledgements:} The authors would like to thank Manos Grillakis, Sergiu Klainerman,
 Joachim Krieger, Matei Machedon,  Igor Rodnianski, and Wilhelm Schlag
for many stimulating discussions over the years regarding the wave-map
problem. We would also especially like to thank Terry Tao
for several key discussions on the nature of induction-on-energy
type proofs.

\ret
%-------------------------------------------------------------------------

\subsection{A guide to reading the paper}

The paper has a ``two tier'' structure, whose aim is to enable the
reader to get quickly to the proof of the main result in
Section~\ref{s:main_proof}.  The first tier consists of
Sections~\ref{old_est_section},~\ref{new_est_section} and ~\ref{s:main_proof},
which play the following roles: \ret

{\bf \underline{Section~\ref{old_est_section}}} is where the notations are set-up.
In addition, in Proposition~\ref{standard_prop1} we review the linear,
bilinear, trilinear and Moser estimates concerning the $S$ and $N$
spaces, as proved in \cite{Tao_WM}, \cite{Tataru_WM1}. The $N$ space
we use is the same as in \cite{Tataru_WM2},\cite{Tao_WM}. For the $S$
space we begin with the definition in \cite{Tao_WM} and add to it the
Strichartz norm $\uS$ defined later in \eqref{phys_str}. This
modification costs almost nothing, but saves a considerable amount of
work in several places. \ret

{\bf \underline{Section~\ref{new_est_section}}} contains new contributions,
reaching in several directions:
\begin{itemize}
\item  {\em Renormalization.} A main difficulty in the study of wave maps is
that the nonlinearity is non-perturbative at the critical energy
level. A key breakthrough in the work of Tao~\cite{Tao_WM} was a
renormalization procedure whose aim is to remove the nonperturbative
part of the nonlinearity. However, despite subsequent improvements in
\cite{Tataru_WM1}, this procedure only applies to the small data
problem. We remedy this in Proposition~\ref{propphiu}, introducing a
large data version of the renormalization procedure. This applies
without any reference to the energy dispersion bounds.

\item {\em $S$ bounds for the paradifferential evolution with a large
  connection.} After peeling off the perturbative part of the
nonlinearity in the wave map equation, one is left with a family of
frequency localized linear paradifferential evolutions as in
\eqref{reduced_lin_eq}.  In the case of the small data problem, by
renormalization this turns into a small perturbation of the linear
wave equation. Here this is no longer possible, as the connection coefficients
$A^\alpha$ are large, and this cannot be improved using the energy dispersion.
However, what the energy dispersion allows us to do is to produce
a large frequency gap $m$ in \eqref{reduced_lin_eq}. As it turns out, this is all that is needed in order to have good estimates for
the equation \eqref{reduced_lin_eq}.

\item {\em New bilinear and trilinear estimates} which take advantage of the
energy dispersion. The main bilinear bound is the $L^2$ estimate in
Proposition~\ref{p:matchfreq}. Ideally one would like to have such
estimates for functions in $S$, but that is too much to ask. Instead
we introduce a narrower class $\W$ of ``renormalizable'' functions
$\phi$ of the form $\phi = U^\dagger w$, where $U \in S$ is a gauge transformation,
while for $w$ we control both $\|w\|_S$ and $\|\Box w\|_N$.
As a consequence of  Proposition~\ref{p:matchfreq} and the
more standard bounds in  Proposition~\ref{standard_prop1}, we later
derive the trilinear estimates in Proposition~\ref{big_tri_prop}, which
are easy to apply subsequently in the proof of our main Theorem.
\end{itemize} \ret

{\bf \underline{Section~\ref{s:main_proof}}} contains the proof of
Theorem~\ref{main_thm_FE}, which is a stronger frequency envelope
version of Theorem~\ref{main_thm}. This is done via an induction on
energy argument. The non-inductive part of the proof is separated into
Propositions~\ref{p:phi_est_low} and \ref{p:phi_est_high}, whose aim
is to bound in two steps the difference between a wave-map $\phi$ and
a lower energy wave map $\tphi$ whose initial data is essentially
obtained by truncating in frequency the initial data for $\phi$.
The arguments in this section use exclusively the results
in Sections~\ref{old_est_section},~\ref{new_est_section}.
\ret

The second tier of the article contains the proofs of all the results
stated in Sections~\ref{old_est_section},~\ref{new_est_section}, with the
exception of those already proved in \cite{Tao_WM} and \cite{Tataru_WM1}.
These are organized as follows:\ret

{\bf \underline{Section~\ref{ext_sect}}}'s content is as follows:
\begin{itemize}
\item{\em  A full description of the
$S$ and $N$ spaces.} Some further properties of these spaces are
detailed in Proposition~\ref{standard_prop2}; most of these are from
\cite{Tao_WM} and \cite{Tataru_WM1}, with the notable exception of the
fungibility estimate \eqref{N_fung}. The bound \eqref{N_fung} is proved using only the definition of $N$.
\item{\em Extension properties for the $S$ space.}
In most of our analysis we do not work with the spaces $S$ and
$N$ globally, instead we use their restrictions to time intervals,
$S[I]$ and $N[I]$. This is not important for $N$, since the
multiplication by a characteristic function of an interval is bounded
on $N$. However, that is not the case for $S$. One can define the
$S[I]$ norm using minimal extensions.  But in our case, we also
need good control of the energy dispersion and of the high modulation
bounds for the extensions. To address this, in Proposition~\ref{ext_prop}
we introduce a canonical way to define the extensions which
obey the appropriate bounds, and which also produce an equivalent $S[I]$ norm.
\item{\em Strichartz and $L^2$ bilinear estimates.} Using the $U^p$
  and $V^p$ spaces\footnote{For further information on the $U^p$ and
    $V^p$ spaces we refer the reader to \cite{KTucp}, \cite{KTnls},
    \cite{HHK}} associated to the half-wave evolutions, we first show
  that  solutions to the wave equation $\Box \phi = F$ with a right
  hand side $F \in N$ satisfy the full Strichartz estimates.  The
  fungibility estimate \eqref{N_fung} plays a significant role here,
  as it allows us to place the solution $\phi$ in a $V^2$ type space,
  see \eqref{SNV2}. A second goal is to prove $L^2$ bilinear bounds
  for products of two such inhomogeneous waves with frequency
  localization and angular frequency separation, see
  Lemma~\ref{wolff_lem}.  This is accomplished using the Wolff
  \cite{Wolff}-Tao \cite{Taobi} type $L^p$ bilinear estimates with $p
  < 2$.
\end{itemize}\ret

{\bf \underline{Section~\ref{s:matchfreq}}} is devoted to the proof of the bilinear
null form estimates in Proposition~\ref{p:matchfreq}. A preliminary
step, achieved in Lemma~\ref{lnore}, is to establish  the counterpart
of the bounds \eqref{ednullw} and \eqref{bal_l2_N} in the absence of the
renormalization factor. The proofs here use only Lemma~\ref{wolff_lem}
and the estimates in Propositions~\ref{standard_prop1},~\ref{standard_prop2}.
\ret

{\bf \underline{Section~\ref{s:der_tri}}} contains the proof of the trilinear
estimates in Proposition~\ref{big_tri_prop}. There are a number of
dyadic decompositions and multiple cases to consider, but this is
largely routine, using either Proposition~\ref{p:matchfreq} or the
estimates in Propositions \ref{standard_prop1} and \ref{standard_prop2}.\ret

{\bf \underline{Section~\ref{s:renormalize}}} is concerned with the
construction of the gauge transformation in
Proposition~\ref{propphiu}. The discrete inductive construction in
\cite{Tao_WM}, \cite{Tataru_WM1} is replaced with a continuous version
which serves to insure that the renormalization matrices $U_{,<k}$ are
exactly orthogonal. To allow for wave-maps which are large in $S$, we
need to forego the simpler inductive way of proving $S$ estimates for
$U_{,k}$ and instead build them up in a less direct fashion using
iterated paradifferential type expansions. On the positive side, this
proof uses only the estimates in
Propositions~\ref{standard_prop1},~\ref{standard_prop2}.  \ret

{\bf  \underline{Section~\ref{s:lin_sect}}} is devoted to the proof of the linear
bounds for the paradifferential equation in
Proposition~\ref{p:para}. A key element in this proof is the gauge
transformation in Proposition~\ref{propphiu}, combined with the
trilinear estimate \eqref{standard_est_tri}.  This would suffice for
connections $A_\alpha$ arising from wave maps $\phi$ which are small in
$S$. However, in our case we need to handle large wave maps, and a
different source of smallness is required.  This is provided by the
large size of the frequency gap $m$, which leads to energy
conservation with small $O(2^{-cm})$ errors. Feeding these almost
apriori energy and characteristic energy bounds back in the bilinear
and trilinear null form estimates turns out to suffice to estimate the
large trilinear contributions, again modulo terms which are small,
i.e. $O(2^{-cm})$.\ret

{\bf \underline{Section~\ref{s:wm_struct}}}'s goal is to provide the description
of finite $S$ norm wave-maps in Proposition~\ref{p:wm_struct}.  The
renormalization bound is a direct consequence of
Proposition~\ref{p:para} and Proposition~\ref{propphiu}. The partial
fungibility of the $S$ norm is tied to the fungibility of the $N$ norm
in the renormalized setting, although the proof is somewhat more technical.

Under the assumption of small energy dispersion, the smallness of high
modulations is given by the trilinear estimate~\eqref{core_L2}. After
that, the bound \eqref{core_N} combined with
Proposition~\ref{propphiu} lead quickly to the frequency envelope
control in \eqref{freq_env_control} via a bootstrap argument.\ret

{\bf \underline{Section~\ref{s:cut}}} contains the proof of the data truncation
result in Proposition~\ref{p:cut}. The argument is self-contained and uses
only $L^2$ type methods.

%\ret
%-------------------------------------------------------------------------

%\subsection{Brief Overview of the Proof}
%A HEURISTIC VERSION OF THE PROOF, AND
%COMPARISON WITH PREVIOUS INDUCTION ON ENERGY RESULTS (NEED TO ADD A LOT OF REFERENCES HERE).
%INCLUDE STUFF LIKE:\ret

%By the small data result in Theorem~\ref{small_data}, we know that
%\eqref{main_eq} has a global regular solution as long as
%\eqref{energy} is such that $E_0\leqslant \epsilon_0$ for some small
%$\epsilon_0 > 0$. We now try to see what is the major obstruction to
%expanding the range of $E_0$ such that regularity with good estimates
%hold. We will see that the problem eventually reduces to a certain
%interplay between linearized stability and certain ``semi-linear''
%interactions which together govern the ability of the non-linear
%solution $\phi$ to spontaneously create large amounts of energy on
%higher and higher frequency blocks as time progresses.

\ret
%-------------------------------------------------------------------------
%%%%%%%%%%%%%%%%%%%%%%%%%%%%%%%%%%%%%%%%%%%%%%%%%%%%%%%%%%%%%%%%%%%%%%%%%%
%-------------------------------------------------------------------------

\section{Standard Constructions, Function Spaces, and
Estimates}\label{old_est_section}

In this section we record the standard portion of the framework we
shall use in our primary demonstration of Theorem \ref{main_thm}.
While we aim to keep our account of things self contained, we also
assume that the reader is thoroughly familiar with the content of
the two papers \cite{Tao_WM} and \cite{Tataru_WM1}. For the sake of
completeness, in Section \ref{ext_sect} below we also include proofs
of several results not contained in these two works, but which are
needed for the more detailed analysis of this paper. Further
notation and estimates that are not needed for the ``first
tier'' of our demonstration of Theorem \ref{main_thm}, but are
needed in later technical sections are also given
in Section \ref{ext_sect} below.\ret

The symbols $\lesssim$, $\gtrsim$, $\sim$, $\ll$,
and $\gg$ are defined with their usual meanings. The constants in these
notations are allowed to vary from line to line.

\ret
%--------------------------------------------------------------------------------

\subsection{Constants}
There will be a number of large and small constants in the present
work. For the most part these are flexible, although the specific
construction of  $F(E)$ and $\epsilon(E)$ from Theorem
\ref{main_thm} will be sensitive to each other as well as to choices
of other constants. Lower case Greek letters such as $\delta$,
$\epsilon$, and $\eta$ will always denote small quantities. We shall
employ a \emph{globally} defined string of small constants:
\begin{equation}
        \delta_0 \ \ll \ \delta_1 \  \ll \ \delta_2 \ \ll \ 1 \ ,
    \qquad \qquad \delta_{i} \ \leqslant  \ \delta^{100}_{i+1} \ .
        \label{delta_list}
\end{equation}
As it occurs often in the sequel, we will set $\delta=\delta_2$ throughout.  For the
convenience of the reader we list here the purposes of these
constants, which all measure various fractional frequency gains in
our dyadic estimates:
\begin{itemize}
        \item The base constant $\delta$ enters our proof through
        the various multilinear estimates for the $S$ and $N$ spaces
        listed below (e.g. in the current section).
        It also influences any portion of our argument
        which is a direct consequence of these estimates, but has
        nothing to do with directly bootstrapping large data Wave-Maps. For
        example, $\delta$ also represents various dyadic gains in our gauge
        construction (see Proposition \ref{propphiu}).
        \item The constant $\delta_1$  measures a small fractional gain
        coming from energy dispersion in $L^2$ and $N$-norm
    null form estimates.
        It enters our proof through estimates \eqref{core_L2} and
        \eqref{ext_L2}, and variations thereof.
        \item The constant $\delta_0$ is reserved for  slowly
        varying frequency envelopes, and for the smallest fractional
        quantities built from the energy dispersion constant
        $\epsilon$. It enters in the core part of our proof of
        Theorem \ref{main_thm}, and is the assumption on the
        frequency envelopes of Proposition \eqref{big_tri_prop}.
\end{itemize}\ret

Large quantities, for example $C$, $F$, $K$, and $M$ will be used in
various contexts as constants in estimates and the size of norms
which are \emph{not} globally defined. We will also often use $m$ to
denote a (possibly) large integer which represents various gaps in
frequency truncations.\ret

To denote  growth and dependence of various estimates on that growth
we employ the following notation in the sequel:

\begin{defn}[Complexity Notation]
We say that a positive function $f(y)$ is of ``polynomial type'' if
$f(y) \leqslant y^M$ for some constant $M$ as $y\to\infty$.
%We say
%that $f(y)$ is of ``exponential type'' if $\ln(f)$ is of polynomial
%type.
We use the notation:
\begin{equation}
    A \ \lesssim_F \ B \ , \notag
\end{equation}
if $A\leqslant K(F) B$ for some function $K$ of polynomial type.
This notation does \emph{not} fix $K$ from line to line, although
$K$ is fixed on any single line where it occurs.
\end{defn}

\ret
%--------------------------------------------------------------------------------

\subsection{Basic harmonic analysis}
As usual we denote by $\xi$ and $\tau$ the spatial and temporal
Fourier variables (resp). We set up both discrete and continuous
spatial Littlewood-Paley (LP) multipliers:
\begin{equation}
        I -P_{-\infty} \ = \ \sum_k P_k \ ,
        \qquad\qquad\qquad I -P_{-\infty}\ = \ \int_{-\infty}^\infty P_k dk \ . \notag
\end{equation}
For the purposes of trichotomy, these two sets of multipliers are
interchangeable, and we will only distinguish them by the use of
$\sum$ or $\int$ in identities. However, for the purposes of proving
Moser type estimates or constructing gauge transformations, the
integral definition of LP projections is essential. We refer the
reader to \cite{Tataru_WM1} for an earlier use of continuous LP
multipliers, and further information. We often denote by
$\phi_k=P_k\phi$. If $\phi$ is any affinely Schwartz function, the
above notation means that we have the identities:
\begin{equation}
        \phi - \lim_{|x|\to\infty}\phi(x) \ = \ \sum_k \phi_k \
        = \ \int_{-\infty}^\infty \phi_k  dk \ . \notag
\end{equation}
Therefore, care must be taken to add constants back into certain
estimates involving very low frequencies.\ret

Many times in the sequel we shall have use for the inequality:
\begin{equation}
    \lp{P_{\mathcal{B}}\phi}{L^q_x} \ \lesssim \ |\mathcal{B}|^{(\frac{1}{r}-\frac{1}{q})}
    \lp{\phi}{L^r_x} \ , \label{Bernstein}
\end{equation}
where $\mathcal{B}\subseteq \R^2_\xi$ is a frequency box.
Furthermore, the $P_k$ multipliers enjoy a commutator structure as
follows:
\begin{equation}
    P_k(\phi\psi) \ = \ \phi \psi_k + L(\nabla_x \phi,2^{-k}\psi) \ , \label{commutator}
\end{equation}
where the bilinear form $L$ is translation invariant and bounded on
all Lebesgue type spaces. Such multilinear expressions occur often
in the sequel. We call an expression of the from:
\begin{equation}
        L(\phi^{(1)},\ldots,\phi^{(k)})(x) \ = \ \int \phi^{(1)}(x +
        y_1)\cdot\cdot\cdot \phi^{(k)}(x+y_k) d\mu(y_1,\ldots,y_k)
        \ , \notag
\end{equation}
where:
\begin{equation}
         \int |d\mu| \ \lesssim \ 1 \ , \notag
\end{equation}
``disposable''. Any disposable operator generates a family of
estimates from any single product estimate involving translation
invariant norms in the usual way (see \cite{Tao_WM}).\ret

We will also use the  variable notation for frequency envelopes from
\cite{Tataru_WM1} (see \cite{Tao_WM} for another definition):

\begin{defn}[Frequency Envelopes]
A frequency envelope $\{c_k\}$ is called
``$(\sigma,\Delta)$-admissible'' if it obeys the bounds:
\begin{equation}
    2^{\sigma(j-k)}c_{k} \ \leqslant \ c_{j} \ \leqslant \ 2^{\Delta(k-j)}c_{k}
    \ , \label{freq_env_defn}
\end{equation}
for any  $j < k$, where $0< \sigma \leqslant \Delta$. If
$\lp{\phi}{Y}$ is any non-negative real valued functional, and
$\{c_k\}$ is a frequency envelope, we define:
\begin{equation}
        \lp{\phi}{Y_c} \ := \
        \sup_k  c_k^{-1}\lp{P_k\phi }{Y} \ . \label{X_c_def}
\end{equation}
There is an exception to this notation for the norm $S[I]$
introduced below, in which case we set:
\begin{equation}
        \lp{\phi}{S_c[I]} \ := \ \lp{\phi}{L^\infty_t(L^\infty_x)[I]} +
        \sup_k c_k^{-1}\lp{\phi}{S_k[I]} \ .  \notag
\end{equation}
\end{defn}\ret

Frequency envelopes may be defined in either the discrete or
continuous settings. It is easy to see that for any such frequency
envelope we have the pair of sum rules (uniformly):
\begin{align}
    \sum_{k'\leqslant k} 2^{Ak'} c_{k'} \ &\lesssim \ (A-\Delta)^{-1} 2^{Ak} c_k \ ,
    &A > \Delta \ , \label{sum_rule1}\\
    \sum_{ k \leqslant k'} 2^{-ak'} c_{k'} \ &\lesssim \ (a-\sigma)^{-1} 2^{-ak} c_k \ ,
    &a > \sigma \ , \label{sum_rule2}
\end{align}
with similar bounds for integrals. These two inequalities capture
the essence of every use we have for the $\{c_k\}$ notation, which
is simply to bookkeep (resp.) $Low \times Low\Rightarrow High$ and
$High \times High \Rightarrow Low$ frequency cascades.

\ret
%--------------------------------------------------------------------------------

\subsection{Function spaces and standard estimates}

We use the function spaces $S$ and $N$ from \cite{Tataru_WM1}--\cite{Tataru_WM2} and
\cite{Tao_WM} with only a few minor modifications.
The spaces of restrictions of $S$ and $N$ functions
to a time interval $I$ are denoted by $S[I]$, respectively $N[I]$,
with the induced norms. The first part of
our proof does not use the precise structure of these
spaces, only the following statement:\ret

\begin{prop}[Standard Estimates and Relations: Part
I]\label{standard_prop1} Let $F$, $\phi$, and $\phi^{(i)}$ be a
collection of test functions, $I\subseteq \mathbb{R}$ any
subinterval (including $\R$ itself). Then there exists
function spaces $S[I]$ and $N[I]$ with the following properties:
\begin{itemize}
    \item (Triangle Inequality for $S$) Let $I=\cup_i^K I_i$ be a decomposition of $I$ into consecutive
    intervals, then the following bounds hold (uniform in $K$):
    \begin{equation}
        \lp{\phi}{S[I]} \ \lesssim \ \sum_i \lp{\phi}{S[I_i]} \ . \label{S_decomp}
    \end{equation}
    \item (Frequency Orthogonality) The spaces $S[I]$ and $N[I]$ are made up
    of dyadic pieces in the sense that:
    \begin{align}
            \lp{\phi}{S[I]}^2 \ &= \ \lp{\phi}{L^\infty_t(L^\infty_x)[I]}^2
            + \sum_k \lp{P_k \phi}{S[I]}^2 \ , \label{S_orth}\\
            \lp{\phi}{N[I]}^2 \ &= \
            \sum_k \lp{P_k\phi}{N[I]}^2 \ . \label{N_orth}
    \end{align}
        \item (Energy Estimates) We have that $L^1_t(L^2_x)[I]\subseteq N[I]$, and also the estimate:
        \begin{align}
            \lp{\phi_k}{S[I]} \ &\lesssim \ \lp{\Box \phi_k}{N[I]} +
                \lp{\phi_k[0]}
                {\dot{H}^1\times L^2} \ . \label{energy_est}
        \end{align}
        \item (Core Product Estimates) We have that:
        \begin{align}
            \lp{\phi^{(1)}_{< k+O(1)}\cdot\phi^{(2)}_{k} }{S[I]} \
            &\lesssim \ \lp{\phi^{(1)}_{k_1}}{S[I]}\cdot\lp{\phi^{(2)}_{k_2}}{S[I]} \ , \label{S_prod2}\\
            \lp{P_k(\phi^{(1)}_{k_1}\cdot\phi^{(2)}_{k_2}) }{S[I]} \
            &\lesssim \ 2^{-(\max\{k_i\}-k)}\lp{\phi^{(1)}_{k_1}}
            {S[I]}\cdot\lp{\phi^{(2)}_{k_2}}{S[I]} \ , \label{S_prod3}\\
            \lp{P_k(\phi_{<k+O(1)}\cdot F_k)}{N[I]} \ &\lesssim \
            \lp{\phi}{S[I]}\cdot\lp{F_{k}}{N[I]} \ , \label{N_prod_est1}\\
            \lp{P_k(\phi_{k_1}\cdot F_{k_2})}{N[I]} \ &\lesssim \ 2^{-\delta(k-k_2)_+}
            \lp{\phi_{k_1}}{S[I]}\cdot\lp{F_{k_2}}{N[I]} \ , \label{N_prod_est2}
        \end{align}
       \item (Bilinear Null Form Estimates) We have that:
        \begin{align}
        \lp{\!P_k\big(\partial^\alpha \phi^{(1)}_{k_1} \cdot\partial_\alpha
                \phi^{(2)}_{k_2}\big)\! }{L^2_t(L^2_x)[I]}
                &\lesssim  2^{\frac{1}{2}\min\{k_i\}}2^{-(\frac{1}{2}+\delta)(\max\{k_i\}-k)}
                \prod_i \lp{\!\phi^{(i)}_{k_i}\!}{S[I]}
                \ , \label{standard_L2_est_bi}\\
                \lp{P_k (\partial^\alpha \phi^{(1)}_{k_1}\cdot\partial_\alpha
                \phi^{(2)}_{k_2})}{N[I]}
                &\lesssim  2^{-\delta(\max\{k_i\}-k)} \prod_i \lp{\phi^{(i)}_{k_i}}{S[I]}
                \ . \label{standard_est_bi}
\end{align}
 \item (Trilinear Null Form Estimate) We have that:
\begin{align}
                 \lp{P_k (\phi^{(1)}_{k_1}\! \cdot \partial^\alpha
                 \phi^{(2)}_{k_2} \!
                 \cdot \partial_\alpha
                \phi^{(3)}_{k_3})}{N[I]}
                &\lesssim  2^{-\delta(\max\{k_i\}-k)}2^{-\delta(k_1-\min\{k_2,k_3\})_+}
                \prod_i \lp{\phi^{(i)}_{k_i}}{S[I]}
                 . \label{standard_est_tri}
        \end{align}
        \item (Moser Estimates) Let $G$ be any bounded function with uniformly bounded derivatives,
        and $\{c_k\}$ a $(\delta,\Delta)$-admissible frequency envelope.
    Then there exists a universal $K>0$ such that:
        \begin{align}
            \lp{G(\phi)}{S[I]} \ &\lesssim \ \lp{\phi}{S[I]}(1+\lp{\phi}{S[I]}^K)
            \ , \label{basic_moser}\\
        \lp{G(\phi)}{S_c[I]} \ &\lesssim \ \lp{\phi}{S_c[I]}(1+\lp{\phi}{S[I]}^K)
            \ . \label{basic_moser_c}
        \end{align}
\end{itemize}
\end{prop}
The space $N$ is the same one as used in  \cite{Tao_WM}, \cite{Tataru_WM1}. To obtain
the space $S$ we start with the one used in \cite{Tao_WM}, \cite{Tataru_WM1} and add
the control of the Strichartz norms $\uS$
defined in \eqref{phys_str}. The bound \eqref{S_decomp} is relatively
straightforward; we prove it in Section~\ref{ext_sect}.
The relations \eqref{S_orth} and \eqref{N_orth}
can be thought of as a part of the definition of the spaces $S$, $N$ starting from their
dyadic versions. The linear estimate \eqref{energy_est}
was proved in  \cite{Tao_WM}; here we show that we can add the Strichartz
component $\uS$ in Corollary~\ref{snstcor}.
The bounds \eqref{S_prod2}--\eqref{N_prod_est2}
as well as \eqref{standard_est_bi},\eqref{standard_est_tri}
were proved in \cite{Tao_WM}. In our context the proofs
of \eqref{S_prod2}, \eqref{S_prod3} need to be augmented to add the
control over the Strichartz norm $\uS$; this is a straightforward
matter which is left for the reader. The bound \eqref{standard_L2_est_bi}
is implicit in \cite{Tao_WM}, but for the reader's convenience
we prove it in Section~\ref{s:matchfreq}.

 The Moser estimates
\eqref{basic_moser} and \eqref{basic_moser_c} were proved in
\cite{Tataru_WM1}. Adding in the $\uS$ norm is again straightforward.
An interesting side remark is that in effect the addition of the
 $\uS$ norm to $S$ can be taken advantage of to simplify considerably
the proof of the Moser estimates in \cite{Tataru_WM1}.  In particular,
one can show that it is possible to take $K=2$. Since it does not
lead to significant improvements in the present article, we leave
this as an exercise for the reader.\ret

At several places in our argument, it will be necessary for us to
introduce some auxiliary norms. We choose to keep these separate from
$S$ defined above for notational purposes:

\begin{defn}[Auxiliary Energy and $X^{s,b}$ Type Norms]
We define:
\begin{align}
        \lp{\phi}{\underline{E}[I]} \ &:= \ \lp{\nabla_{t,x}\phi}{L^\infty_t (L^2_x)[I]}
        + \sup_{\omega}\lp{\snabla_{t,x}^\omega
        \phi}{L^\infty_{t_\omega}(L^2_{x_\omega})[I]} \ ,
        \label{uE_norm}\\
         \lp{\phi}{\underline{X}_k[I]} \ &:= \ 2^{-\frac{1}{2}k}\lp{\Box P_k \phi}{L^2_t(L^2_x)[I]}
         \ . \label{uX_norm}
\end{align}
Here the second term in the RHS of \eqref{uE_norm} represents the
energy of $\phi$ on characteristic hyperplanes, see \cite{Tataru_WM2},
\cite{Tao_WM}.  We also define $\underline{X}[I]$ as the square sum of
$\uX_k[I]$, and $\underline{X}_c[I]$ according to \eqref{X_c_def}.
Notice that there are no square sums or frequency localizations in the
norm $\underline{E}$. The size of this norm depends only on the
initial energy of any (global) classical solution to \eqref{main_eq}.
\end{defn}\ret

In the sequel, it will be also be notationally convenient for us to
work with the following definition which one should think of as a
variant of the $S[I]$ space introduced above. The reader should keep
in mind that this not even a  quasinorm due to the lack of any good
additivity property:

\begin{defn}(Renormalizable Functions)
Let $C> 0$ be a large parameter.  We define a
non-linear functional $\W_k$ on $S$ as follows:
\begin{multline}
    \llp{\phi}{\W_k } \
    := \ \inf_{U\in SO(d)} \Big[
    \big(\lp{U}{S\cap\uX} + \sup_{j\geqslant k}2^{C(j-k)}
    \lp{P_jU}{S\cap\uX} \big)\\
    \cdot\  \sup_{k'} 2^{|k'-k|}
    \big( \lp{P_{k'}(U \phi_k)[0]}{\dot{H}^1\times L^2}
    + \lp{P_{k'}\Box(U \phi_k)}{N} \big) \Big] \ . \label{big_S_def}
\end{multline}
The functionals $\W[I]$, $\W_c[I]$ are also defined as  above.
\end{defn}\ret

Notice that while the definition of $\W$ is  nonlinear, one still
has the scaling relation $\llp{\lambda\phi}{\W[I]}=\lambda\llp{\phi}{\W[I]}$.
The reader should note that while these bounds are cumbersome to state,
they are all natural in light of Propositions \ref{propphiu}--\ref{p:para} below.

\ret
%-------------------------------------------------------------------------
%%%%%%%%%%%%%%%%%%%%%%%%%%%%%%%%%%%%%%%%%%%%%%%%%%%%%%%%%%%%%%%%%%%%%%%%%%
%-------------------------------------------------------------------------

\section{New Estimates and Intermediate
Constructions}\label{new_est_section}

In this section we introduce the main technical components of the
paper. We begin with the core underlying tools that allow us to
handle more complicated constructions. In a later sub-section we derive
some further useful results that encapsulate many of the repetitive
computations in the sequel.\ret

\subsection{Core technical estimates and constructions}

The right hand side in the equation~\eqref{main_eq} is
nonperturbative even when the energy is small. In the case of larger
energies, it becomes quite a bit more difficult to handle things in
a perturbative manner. Therefore, we introduce a set of tools which
are general enough to handle large data situations. The first two of
these work without any additional properties (e.g. energy
dispersion), and form the technical heart of the paper. The first is
a novel gauge construction that should be of more general use. It
should be noted that this construction is stable regardless of the
size of the energy or the convexity properties of the target, as its
key properties depend only on the compactness of the underlying
gauge group.\ret

\begin{prop}[The ``Diffusion Gauge'']\label{propphiu}
Let $\phi$ be a wave-map in a time interval $I$ with energy $E$,
$S[I]$ norm $F$, and $S[I]$ norm $(\delta,\Delta)$-admissible
 envelope $\{c_k\}$. Let the antisymmetric $B$ be defined by:
\begin{equation}
        (B^a_b)_{<k}  \ =  \ \int_{-\infty}^k  \big(\mathcal{S}^a_{bc}(\phi) -
        \mathcal{S}^b_{ac}(\phi)\big)_{<k'-10}
        \phi^c_{k'}\ dk' \ , \label{B_def}
\end{equation}
where $\mathcal{S}^a_{bc}$ is a  smoothly bounded $(a,b)$ symmetric
matrix valued vector. We denote the integrand by $B_k$. Then for
each real number $k$ there exists an orthogonal matrix $U_{,<k}$
defined on all of $\mathbb{R}^{2+1}$ with the following
properties:\ret

\begin{itemize}
    \item ($U_{,<k}$ is a Sum of Frequency Localized Pieces in $S_k$)
    For each real number $k$ there exists a matrix $U_{,k}$ such
    that:
    \begin{equation}
            U_{,<k} \ = \ \int_{-\infty}^k U_{,k'} dk' \ , \notag
    \end{equation}
    where each $U_{,k} = U_{,<k} B_{k}$, and each $U_{,k}$ obeys the bounds:
    \begin{align}
            \lp{P_{k'} U_{,k}}{S \cap \uX}  &\lesssim_F
            2^{-\delta|k-k'|}2^{-C(k'-k)_+} c_{k} \ , \label{env_est1}\\
            %\sum_{|J|=0,1,2} 2^{3k-|J|}
            \| P_{k'} \nabla^J_{t,x} U_{,k}\|_{L^1_t(L^1_x)}
             &\lesssim_F 2^{(|J|-3)k}2^{ - C(k'-k)} c_k \ , \quad k' >
            k+10  , \ \ \ |J|\leqslant 2 \ ,
            \label{env_est1h}\\
        \lp{ P_{k'} \big(U_{,<k-20}\cdot G_k\big)}{N}  &\lesssim_F
            2^{-|k'-k|}\lp{G_k}{N} \ ,  \label{env_est3}\\
            \lp{P_k \big(\Box U_{,k_1}\cdot\psi_{k_2}\big)}{N}  &\lesssim_F
            2^{-|k-k_2|}
            2^{-\delta(k_2-k_1)}
            c_{k_1}\lp{\psi_{k_2}}{S} \ , \ \ k_1 < k_2-10
            \ . \label{env_est2}
    \end{align}
    In addition, if $\tc_k$ is a $(\delta_0,\Delta)$-admissible frequency envelope
    for the energy $\lp{\nabla_{t,x} \phi_k }{L^\infty_t(L^2_x)[I]}$ then we have a similar
    bound for $U_{,k}$:
        \begin{equation}
            \lp{ P_{k'} \nabla_{t,x} U_{,k}}{L^\infty_t(L^2_x)}   \
            \lse   \ 2^{-|k-k'|- C (k'-k)_+}\tc_{k}  \ .
            \label{env_est0}
        \end{equation}
    Here $C\geqslant 0$ is any constant.\ret
    \item (The Matrix $U$ Approximately Renormalizes $A_\alpha = \nabla_\alpha B$)
    We have the formula:
    \begin{equation}
            U^\dagger_{,<k} \nabla_\alpha U_{,<k} \ = \ \nabla_\alpha B_{<k} -
            \int_{-\infty}^k \big[B_{k'},U^\dagger_{,< k'} \nabla_\alpha U_{,<k'}\big]
            \, dk' \ . \label{renorm_form}
    \end{equation}
\end{itemize}

\end{prop}\ret

This result is proved in Section~\ref{s:renormalize}.
Next, we state a
technical proposition that will help us to deal with the
non-fungibility of the $S$ norm.
The wave map nonlinearity is nonperturbative. However,
due to the small energy dispersion, at fixed frequency we are
able to perturbatively replace the nonlinearity in the wave
map equation with a paradifferential term, i.e. a linear term
involving the lower frequencies of the wave map.
This term is large, and due to the non-fungibility of the $S$ norm,
it cannot be made small on small time intervals. Fortunately, it has another
redeeming feature, namely a large frequency gap (see $m$ below).
We take advantage of this in Section~\ref{s:lin_sect} to prove that:\ret

\begin{prop}[Gauge Covariant ${S[I]}$ Estimate] \label{p:para}
Let $\psi_k = P_k \psi$ be a solution to the linear problem:
\begin{equation}
    \Box \psi_k \ = \ -2A^\alpha_{<k-m}\partial_\alpha\psi_k  +  G
    \ , \label{reduced_lin_eq}
\end{equation}
where $A^\alpha_{<k-m}$ is the $\mathfrak{so}(N)$ matrix:
\begin{equation}
    (A^\alpha_{<k-m})^a_b \ = \
   \big(\mathcal{S}^a_{bc}(\phi) - \mathcal{S}^b_{ac}(\phi)\big)_{<k-m}\partial^\alpha
    \phi^c_{<k-m} \ . \label{red_con}
\end{equation}
Assume that ${\phi}$ is a classical Wave-Map on $I$ with the bounds:
\begin{equation}
    \lp{{\phi}}{\underline{E}[I]} +
    \lp{{\phi}}{\uX[I]} +
    \lp{{\phi}}{S[I]}\ \leqslant \ F \ . \label{A_tdphi_est}
\end{equation}
Furthermore, assume that $m \geqslant m(F)>20$, for a certain function
$m(F)\sim\ln(F)$  (to be defined in the proof). Then we have the
estimate:
\begin{equation}
    \llp{\psi_k}{\W[I]} \ \lesssim_F \
    \lp{\psi_k[0]}{\dot{H}^1\times L^2} +  \lp{ G}{N[I]}
    \ . \label{linearized_est}
\end{equation}
\end{prop}\ret

\begin{rem}\label{direct_para_rem}
As will become apparent in the proof of estimate
\eqref{linearized_est}, the \emph{only} use of the large frequency
gap parameter $m$ is to be able to bootstrap the RHS involving
$\psi_k$. In the sequel, there will be situations where one already
has good $S[I]$ norm bounds on $\psi_k$, and the task is to provide
a renormalization $w_{,k}$ such that $\Box w_{,k}$ has good $N$
norm bounds. Therefore,
we state the following:
\begin{itemize}
        \item Let $\psi_k$, $A^\alpha_{<k-m}$, and $G$ be defined as
        in Proposition \ref{p:para}. Then by simply assuming that
        $m>20$ we have the following estimate for $\psi_k$:
        \begin{equation}
            \llp{\psi_k}{\W[I]} \ \lesssim_F \
            \lp{\psi_k}{S[I]} +  \lp{P_k G}{N[I]}
            \ . \label{easy_linearized_est}
        \end{equation}
        \item Furthermore, in the above situation, the
        renormalization on the LHS of estimate
        \eqref{easy_linearized_est} is given by a matrix as in
        Proposition \ref{propphiu} where the pieces $B_j$ are defined from
        $A^\alpha_{<k-m}$ in the obvious way (this is of course true
        for estimate \eqref{linearized_est} as well).
\end{itemize}
See Remark \ref{direct_para_rem_proof} in Section \ref{s:lin_sect}
for more details.
\end{rem}\ret

Next, we state a gauged version of certain improved multilinear
estimates for the wave equation. Roughly speaking, these estimates
imply that matched frequency interactions in the RHS of
\eqref{main_eq} behave in a perturbative fashion in the presence of
energy dispersion. The heart of these estimates lies in the
Wolff-Tao bilinear estimates (see \cite{Wolff} and \cite{Taobi}) for
non-parallel waves, and the parallel wave cancelation property of
the ``$Q_0$ null structure'' which was originally
investigated in \cite{MR94h:35137}:\ret

\begin{prop}[Improved Matched Frequency Estimates]\label{p:matchfreq}
Let $\phi^{(i)}_{k_i}$ be functions localized at frequency $k_i$. Assume that
these functions are normalized as follows:
\begin{align}
    \llp{\phi^{(i)}_{k_i}}{\W[I]} \ &\leqslant \  1 \ ,
    & \lp{\phi^{(1)}_{k_1}}{L^\infty_t(L^\infty_x)[I]} \ &\leqslant \ \eta \ . \label{e_disp_ass}
\end{align}
Then the following estimates  hold:\ret

\begin{itemize}
    \item (Bilinear $L^2$ Estimate) We have that:
    \begin{equation}
        \lp{ \partial^\alpha \phi^{(1)}_{k_1} \partial_\alpha \phi^{(2)}_{k_2}}{L^2_t(L^2_x)[I]}
         \ \lesssim  \ 2^{\frac12 \max\{k_1,k_2\}} \eta^\delta \ .
        \label{ednullw}
    \end{equation}\ret
    \item (Bilinear $N$ Estimate) Assume that in addition to \eqref{e_disp_ass}
    we also have the high modulation bounds:
    \begin{align}
         \lp{ \Box \phi^{(1)}_{k_1}}{L^2_t(L^2_x)[I]}  \
         &\leqslant \ 2^{\frac{k_1}2} \eta \ ,
        &\lp{ \Box \phi^{(2)}_{k_2}}{L^2_t(L^2_x)[I]}  \
        &\leqslant \ 2^{\frac{k_2}2} \eta
         \ . \label{high_mod_ass}
    \end{align}
    Then the following estimate holds:
    \begin{equation}
        \lp{ \partial^\alpha \phi^{(1)}_{k_1} \partial_\alpha  \phi^{(2)}_{k_2} }{N[I]}
        \ \lesssim \ 2^{C|k_1-k_2|}\eta^\delta  , \label{bal_l2_N}
    \end{equation}
\end{itemize}
\end{prop}\ret

This is proved in Section~\ref{s:matchfreq}.
Finally, we list a technical a result concerning initial data
frequency truncation.  This does not preserve the space of functions
with values in $T\mathcal{M}$, so it has to be followed by a
non-linear physical space projection $\Pi$ back onto $T\mathcal{M}$.
We will show that in the energy dispersed case, this operation is
very well behaved in the energy norm. Theorems of this type may be
useful for other problems involving the need for a ``non-linear
Littlewood-Paley theory'' of functions with values in a manifold:\ret

\begin{prop}\label{p:cut}
For each $E > 0$ there exists $\epsilon_0 > 0$ so that for each
initial data set  $\phi[0]$ for \eqref{main_eq} with energy $E$ and
energy dispersion $\epsilon \leqslant \epsilon_0$  and
 $k,k_* \in \Z$ we have
\begin{equation}
    \lp{ P_{k}  \big(P_{< k_*} \phi[0] - \Pi ( P_{< k_*}
    \phi[0])\big)}{\dot H^1
    \times L^2}   \ \lesssim_E \ \epsilon^\frac{1}{4} 2^{-\frac{1}{2} |k-k_*|} \ . \label{cutphi}
\end{equation}
\end{prop}
This is proved in Section~\ref{s:cut} using Moser estimates and some
integral identities involving the continuous Littlewood-Paley theory
developed in \cite{Tataru_WM1}. \ret

%-----------------------------------------------------------------------------------

\subsection{Derived estimates and intermediate constructions}

A corollary of the above Propositions is the following, which will
be needed for the proof of our Main Theorem. The reader should keep
in mind that this Proposition is merely a bookkeeping device that
will allow us to avoid many repetitive calculations in the sequel:\ret

\begin{prop}[Improved Multilinear Estimates]\label{big_tri_prop}
Let $\phi^{(i)}$ be three test functions defined on a time interval $I$ normalized
so that:
\begin{align}
    \lp{\phi^{(1)}}{S[I]} \ &\leqslant \ 1 \ ,
    &\sup_{i=2,3} \llp{\phi^{(i)}}{\W[I]} \
    &\leqslant \ 1 \ ,
     \label{normalization}
\end{align}
Suppose in addition that  $\phi^{(2)}$ has the improved energy
dispersion bound  on $I$:
\begin{equation}
    \sup_k \lp{P_k \phi^{(2)}}{L^\infty[I]}
    \ \leqslant \ \eta \ . \label{addtl_bnd}
\end{equation}
Finally, let $\{c_k\}$ be any $(\delta_0,\delta_0)$-admissible
frequency envelope, and $0\leqslant m$ an additional integer subject
to the condition:
\begin{equation}
        m \ \leqslant \ \sqrt{\delta_1} |\ln(\eta)| \ . \label{trilin_delta_conds}
\end{equation}
Then one has  the following multilinear bounds:\ret

\begin{enumerate}[i)]
    \item (Core Trilinear $L^2$ Estimate) Suppose  along with the above assumptions
    that $\phi^{(3)}$ has unit $\W_c[I]$ norm for the
    frequency envelope $\{c_k\}$. Then for any disposable
    trilinear form $L$ we have the bound:
    \begin{equation}
        \lp{L(\phi^{(1)},\partial^\alpha\phi^{(2)}, \partial_\alpha\phi^{(3)}) }{L^2_t(\dot{H}^{-\frac{1}{2}})_c[I]}
        \lesssim \ \eta^{\delta_1}  \ . \label{core_L2}
    \end{equation}\ret
    \item (Additional Trilinear $L^2$ Estimate) Suppose again that we have the conditions
    \eqref{normalization}--\eqref{addtl_bnd}, and that this time
    $\phi^{(1)}$ has unit $S_c[I]$ norm for the frequency envelope $\{c_k\}$. Then for any disposable trilinear form
    $L$ we have the bound:
    \begin{equation}
        \lp{P_k L(\phi^{(1)},\partial^\alpha\phi^{(2)}, \partial_\alpha\phi^{(3)})}{L_t^2(L^2_x)[I]}
        \lesssim \ 2^\frac{k}{2} \eta^{\delta_1}
        \big(c_k +\lp{P_{<k}\phi^{(1)}}{S[I]} \big)
        \ . \label{ext_L2}
    \end{equation}\ret
    \item (Core Trilinear $N$ Estimate) For a positive integer $m$ and integer $k$
    and disposable trilinear form $L$, define the following trilinear form:
    \begin{multline}
        T^{m}_k(\phi^{(1)},\phi^{(2)},\phi^{(3)}) \ := \
        P_k L(\phi^{(1)},\partial^\alpha\phi^{(2)}, \partial_\alpha\phi^{(3)})\\
        - L(\phi^{(1)}_{<k-m},\partial^\alpha\phi^{(2)}_{<k-m}, \partial_\alpha\phi^{(3)}_k)
        - L(\phi^{(1)}_{<k-m},\partial^\alpha\phi^{(2)}_{k}, \partial_\alpha\phi^{(3)}_{<k-m})
        \ . \label{basic_trilin}
    \end{multline}
    Suppose in addition to the \eqref{normalization}--\eqref{addtl_bnd} we also
    have unit $\W_c[I]$ norm
    of $\phi^{(3)}$, and furthermore the  high modulation
    bounds:
    \begin{align}
        \lp{\phi^{(2)}}{\uX[I]} \ &\leqslant \ \eta \ ,
        &\lp{\phi^{(3)}}{\uX_c[I]} \ &\leqslant \ \eta \ .
         \label{addtl_X_bnds}
    \end{align}
    %where $\{\underline{c}_k\}$ is any other $(\delta_0,\delta_0)$-admissible frequency envelope.
    Then  the following trilinear estimate holds:
    \begin{equation}
        \lp{T^m_k(\phi^{(1)},\phi^{(2)},\phi^{(3)})}{N[I]} \ \lesssim \
        \eta^{\delta_1} c_k
    %+\eta^{-C\delta_1} \underline{c}_k^\delta c_k^{1-\delta}
    \ . \label{core_N}
    \end{equation}\ret

    \item (Additional Trilinear $N$ Estimate) Suppose in addition to
    \eqref{normalization}--\eqref{addtl_bnd} we have unit $S_c[I]$ norm
    of $\phi^{(1)}$, and in addition the  high modulation bounds:
    \begin{equation}
        \sup_{i=2,3} \lp{\phi^{(i)}}{\uX[I]} \ \leqslant \ \eta \
          . \label{addtl_X_bnds2}
    \end{equation}
    Then if $T^m_k$ is defined as on line \eqref{basic_trilin} we have the bound:
    \begin{equation}
        \lp{T^m_k(\phi^{(1)},\phi^{(2)},\phi^{(3)})}{N[I]} \ \lesssim \
        \eta^{\delta_1}\big(c_k + \lp{P_{<k}\phi^{(1)}}{S[I]}\big)  \ . \label{ext_N}
    \end{equation}
\end{enumerate}
\end{prop}\ret

\begin{rem}
If the functions $\phi^{(i)}$ admit a common frequency envelope
$\{c_k\}$ then we can relax the admissibility condition on $\{c_k\}$ and work with
$(\delta_0,\Delta)$ frequency envelopes. Precisely, for any
$(\delta_0,\Delta)$-admissible frequency envelope $\{c_k\}$
we have the following:
\begin{itemize}
  \item If \eqref{normalization} is replaced by
\begin{align}
    \lp{\phi^{(1)}}{S_c[I]} \ &\leqslant \ 1 \ ,
    &\sup_{i=2,3} \llp{\phi^{(i)}}{\W_c[I]} \
    &\leqslant \ 1 \ ,
     \label{normalization-eq}
\end{align}
then \eqref{core_L2} follows.
\item If in addition \eqref{addtl_X_bnds} is replaced by
 \begin{align}
        \lp{\phi^{(2)}}{\uX_c[I]} \ &\leqslant \ \eta \ ,
        &\lp{\phi^{(3)}}{\uX_{{c}}[I]} \ &\leqslant \ \eta \ ,
         \label{addtl_X_bnds-eq}
 \end{align}
then the following version of \eqref{core_N} holds:
  \begin{equation}
        \lp{T^m_k(\phi^{(1)},\phi^{(2)},\phi^{(3)})}{N[I]} \ \lesssim \
        \eta^{\delta_1} c_k \ . \label{core_N-eq}
    \end{equation}
\end{itemize}
\end{rem}\ret

\begin{rem}\label{trilinear_rem}
As will become apparent in the sequel, the \emph{only} use of the
renormalized norms $\W[I]$ and the high modulation bounds $\uX[I]$
in the estimates of Proposition \ref{big_tri_prop} is to ensure the
smallness coming from the parameter $\eta$. Thus, under the simpler
assumption that the $\phi^{(i)}$ are only normalized so that
$\lp{\phi^{(i)}}{S[I]}\leqslant 1$ we have the following:
\begin{itemize}
        \item If $\phi^{(3)}$ has $(\delta_0,\delta_0)$-admissible $S[I]$ norm frequency envelope
        $\{c_k\}$, then estimate \eqref{core_L2} holds with
        $\eta=1$.
%        \item If $\phi^{(1)}$ as $(\delta,\delta)$-admissible
%        $S[I]$ norm frequency envelope
%       $\{c_k\}$, then estimate \eqref{ext_L2} holds with
 %      $\eta=1$.
        \item If $\phi^{(3)}$ has a $(\delta_0,\delta_0)$-admissible
        $S[I]$ norm frequency envelope
        $\{c_k\}$, and if we let $m\geqslant 10$ be any integer, then
        we may replace estimate \eqref{core_N} with the bound:
        \begin{equation}
            \lp{T^m_k(\phi^{(1)},\phi^{(2)},\phi^{(3)})}{N[I]} \ \lesssim \
            2^{4\delta_0 m} c_k  \ . \label{large_core_N}
        \end{equation}
       \item If $\phi^{(1)}$ has $(\delta_0,\delta_0)$-admissible
       $S[I]$ norm frequency envelope
        $\{c_k\}$, and if we let $m\geqslant 0$ be any integer, then
        we may replace estimate \eqref{ext_N} with the bound:
       \begin{equation}
           \lp{T^m_k(\phi^{(1)},\phi^{(2)},\phi^{(3)})}{N[I]} \ \lesssim \
            2^{4\delta_0 m}  \big(c_k +\lp{P_{<k}\phi^{(1)}}{S[I]} \big)  \ . \label{large_ext_N}
       \end{equation}
\end{itemize}
For further details, see Remarks \ref{large_core_L2_rem},
\ref{large_core_N_rem} and \ref{large_ext_N_rem}
in Section \ref{s:der_tri} below.
\end{rem}

Next, we state a result that ties together many of the previous
Propositions. This is a structure theorem for large data wave-maps
with says that in the presence of good $S[I]$ norm bounds one has
some additional regularity properties, as well as a crucial
``fungibility'' property that is central to energy norm
inductions.\ret

\begin{prop}[Structure of Finite $S$ Norm Wave-Maps]\label{p:wm_struct}
Let $\phi$ be a wave-map defined on the interval $I$ with energy $E$ and $S$ norm $F$. Then the following is true:\ret

\begin{itemize}
    \item (Additional Norm Control) We have the  bounds:
    \begin{equation}
        \lp{\phi}{\uX[I]} +
        \lp{\phi}{\underline{E}[I]} \ \lesssim_F \ 1 \ . \label{ext__phi_bnds}
    \end{equation}\ret
 \item
    (Renormalization) If $\{c_k\}$ is a $(\delta_0,\Delta)$-admissible frequency
    envelope for $\lp{\phi}{S[I]}$, then we may renormalize our wave-map as follows:
    \begin{equation}
            \llp{\phi}{\W_c[I]} \ \lesssim_F \
            1 \ , \label{renorm_phi_bnds}
    \end{equation} \ret
    \item (Partial Fungibility) If $\lp{\phi}{S[I]} = F$, then there exists a collection of
    subintervals $I=\cup_{i=1}^K I_i$, such that $K=K(F)$ depends only on $F$, and
    such that the following bound holds on each $I_i$:
    \begin{equation}
        \lp{\phi}{S[I_i]} \ \lesssim_E \ 1 \ . \label{fung_bnd}
    \end{equation}\ret
       \item (Smallness of High Modulations) Suppose in addition that we have
    energy dispersion $\sup_k \lp{P_k \phi}{L^\infty_t(L^\infty_x)[I]} \leqslant \epsilon$. Then
    we also have the  estimate:
    \begin{equation}
        \lp{\phi}{\uX[I]} \ \lesssim_F \ \epsilon^{\delta_1}
        \ . \label{small_mod_est}
    \end{equation}\ret
%    where the extension of $\phi$ is compatible with the extensions in all other
%   estimates on $I$ listed above and below.\\
%   \item (Improved Bilinear Estimates) Let $\phi^{(i)}$ be a collection of WM that satisfy both
%    of the above $S$ and energy dispersion bounds. Then we have the estimate:
%    \begin{equation}
%        \lp{P_k (P_{k_1}\phi^1,P_{k_2}\partial^\alpha
%            \phi^2,P_{k_3}\partial_\alpha \phi^3)}{N[I]}
%        \ \lesssim_F \ 2^{C\sum_j |k-k_j|} \epsilon^\delta \ . \label{w_imp_est}
%    \end{equation}
%    as well as:
%    \begin{equation}
%        \lp{P_k (\partial^\alpha
%            \phi^{(1)}_{k_1}\cdot \partial_\alpha \phi^{(2)}_{k_2})}{N[I]}
%        \ \lesssim_F \ 2^{C\sum_j |k-k_j|} \epsilon^\delta \ . \label{w_imp_est2}
%    \end{equation}\\
    \item (Frequency Envelope Control) Suppose that $\phi$  has
    sufficiently small energy dispersion $\epsilon < \epsilon(F)$. Then if $\{c_k\}$ is a $(\delta_0,\Delta)$-admissible
    $\dot{H}^1\times L^2$ frequency envelope for $\phi[0]$ we have:
    \begin{equation}
            \lp{\phi_k}{S[I]} \ \lesssim_F \ c_k  \ .
            \label{freq_env_control}
    \end{equation}%\\
%    \item (Improved High Frequency  Control) Suppose that $\phi$ again has
%    $\epsilon$ energy dispersion. If we have the improved high frequency bound:
%   \begin{align}
%        \lp{P_k \phi[0] }{\dot{H}^1\times L^2} \ &\leqslant \
%     2^{(k_*-k)} \eta \ ,
%    &\hbox{for}\ k\geqslant k_* \ , \label{high_freq_assumpt}
%   \end{align}
%   then we also have:
%   \begin{equation}
%           \lp{P_k \phi}{S[I]} \ \lesssim_F \ 2^{\delta (k_*-k)}(\eta + \epsilon^\delta)  \ ,
%           \label{non_prod_est}
%    \end{equation}
 %   where $P_k$ is the \textbf{same} multiplier on both lines above.
\end{itemize}
\end{prop}\ret

Finally, for the reader's convenience we group together
the results which enable us to carry out our bootstrapping
arguments:

\begin{prop}[Bootstrapping Tool] \label{propboot}
Let $I=[a,b]$ be an interval and $c$ a $(\delta_0,\Delta)$ frequency
envelope. Then for each affinely Schwartz function $\phi$ in $I$ the following properties hold:\ret

\begin{itemize}
    \item (Seed $S$ bound)  Let $I_n \subset I$ be a decreasing sequence of intervals
    which converges to the point $t=0$. Then:
    \begin{equation}\label{seed_S}
        \lim_{n \to \infty} \lp{\phi}{S[I_n]}  \ \lesssim  \ \lp{\phi[0]}{\dot{H}^1\times L^2} \ ,
        \qquad
        \lim_{n \to \infty} \lp{\phi}{S_c[I_n]}  \
        \lesssim  \ \lp{\phi[0]}{(\dot{H}^1\times L^2)_c} \ .
    \end{equation}\ret

    \item (Continuity Properties)
    For each subinterval $J \subset I$ we have $\phi \in S[J] \cap S_c[J]$, and
     its $S$  norm $\lp{\phi}{S[J]}$, its $S_c$ norm
    $\lp{\phi}{S_c[J]}$, and its energy dispersion  norm
    $\sup_k \lp{P_k \phi}{L^\infty_t(L^\infty_x)[J]}$ all depend continuously
    on the endpoints of $J$. \ret
    %In addition, if $\phi$ is a classical Wave-Map on $I$,
    %then it extends to a classical Wave-Map
    %$\phi \in S[I_1]$ for a larger interval
    %$I_1=[a_1,b_1]$ with $a_1 <a < b < b_1$.

    \item (Closure and Extension Property) Let $I_n$ be an increasing sequence of intervals
    and $\cup I_n = I = (a,b)$. Let $\phi$ be a classical Wave-Map in $I$
    which satisfies the uniform bounds:
    \begin{equation}
        \lp{ \phi}{S[I_n]} \ \leqslant \ F \ ,
        \qquad \sup_k \lp{P_k \phi}{L^\infty_t(L^\infty_x)[I_n]}
        \ \leqslant \ \epsilon \ . \notag
    \end{equation}
    with $\epsilon \leqslant \epsilon(F)$. Then $\phi \in S[I]$,
    and furthermore it can be extended to a  classical Wave-Map in a
    larger interval $I_1=[a_1,b_1]$ with $a_1 <a < b < b_1$.
\end{itemize}
\end{prop}

\begin{proof}%[proof of Proposition \ref{propboot}]
The first part a direct consequence of the solvability
bound \eqref{energy_est} since $\Box \phi \in L^1_t L^2_x[I]$ as well as
$\Box \phi \in (L^1_t L^2_x)_c[I]$. \ret

For the second part we first consider the $S$ norm. Let $J_n \subset I$ be a sequence of intervals converging to $J$. We consider a sequence of rescalings mapping $J$ to $J_n$,
\[
(t,x) \to (\lambda_n t + t_n^0, \lambda_n x), \qquad \lambda_n \to 1, \quad t_n \to 0
\]
This allows us to map functions in $J_n$ to functions in $J$,
\[
 \phi \to \phi_n(t,x) = \phi(\lambda_n t + t_n^0, \lambda_n x)
\]
Hence using the scale invariance of the $S$ norm, we have
\[
 \| \phi\|_{S[J_n]}= \|\phi_n\|_{S[J]} \to \|\phi\|_{S[J]}
\]
where in the last step we simply use the fact that
convergence in the Schwartz space implies the convergence in $S[J]$.

For $S_c$ norms the proof is similar. The dyadic convergence
 $\| \phi_k\|_{S[J_n]} \to \| \phi_k\|_{S[J]}$ follows by the same rescaling
argument. This implies the $S_c$ convergence since the tails are small,
\[
 \lim_{k \to \pm \infty} c_k^{-1} \| \phi_k\|_{S[I]} = 0
\]
which is due to the Schwartz regularity of $\phi$. A similar decay of the
tails yields the continuity of the energy dispersion norm. \ret

For the last part we observe that by \eqref{freq_env_control},
for each $(\delta_0,\Delta)$ frequency envelope $c$
 we obtain a uniform bound for $\lp{P_k \phi}{S_c[I_n]}+\lp{P_k \phi}{\uX_c[I_n]}$.
Letting $n \to \infty$ we directly obtain
$P_k \phi \in {\uX_c[I]}$, which shows that for each $k$ we have
$P_k \phi \in S[I]$ and $\lp{ P_k \phi}{S[I_n]} \to \lp{ P_k \phi}{S[I]} \lesssim c_k$. Hence $\phi$ is a Schwartz wave map in $[a,b]$,
therefore by the local well-posedness result it admits a   Schwartz
extension to a larger interval.
\end{proof}

\ret
%-------------------------------------------------------------------------
%%%%%%%%%%%%%%%%%%%%%%%%%%%%%%%%%%%%%%%%%%%%%%%%%%%%%%%%%%%%%%%%%%%%%%%%%%
%-------------------------------------------------------------------------

\ret
%-----------------------------------------------------------------------------------
%%%%%%%%%%%%%%%%%%%%%%%%%%%%%%%%%%%%%%%%%%%%%%%%%%%%%%%%%%%%%%%%%%%%%%%%%%%%%%%%%%%%
%%%%%%%%%%%%%%%%%%%%%%%%%%%%%%%%%%%%%%%%%%%%%%%%%%%%%%%%%%%%%%%%%%%%%%%%%%%%%%%%%%%%
%-----------------------------------------------------------------------------------

\section{Proof of the Main Result}\label{s:main_proof}

The purpose of this Section is to use the setup of the previous two
Sections to prove the following result, which easily implies our
main Theorem \ref{main_thm} as well as Corollary~\ref{maint_nonpert}.\ret

\begin{thm}[Frequency Envelope Version of the Main Theorem]\label{main_thm_FE}
\ \  There exist two functions $1\ll F(E)$ and $0 < \epsilon(E)\ll 1$ of
the energy \eqref{energy} such that if $\phi$ is a finite
energy solution to \eqref{main_eq} in a closed interval $I\times\R^2$, where
$I=[a,b]$, with energy $E$ and dispersion \eqref{energy_disp}, then
estimate \eqref{S_est} holds in $S[I]$. In addition, there exists a
universal polynomial $K(F)$ such that if $\{c_k\}$ is any
$(\delta_0,\Delta)$-admissible frequency envelope for $\phi[0]$,
we have the bound:
\begin{equation}
    \lp{\phi}{S_c[I]} \ \leqslant \ K(F) \ . \label{Sc_est}
\end{equation}
In particular, one may extend $\phi$ to a finite energy Wave-Map on open neighborhood $I\subseteq (a-i_0,b+i_0)$
whose additional length $i_0$ depends \textbf{only} on $E$, $\{c_k\}$, and $\epsilon$.
\end{thm}\ret

We immediately observe that it suffices to prove the result for classical
wave-maps. This is due to the small data result in Theorem~\ref{Tat_Thm},
which implies that any finite energy wave map in a closed interval
can be approximated in $S$ by classical wave maps. In addition, the $S$
convergence easily implies the convergence of the energy dispersion norm
 \eqref{energy_disp}.

In the sequel we simply focus on proving \eqref{S_est}. The estimate
\eqref{Sc_est} is an immediate consequence of
\eqref{freq_env_control}. In fact, it would be tempting to use the
more direct analysis employed in the proof of \eqref{freq_env_control}
to establish \eqref{S_est} as well in a single go. Such a strategy
seems to fail basically due to linearized $Low\times High \Rightarrow
High$ frequency interactions. These interactions need to be handled
via Proposition \ref{p:para}, which in turn requires one to already
control $S$ type norms (e.g. in assumption \eqref{A_tdphi_est}). To
avoid this dilemma, we employ a simple induction scheme to reduce
things to estimates for Wave-Maps of (slightly) smaller energy. The
reader should keep in mind however that modulo this single $Low\times
High$ obstruction, our analysis would work to prove \eqref{S_est} and
\eqref{freq_env_control} simultaneously. More specifically, the
remaining estimates basically boil down to using
\eqref{ednullw}--\eqref{bal_l2_N} to eliminate matched frequency
``semilinear'' type interactions (this is the \emph{only} place where
energy the dispersion \eqref{energy_disp} really comes in), and
\eqref{standard_est_bi}--\eqref{standard_est_tri} to kill off $High
\times High \Rightarrow Low$ frequency cascades.\ret

We now construct the functions $F(E)$ and $\epsilon(E)$ such that
\eqref{energy_disp} and \eqref{S_est} hold.  Precisely, we will show
that there exists a strictly positive nonincreasing function defined
for \emph{all} values of $E$, $c_0 = c_0(E) \ll 1$, so that if the
conclusion of the Theorem holds up to energy $E$ then it also holds
up to energy $E+c_0$. It is important here that $c_0$ depends
\emph{only} on $E$ and not on the size of $F(E)$ or $\epsilon(E)$,
as otherwise we would only be able to conclude the usual first step
in an induction on energy proof which is establishing that the set
of regular energies is open.\footnote{In this latter setup, one is
then left with the arduous task of eliminating minimal energy blowup
solutions. Our strategy is a bit more direct because we accomplish
this as well in our construction of $c_0$, so we are able to avoid a
good deal of repetitive analysis.} Also, we note here the
monotonicity of $c_0$ is only used to conclude that $c_0$ admits a
positive lower bound on any
compact set.\ret

According to Theorem \ref{Tat_Thm} we know that $\epsilon(E)$ and
$F(E)$ can be constructed up to some $E_0\ll 1$. We now assume that
$E_0$ is fixed by induction, and to increase its range we consider a
solution $\phi$ defined on an interval $I$ with energy $E[\phi]= E_0 +
c$, $c \leqslant c_0(E_0)$ and with energy dispersion $\leqslant \epsilon$ (at
first this is a free parameter which may take as small as we like).
We will compare $\phi$ with a solution $\tphi$ with energy $E_0$.  To
construct $\tphi$ we reduce the initial data energy of $\phi[0]$ by
truncation in frequency.  We define the ``cut frequency'' $k_* \in \R$
according to (this can be done by adjusting the definition of the
$P_{<k}$ continuously if necessary):
\begin{equation}
    E[\Pi P_{\leqslant  k_*}\phi[0]] \ = \ E_0 \ . \notag
\end{equation}
We consider the Wave-Map $\td\phi$ with this initial data
$\td\phi[0]=\Pi P_{\leqslant k_*}\phi[0]$.  This Wave-Map exists
classically for at least a short amount of time according to Cauchy
stability, and where it exists classically we have:
\begin{equation}
    E[\td{\phi}(t)] \ = \ E_0 \ . \label{td_phi_eng}
\end{equation}
Since $\phi$ has energy dispersion $\leqslant \epsilon$, by \eqref{cutphi}
it follows that $\tphi$ has energy dispersion $ \lesssim_{E_0}
\epsilon^\frac12$ at time $t=0$. Again by the usual Cauchy stability
theory, if $\epsilon$ is chosen small enough in comparison to the
inductive defined parameter $\epsilon(E_0)$ it follows that there
exists a non-empty interval $J$ where $\tphi$ satisfies:
\begin{equation}
    \sup_k \lp{P_k \tphi}{L^\infty_t(L^\infty_x)[J_0]} \ \leqslant \  \epsilon(E_0) \ . \label{en_disp_t}
\end{equation}
Then our induction hypothesis guarantees that we have the dispersive bounds:
\begin{equation}
    \lp{\td{\phi}}{S[J_0]}  \ \leqslant  \
    F(E_0) \ . \label{td_phi_ind_bnd}
\end{equation}
The plan is now very simple. On one hand, we try to pass the
space-time control of $\td{\phi}$ up to $\phi$ via linearization
around $\td{\phi}$ to control the low frequencies, and conservation of
energy and perturbation theory to control the high frequencies. On the
other hand, we need to pass the good energy dispersion bounds from
$\phi$ back down to $\tphi$ in order to increase the size of
$J\subseteq I$ on which \eqref{en_disp_t} holds until it eventually
fills up all of $I$.  To achieve all of this, we proceed via two core
estimates:\ret

\begin{prop}[Evolution of Low Frequency Errors]\label{p:phi_est_low}
Let $\phi$ be a Wave-Map defined on an interval $J$ with energy $E+c$
with $ 0 < c \lesssim 1$ and bounds:
\begin{align}
    \sup_k \lp{P_k \phi}{L^\infty_t(L^\infty_x)[J]} \ &\leqslant \ \epsilon \ ,
    &\lp{\phi}{S[J]} \ &\leqslant \ F \ . \label{phi_low_freq_assumt}
\end{align}
Suppose in addition that $\td{\phi}$ is the Wave-Map with energy $E$ defined by
$\tphi[0] = \Pi P_{\leqslant k_*}\phi[0]$, and that $\td{\phi}$ is classical on $J$
with bounds:
\begin{align}
    \sup_k \lp{P_k \td\phi}{L^\infty_t(L^\infty_x)[J]} \ &\leqslant \ \td{\epsilon} \ ,
    &\lp{\td\phi}{S[J]} \ &\leqslant \ \td{F}  \ , \label{tdphi_low_freq_assumt}
\end{align}
Assume also that the two energy dispersion constants are chosen so that:
\begin{align}
    \epsilon \ &\leqslant \ \td{\epsilon} \ ,
    &\epsilon\ &\leqslant \ (C F)^{-\delta^{-10}_0} \ ,
    &\td\epsilon\ &\leqslant \ (C \td{F})^{-\delta^{-10}_0} \ , \label{epsilon_bnds}
\end{align}
where we may assume that  $F\geqslant \td{F}\geqslant E\geqslant  C^{-\frac{1}{2}}$ and $C$ is
a sufficiently large constant. Then in addition we have the bound:
\begin{equation}
    \lp{\tphi - P_{\leqslant k_*} \phi}{S[J]}
    \ \lesssim_F \ \epsilon^{\delta_0} \ . \label{low_zip_est}
\end{equation}
\end{prop}\ret

\begin{prop}[High Frequency Evolution Estimates]\label{p:phi_est_high}
Let $\phi$ and $\td\phi$ be defined as in the last Proposition, in
particular with the bounds \eqref{phi_low_freq_assumt} and
\eqref{tdphi_low_freq_assumt} (resp), and that the dispersion
constants obey \eqref{epsilon_bnds}. Then there exists a universal
function $c_0(E)$ with $c_0^{-1}\lesssim_E 1$ such that if we
assume $c_0= c_0(E)$ in the
definition of $\tphi$ we in addition have the bound:
\begin{equation}
    \lp{\td\phi -  \phi}{S[J]} \ \lesssim_{\td{F}} \ 1  \ ,
    \label{high_est}
\end{equation}
\end{prop}\ret

We postpone the proof of the above Propositions
to show how to use them to conclude our induction.
By the seed bound \eqref{seed_S} we may assume that in addition to
\eqref{en_disp_t} and \eqref{td_phi_ind_bnd} above we also have:
\begin{equation}
    \lp{{\phi}}{S[J]}  \ \leqslant  \
    2F(E_0+c) \ , \notag
\end{equation}
on some interval $J$. With this setup, and by an application of the continuity property in Proposition~\ref{propboot} it suffices
to combine Propositions \ref{p:phi_est_low}--\ref{p:phi_est_high} to
show the following:\ret

\begin{cor} \label{cor:maint}
  Assume there exists functions $\epsilon(E)$ and $F(E)$ defined up to
  $E_0$ such that \eqref{energy_disp} implies \eqref{S_est}. Choose
  $c_0(E_0)$  according to Proposition \ref{p:phi_est_high}.
Then  there exists extensions of $\epsilon(E)$ and $F(E)$ such that for each
 $0  < c \leqslant c_0(E_0)$ and each classical  Wave-Map $\phi$  in a time
  interval $J$ with energy $E_0+c$ and the bounds:
\begin{align}
    \sup_k \lp{P_k\phi}{L^\infty_t(L^\infty_x)[J]} \ &\leqslant \ \epsilon(E_0+c) \ ,
    &\lp{\phi}{S[J]} \ &\leqslant \ 2F(E_0+c) \ , \label{cor_ass_bnds}
\end{align}
 we  have:
\begin{equation}
    \lp{\phi}{S[J]} \ \leqslant \ F(E_0+c) \ . \label{cor_get_bnds}
\end{equation}
\end{cor}\ret

\begin{proof}
In addition to \eqref{cor_ass_bnds} $\Rightarrow$
  \eqref{cor_get_bnds}, we will make the additional assumption
that $\tphi$ is defined as a Schwartz wave map in $J$ and satisfies:
\begin{equation}
    \sup_k \lp{P_k \tphi}{L^\infty_t(L^\infty_x)[J]} \ \leqslant \  \epsilon(E_0) \ , \label{en_disp_t_imo}
\end{equation}
and  show that if the extensions to
  $\epsilon(E)$ and $F(E)$ are chosen correctly then we in addition have
  the following improvement to \eqref{en_disp_t_imo}:
\begin{equation}
    \sup_k \lp{P_k \tphi}{L^\infty_t(L^\infty_x)[J]} \ \leqslant \  \frac{1}{2}\epsilon(E_0) \ . \label{en_disp_t_imp}
\end{equation}
To see that this is sufficient, we first note that by \eqref{en_disp_t}
the bound \eqref{en_disp_t_imo} holds in a smaller interval $J_0
\subset J$. Extending $J_0$ to a maximal interval
in $J$, denoted still $J_0$, so that
\eqref{en_disp_t_imo} holds, by the closure property in
Proposition~\ref{propboot} it follows that $J_0$ must be closed.
The same part of Proposition~\ref{propboot} shows that
$\tphi$ has a Schwartz extension to a neighborhood of $J_0$.
Then by  \eqref{en_disp_t_imp} applied in $J_0$ and the continuity
property in Proposition~\ref{propboot} it follows that \eqref{en_disp_t_imo}
holds in a larger interval. Hence $J_0$ must be both closed and open in $J$, and
therefore $J_0=J$.\ret

It remains to find extensions $\epsilon(E)$ and $F(E)$ so that
\eqref{cor_ass_bnds} together with \eqref{en_disp_t_imo} imply
\eqref{cor_get_bnds} and \eqref{en_disp_t_imp}.
Our extensions of $\epsilon(E)$ and $F(E)$ in $(E_0,E_0+c_0]$ are constant:
\begin{equation}
    \epsilon(E)  \ =  \ \epsilon \ ,
    \qquad F(E)  \ =  \ F \ ,
    \qquad E \in (E_0,E_0+c_0] \ . \notag
\end{equation}
Let $K_1({F})$ and $K_2(\td{F})$ be the implicit polynomials from
lines \eqref{low_zip_est} and \eqref{high_est} (resp). In order to get the improvement
\eqref{en_disp_t_imp} we need that:
\begin{equation}
    \epsilon^{\delta_0} \cdot K_1({F}) \ \ll \ \epsilon(E_0) \ . \notag
\end{equation}
In order to conclude \eqref{cor_get_bnds} we need that:
\begin{equation}
     K_2(F(E_0)) \
    \ll \ F \ . \notag
\end{equation}
Finally, we also need to choose $\epsilon$ and $F$ so that
\eqref{epsilon_bnds} holds, and so that (which is of course
redundant):
\begin{equation}
    E_0^\frac{1}{2} \epsilon^\frac{1}{2} \ \ll \ \epsilon(E_0) \ , \notag
\end{equation}
which was used right before line \eqref{en_disp_t} to get things
started. All of these goals can easily be satisfied as long as we
choose $\sigma=\delta_0^{20}$, with $\delta_0\ll 1$ sufficiently small,
and then first choose $F$, followed by $\epsilon$, such that:
\begin{align}
    F(E_0) \ &\ll \ F^{\sigma} \ ,
    &\epsilon^\sigma \ &\ll \ \min\{\epsilon(E_0),{F}^{-1}\} \ . \notag
\end{align}
Notice that this process can be carried on indefinitely, regardless of
the size of $E$, because we have taken care to decouple the step size
$c_0$ from the growth and decay properties of $F$ and $\epsilon$.

We remark that the above proof allows us to estimate the size of $E(F)$
and $\epsilon(F)$. Indeed, what we have obtained are
piecewise constant functions $c_0(E)$, $\epsilon(E)$, and $F(E)$ which at the
jump points $E_n$ are given by the recurrence relation:
\begin{equation}
    E_{n+1}  \ =  \ E_n + c_0(E_n) \ , \notag
\end{equation}
and which satisfy:
\begin{equation}
    c_0(E_n)  \ =  \ c E_n^{-\sigma^{-1}} \ ,  \quad
    F(E_{n+1})  \ = \ C F(E_n)^{\sigma^{-1}} \ ,
    \quad \epsilon(E_{n+1})  \ = \
    c F(E_n)^{\sigma^{-2}} \ , \notag
\end{equation}
with sufficiently small $\sigma,c$ and sufficiently large $C$.
The first relation shows that:
\begin{equation}
    E_n  \ \approx  \ n^\frac{1}{\sigma^{-1}+1} \ , \notag
\end{equation}
while the next two give relations of the form:
\begin{equation}
    F(E_n)  \ \leqslant   \ C_1^{\sigma^{-n}},\qquad  \epsilon(E_n) \geqslant
    c_1^{\sigma^{-n}} \ . \notag
\end{equation}
Together the last two bounds yield estimates for $F$ and $\epsilon$ of the form:
\begin{equation}
    F(E)  \ \leqslant  \ e^{C e^{E^M}} \ ,
    \qquad \epsilon(E)  \ \geqslant  \ e^{-C e^{E^M}} \ , \notag
\end{equation}
again with $C$ and $M$ sufficiently large.

\end{proof}\ret

The remainder of this section is devoted to the proof of Propositions
\ref{p:phi_est_low}--\ref{p:phi_est_high}. This will be done in order because we will use
some of the estimates of Proposition \ref{p:phi_est_low}
in our demonstration of Proposition \ref{p:phi_est_high}.\ret

%----------------------------------------------------------------------------

\begin{proof}[Proof of Proposition \ref{p:phi_est_low}]
Denoting:
\begin{equation}
    \psi \ = \ P_{\leqslant k_*} \phi - \tphi \ , \label{psi_def}
\end{equation}
we will prove the stronger bound:
\begin{equation}
        \lp{ \psi}{S_c[J]}
        \ \leqslant 1 \ , \label{l1_low_zip_est}
\end{equation}
where $\{c_k\}$ is the $(\delta_0,\delta_0)$-admissible frequency envelope
$c_k=2^{-\delta_0|k - k_* |}\epsilon^{\delta_0}$.
We first consider the initial data for $\psi$. By an immediate
application of Proposition \ref{p:cuta} and the energy dispersion
bound \eqref{p:phi_est_low} we have:
\begin{equation}
    \lp{P_k \psi[0]}{\dot H^1 \times L^2}
    \ \lesssim_E \ \epsilon^\frac14   2^{- \frac{1}{2}|k-k_*|} \ . \label{psidata}
\end{equation}
Since $\psi$ is a Schwartz function, this implies that
for a small interval $I \subset J$ containing $t=0$ we have:
\begin{equation}
        \lp{\psi}{S_c[I]}
        \ \leqslant 1 \ , \label{l1_low_zip_estseed}
\end{equation}
Using this as a seed bound, by the continuity property in Proposition~\ref{propboot} it suffices to prove
that \eqref{l1_low_zip_est} holds under
the bootstrap assumption:
\begin{equation}
        \lp{\psi}{S_c[I]}
        \ \leqslant \  2 \ . \label{l1_low_zip_est-boot}
\end{equation}
As a preliminary step we use the general
renormalization bound \eqref{renorm_phi_bnds} as well as the high
modulation bound \eqref{small_mod_est}, which in light of the
estimates on each of lines \eqref{phi_low_freq_assumt} and
\eqref{tdphi_low_freq_assumt} imply the set of inequalities:
\begin{align}
    \llp{\phi}{\W[J]} \ &\lesssim_F \ 1 \ ,
    &\llp{\td{\phi}}{\W[J]} \ &\lesssim_{\td{F}} \ 1
    \ , \label{phi_tdphi_renorm_bnds}\\
    \lp{\phi}{\uX[J]} \ &\lesssim_F \ \epsilon^{\delta_1} \ ,
    &\lp{\td{\phi}}{\uX[J]} \ &\lesssim_{\td{F}} \ \td{\epsilon}^{\delta_1}
    \ . \label{phi_tdphi_mod_bnds}
\end{align}
The proof is deduced in a series of steps:\ret

%----------------------------------------------------------------------

\step{1}{Outline of the proof}
The equation for $\psi$ has the form:
\begin{equation}
    \Box \psi \ = \ - P_{\leqslant k_*} \left( \mathcal{S}(\phi)\partial^\alpha \phi \partial_\alpha \phi\right)
    + \mathcal{S}(\tphi)\partial^\alpha \tphi \partial_\alpha \tphi \ . \label{psi_eq}
\end{equation}
This may be rewritten as follows:
\begin{equation}
    \Box \psi \ = \ -   \mathcal{D}(\tphi,\psi) + \mathcal{C}(\phi)
    \ , \label{psi_eq_comm}
\end{equation}
where the difference $\mathcal{D}$ and the generalized commutator $\mathcal{C}$ are defined as follows:
\begin{align}
    \mathcal{D}(\tphi,\psi) \ &= \ \mathcal{S}(\tphi+\psi)\partial^\alpha( \tphi+\psi) \partial_\alpha( \tphi+\psi)
        - \mathcal{S}(\tphi)\partial^\alpha \tphi \partial_\alpha \tphi \ , \label{general_phi_diff_def}\\
    \mathcal{C}(\phi) \ &= \
     \mathcal{S}(\phi_{\leqslant k_*})\partial^\alpha
    \phi_{\leqslant k_*} \partial_\alpha \phi_{\leqslant k_*}
    - P_{\leqslant k_*} \left( \mathcal{S}(\phi)\partial^\alpha \phi
    \partial_\alpha \phi\right) \ . \label{general_phi_comm_def}
\end{align}
This form of the equation will be used for proving pure $L^2$ estimates.

Alternatively, freezing the spatial frequency $k$ and introducing
a frequency gap parameter $m\geqslant 20$, we will write \eqref{psi_eq}
in the following paradifferential form:
\begin{equation}
    \Box\psi_k + 2\td{A}^\alpha_{<k-m}\partial_\alpha\psi_k
    \ = \ \mathcal{D}^m_{k}(\tphi,\psi) + \mathcal{L}^m_{k}(\tphi,\psi) + \mathcal{C}^m_{k}(\phi)
    \ , \label{more_crap}
\end{equation}
which will be useful for establishing $N$ estimates.
Here we are writing:
\begin{equation}
    \td{A}^\alpha_{<k-{{m}}} \ =\ {A}^\alpha_{<k-{{m}}}(\tphi) \ := \
    \big(\mathcal{S}(\tphi)_{<k-{{m}}} - \mathcal{S}^\dagger(\tphi )_{<k-{{m}}}
     \big)\partial^\alpha\tphi_{<k-{{m}}}
     \ . \label{A_phi_connection}
\end{equation}
These terms are chosen roughly as follows. The term $\mathcal{D}^m_{k}$ denotes differences of the
form \eqref{general_phi_diff_def} between $\tphi$ and $\psi$ which are frequency localized according to the
general $T^m$ structure defined on line \eqref{basic_trilin}.
In particular, these never contain $Low\times Low\times High$
or $Low\times High\times Low$ interactions.
The term $\mathcal{L}^m_{k}$ contains certain $Low\times Low\times High$ and $Low\times High\times Low$
interactions in $\tphi$ and $\psi$ differences, with the additional structure that $\psi$ is always at
$Low$ frequency with a (possibly large) $m$ dependent gap. Finally,
the expression $\mathcal{C}^m_{k}(\phi)$ contains $\phi$ dependent commutators of the form \eqref{general_phi_comm_def}.

With this setup, we prove the following estimates. First, we show that
the commutators are always favorable, regardless of $m$:
\begin{equation}
    \lp{(\mathcal{C},\mathcal{C}^m)}{(L^2_t(\dot{H}^{-\frac{1}{2}})\cap N)_c[I]} \
    \lesssim_F \ \epsilon^{\frac{1}{2} \delta_1^2} \ . \label{general_comm_est}
\end{equation}
Second, under the bootstrapping assumption
\eqref{l1_low_zip_est-boot}, we will show the first two terms on the
RHS of \eqref{more_crap} may be estimated as follows:
\begin{align}
    \lp{\mathcal{D}^m_{k}(\tphi,\psi)}{N_c[I]} \ &\lesssim_{\td{F}} \ 2^{4\delta_0m} \ , \label{bad_N_bound} \\
    \lp{\mathcal{L}^m_{k}(\tphi,\psi)}{N_c[I]} \ &\lesssim_{\td{F}}  \ \td{\epsilon}^{\delta_0} + 2^{-\delta_0 m}
    \ . \label{L_bound}
\end{align}
While the second of these last two estimates is favorable for closing
a bootstrap via Proposition \ref{p:para}, the first is not. However,
via Remark \ref{direct_para_rem} the above estimates with $m=20$ allow
us to gain renormalization control of $\psi$, namely:
\begin{equation}
    \lp{\psi}{\mathcal{W}_c[I]} \ \lesssim_{\td{F}} \ 1 \ . \label{psi_renorm_control}
\end{equation}
To close the bootstrap, we now use two additional estimates. The first
shows that with \eqref{phi_tdphi_renorm_bnds} and
\eqref{psi_renorm_control}, we have improved $L^2$ control:
\begin{equation}
    \lp{\mathcal{D}(\tphi,\psi)}{L^2_t(\dot{H}^{-\frac{1}{2}})_c[I]} \ \lesssim_{\td{F}}
    \ \td{\epsilon}^{\delta_1} \ . \label{psi_L2_bound}
\end{equation}
In particular by this, \eqref{general_comm_est}, and the gap condition
\eqref{epsilon_bnds} we have:
\begin{equation}
    \lp{\psi}{\uX_c[I]} \ \lesssim_{\td{F}}
    \ \td{\epsilon}^{\delta_1} \ . \label{psi_X_bound}
\end{equation}
Finally, we show that this last estimate, \eqref{psi_renorm_control},
and \eqref{phi_tdphi_renorm_bnds}--\eqref{phi_tdphi_mod_bnds} allow
the following drastic improvement to \eqref{bad_N_bound}:
\begin{equation}
    \lp{\mathcal{D}^m_{k}(\tphi,\psi)}{N_c[I]} \
    \lesssim_{\td{F}} \ \td{\epsilon}^{\delta_1^2}
    \ . \label{good_N_bound}
\end{equation}
The bootstrap is therefore concluded by choosing
$m=\delta_1|\ln(\td{\epsilon})|$ in estimate \eqref{L_bound}, and
applying the linear bound \eqref{linearized_est}
for the paradifferential flow, with the estimates
\eqref{good_N_bound}, \eqref{general_comm_est} for the right hand side and
\eqref{psidata} for the initial data.\ret

%------------------------------------------------------------------------------------
\step{2}{The algebraic decomposition}
Here we derive
the form of the RHS of \eqref{more_crap}.
To uncover this,  we shall employ the following generic notation.
We let $T$ be a trilinear expression of the form:
\begin{equation}
    T\big(\mathcal{S}(\phi^{(1)}),\phi^{(2)},\phi^{(3)}\big)
    \ = \ L\big(\mathcal{S}(\phi^{(1)}),
    \partial^\alpha\phi^{(2)},\partial_\alpha\phi^{(3)}\big) \ ,
    \label{T_form_def}
\end{equation}
with $L$ disposable, and $\mathcal{S}$ is a smooth function with
uniformly bounded derivatives.  From this we may define the $T$
dependent expressions $\mathcal{D}$ and $\mathcal{C}$ as on lines
\eqref{general_phi_diff_def}--\eqref{general_phi_comm_def}.

The frequency localized
equation for $\psi$ is:
\begin{equation}
     \Box \psi_k \ = \ - P_k P_{\leqslant k_*} \big( \mathcal{S}(\phi)\partial^\alpha \phi \partial_\alpha \phi\big)
        + P_k\big(\mathcal{S}(\tphi)\partial^\alpha \tphi \partial_\alpha \tphi\big) \ , \label{psik_equation}
\end{equation}
which may be written in the form:
\begin{equation}
     \Box \psi_k \! = \!
     2 \mathcal{S}(\tphi)_{<k-{{m}}}\partial^\alpha \tphi_{<k-{{m}}}\partial_\alpha \tphi_k
     -
     2 P_{\leqslant k_*} \big[\mathcal{S}(\phi)_{<k-{{m}}}
     \partial^\alpha \phi_{<k-{{m}}}\partial_\alpha \phi_k \big] + T_{1;k}^m
     \ , \label{psik_bulk_decomp}
\end{equation}
where we are writing:
\begin{equation}
    T_{1;k}^m \ = \ T^{{m}}_k\big( \mathcal{S}(\tphi), \tphi, \ \tphi \big)
     - P_{\leqslant k_*} T^{{m}}_k\big(
     \mathcal{S}(\phi), \phi, \ \phi \big)  \ = \ \mathcal{D}^m_{1;k} +
     \mathcal{C}^m_{1;k} \ . \label{T1_line}
\end{equation}
with $T^m$ defined as on line \eqref{basic_trilin}. We now employ the geometric
identity for the second
fundamental form:
\begin{equation}
    \sum_c \mathcal{S}_{ab}^c(\phi)\nabla_{t,x}\phi^c \ \equiv \ 0 \ , \label{sff_iden}
\end{equation}
which follows simply because the constraint on the image of $\phi$ to
lie in $\mathcal{M}$ implies that $\nabla_{t,x}\phi$ lies in
$T_\phi\mathcal{M}$. This is valid for $\tphi$ as well, because it is
an exact wave-map. Therefore, we have the zero expression:
\begin{equation}
     P_k\left( 2 \mathcal{S}(\tphi)^\dagger\partial^\alpha \tphi\partial_\alpha \tphi_{<k-m+2}
        - 2 P_{\leqslant k_*} \big[\mathcal{S}(\phi)^\dagger
        \partial^\alpha \phi\partial_\alpha \phi_{<k-m+2} \big]\right) \ = \ 0 \ , \label{zero_exp}
\end{equation}
which if added to the first two terms on RHS of \eqref{psik_bulk_decomp}
produces:
\begin{equation}
     \hbox{(First two  R.H.S.)}\eqref{psik_bulk_decomp}
     = \ 2\td{A}^\alpha_{<k-{{m}}} \partial_\alpha\tphi_k -
     2P_{\leqslant k_*} \big[ A^\alpha_{<k-{{m}}} \partial_\alpha\phi_k \big]
     + T_{2;k}^m  \ , \label{A_diff}
\end{equation}
where both $\td{A}^\alpha$ and  ${A}^\alpha_{<k-{{m}}}
={A}^\alpha_{<k-{{m}}}(\phi)$ are defined as on line \eqref{A_phi_connection}.
Here the trilinear form $T_{2;k}^m$ is a difference:
\begin{equation}
    T_{2;k}^m \ =\  T_{2;k}^m(\tphi) - P_{\leqslant k_*}T_{2;k}^m(\phi)
    \ = \ \mathcal{D}^m_{2;k} + \mathcal{C}^m_{2;k}
    \ , \label{T2_line}
\end{equation}
where each individual form is
defined as a $T^m$ from line
\eqref{basic_trilin} applied separately to the two trilinear expressions on the
LHS of \eqref{zero_exp}.

We now assign the generalized  difference labels on the RHS of \eqref{more_crap}
by setting $\mathcal{D}^m_k=\sum_i \mathcal{D}^m_{i;k}$, where the two summands
were defined on lines \eqref{T1_line} and \eqref{T2_line}.

To assign $\mathcal{C}^m_k$, we further denote by
$\mathcal{C}^m_{3;k}$ the corresponding expression which results from
commuting the $P_{\leqslant k_*}$ in the second term on the RHS of line \eqref{A_diff}.
We then set $\mathcal{C}^m_k=\sum_i \mathcal{C}^m_{i;k}$.

With these choices, the equation
\eqref{psik_equation} may be written in the form:
\begin{equation}
    \Box\psi_k \ = \ 2{A}^\alpha_{<k-{{m}}}(\tphi) \partial_\alpha\tphi_k-
    2{A}^\alpha_{<k-{{m}}}(\tphi+\psi) \partial_\alpha(\tphi_k + \psi_k) +
    \mathcal{D}_{k}^m + \mathcal{C}_{k}^m \ . \notag
\end{equation}
As a final step we assign:
\begin{equation}
     \mathcal{L}_{k}^m \ = \  2\big( {A}^\alpha_{<k-{{m}}}(\tphi)
     - 2{A}^\alpha_{<k-{{m}}}(\tphi+\psi) \big)\partial_\alpha(\tphi_k+\psi_k) \ . \label{Lpsi_def}
\end{equation}
and the form of \eqref{more_crap} is achieved.\ret

%--------------------------------------------------------------------

The remainder of the proof shows estimates \eqref{general_comm_est}, \eqref{bad_N_bound}, \eqref{L_bound},
\eqref{psi_L2_bound}, and \eqref{good_N_bound}.\ret

%--------------------------------------------------------------------
\step{3}{Estimates for commutators} Here we demonstrate \eqref{general_comm_est}.
Let $\mathcal{C}$ be any expression of the form:
\begin{equation}
     \mathcal{C} \ = \ T\big(\mathcal{S}(P_{\leqslant k_*}\phi),P_{\leqslant k_*}\phi,P_{\leqslant k_*}\phi\big)
    - P_{\leqslant k_*} T\big(\mathcal{S}(\phi),\phi,\phi\big) \ . \label{G_def}
\end{equation}
We will prove the general pair of bounds:
\begin{align}
        \lp{\mathcal{C}}{L^2_t(\dot{H}^{-\frac{1}{2}})_c[I]} \ &\lesssim_F \
        \epsilon^{\frac{1}{2}\delta_1} \ ,
        & \lp{\mathcal{C}}{N_c[I]}  \ &\lesssim_F \ \epsilon^{\frac{1}{2}\delta_1^2} \ .
        \label{basic_G_bound}
\end{align}
As a preliminary step we decompose $\mathcal{C}=\mathcal{C}_1+\mathcal{C}_2$ where:
\begin{align}
    \mathcal{C}_1 \ &= \ T\big(\mathcal{S}(\phi_{\leqslant k_*})- \mathcal{S}(\phi)_{\leqslant k_*}
    ,\phi_{\leqslant k_*},\phi_{\leqslant k_*}\big) \ , \notag\\
    \mathcal{C}_2 \ &= \
     T\big(\mathcal{S}(\phi)_{\leqslant k_*},
    \phi_{\leqslant k_*},\phi_{\leqslant k_*}\big)
    - P_{\leqslant k_*} T\big(\mathcal{S}(\phi),\phi,\phi\big) \ . \notag
\end{align}
These terms are handled separately:\ret

%----------------------------------------------------------------------------
\step{3A}{Estimates for $\mathcal{C}_1$} This is based on the Moser type
estimate:
\begin{equation}
    \lp{\mathcal{S}(\phi_{\leqslant k_*})- \mathcal{S}(\phi)_{\leqslant k_*}}{S_k[I]}
    \ \lesssim_F \ 2^{-\delta |k-k_*|} \ . \label{Sdiff_FE}
\end{equation}
To prove this, we further decompose the difference as:
\begin{equation}
    \mathcal{S}(\phi_{\leqslant k_*})- \mathcal{S}(\phi)_{\leqslant k_*}
    \ = \ P_{> k_*}\mathcal{S}(\phi_{\leqslant k_*}) + P_{\leqslant k_*}\big(
    \mathcal{S}'(\phi_{\leqslant k_*},\phi_{>k_*})\phi_{>k_*}
    \big)  \ , \notag
\end{equation}
where here $\mathcal{S}'$ a bounded and smooth function or its arguments which results
from the difference
$\mathcal{S}(\phi_{\leqslant k_*} + \phi_{>k_*}) -
\mathcal{S}(\phi_{\leqslant k_*})$. The bound \eqref{Sdiff_FE} now
follows by directly applying the Moser estimate
\eqref{basic_moser_c} to the first term on the last line above, and
by applying a combination of the product estimate \eqref{S_prod3}
and the Moser estimate \eqref{basic_moser} to the second.

To conclude the proof of \eqref{basic_G_bound} for the term $\mathcal{C}_1$ we
need to split the output frequency into two cases: $k\leqslant
k_*+10$ or $k > k_*+10$. In the first case, we directly use
\eqref{ext_L2} and \eqref{ext_N}, which together provide
\eqref{basic_G_bound} in light of
\eqref{phi_tdphi_renorm_bnds}--\eqref{phi_tdphi_mod_bnds} and the
additional $L^\infty$ estimate:
\begin{align}
    \lp{P_{<k}\big[\mathcal{S}(\phi_{\leqslant k_*})- \mathcal{S}(\phi)_{\leqslant k_*}\big]}
    {L^\infty_t(L^\infty_x)[I]} \ &\lesssim_F \ 2^{-\delta |k-k_*|} \ ,
    &\hbox{for}\ &k\leqslant k_* +10 \ . \notag
\end{align}
This last inequality follows  from \eqref{Sdiff_FE}, and the fact
that the difference $\mathcal{S}(\phi_{\leqslant k_*})-
\mathcal{S}(\phi)_{\leqslant k_*}$ is  rapidly decaying outside of a
compact set, so in particular one can control the $L^\infty$ norm by
summing dyadically from $k=-\infty$. Note that  the use of
\eqref{ext_L2} and \eqref{ext_N} costs a power of $F$ because these
estimates are in normalized form.

In the second case ($k>k_*+10$), we establish \eqref{basic_G_bound}
by directly appealing to estimates \eqref{core_L2} and
\eqref{core_N}, which suffice because of \eqref{Sdiff_FE} and the
observation that due to the fact $T$ is translation invariant we
have the identity:
\begin{equation}
    P_k \mathcal{C}_1 \ = \ P_k T\big(P_{k+O(1)}\big[
    \mathcal{S}(\phi_{\leqslant k_*})- \mathcal{S}(\phi)_{\leqslant k_*}\big] ,
    \phi_{\leqslant k_*} , \phi_{\leqslant k_*}\big) \ . \notag
\end{equation}\ret

%--------------------------------------------------------------------

\step{3B}{Estimating the Term $\mathcal{C}_2$} We first observe that from the
definition we have $P_k \mathcal{C}_2\equiv 0$ whenever $k>k_*+10$. Thus, we
only need to deal with the frequency range $k\leqslant k_*+10$. We
split this range into two regions: either $k_*-m\leqslant k
\leqslant k_*+10$ or $k < k_*-m$. Here $m$ is defined as follows:
\begin{equation}
        2^{- m} \ = \ \epsilon^{\delta_1} \ .
        \label{exp_m_bnd}
\end{equation}
Note that this definition has nothing to do with the $m$ in the
decomposition \eqref{more_crap},
and is only local to this step. We now estimate separately:\ret

%------------------------------------------------------------------
\step{3B.1}{The range $k_*-m\leqslant k \leqslant k_*+10$} We may
write:
\begin{multline}
    P_k \mathcal{C}_2 \ = \ T^m_k\big(\mathcal{S}(\phi)_{\leqslant k_*},\phi_{\leqslant k_*},\phi_{\leqslant k_*}\big)
    -  T^m_k\big(\mathcal{S}(\phi),\phi,\phi\big)\\
    + 2^{-k_*}\Big( \td{L}_1\big(\nabla_x \mathcal{S}(\phi)_{<k-m},\partial_\alpha\phi_{<k-m},\partial^\alpha\phi_k\big)
    + \td{L}_2\big(\mathcal{S}(\phi)_{<k-m},\! \nabla_x\partial_\alpha\phi_{<k-m},\partial^\alpha\phi_k\big) \\
    + \td{L}_3\big(\nabla_x\mathcal{S}(\phi)_{<k-m},\partial_\alpha\phi_k,\partial^\alpha\phi_{<k-m}\big)
    + \td{L}_4\big(\mathcal{S}(\phi)_{<k-m},\partial_\alpha\phi_k,\nabla_x\partial^\alpha\phi_{<k-m}\big)\Big)
     \label{high_G2}
\end{multline}
where the $T^m_k$ are defined as on line \eqref{basic_trilin} with
the additional structure and frequency localizations from the
definition of $\mathcal{C}_2$. The $\td{L}_i$ are an additional collection of
translation invariant and disposable trilinear forms resulting from
the commutator rule \eqref{commutator} applied to the second and
third terms on the RHS of line \eqref{basic_trilin}. In particular,
this commutator is trivial unless $k_*-10 \leqslant k \leqslant
k_*+10$,  so $\td{L}_i\equiv 0$ without this further restriction.

For the first two terms on the RHS \eqref{high_G2}, we use
\eqref{phi_tdphi_renorm_bnds}--\eqref{phi_tdphi_mod_bnds} which allows us to apply
\eqref{core_L2} or \eqref{core_N}, and these suffice to give
\eqref{basic_G_bound} in this case because of the frequency gap
\eqref{exp_m_bnd} and the conditions \eqref{delta_list} on the $\delta_i$.

It remains to estimate the commutators. From the version of estimate
\eqref{core_L2} in Remark \ref{trilinear_rem}, and the fact that:
\begin{equation}
        \lp{\nabla_x \mathcal{S}(\phi)_{<k-m}}{S[I]}
        + \lp{\nabla_x\phi_{<k-m}}{S[I]}
        \ \lesssim_F \ 2^{k-m} \ , \notag
\end{equation}
we directly have the $L^2$ bound from line \eqref{basic_G_bound}
via \eqref{exp_m_bnd} and the range
restriction $k_*-10 \leqslant k \leqslant k_*+10$. To
prove the $N$ estimate,  we similarly only need to
show:
\begin{equation}
       \lp{2^{-k_*}\td{L}_i}{N[I]} \ \lesssim_F \ 2^{- m}
        \ . \notag
\end{equation}
To estimate $2^{-k_*}{L}_1$ in $N[I]$ we use
\eqref{standard_est_tri} as follows (again using  $k=k_*+O(1)$):
\begin{align}
        &\lp{2^{-k_*}{L}_1
%   \big(\nabla_x\mathcal{S}
%         (\phi)_{<k-m},\partial_\alpha\phi_k,\partial^\alpha\phi_{<k-m}\big)
    }{N[I]} \
         \lesssim_F \ 2^{-k_*} \sum_{k_1,k_2\leqslant k-m}
%         \substack{ k_1:\ k_1\leqslant k_2\\
%         k_2:\ k_2<k-m }}
     2^{k_1}2^{-\delta(k_1-k_2)_+}
%         \  +  \ 2^{-k_*} \sum_{
%         \substack{ k_1:\ k_1<k-m\\
%         k_2:\ k_2<k_1 }} 2^{k_1}2^{-\delta(k_1-k_2)} \ , \notag\\
         \  \lesssim_F   \ 2^{-m} \ . \notag
\end{align}
The details of these calculations for  other ${L}_i$ are similar and
left to the reader.\ret

%--------------------------------------------------------------------------

\step{3B.2}{The range $k < k_*-m$} Here we simply decompose:
\begin{equation}
    P_k \mathcal{C}_2 \ = \ -
    \sum_{  k_i:\  \text{max} \{k_i\}> k_*}
    P_k  L\big(\mathcal{S}(\phi)_{k_1},\partial^\alpha\phi_{k_2},
    \partial_\alpha\phi_{k_3}\big) \ , \notag
\end{equation}
so in particular at least one of the second two factors must be in
the range $k_i>k_*-10$. We remark that this sum has a $T^m$
structure of the form \eqref{basic_trilin}, so smallness is
guaranteed. The main issue is to also recover the exponential
falloff in the definition of $\{c_k\}$. This may be achieved via a
direct application of estimates \eqref{core_L2} and \eqref{core_N}
by first introducing a high frequency
$(\delta_0,\delta_0)$-admissible $\mathcal{W}[I]$ envelope for
$P_{>k_*-10}\phi$ which we denote by $\{d_k\}$. In particular, we
have $d_k\lesssim_F 2^{\delta_0(k-k_*)}$, so we directly have
\eqref{basic_G_bound} for the above sum.\ret

%------------------------------------------------------------------------------------------------

\step{4}{Estimates for matched frequency differences}
Here we prove  \eqref{bad_N_bound},
\eqref{psi_L2_bound}, and \eqref{good_N_bound}.
To do this, it is enough to demonstrate the bound:
\begin{equation}
    \lp{\mathcal{D}^m(\tphi,\psi)}{N_c[I]} \ \lesssim_{\td{F}}
    \ 2^{4\delta_0 m} \ , \label{gen_bad_N_bound}
\end{equation}
under the assumptions \eqref{l1_low_zip_est-boot} and \eqref{phi_tdphi_renorm_bnds}, the bound:
\begin{equation}
    \lp{\mathcal{D}(\tphi,\psi)}{L^2_t(\dot{H}^{-\frac{1}{2}})_c[I]} \ \lesssim_{\td{F}}  \ \td{\epsilon}^{\delta_1}
    \ , \label{gen_psi_L2_bound}\\
\end{equation}
under the additional assumption \eqref{psi_renorm_control}, and finally the improved estimate:
\begin{equation}
    \lp{\mathcal{D}^m(\tphi,\psi)}{N_c[I]} \ \lesssim_{\td{F}} \ \td{\epsilon}^{\delta_1^2}
     \ , \label{gen_good_N_bound}
\end{equation}
assuming all of the above and also using \eqref{phi_tdphi_mod_bnds}
and \eqref{psi_X_bound}.  Here $\mathcal{D}$ is any expression of the
form \eqref{general_phi_diff_def} for a general trilinear form $T$ as
on line \eqref{T_form_def}, and $\mathcal{D}^m$ denotes a similar
expression in terms of $T^m$ from line \eqref{basic_trilin}.  To a
some extent these tasks are redundant, so we will make some effort to
collapse cases.\footnote{We note that one can eliminate this
  redundancy and also simplify some of the case analysis by simply
  replacing the $S_c[I]$ bootstrap of the current proof with a direct
  bootstrap with respect to the $\mathcal{W}_c[I]$ space. We will not
  pursue this here.}

First we introduce the following decomposition of differences
of trilinear expressions of the form \eqref{T_form_def}:
\begin{align}
    T\big(\mathcal{S}(\tphi),\tphi,\tphi\big) - T\big(\mathcal{S}(\tphi+ \psi),\tphi+ \psi,\tphi+ \psi\big)  \
    &= \ R_0 + R_1 + R_2 + R_3 + R_4 \ . \notag
\end{align}
Here we have set:
\begin{subequations}\label{general_decomp}
\begin{align}
    R_0 \ &= \ T\big(\mathcal{S}'(\tphi,\psi)\psi,\tphi,\tphi\big) \ , \label{R0_def}\\
    R_1 \ &= \ T\big(\mathcal{S}'(\tphi,\psi)\psi,\psi,\tphi\big)
    + T\big(\mathcal{S}'(\tphi,\psi)\psi,\tphi,\psi\big)
    + T\big(\mathcal{S}'(\tphi,\psi)\psi,\psi,\psi\big) \ , \label{R1_def}\\
    R_2 \ &= \ -T\big(\mathcal{S}(\tphi),\psi,\tphi\big)  \ , \label{R2_def}\\
    R_3 \ &= \ -T\big(\mathcal{S}(\tphi),\tphi,\psi\big)  \ , \label{R3_def}\\
    R_4 \ &= \ -T\big(\mathcal{S}(\tphi),\psi,\psi\big)  \ . \label{R4_def}
\end{align}
\end{subequations}
Here we are using $\mathcal{S}'$ as shorthand for the formula:
\begin{equation}
     \mathcal{S}(\tphi)-\mathcal{S}(\tphi+\psi) \ := \
     \mathcal{S}'(\tphi,\psi)\psi \ , \label{diff_line}
\end{equation}
so that it symbolically represents some additional set of smooth
functions obeying the same bounds as the original second fundamental
form $\mathcal{S}$. We proceed to estimate the above terms via two subcases:

%------------------------------------------------------------------------------------------------
\ret \step{4A}{Estimates
  \eqref{gen_bad_N_bound}--\eqref{gen_good_N_bound} for the terms
  $R_1$--$R_4$} The bound \eqref{gen_bad_N_bound} for the terms
$R_2$--$R_4$ follows immediately from the estimate \eqref{large_core_N} of
Remark \ref{trilinear_rem} in view of  the bounds \eqref{l1_low_zip_est-boot}
and \eqref{phi_tdphi_renorm_bnds}.

The estimates \eqref{gen_psi_L2_bound} and \eqref{gen_good_N_bound}
for terms $R_2$--$R_4$ under the additional assumptions
\eqref{psi_renorm_control} and then \eqref{psi_X_bound} result
directly from estimates \eqref{core_L2} and \eqref{core_N}.

It remains to prove  these estimates for $R_1$.  This will be
accomplished with the aid of the following three bounds, which we
state in more detail here for their use in the next step:
\begin{align}
    \lp{P_{>k_* + 10}\tphi}{S_c[I]} \ &\lesssim \ 1
    \ , \label{high_tphi} \\
    \lp{\mathcal{S}'(\tphi,\psi)\psi }{S_k[I]} \ &\lesssim_{\td{F}}\ c_k \ , \label{Sprime_Sk}\\
     \lp{P_{<k}\big[ \mathcal{S}'(\tphi,\psi)\psi \big]}{S[I]} \ &\lesssim_{\td{F}} \
     c_k \ , &\hbox{for}\ k\leqslant k_*+20 \ . \label{Sprime_linfty}
\end{align}
The proof of the first bound follows immediately from the
bootstrapping assumption \eqref{l1_low_zip_est-boot}. The second
bound follows from the first and estimates
\eqref{l1_low_zip_est-boot} and \eqref{phi_tdphi_renorm_bnds} after an application of the product
bounds \eqref{S_prod2}--\eqref{S_prod3} and the Moser estimate
\eqref{basic_moser_c}. The last estimate above follows by summing
the second and using the explicit form of the frequency envelope
$\{c_k\}$, and also using the fact that the product vanishes for
$k=-\infty$.

The proof of \eqref{gen_bad_N_bound}--\eqref{gen_good_N_bound} for $R_1$
is a direct application of estimates \eqref{large_core_N}, \eqref{core_L2}, and \eqref{core_N} in
conjunction with the bounds \eqref{Sprime_Sk}--\eqref{Sprime_linfty}
above.\ret

%------------------------------------------------------------------------------------------------

\ret \step{4B}{Estimates \eqref{gen_bad_N_bound}--\eqref{gen_good_N_bound}
for  $R_0$} We first demonstrate \eqref{gen_psi_L2_bound}--\eqref{gen_good_N_bound}.
We split things into two output frequency cases:
$k\leqslant k_*+20$ and $k > k_*+20$. In the first, we combine the
bounds \eqref{ext_L2} and \eqref{ext_N} with both of \eqref{Sprime_Sk}
and \eqref{Sprime_linfty} above. Notice that the condition
$k\leqslant k_*+20$ and the specific form of the frequency envelope
$\{c_k\}$ below $k_*$ gives the desired result.

To deal with $R_0$ in the output range $k>k_*+20$ requires
additional work. Note that any estimate of the form
\eqref{Sprime_linfty} is false for $k>k_*+20$, so the needed
frequency envelope control needs to come from the second and third
$\td{\phi}$ factors. For this we employ the following version of \eqref{high_tphi}:
\begin{equation}
     \lp{P_{>k_* + 10}\tphi}{\W_c[I]} \ \lesssim_{\td{F}} \ 1
        \ , \label{Whigh_tphi}
\end{equation}
which follows from the assumption \eqref{psi_renorm_control}. To
use it, we split $R_0$ into a sum of pieces (we may drop $L$ from
the picture):
\begin{align}
    R_0 \ = \ &\mathcal{S}'(\tphi,\psi)\psi\cdot\big[
    \partial^\alpha P_{<k_*+10}\tphi \partial_\alpha P_{<k_*+10}\tphi
    +2 \partial^\alpha P_{< k_*+10}\tphi \partial_\alpha P_{\geqslant k_*+10}\tphi \notag\\
    &+ \partial^\alpha P_{\geqslant k_*+10}\tphi \partial_\alpha P_{\geqslant k_*+10}\tphi \big]
     \ = \ R_{0,1} + R_{0,2} + R_{0,3} \ .  \label{R0_decomp}
\end{align}
The estimates \eqref{gen_psi_L2_bound}--\eqref{gen_good_N_bound} for
the pieces $P_{>k_*+20}R_{0,2}$ and $P_{>k_*+20} R_{0,3}$ are an
immediate consequence of estimate \eqref{core_L2} and \eqref{core_N}
because in either case we can use $\{c_k\}$ frequency envelope control
of the high frequencies of at least one of the $\td{\phi}$ factors. To
handle the term $P_{>k_*+20}R_{0,1}$ we again use \eqref{core_L2} and
\eqref{core_N} along with the bound \eqref{Sprime_Sk}, which provides
the needed $\{c_k\}$ factor on account of the forced $P_k$ frequency
localization of the first factor in $P_k P_{>k_*+20}R_{0,1}$.

The proof of estimate \eqref{gen_bad_N_bound} is similar to \eqref{gen_bad_N_bound}
above except that we use \eqref{high_tphi} instead of \eqref{Whigh_tphi},
and \eqref{large_core_N}--\eqref{large_ext_N} instead of \eqref{core_N} and \eqref{ext_N}.\ret

%------------------------------------------------------------------------------------------------

\step{5}{Estimates for $\psi$ at low frequency}
Here we prove \eqref{L_bound}. From the definition of $\mathcal{L}_k^m$ on line \eqref{Lpsi_def},
we see that is suffices to prove
the two general bounds:
\begin{align}
        \lp{
        \td{\phi}^{(1)}_{<k-m}\partial^\alpha\psi_{<k-m}\partial_\alpha\tphi^{(2)}_k
        }{N[I]} \
        &\lesssim_{\td{F}} \ (\td{\epsilon}^{\delta_0} + 2^{-\delta_0 m})c_k \ , \label{ext_bnd1}\\
        \lp{\big(\mathcal{S}'(\tphi,\psi)\psi\big)_{<k-m}
        \partial^\alpha \tphi^{(1)}_{<k-m}\partial_\alpha\tphi^{(2)}_k
        }{N[I]} \ &\lesssim_{\td{F}} \
        (\td{\epsilon}^{\delta_0} + 2^{-\delta_0 m})c_k
        \ , \ \label{ext_bnd2}
\end{align}
where $\lp{\td{\phi}^{(i)}}{S[I]}\lesssim_{\td{F}}1$, and where the
$\td{\phi}^{(i)}$ also has high frequency improvement
\eqref{high_tphi}.  We split into two cases:\ret

%-------------------------------------------------------------------
\step{5A}{The range $k<k_*+10$} For estimate
\eqref{ext_bnd1} we use \eqref{N_prod_est1} and
\eqref{standard_est_bi} which together give:
\begin{equation}
        \hbox{(L.H.S.)}\eqref{ext_bnd1} \ \lesssim_{\td{F}} \
        \sum_{j<k-m} c_j \
        \lesssim_{\td{F}} \ c_{k-m} \ \lesssim_{\td{F}} \
        2^{-\delta_0 m}c_k \ . \notag
\end{equation}
For estimate
\eqref{ext_bnd2} we use \eqref{standard_est_tri} and
\eqref{Sprime_Sk}:
\begin{equation}
        \hbox{(L.H.S.)}\eqref{ext_bnd2} \ \lesssim_{\td{F}}\
        \sum_{j,k_1<k-m}\!\!\!\!\!\! 2^{-\delta(j-k_1)_+}c_j
        \ \lesssim_{\td{F}}\
        \sum_{j<k-m}\!\!\! (k-m-j)c_j
         \ \lesssim_{\td{F}}\
        2^{-\delta_0 m}c_k \ . \notag
\end{equation}\ret

%-------------------------------------------------------------------
\step{5B}{The range $k>k_*+10$} Here we use the fact that
the $\td{\phi}^{(i)}$ also have high frequency improvement
\eqref{high_tphi}, which already incorporates $\{c_k\}$. We only
need to gain a small factor. By \eqref{N_prod_est1} and
\eqref{standard_est_bi}, and the fact that
$\sum_k\lp{\psi}{S_k[I]}\lesssim_{\td{F}}\td{\epsilon}^{\delta_0}$
we immediately have \eqref{ext_bnd1}.

The proof of \eqref{ext_bnd2} follows from \eqref{Sprime_Sk} and
\eqref{standard_est_tri}, and the summation:
\begin{equation}
        \hbox{(L.H.S.)}\eqref{ext_bnd2}   \ \lesssim_{\td{F}} \
        c_k
        \sum_{j<k-m} |k_*-j|c_j
           \ \lesssim_{\td{F}}\
        \td{\epsilon}^{\delta_0 }c_k \ . \notag
\end{equation}
Notice that the exponential falloff of $\td{\phi}^{(1)}$ was used in the
range $k_1>k_*$ to reduce the factor of $(k-m-j)$ from the previous
step to $|k_*-j|$ here, and thus avoid a logarithmic divergence.
This concludes our demonstration of Proposition \ref{p:phi_est_low}.
\end{proof}\ret\ret

%----------------------------------------------------------------

\begin{proof}[Proof of Proposition~\ref{p:phi_est_high}]
We denote the difference to be estimated by:
\begin{equation}
    \phi^{high} \ = \ \phi - \td{\phi} \ , \notag
\end{equation}
which represents the evolution of high frequencies in $\phi$. This
solves the equation:
\begin{equation}
    \Box \phi^{high}  \ =  \ - \mathcal{S}(\phi) \partial^\alpha \phi \partial_\alpha \phi
    + \mathcal{S}(\tphi) \partial^\alpha \tphi \partial_\alpha \tphi
    \ . \label{high_evolution}
\end{equation}
A natural attempt is to argue directly as in the preceding  proof,
namely to replace this nonlinear equation for $\phi^{high}$ with a linear
paradifferential equation plus a nonlinear perturbative term.
However, if we do that directly we encounter some difficulties:

Precisely, the initial data $  \phi^{high}[0]$ has size on the order
of $c_0 = c_0(E)$. Solving the linear paradifferential equation we
loose a constant which depends on the $S[I]$ size of the
coefficients, namely at least $K(\td{F}(E))$. Thus the solution for
the approximate linear flow will have size $c_0(E) K(\td{F}(E))$,
and the key point is that we cannot expect this to be small. This
would cause the nonlinear effects to be truly non-perturbative,
and therefore outside the scope of the current paper.

One fix to this would be to allow $c_0$ to depend on $\td{F}(E)$
instead of $E$. However, this would weaken our conclusion to the
point where the induction on energy argument only works to show that
the set of energies where one has regularity is \emph{open}. While
this is the usual first step in an induction on energy strategy, it
still leaves one to deal with the heart of the matter which is the
task of showing that there is no finite upper bound to the set of
regular energies. Our path here will be a to establish this
latter claim more directly.

As a first step in our argument, we subdivide the time interval $I$
into consecutive subintervals $I_k$, and we can insure that on each such
subinterval we have the partial fungibility:
\begin{equation}
        \lp{ \tphi}{S[I_k]}  \ \lesssim  \ E \ .
        \label{local_tphi_bnd}
\end{equation}
This is possible by estimate \eqref{fung_bnd}. This estimate
remedies the bootstrapping argument within the first interval
because by design $\phi^{high}$ has small initial energy (see
\eqref{phi_high_eng_bnd} below). However, one might expect that in each
subinterval $I_i$ the energy of $\phi^{high}$ may grow by a factor
$K(E)$ factor, and the number of intervals where
\eqref{local_tphi_bnd} is true unfortunately depends on $\td{F}(E)$.
Thus a brute force bound would allow the energy of $\phi^{high}$ to
grow by a $K(F(E))$ factor, which would bring us back to the core
difficulty. However, such a brute force approach does not make any
good use of the fact that both $\tphi$ and $\phi$ are true
Wave-Maps, and therefore exactly conserve their energy.

In order to take advantage of this last observation, we compute that
for each fixed $t\in I$:
\begin{equation}
        E(\phi)  \ =  \ E(P_{>k_*} \phi) + E(P_{\leqslant k_*}\phi)
        + 2 \langle P_{>k_*} \phi,
        P_{\leqslant k_*} \phi\rangle \ \geqslant \
         E(P_{>k_*} \phi) + E(P_{\leqslant k_*} \phi) \ , \notag
\end{equation}
where $\langle \cdot,\cdot \rangle$ denotes the $\dot H^1 \times
L^2$ inner product. The last inequality holds because $P_{\leqslant
k_*} P_{>k_*}$ is a nonnegative operator because it has a
nonnegative symbol. By Proposition \ref{p:phi_est_low} just proved,
we have the fixed time bound in $I$:
\begin{equation}
        \lp{(P_{\leqslant k_*} \phi - \tphi)[t]}{ \dot H^1 \times L^2}  \
        \lesssim_F \ \epsilon^{\delta_0} \ . \notag
\end{equation}
Hence we have:
\begin{equation}
        | E(P_{\leqslant k_*} \phi) - E(\tphi)|+| E(\phi-\tphi)-E(P_{>k_*} \phi)|
        \ \lesssim_F \  \epsilon^{\delta_0}\ . \notag
\end{equation}
Thus
\begin{align}
        E(\phi^{high}) \ &\leqslant \ E(P_{>k_*} \phi) + \epsilon^{\delta_0} K(F(E))
         \ , \notag\\
        &\leqslant \ E(\phi) - E(P_{\leqslant k_*} \phi) + \epsilon^{\delta_0} K(F(E))
         \ , \notag\\
       &\leqslant \  E(\phi) - E( \tphi) + 2\epsilon^{\delta_0} K(F(E))
        \ , \notag\\
       &\leqslant \ c_0(E) + 2\epsilon^{\delta_0} K(F(E)) \ .
        \label{phi_high_eng_bnd}
\end{align}
This calculation shows that if $\epsilon$ is small enough with
respect to the function $K(F)$ which appears implicity on the RHS of
estimate \eqref{low_zip_est} then we have a good uniform bound on
the energy of $\phi^{high}$. The argument now proceeds in
a series of steps:\ret

%-----------------------------------------------------------------------------

\step{1}{The bootstrapping construction, and the main estimates} We
now fix the interval $I_i\subseteq I$ and consider an $S$-norm
bootstrap for the $\phi^{high}$ on subintervals $J\subseteq I_i$,
where we may assume $J$ is centered about $t=0$. We seek to prove the
bound:
\begin{equation}
    \lp{\phi^{high}}{S[J]}\ \leqslant \ 1 \ .
    \label{phihigh_s}
\end{equation}
Due to \eqref{phi_high_eng_bnd} we have:
\begin{equation}
    \lp{\phi^{high}[0]}{\dot H^1 \times L^2}  \ \lesssim  \ c_0 \ .
    \label{phihdata}
\end{equation}
Hence by the second part of Proposition~\ref{propboot}
we obtain the seed bound:
\begin{equation}
    \lp{\phi^{high}}{S[J_0]}  \ \lesssim  \ c_0 \ , \notag
\end{equation}
for a small enough interval $J_0 \subset J$. Taking this into account
and also the continuity of the $S$ norm in Proposition~\ref{propboot},
it suffices to prove \eqref{phihigh_s} under the additional
bootstrap assumption:
\begin{equation}
    \lp{\phi^{high}}{S[J]}\ \leqslant \ 2  \ .
    \label{phihighboot}
\end{equation}

Combining \eqref{phihighboot} with \eqref{local_tphi_bnd} we obtain:
\begin{equation}
    \lp{\phi}{S[J]} + \lp{\tphi}{S[J]}\ \lesssim_E 1 \ .
    \label{phitphi_s}
\end{equation}
By Proposition~\ref{p:wm_struct} this gives:
\begin{equation}
    \lp{ \phi}{\W[J]} + \lp{ \tphi}{\W[J]}  \ \lesssim_E  \ 1 \ ,
    \qquad \lp{\phi}{\uX[J]} + \lp{\tphi}{\uX[J]}
     \ \lesssim_E  \ \td{\epsilon}^{\delta_1} \ .
    \label{phitphie}
\end{equation}
We rewrite the bounds \eqref{phihdata}, \eqref{phihighboot} and \eqref{phitphie} using frequency
envelopes. Precisely, we can find a common
$(\delta_0/2,\delta_0/2)$-admissible normalized    frequency envelope
$c_k$ so that $c_{k_*} = 1$ and the following bounds hold:
\begin{subequations}\label{samec}
\begin{align}
    \lp{\phi^{high}[0]}{(\dot H^1 \times L^2)_c} \ &\lesssim_E  \ c_0 \ , \label{samec1} \\
    \lp{\phi^{high}}{S_c[J]}\ &\lesssim_E  \   1 \ , \label{samec2} \\
    \lp{\phi}{\W_c[J]} + \lp{\tphi}{\W_c[J]}\ &\lesssim_E  \  1 \ , \label{samec3} \\
    \lp{ \phi}{\uX_c[J]} + \lp{\tphi}{\uX_c[J]}\ &\lesssim_E \ \td{\epsilon}^{\delta_1} \ . \label{samec4}
\end{align}
\end{subequations}
From these four bounds, together with the energy dispersion on lines
\eqref{phi_low_freq_assumt}--\eqref{tdphi_low_freq_assumt}, we will obtain the following
vastly improved frequency envelope $S$ bound for $\phi^{high}$:
\begin{equation}
    \lp{\phi^{high}}{S_{\lambda c}[J]}  \ \lesssim_E\ 1 \ ,
    \qquad \lambda = c_0  +  \td{\epsilon}^{\delta_0 \delta_1^2}
    \ . \label{phihc}
\end{equation}
The second term on the right is small ($\ll 1$) due to \eqref{epsilon_bnds},
therefore the desired conclusion \eqref{phihigh_s} follows if $c_0$
is chosen appropriately small, $c_0 \ll_E 1$.
\ret

It remains to show that \eqref{phi_low_freq_assumt}--\eqref{tdphi_low_freq_assumt}
together with \eqref{samec} imply \eqref{phihc}.  By estimate
\eqref{l1_low_zip_est}, we may reduce this demonstration to the
frequency range $k>k_*-10$.
The mechanics of our argument is to decompose the $P_k$ frequency
localized version of \eqref{high_evolution} as follows:
\begin{equation}
    \Box\phi^{high}_k + 2A^\alpha(\phi)_{<k-m}\partial_\alpha \phi^{high}_k
    \ = \ \mathcal{T}^m_k(\tphi) + \mathcal{T}^m_k(\phi) + \mathcal{L}^m_k(\tphi,\phi^{high})
    \ , \label{para_high}
\end{equation}
where the large gap
parameter $m$ is consistent with Proposition~\ref{big_tri_prop}:
\begin{equation}
    2^{-m}  \ =  \ \td{\epsilon}^{\delta_1} \ .
    \label{choosem}
\end{equation}
Here the terms $\mathcal{T}^m_k$ are matched frequency
trilinear expressions of the form \eqref{basic_trilin}, while the term $\mathcal{L}^m_k$
denotes certain trilinear expressions between $\tphi$ and $\phi^{high}$ which contain
at least one  ($m$ dependent) low or high frequency factor with improved exponential
bounds.

Our first round of estimates shows  that:
\begin{equation}
    \lp{\mathcal{T}^m_k(\tphi)}{N[J]} \ \lesssim_{E} \ \td{\epsilon}^{\delta_1^2} c_k \ ,
    \qquad \lp{\mathcal{T}^m_k(\phi)}{N[J]} \ \lesssim_E \ \epsilon^{\delta_1^2}c_k \ . \label{high_Tm_bounds}
\end{equation}
Our second round of estimates gives the exponential control:
\begin{equation}
    \lp{\mathcal{L}^m_k(\tphi,\phi^{high})}{N[J]} \ \lesssim_E \
      2^{-\frac14 \delta_0 m}  2^{-\frac{1}{2}\delta_0|k-k_*|} \ = \
    \td{\epsilon}^{\frac14 \delta_0 \delta_1} 2^{-\frac{1}{2}\delta_0|k-k_*|}
     \ . \label{phi_high_low_ests}
\end{equation}
An application of \eqref{linearized_est} using \eqref{samec1}  and \eqref{high_Tm_bounds}--\eqref{phi_high_low_ests}
then implies \eqref{phihc}.\ret

%--------------------------------------------------------------------------------
\step{2}{Algebraic derivation of \eqref{para_high}}
We first write the frequency localized equation for $\phi^{high}$
as follows:
\[
\begin{split}
    \Box \phi^{high}_k   = &  - 2\mathcal{S}(\phi)_{<k-m} \partial^\alpha
    \phi_{<k-m} \partial_\alpha \phi_k
    + 2\mathcal{S}(\tphi)_{<k-m} \partial^\alpha \tphi_{<k-m} \partial_\alpha \tphi_k \\
  &  + \mathcal{T}^m_{1;k}(\phi) +  \mathcal{T}^m_{1;k}(\tphi)
    \ , \notag
\end{split}
\]
where the $\mathcal{T}^m_{1;k}$ are trilinear forms as defined on line \eqref{basic_trilin}.
Adding to this a zero expression similar to \eqref{zero_exp} (i.e. without $P_{\leqslant k_*}$
on the second factor), and further decomposing
the result into principle terms and $T^m$ interactions, we have:
\[
    \Box \phi^{high}_k   =   - 2A^\alpha(\phi)_{<k-m} \partial_\alpha \phi_k
    + 2A^\alpha(\tphi)_{<k-m} \partial_\alpha \tphi_k
    + \mathcal{T}^m_{k}(\phi) +  \mathcal{T}^m_{k}(\tphi)
    \ , \notag
\]
Then equation \eqref{para_high} is achieved by setting:
\begin{equation}
    \mathcal{L}^m_k \ = \ -\big(A^\alpha(\tphi+\phi^{high})_{<k-m}
    -A^\alpha(\tphi)_{<k-m}\big)\partial_\alpha \tphi_k \ . \label{phi_high_L_form}
\end{equation}\ret

%--------------------------------------------------------------------------------
\step{3}{Control of matched frequency interactions}
Here we prove \eqref{high_Tm_bounds}. This is an immediate consequence
of \eqref{core_N} using  \eqref{samec3}--\eqref{samec4}.
\ret

%--------------------------------------------------------------------------------
\step{4}{Control of separated frequency interactions}
Here we prove the estimate \eqref{phi_high_low_ests}. We decompose line \eqref{phi_high_L_form}
as a sum of two terms $\mathcal{L}^m_k=R_1+R_2$ where:
\begin{align}
    R_{1;k} \ &= \ -\mathcal{S}(\tphi)_{<k-m}\partial^\alpha\phi^{high}_{<k-m}\partial_\alpha \tphi_k
    \ , %\label{high_R1_def}
    \notag\\
    R_{2;k} \ &= \ \big[ \mathcal{S}'(\tphi,\phi^{high})\phi^{high}\big]_{<k-m}
    \partial^\alpha\phi_{<k-m}\partial_\alpha \tphi_k
    \ , %\label{high_R2_def}
    \notag
\end{align}
where  $\mathcal{S}$ now denotes the antisymmetrization of the original second fundamental form,
and $\mathcal{S}'$ is defined as on line \eqref{diff_line}. Recall that we are restricted to the
conditions $k\geqslant k_*-10$ and $m\geqslant 20$. We proceed to estimate
each term separately. In doing so, we repeatedly use the following estimates:
\begin{align}
    \lp{P_k \tphi}{S[J]} \ &\lesssim \ 2^{-\delta_0(k-k_*)} \ ,
    &k\geqslant k_* \ , \label{high_tdphi}\\
    \lp{P_k \phi^{high}}{S[J]} \ &\lesssim \ 2^{-\delta_0(k_*-k)} \ ,
    &k\leqslant k_* \ ,  \label{low_phihigh}\\
    \lp{P_k\big[ \mathcal{S}'(\tphi,\phi^{high})\phi^{high}\big]}{S[J]} \ &\lesssim_E \
    2^{-\delta_0(k_*-k)} \ ,
    &k\leqslant k_* \ ,  \label{Sprod_phihigh1}\\
    \lp{ \mathcal{S}'(\tphi,\phi^{high})\phi^{high}}{S[J]} \ &\lesssim_E \
    1 \ .  \label{Sprod_phihigh2}
\end{align}
Estimates \eqref{high_tdphi}--\eqref{low_phihigh} are simply a
weaker restatement of \eqref{l1_low_zip_est}
for the convenience of the reader. Estimates \eqref{Sprod_phihigh1}--\eqref{Sprod_phihigh2}
follow from \eqref{low_phihigh}, \eqref{phihighboot}--\eqref{phitphi_s},
and the Moser and product estimates
\eqref{S_prod2}--\eqref{S_prod3} and \eqref{basic_moser}--\eqref{basic_moser_c}
after a standard summation argument.\ret

%--------------------------------------------------------------------------------
\step{4A}{Estimating $R_{1;k}$}
After an application of \eqref{N_prod_est1}--\eqref{N_prod_est2} to peel off the first factor,
it suffices to show the bound:
\begin{equation}
    \lp{\partial^\alpha\phi^{high}_{<k-m}\partial_\alpha \tphi_k}{N[J]} \ \lesssim_E \
    2^{-\frac{1}{4}\delta_0 m}
    2^{\frac{1}{2}\delta_0(k_*-k)} \ . \notag%\label{easy_sep_bnd}
\end{equation}
If $k <k_*$ this follows at once from \eqref{phitphi_s} and
\eqref{low_phihigh}, and
summing over \eqref{standard_est_bi}.

If $k>k_*$, we use   \eqref{high_tdphi}, and \eqref{low_phihigh} or \eqref{phihighboot}
in \eqref{standard_est_bi} to obtain after summation:
\begin{equation}
    \lp{\partial^\alpha\phi^{high}_{<k-m}\partial_\alpha \tphi_k}{N[J]}
    \ \lesssim_E\ 2^{-\delta_0(k-k_*)} \sum_{j < k-m} 2^{-\delta_0(k_*-j)_+} \ . \notag
\end{equation}
If $k_* < k < k_*+m$ then the expression on the right gives
$ 2^{-\delta_0 m}$ which suffices. If $k>k_*+m$
then the expression on the right gives $|k-k_*| 2^{-\delta_0(k-k_*)}$
which is again sufficient for \eqref{phi_high_low_ests}.
\ret

%--------------------------------------------------------------------------------
\step{4B}{Estimating $R_{2;k}$}
In this final step we show the estimate:
\begin{equation}
    \lp{ R_{2;k}}{N[J]} \ \lesssim_E \ 2^{-\frac14 \delta_0 m}
    2^{\frac{1}{2}\delta_0(k_*-k)} \ . \notag
\end{equation}
In the case when $k \leqslant k^*$, using \eqref{Sprod_phihigh1},
\eqref{phitphi_s}, and the trilinear estimate \eqref{standard_est_tri} we
have the sum:
\begin{equation}
    \lp{ R_{2;k}}{N[J]} \ \lesssim_E \  \sum_{j,k_1< k-m} 2^{-\delta(j-k_1)_+}2^{-\delta_0(k_*-j)}
    \ \lesssim_E \  2^{-\delta_0 m} 2^{-\delta_0|k-k_*|} \ . \notag
\end{equation}
Finally, in the case where $k>k_*$ we use \eqref{phitphi_s}, \eqref{high_tdphi}, and
\eqref{Sprod_phihigh1}--\eqref{Sprod_phihigh2} in conjunction with \eqref{standard_est_tri}
to achieve the sum:
\begin{equation}
    \lp{ R_{2;k}}{N[J]} \ \lesssim_E \ 2^{-\delta_0(k-k_*)}
    \sum_{j,k_1<k-m} 2^{-\delta(j-k_1)_+}2^{-\delta_0(k_*-j)_+} \ . \notag
\end{equation}
If $k_* < k < k_*+m$ then the expression on the right gives $
2^{-\delta_0 m}$ which is enough for \eqref{phi_high_low_ests}. If
$k>k_*+m$ then the expression on the right gives $|k-k_*|
2^{-\delta_0(k-k_*)}$ which is again sufficient. This concludes our
demonstration of Proposition \ref{p:phi_est_high}.
\end{proof}

\ret
%-------------------------------------------------------------------------
%%%%%%%%%%%%%%%%%%%%%%%%%%%%%%%%%%%%%%%%%%%%%%%%%%%%%%%%%%%%%%%%%%%%%%%%%%
%%%%%%%%%%%%%%%%%%%%%%%%%%%%%%%%%%%%%%%%%%%%%%%%%%%%%%%%%%%%%%%%%%%%%%%%%%
%-------------------------------------------------------------------------

\section{The Iteration Spaces: Basic Tools and Estimates}\label{ext_sect}
This is a continuation of Section \ref{old_est_section}, and our
purpose is to fill in any gap between the notation and additional
structure of basic function spaces used in this paper and the spaces
developed in \cite{Tataru_WM2}--\cite{Tataru_WM1} and \cite{Tao_WM}.\ret

%----------------------------------------------------------------------------

\subsection{Space-time and angular frequency cutoffs}
As usual we denote by $Q_j$ the multiplier with symbol:
\begin{equation}
        q_j(\tau,\xi) \ = \ \varphi\big(2^{-j} \big| |\tau|-|\xi|\big| \big)
         \ , \notag
\end{equation}
where $\varphi$ truncates smoothly on a unit annulus. We denote by
$Q_j^\pm$ the restriction of this multiplier to the upper or lower
time frequency space. At times we also denote by $Q_{|\tau|\lesssim
|\xi|}=Q_{<C}$ for some $C>0$.\ret

We denote by $\kappa\in K_l$ a collection of caps of diameter $\sim
2^{-l}$ providing a finitely overlapping cover of the unit sphere.
According to this decomposition, we cut up the spatial frequency
domain according to:
\begin{equation}
        P_k \ = \ \sum_{\kappa\in K_l} P_{k,\kappa} \ . \notag
\end{equation}
These decompositions often occur in conjunction with modulation
cutoffs on the order of $j=k-2l$, and a central principle is that
the corresponding multipliers $Q^\pm_{<k-2l} P_{k,\pm\kappa}$ are
uniformly disposable.\ret

%----------------------------------------------------------------------------

\subsection{The $S$ and $N$ function spaces}

\begin{defn}[Dyadic Iteration Space]
For each integer $k$ we define the following frequency localized
norm:
\begin{equation}
        \lp{\phi}{S_k} \ := \ \lp{\nabla_{t,x}\phi_k }{L^\infty_t(L^2_x)}
        +  \lp{\nabla_{t,x}\phi_k}{X_\infty^{0,\frac{1}{2}}}
        + \lp{\phi_k}{\underline{S}}
        + \sup_{j < k-20} \lp{\phi}{S[k;j]} \ .
        \label{S_norm}
\end{equation}
In general, the fixed frequency space $X_p^{s,b}$ is defined as:
\begin{equation}
        \lp{P_k \phi}{X_p^{s,b}}^p \ := \ 2^{psk}
        \sum_{j}2^{pbj}\lp{Q_jP_k\phi}{L^2_t(L^2_x)}^p \ , \notag
\end{equation}
with the obvious definition for $X_\infty^{s,b}$. Here we define the
``physical space Strichartz'' norms:
\begin{equation}
        \lp{\phi_k}{\underline{S}} \ :=\  \sup_{(q,r):\ \frac{2}{q}+ \frac{1}{r}\leqslant
        \frac{1}{2} }2^{(\frac{1}{q}+ \frac{2}{r}-1)k}
        \lp{\nabla_{t,x}\phi_k}{L^q_t(L^r_x)}
        \ , \label{phys_str}
\end{equation}
 the ``modulational Strichartz'' norms:
\begin{align}
    \lp{\phi}{S[k;j]} \ &:= \ \sup_\pm \big(\sum_{\kappa\in K_{l}} \lp{Q^\pm_{<k-2l}
        P_{k,\pm\kappa}\phi}{S[k,\kappa]}^2\big)^\frac{1}{2} \ ,
\qquad  l=\frac{k-j}{2} > 10 \ , \label{mod_Str_space}
\end{align}
and the ``angular Strichartz'' space  in terms of the
three components:
\begin{multline}
        \lp{\phi}{S[k,\kappa]} \ := \ 2^k \sup_{\omega\notin 2\kappa}
        \hbox{dist}(\omega,\kappa)\lp{\phi}{L^\infty_{t_\omega}(L^2_{w_\omega})}
        + 2^k\lp{\phi}{L^\infty_t(L^2_x)}\\
        + 2^{\frac{1}{2}k}|\kappa|^{-\frac{1}{2}}
        \inf_{\sum_\omega \phi^\omega = \phi} \sum_\omega \lp{\phi^\omega }
        {L^2_{t_\omega}(L^\infty_{x_\omega})}\ . \label{str_norm}
\end{multline}
The first component on the RHS above will often be referred to as
$N\!\!F\!\!A^*$. We define $S$ as the space of functions $\phi$ in $\R^{2+1}$
with $\nabla_{x,t} \phi \in C(\R;L^2_x)$ and finite norm:
\begin{equation}
    \lp{\phi}{S}^2  \ =  \ \lp{\phi}{L^\infty_t(L^\infty_x)}^2
    + \sum_k \lp{\phi}{S_k}^2 \ , \notag
\end{equation}
and also use the frequency envelope convention from Section \ref{old_est_section}
to define $S_c$.
\end{defn}\ret

To measure the derivatives of functions in $S$ we  introduce
a related space $DS$:

\begin{defn}[Differentiated $S$ functions]
We define the norm:
\begin{equation}
        \lp{\phi}{DS_k} \ := \ \lp{\phi_k }{L^\infty(L^2_x)}
        +  \lp{\phi_k}{X_\infty^{0,\frac{1}{2}}}
        + \lp{\phi_k}{D\underline{S}_k}
        + 2^k \sup_{j < k-20} \lp{\phi}{S[k;j]} \ .
        \label{dtS_norm}
\end{equation}
where the $D\uS$ norm is as in \eqref{phys_str} but without the
gradient:
\begin{equation}
        \lp{\phi}{D\underline{S}_k} \ :=\  \sup_{(q,r):\ \frac{2}{q}+ \frac{1}{r}\leqslant
        \frac{1}{2} }2^{(\frac{1}{q}+ \frac{2}{r}-1)k}
        \lp{\phi}{L^q_t(L^r_x)}
        \ . \label{dphys_str}
\end{equation}
The $DS$ space is defined as the space of functions for which
the square sum of the $DS_k$ norms is finite:
\begin{equation}
    \lp{\phi}{DS}^2  \ =  \ \sum_k \lp{\phi}{DS_k}^2 \ . \notag
\end{equation}
\end{defn}

We remark that by definition we have:
\begin{equation}
    \lp{\phi}{S_k}  \ \approx  \ \lp{\nabla_{x,t} \phi}{DS_k} \ .
    \label{s=ds}
\end{equation}\ret

% \begin{rem}\label{weak_high_mod_rem}
% We remark here that the product estimates
% \eqref{S_prod2}--\eqref{S_prod3} are still valid for
% \eqref{dtS_norm}, as can easily be shown by estimating the high
% modulation output of products in $L^2_t(L^2_x)$ with the additional
% assumption that at least one term in the product is also at high
% modulation.
% \end{rem}\ret

\begin{defn}[Dyadic Source Term Space]
For each integer $k$ we define the following frequency localized
norm:
\begin{multline}
        \lp{F}{N_k} \ := \ \inf_{F_A + F_B + \sum_{l,\kappa}
        F_{C}^{l,\kappa} = F}\Big(
        \lp{P_k F_A}{L^1_t(L^2_x)} + \lp{P_k
        F_B}{X_1^{0,-\frac{1}{2}}}\\
        + \sum_{\pm} \sum_{l>10}
        \big(  \sum_{\kappa}  \inf_{\omega\notin 2\kappa}
        \hbox{dist}(\omega,\kappa)^{-2}\lp{Q^\pm_{<k-2l}
        P_{k,\pm\kappa}F_C^{l,\kappa} }{L^1_{t_\omega}(L^2_{x_\omega})}^2\big)^\frac{1}{2}
        \Big) \ . \label{N_norm}
\end{multline}
We will often refer to the last component on the RHS above as
$N\!\!F\!\!A$, and the norm applied to a fixed $Q^\pm
F_C^{l,\kappa}$ as $N\!\!F\!\!A[\pm\kappa]$.
\end{defn}\ret

%These ``spaces'' just defined\footnote{Notice also the
%notational difference with \cite{Tao_WM},
%which made $S_k$ and $N_k$ into Banach spaces by allowing frequency leakage, and
%also implicitly restricted them to time slabs $I\times\R^2$.
%Our definitions of $S_k$ and $N_k$ correspond to the labels
%$S[k]$ and $N[k]$ defined in Section 10 of  \cite{Tao_WM}.}  are only
%semi-norms. We construct full
%$\ell^2$ and frequency envelope based norms as follows:\\

%\begin{defn}[Full Norms]
%For functions $\phi$ and $F$, and a we form the following Banach
%space norms:
%\begin{align}
%        \lp{\phi}{S} \ &:= \ \lp{\phi}{L^\infty_t(L^\infty_x)}
%        + \big(\sum_k \lp{P_k \phi}{S_k}^2\big)^\frac{1}{2} \ ,
%        \label{full_S_norm}\\
%        \lp{F}{N} \ &:= \ \big(\sum_k \lp{P_k F}{N_k}^2\big)^\frac{1}{2} \ .
%        \label{full_N_norm}
%\end{align}
%Frequency enveloped versions $S_c$ and $N_c$ are defined as on line
%\eqref{X_c_def} with the exception that for $S_c$ we add the
%$L^\infty$ norm as follows:
%\begin{equation}
%        \lp{\phi}{S_c} \ := \ \lp{\phi}{L^\infty_t(L^\infty_x)}
%        + \sup_k \max\{ c_k^{-1}\lp{P_k \phi}{S_k},1\} \ .
%        \label{Sc_norm}
%\end{equation}
%We also use the notation $S\setminus L^\infty$ to denote the second member
%on the RHS of line \eqref{full_S_norm}.
%\end{defn}\ret

For any \emph{closed} interval $I=[t_0,t_1]$ we define spaces
$\uX[I]$, $\uX_c[I]$, $\underline{E}[I]$,
$L^\infty_t(L^\infty_x)[I]$, etc. as the restriction of these
classical $L^p$ based norms to the time slab $I\times
\mathbb{R}^2_x$. We also need a similar procedure for the non-local
$S$ and $N$ spaces. As usual we define $S_k[I]$, $S[I]$, $S_c[I]$,
$N[I]$, etc. in terms of minimal extension.\footnote{We will modify
this procedure somewhat below by an equivalent norm, but for the
most part they are interchangeable.} For example:
\begin{equation}
    \lp{\phi}{S_k[I]} \ = \ \inf_{\Phi} \big\{\lp{\Phi}{S_k} \, \big| \ \ \Phi\equiv \phi \ \hbox{on} \
    I\times \R^2 \big\} \ . \label{ext_norm}
\end{equation}
On an \emph{open} time interval $(t_0,t_1)$ we may also define
localized norms by taking $\lp{\phi}{S(t_0,t_1)}=\sup_{I\subseteq
(t_0,t_1)}\lp{\phi}{S[I]}$. This  definition will only be important
for us  as a convenience when stating results like Theorem
\ref{main_thm}, so the reader is safe to ignore the distinction and
always assume that $I$ denotes
a \emph{closed} time interval. We now state a continuation of Proposition \ref{standard_prop1}:\ret

\begin{prop}[Standard Estimates and Relations: Part
2]\label{standard_prop2} Let $F$, $\phi$, and $\phi^{(i)}$ be a
collection of test functions, $I\subseteq \mathbb{R}$ any
subinterval (including $\R$ itself). Then the following  list of
properties for the $S[I]$ and $N[I]$ spaces hold:
\begin{itemize}
\item (Time Truncation of S) Let $\chi_I$ be the characteristic function
of $I$. Then
    \begin{align}
        \lp{\phi}{DS_k[I]}   \ &\approx   \ \lp{ \chi_I \phi}{DS_k}   \ \lesssim  \ \lp{\phi}{DS_k}
        \ , \label{dscut}\\
        \lp{\phi}{S_k[I]}  \ &\approx  \ \lp{\chi_I\nabla_{x,t} \phi}{DS_k} \ .
        \label{s-ds}
        \end{align}
       \item (Time Truncation of $N$) Let $I=\cup_i^K I_i$ be a decomposition of $I$ into consecutive
        intervals, and let $\chi_{I},\chi_{I_i}$ be the corresponding sharp time cutoffs.
        Then the following bounds hold (uniform in $K$):
        \begin{align}
         \lp{\chi_I F}{N} \ &\lesssim \ \lp{F}{N} \ ,
         \label{time_cutoff}\\
        \sum_i \lp{\chi_{I_i}F}{N}^2 \ &\lesssim \ \lp{F}{N}^2 \ . \label{N_fung}
        \end{align}
        Furthermore, for any Schwartz function $F$ the quantity $\lp{\chi_I F}{N}$
        is continuous in the endpoints of $I$.
        \item (Basic  $S$ and $N$ Relations) We have that:
        \begin{align}
                \lp{\phi^{(1)}_{k_1}\cdot\phi^{(2)}_{k_2}}{DS[I]} \
                &\lesssim \ 2^{(k_1-k_2)}\lp{\phi^{(1)}_{k_1}}{DS[I]}
                \cdot\lp{\phi^{(2)}_{k_2}}{S[I]}\ , \quad k_1<k_2-10 \ ,
                \label{S_DS_prod}\\
%                \lp{\phi}{\underline{E}[I]}
%                \ &\lesssim \ \lp{\Box \phi}{N[I]} +
%               \lp{\phi[0]}
%               {\dot{H}^1\times L^2} \ , \label{energy_est}\\
                \langle \phi , F_k \rangle  \ & \lesssim \ \lp{\phi}{DS}\cdot \lp{F_k}{N} \ ,
                \label{duality_est}\\
%                \lp{\partial_t \phi_k}{S_k[I]} \ &\lesssim \
%               2^k \lp{\phi_k}{S_k \cap \uX_k [I]} \ ,
%               \label{SX_est}\\
                \lp{ \phi_k}{S} \ &\lesssim \
                \lp{\nabla_{t,x}\phi_k}{X_{1}^{0,\frac{1}{2}}} \ , \label{XS_est}\\
                \lp{ F_k}{X_\infty^{0,-\frac{1}{2}}} \ &\lesssim \
                \lp{F_k}{N} \ . \label{XN_est}
        \end{align}
        \item ($L^q_t(L^r_x)$ and Disposability Estimates) We have that:
        \begin{align}
                \lp{Q_j \phi_k}{L^2_t(L^\infty_x)} \! &\lesssim \!
                2^{-(j-k)_+}
                2^{-\frac{1}{2}j}\lp{\phi_k}{S} ,
                \label{L2Linfty_est}\\
                \lp{\phi_k}{L^\infty_t(L^2_x)} +
                \lp{Q_{\leqslant j}\phi_k}{L^\infty_t(L^2_x)} +
                \lp{Q_j \phi_k}{L^\infty_t(L^2_x)} \! &\lesssim \!
                2^{-k}\lp{\phi_k}{S} \ , \label{L2Linfty_Q_est}\\
                \lp{Q_{\leqslant j}\phi_k}{L^\infty_t(L^\infty_x)} +
                \lp{Q_j \phi_k}{L^\infty_t(L^\infty_x)} \! &\lesssim
                \!
                \lp{\phi_k}{S} \ . \label{Linfty_Q_est}
        \end{align}
        \item (Fine Product Estimates) We have that:
        \begin{align}
            &\lp{\phi^{(1)}_{k_1}\cdot\phi^{(2)}_{k_2}}{L^2_t(L^2_x)} \ \lesssim \
            |\kappa|^\frac{1}{2}2^{-\frac{1}{2}k_1}
            \lp{\phi^{(1)}_{k_1}}{S[k_1,\kappa]}\cdot\sup_{\omega\in\kappa}
            \lp{\phi^{(2)}_{k_2}}{L^\infty_{t_\omega}(L^2_{x_\omega})}
             , \label{PW_product}\\
            &\lp{P_{<j-10}Q_{<j-10} \phi^{(1)}  \cdot\phi^{(2)}_{k_2}}{S[k_2;j]} \
            \lesssim \ \lp{\phi^{(1)}}{L^\infty_t(L^\infty_x)}
            \cdot\lp{\phi^{(2)}_{k_2}}{S} \ , \label{S_prod4}\\
            &\lp{\nabla_{t,x} P_k Q_j(
            \phi^{(1)}_{k_1} \cdot \phi_{k_2}^{(2)})
            }{X^{0,\frac{1}{2}}_1} \!  \lesssim \!
            2^{\delta(k-\max\{ k_i\})}2^{\delta(j-\min\{ k_i\})}
            \lp{\phi^{(1)}_{k_1}}{S} \lp{\phi^{(2)}_{k_2}}{S} \ ,
            \label{XS_X_est}\\
            &\lp{ P_k(Q_j F_{k_1}\cdot\phi_{k_2})}{N} \ \lesssim   \
            2^{\delta(k-\max\{ k_i\} )}2^{\delta(j-\min \{k_i\})}
            \lp{F_{k_1}}{ X^{0,-\frac{1}{2}}_\infty}\lp{\phi_{k_2}}{S}
            %&j\leqslant \min(k_i)
        \ , \label{XS_N_est}
        \end{align}
        where in estimates \eqref{XS_X_est}--\eqref{XS_N_est}
        we are assuming $j\leqslant \min\{k_i\}$.
\end{itemize}
\end{prop}
Estimates \eqref{dscut}--\eqref{S_DS_prod} are proved next. The rest
of the above bounds are standard, and with the exception of
\eqref{S_prod4} which is Lemma 9.1 in \cite{Tataru_WM1}, may be
found in \cite{Tao_WM}. For the convenience of the reader we give
the detailed citations (CMP copy). Estimate \eqref{duality_est} is
estimate (94) on p. 487. Estimate \eqref{XS_est} is Lemma 8 on p.
483. Estimate \eqref{XN_est} is Lemma 10 on p. 487. Estimates
\eqref{L2Linfty_est}--\eqref{Linfty_Q_est} are listed in estimates
(81)--(84) on p. 483. Estimate \eqref{PW_product} is (by duality)
the estimate in Step 2 on p. 479, and it also follows more or less
immediately by inspecting the third term on line \eqref{str_norm}
above. Estimate \eqref{XS_X_est} is Lemma 13 on p. 515. Estimate
\eqref{XS_N_est} is Lemma 12 on p. 501.

\ret

\begin{proof}[Proof of estimates \eqref{dscut}--\eqref{s-ds} and \eqref{S_decomp}]
Without any loss of generality we replace $\chi_I$ by
$\chi=\chi_{t<0}$. Our main observation here will be that the
multipliers $Q_j$ applied to $\chi$ act like time-frequency cutoffs
onto dyadic sets $|\tau|\sim 2^j$. For each of these we have the
Strichartz type estimate:
\begin{equation}
    \lp{Q_j \chi }{L^2_t(L^\infty_x)} \
    \lesssim \ 2^{-\frac{1}{2}j} \ . \label{chi_str}
\end{equation}
Therefore, one  can look upon the estimate \eqref{dscut} as some
version of the product bound \eqref{S_prod2}. We rescale to $k=0$,
and set $\phi_0=P_0\phi$.\ret

We begin with the proof of \eqref{dscut}. The $D\underline{S}$
bound in this estimate is immediate. Therefore we focus on proving the
$X^{s,b}$ and $S[0,\kappa]$ sum portions of the estimate. This is
split into
cases:\ret

\step{1}{Controlling the $X^{0,\frac{1}{2}}_\infty$ norm} Freezing
$Q_j$, our goal is to show that:
\begin{equation}
    \lp{ Q_j\big( \chi\cdot \phi_0 )}{L^2_t(L^2_x)} \ \lesssim \
    2^{-\frac{j}{2}}\lp{\phi_0}{DS} \ . \label{X_est2}
\end{equation}
We now split into subcases.\ret

\step{1.A}{$\chi$ at low modulation} In this case we look at the
contribution of the product $Q_j(Q_{<j-10} \chi\cdot \phi_0)$. We
may freely insert the multiplier $Q_{[j-5,j+5]}$ in front of
$\phi_0$. Then \eqref{X_est2} is immediate from $L^\infty$ control
of $Q_{<j-10} \chi$.\ret

\step{1.B}{$\chi$ at high modulation} In this case we'll rely on the
even stronger $L^2$ bound:
\begin{equation}
         \lp{Q_{\geqslant j-10}\chi\cdot \phi_0}{L^2_t(L^2_x)}
        \ \lesssim \ \lp{Q_{\geqslant j-10}\chi}{L^2_t(L^\infty_x)}
        \lp{\phi_0}{L^\infty_t(L^2_x)} \lesssim \ 2^{-\frac{1}{2}j}
        \lp{\phi_0}{S}
        \ , \label{step_2c_bnd}
\end{equation}
which results from summing over \eqref{chi_str}. In particular,
isolating the LHS of this last line at frequency $Q_j$ we have
\eqref{X_est2} for this term.\ret

\step{2}{Controlling the  $S[0;j]$ norms} Freezing $j < -20$ we need
to demonstrate:
\begin{equation}
    \lp{ Q_{<j} \big( \chi\cdot \phi_0 )}{S[0;j]} \ \lesssim \
    \lp{\phi_0}{DS} \ . \notag%\label{S0j_toshow}
\end{equation}

\step{2.A}{$\chi$ at low modulation} The contribution of
$Q_{<j-10}\chi\cdot Q_{<j }\phi_0$ is  bounded via estimate
\eqref{S_prod4}. Notice that $\chi$ is automatically at zero spatial
frequency.\ret

\step{2.B}{$\chi$ at high modulation} Adding over estimate
\eqref{step_2c_bnd} we have:
\begin{equation}
        \lp{Q_{<j}(Q_{\geqslant j-10}\chi\cdot \phi_0)}{X_1^{0,\frac{1}{2}}}
        \ \lesssim \
        \lp{\phi_0}{S} \ , \notag
\end{equation}
which is sufficient via the differentiated version of
\eqref{XS_est}.\ret

To wrap things up here, we need to demonstrate the bounds
\eqref{s-ds} and \eqref{S_decomp}. Beginning with $\phi \in S[I]$ we
consider an extension $\phi \in S$ with comparable norm.  Then by
\eqref{dscut} and \eqref{s=ds} we have the chain of inequalities:
\begin{equation}
    \lp{ \nabla_{x,t} \phi}{DS_k[I]}  \ \lesssim  \ \lp{\chi_I \nabla_{x,t} \phi_k}{DS_k}
     \ \lesssim  \ \lp{ \nabla_{x,t} \phi_k}{DS_k}  \ \approx  \ \lp{\phi}{S_k} \ . \notag
\end{equation}
It remains to prove the converse. We begin with the energy norm,
observing that for $\phi \in S[I]$, for $I=[-i_0,i_0]$, we have:
\begin{equation}
    \lp{ \phi_k[\pm i_0]}{\dot H^1 \times L^2}
     \ \lesssim  \ \lp{\nabla_{x,t} \phi_k}{DS_k[I]} \ . \notag
\end{equation}
We extend $\phi$ to $I^\pm=[\pm i_0,\pm\infty)$ as a solution to the homogeneous wave
equation with data $\phi[\pm i_0]$ and use \eqref{s=ds} to compute:
\begin{align}
    \lp{\phi}{S_k[I]}  \ &\lesssim \ \lp{\phi}{S_k}  \ \approx\
    \lp{\nabla_{x,t} \phi}{DS_k}  \ \lesssim  \ \lp{\chi_I \nabla_{x,t} \phi}{DS_k}
    + \lp{(1-\chi_I) \nabla_{x,t} \phi}{DS_k}\notag\\
    &\lesssim\  \lp{\nabla_{x,t} \phi}{DS_k[I]} + \|\phi_k[\pm i_0]\|_{\dot H^1 \times L^2}
    \notag\\
    &\lesssim\   \lp{\nabla_{x,t} \phi}{DS_k[I]} \ . \notag
\end{align}
The proof of \eqref{s-ds} is concluded.

Finally, we use \eqref{dscut} and \eqref{s-ds} to prove
\eqref{S_decomp}:
\begin{equation}
    \lp{ \phi}{S_k[I]}  \ \approx \ \lp{\nabla_{x,t} \phi}{DS_k[I]}
    \ \lesssim \ \sum_i \lp{\chi_{I_i} \nabla_{x,t} \phi}{DS_k}
    \ \lesssim \ \sum_i \lp{\phi}{S_k[I_i]} \ . \notag
\end{equation}
The proof is concluded.
\end{proof}\ret

%-------------------------------------------------------------------------
\begin{proof}[Proof of estimate \eqref{N_fung}]
Since:
\begin{equation}
    \lp{ F}{N}^2  \ \approx  \ \sum_k \lp{ P_k F}{N_k}^2 \ , \notag
\end{equation}
it suffices to show that the similar relation holds for the $N_k$
spaces:
\begin{equation}
    \sum_n \lp{ \chi_{I_n} F}{N_k}^2  \ \lesssim  \ \lp{ F }{N_k}^2 \ .
    \label{nksqsum}
\end{equation}
The space $N_k$ is an atomic space,
therefore is suffices to prove \eqref{nksqsum} for each atom.

\ret
%---------------------------------------------------------------------
\step{1}{$L^1_t( L^2_x)$ atoms} For these we directly have the stronger
relation:
\begin{equation}
    \sum_n \lp{ \chi_{I_n} F_k}{L^1_t( L^2_x)}  \ \lesssim\
    \lp{F_k }{L^1_t(L^2_x)} \ . \notag
\end{equation}

\ret
%---------------------------------------------------------------------
\step{2}{$\dot X^{0,-\frac12}_1$ atoms} For $F$ localized at
frequency $2^k$ we will prove the relation:
\begin{equation}
    \sum_n \lp{ \chi_{I_n} F_k }{L^1_t(L^2_x) + \dot X^{0,-\frac12}_1 }^2
     \ \lesssim  \ \lp{F_k}{\dot X^{0,-\frac12}_1}^2 \ . \notag
\end{equation}
Without any restriction in generality we can assume that $F_k$ is
also localized in modulation at $2^j$. By rescaling we can take $j =
0$. At modulation $1$ the $\dot X^{0,-\frac12}_1$ is equivalent to
the $L^2_t(L^2_x)$ norm. Then the last bound would follow from the
stronger estimate:
\begin{equation}
    \sum_n \lp{ Q_{< -4} (\chi_{I_n} Q_0 F_k)}{L^1_t( L^2_x) }^2 \!+\!
    \lp{ Q_{> -4} (\chi_{I_n} Q_0F_k) }{ L^2_t(L^2_x) }^2
      \lesssim   \lp{ Q_0F_k}{L^2_t(L^2_x)}^2 \ .
    \label{l2sqsum}
\end{equation}
We trivially have:
\begin{equation}
    \sum_n  \lp{\chi_{I_n} Q_0 F_k }{ L^2_t(L^2_x) }^2
    \ \lesssim \ \lp{ Q_0 F_k }{L^2_t(L^2_x)}^2 \ , \notag
\end{equation}
therefore it remains to prove the $L^1_t(L^2_x)$ bound on line
\eqref{l2sqsum}. We do this in two cases:\ret

%----------------------------------------------------------
\step{2.A}{Small intervals} We parse the collection of intervals
$I_n$ into two sub-collections, intervals $J_n$ such that $|J_n|\geqslant 1$,
and intervals $K_n$ such that $|K_n|<1$. In the latter can we may
drop the outer $Q_{<-4}$ and simply use H\"older's inequality to estimate:
\begin{equation}
    \sum_n  \lp{\chi_{K_n} Q_0 F_k }{ L^1_t(L^2_x) }^2
    \ \lesssim \ \sum_n \lp{ \chi_{K_n} Q_0 F_k }{L^2_t(L^2_x)}^2 \ , \notag
\end{equation}
so the estimate follows as above.\ret

%----------------------------------------------------------
\step{2.B}{Large intervals} In this case we break the first term on
LHS \eqref{l2sqsum} up as follows:
\begin{align}
    \sum_m \lp{ Q_{< -4} (\chi_{J_m} Q_0 F_k)}{L^1_t( L^2_x) }^2
    \ &= \ \sum_m \lp{ Q_{< -4} (Q_{[-10,10]}\chi_{J_m}\cdot Q_0 F_k)}{L^1_t( L^2_x) }^2 \label{kgp}\\
    \ &\lesssim \ \sum_m \lp{ Q_{[-10,10]}\chi_{J_m}\cdot Q_0 F_k}{L^1_t( L^2_x) }^2.
 \notag
\end{align}
Denoting $J_m = [a_m,b_m]$, for $Q_{[-10,10]}\chi_{J_m}$
we have the pointwise bounds:
\begin{equation}
    |Q_{[-10,10]}\chi_{J_m}(t)| \ \leqslant \  (1+|t-a_m|)^{-N} + (1+|t-b_m|)^{-N} \ . \notag
\end{equation}
Hence by Cauchy-Schwartz we obtain:
\begin{align}
    \hbox{L.H.S.}\eqref{kgp} \! &\lesssim \!
    \sum_m \lp{ (1+|t-a_m|)^{-\frac{N}2}   Q_0 F_k}{L^2_t(L^2_{x}) }^2
    \!+\! \lp{ (1+|t-b_m|)^{-\frac{N}2}   Q_0 F_k}{ L^2_t(L^2_{x})}^2\notag\\
    &\lesssim \ \int_{\R} \sum_{m}\left((1+|t-a_m|)^{-N} + (1+|t-b_m|)^{-{N}}\right)
    \lp{Q_0 F_k(t)}{L^2_x}^2 dt \ . \notag
\end{align}
Since the intervals $J_m$ are disjoint and of size at least $1$,
the last sum above is bounded by $\lesssim 1$,
therefore we obtain:
\begin{equation}
    \hbox{L.H.S.}\eqref{kgp}  \ \lesssim  \ \lp{Q_0 F_k}{L^2_t(L^2_x)}^2 \ . \notag
\end{equation}
\ret

%----------------------------------------------------------------
\step{3}{$N\!F\!A$ atoms} In this case we can express $F_k$ as:
\begin{equation}
    F_k \ =  \ F_k^+ + F_k^- \ = \ \sum_{\pm,\kappa} F_{k,\pm\kappa}^\pm \ , \notag
\end{equation}
where  $F_{k,\pm\kappa}^\pm$ is supported in the wedge carved by the
multiplier $Q_{<k-2j}^\pm P_{k,\pm\kappa}$, with $j > 10$, and
furthermore:
\begin{equation}
    \sum_{\kappa} \lp{F_{k,\kappa}^\pm}{N\!F\!A[\pm\kappa]}^2  \ \leqslant  \ 1 \ . \notag
\end{equation}
Without loss of generality we may assume we are in the $+$ case, and
we rescale to $k=2j$, and so in particular $k>20$.

By summing over \eqref{XN_est} we have the $L^2_t(L^2_x)$ bound:
\begin{equation}
    \lp{ F_{k}^+}{L^2_t(L^2_x)} \ \lesssim \ 1 \ .
    \label{NFA_L2}
\end{equation}
The $N\!F\!A[\kappa]$ norms are translation invariant, and are
defined using characteristic $L^1_{t_\omega}( L^2_{x_\omega})$
norms. Thus they directly satisfy the inequality:
\begin{equation}
    \sum_n     \lp{ \chi_{I_n} F_{k,\kappa}^+}{N\!F\!A[\kappa]}^2 \ \leqslant \
    \lp{F_{k,\kappa}^+}{N\!F\!A[\kappa]}^2  \ . \label{l2nfa}
\end{equation}
We write:
\begin{equation}
    \chi_{I_n} F_k^+  \ =  \ Q_{>0} (\chi_{I_n} F_k^+) + Q_{<0}^- (\chi_{I_n} F_k^+) +
    \sum_\kappa  Q_{<0}^+( \chi_{I_n}
    F_{k,\kappa}^+) \ , \notag
\end{equation}
and estimate the first component in $\dot X^{0,-\frac12}_1$, the
second in $L^1_t(L^2_x)$, and the third in $N\!F\!A\!$. We have from
line \eqref{NFA_L2}:
\begin{equation}
    \sum_n \lp{ Q_{>0}( \chi_{I_n} F_k^+)}{\dot X^{0,-\frac12}_1}^2  \ \lesssim\
    \sum_n  \lp{  \chi_{I_n} F_k^+}{L^2_t(L^2_x)}^2  \ \lesssim\
     \lp{   F_k^+}{L^2_t(L^2_x)}^2 \ \lesssim\ 1 \ . \notag
\end{equation}
Next, using the restriction on the Fourier support of $F^+_k$, we
have for any single interval the bound:
\begin{equation}
        \lp{Q_{<0}^- (\chi_{J_n} F_k^+)}{L^1_t(L^2_x)}
        \ \lesssim \ \lp{Q_{[-20,20]}\chi_{J_n}\cdot
        F_k^+}{L^1_t(L^2_x)} \ . \notag
\end{equation}
Using this and \eqref{NFA_L2}, one may proceed as in \textbf{Step
2.A} and \textbf{Step 2.B} above. On the other hand, by
\eqref{l2nfa} and the disposability of $Q_{<0}^+$ on the Fourier
support of the multiplier $P_{k,\kappa}$ we have:
\begin{equation}
    \sum_{n,\kappa} \lp{ Q_{<0}^+( \chi_{I_n} F_{k,\kappa}^+)}{N\!F\!A[\kappa]}^2
    \lesssim \sum_{n,\kappa}  \lp{\chi_{I_n} F_{k,\kappa}^+}{N\!F\!A[\kappa]}^2
    \lesssim \sum_{\kappa} \lp{ F_{k,\kappa}^+}{N\!F\!A[\kappa]}^2 \lesssim 1 \ . \notag
\end{equation}
\end{proof}\ret

%-------------------------------------------------------------------------------
\begin{proof}[Proof of estimate \eqref{S_DS_prod}]
This is a minor variation of \eqref{N_prod_est1}, and the proof is
similar to that of estimate \eqref{dscut} above. We rescaling to
$k_2=0$, discard $I$, and set
$\lp{\phi^{(1)}_{k_1}}{DS}=\lp{\phi^{(2)}_{0}}{S}=1$. Using the fact
that:
\begin{equation}
        \lp{Q_{<k_1+10}\phi^{(1)}_{k_1}}{S} \ \lesssim\
        2^{k_1}\lp{\phi^{(1)}_{k_1}}{DS} \ , \qquad
        \lp{P_0 f}{DS} \ \lesssim \ \lp{f}{S} \ , \notag
\end{equation}
along with the usual $S$ algebra estimate \eqref{S_prod2}, we
control the low modulation contribution of the first factor.
Furthermore, it is always possible to gain in the $D\underline{S}$
component by using Bernstein's inequality and energy estimates for
the high frequency factor.

We are reduced to bounding the contribution of
$Q_{>k_1+10}\phi^{(1)}_{k_1}\cdot\phi^{(2)}_{0}$.
%\begin{equation}
%            \lp{Q_{>k_1+10}\phi^{(1)}_{k_1}\cdot\phi^{(2)}_{0}}{X_\infty^{0,\frac{1}{2}}}
%            +
%            \sup_{j<-20}\lp{Q_{>k_1+10}\phi^{(1)}_{k_1}\cdot\phi^{(2)}_{0}}{S[0;j]}\!
%            \lesssim \! 2^{k_1}\lp{\phi^{(1)}_{k_1}}{DS}
%            \cdot\lp{\phi^{(2)}_{0}}{S} \ . \notag
%\end{equation}
As a general tool we have the $L^2$ bound:
\begin{equation}
         \lp{ Q_{>j-10}\phi^{(1)}_{k_1}\cdot\phi^{(2)}_{0})}{L^2_t(L^2_x)} \
            \lesssim \
            \lp{Q_{>j-10}\phi^{(1)}_{k_1}}{L^2_t(L^\infty_x)}
            \lp{\phi^{(2)}_{0}}{L^\infty_t(L^2_x)} \ \lesssim  \
            2^{k_1}2^{-\frac{1}{2}j} \ . \notag
\end{equation}
In particular, via the differentiated version of \eqref{XS_est}, if
the output modulation is $j<k_1$ we have both:
\begin{equation}
            \lp{Q_j(Q_{>k_1+10}\phi^{(1)}_{k_1}\cdot\phi^{(2)}_{0})}{X_\infty^{0,\frac{1}{2}}}
            + \lp{Q_{>k_1+10}\phi^{(1)}_{k_1}\cdot\phi^{(2)}_{0}}{S[0;j]}\
            \lesssim \ 2^{k_1} \ . \notag
\end{equation}

On the other hand, if the output modulation is $j>k_1$, then by
again using the above general $L^2$ estimate,  it suffices to show:
\begin{equation}
            \lp{Q_j(Q_{[k_1+10,j-10]}\phi^{(1)}_{k_1}\cdot\phi^{(2)}_{0})}{X_\infty^{0,\frac{1}{2}}}
            +
            \lp{Q_{[k_1+10,j-10]}\phi^{(1)}_{k_1}\cdot\phi^{(2)}_{0}}{S[0;j]}\!
            \lesssim \! \lp{\phi^{(1)}_{k_1}}{L^\infty_t(L^\infty_x)} \ , \notag
\end{equation}
and then conclude via an application of Bernstein's inequality.
%we are reduced to demonstrating this last bound for the localized
%product $Q_{[k_1+10,j-10]}\phi^{(1)}_{k_1}\cdot\phi^{(2)}_{0}$.
%\begin{equation}
%            \lp{Q_j(Q_{[k_1+10,j-10]}\phi^{(1)}_{k_1}\cdot\phi^{(2)}_{0})}{X_\infty^{0,\frac{1}{2}}}
%            + \lp{Q_{[k_1+10,j-10]}\phi^{(1)}_{k_1}\cdot\phi^{(2)}_{0}}{S[0;j]}\
%            \lesssim \ 2^{k_1} \ . \notag
%\end{equation}
For the first term on the LHS, we may freely insert a
$Q_{[j-5,j+5]}$ multiplier in front of the second factor, which
suffices. For the second term, we directly use \eqref{S_prod4}.
\end{proof}

%--------------------------------------------------------------------------------

\subsection{Extension and restriction for $S$ and $N$ functions}
In the sequel we will build up estimates through an iterative
process by which we first prove bounds in weak spaces (such as $P_k
L^\infty$, $\uX$, and $\underline{E}$), and then show that these may
be used in conjunction with bootstrapping to establish uniform
bounds in much stronger spaces (such as $S$). This process
unfortunately leads to some technical difficulties regarding
compatibility of extensions in various norms. To tame this
difficulty, we will make use of a variable but universal extension
process. Because this feature is more of a technicality in our
proof, we state here for the convenience of the reader where such
extensions are necessary in the sequel:
\begin{itemize}
    \item The primary use of compatible extensions is in the proof
    or Proposition \ref{p:matchfreq}, most importantly in the proof
    of estimate \eqref{bal_l2_N}. Here we are forced to use several
    norms simultaneously in a single estimate that involves space-time
    frequency cutoffs. As will become apparent soon, in such a
    situation
    choosing extensions needs to be done carefully because it
    is not immediate that this can be done in a way that retains
    smallness of the various component norms.
    \item Universal extensions are also used in a key way in the
    proof of Proposition \ref{p:para} because we need to know that
    extensions still enjoy good characteristic energy estimates when
    these estimates are only known on a finite interval. This
    extended control needs to be used in conjunction with $S$ norm
    control in estimates requiring space-time frequency cutoffs
    (see Lemma \ref{para_trilin_lem} in Section \ref{s:lin_sect}).
    \item A secondary use of compatible extensions occurs because we
    do not include $\underline{X}$ as a component of $S$ defined
    above. Doing this allows us to quote standard product estimates from
    \cite{Tao_WM} modulo physical space Strichartz components. The
    price one pays is that $\underline{X}$ bounds are established
    separately, and one then needs to included this a-posteriori
    into extension estimates. For example, this feature is used at the
    beginning of the proof of Proposition \ref{propphiu} to extend
    the connection $B$ with good $S$ and $\underline{X}$ bounds.
\end{itemize}\ret

\begin{prop}[Existence of $S$ Extensions/Restrictions]\label{ext_prop}
 Let $\phi$ be any affinely Schwartz function defined on an interval
$I=[-i_0,i_0]$.
\begin{itemize}
%        \item (Seed Bounds) Define $E=\lp{\phi[0]}{\dot{H}^1\times L^2}$.
%Then there exists an interval $I_0\subseteq I$ containing $t=0$, depending on the profile of $\phi$, such that:
%        \begin{align}
%                \lp{\Phi^{I_0,\eta}}{S} \ &\lesssim_E \ 1 \ ,
%                \label{initial_bnds1}\\
%                \lp{\Phi^{I_0,\eta}}{S_k} \ &\lesssim_E \ c_k \ ,
 %               \label{initial_bnds2}
 %       \end{align}
 %where $\{c_k\}$ is any $(\delta,\Delta)$-admissible $\dot{H}^1\times L^2$
 %frequency envelope for $\phi[0]$.\ret
        \item (Canonical extension) For every $0< \eta\leqslant 1$ there
exists a canonical extension
$\Phi^{I,\eta}$ which is compactly supported in time and for which
the following estimates are true:
        \begin{align}
                \lp{P_k\Phi^{I,\eta}}{S} \ &\lesssim \ \lp{P_k\phi}{S[I]} \ , \label{can_S}\\
        \lp{ P_k\Box\Phi^{I,\eta}}{N} \ &\lesssim \ \lp{P_k\phi}{\underline{E}[I]}
        + \lp{P_k\Box\phi}{N[I]}\ , \label{can_N}\\
%       \llp{P_{k'}\Box\Phi^{I,\eta}}{\W} \ &\lesssim \ (1+|\ln(\eta)|)\llp{P_k\phi}{\W[I]} \ , \label{can_W}\\
        \lp{P_k\Phi^{I,\eta}}{L^\infty_t(L^\infty_x)} \ &\lesssim \ \eta^{-\frac{5}{2}}
            \lp{P_k\phi}{L^\infty_t(L^\infty_x)[I]}+ \eta^\frac{1}{2}  \lp{P_k\phi}{\uX[I]}
        \ , \  |I|\geqslant 2^{-k}\eta^2  \ , \label{can_Linfty}\\
        \lp{P_k\Phi^{I,\eta}}{\uX} \ &\lesssim \ \eta^\frac{1}{2}\lp{P_k\phi}{\underline{E}[I]}
            + \lp{P_k\phi}{\uX[I]} \ , \label{can_X}\\
                \lp{P_k \Phi^{I,\eta}}{\underline{E}} \ &\lesssim \ \lp{P_k\phi}{\underline{E}[I]} \ , \label{can_E}\\
%                \lp{P_{<k}\Phi^{I,\eta}}{L^\infty_t(L^\infty_x)} \ &\lesssim \
%           \eta^{-?}\lp{P_{<k}\phi}{L^\infty(I)} \ , \label{canLinfty2}\\
        \lp{ P_k\Box\Phi^{I,\eta}\cdot\psi_j}{N}  &\lesssim  2^{k-j}
        \lp{P_k\phi}{\underline{E}[I]}\lp{\psi_j}{S}
        + \lp{P_k\Box\phi\cdot\psi_j}{N[I]}\ ,  \label{can_N_prod}
        \end{align}
    where the last bound holds under the additional condition that $j>k+10$.\ret
\item (Secondary extension) For every $0< \eta\leqslant 1$ there
exists an  extension $\td{ \Phi}^{I,\eta}$ which is compactly
supported in time and such that \eqref{can_S},\eqref{can_N} hold.
Furthermore, for this extension the following improvement of
\eqref{can_Linfty} is valid:
\begin{equation}
    \lp{P_k\td{\Phi}^{I,\eta}}{L^\infty_t(L^\infty_x)} \ \lesssim \ \eta^{-\frac12}
    \lp{P_k\phi}{L^\infty_t(L^\infty_x)[I]}+
    \eta \lp{P_k\phi}{\uE[I]} \ .
    \label{can_Linftyb}
\end{equation}
% \begin{align}
%                \lp{P_k\tilde \Phi^{I,\eta}}{S} \ &\lesssim \ \lp{P_k\phi}{S[I]} \ , \label{can_Sb}\\
%  \lp{ P_k\Box\tilde \Phi^{I,\eta}}{N} \ &\lesssim \ \lp{P_k\phi}{\underline{E}[I]}
%       + \lp{P_k\Box\phi}{N[I]}\ , \label{can_Nb} \\
% \lp{P_k\tilde \Phi^{I,\eta}}{L^\infty_t(L^\infty_x)} \ &\lesssim \ \eta^{-\frac12}
%        \lp{P_k\phi}{L^\infty_t(L^\infty_x)[I]}+
%\eta \lp{P_k\phi}{\uE[I]}
% \label{can_Linftyb}
%\end{align}
%         \item (Continuity) These extensions are continuous in the sense that:
%         \begin{align}
%%                \lim_s \lp{\Phi^{J_s,\eta}}{S\setminus L^\infty} \
%%               &\to \ \lp{\Phi^{I,\eta}}{S\setminus L^\infty} \ , \label{can_cont1}\\
%                 \lim_s \lp{\Phi^{J_s,\eta}}{S_c} \
%                 &\to \ \lp{\Phi^{I,\eta}}{S_c} \ , \label{can_cont2}
%         \end{align}
%         for any frequency envelope $\{c_k\}$, whenever
%         $\phi$ is defined on $I$ with $I\subseteq J_s$ and $\lim_s J_s=I$.
\end{itemize}
\end{prop}\ret

The canonical extension above will be used most of the time. Its only
disadvantage is that in order to control the $L^\infty_t(L^\infty_x)$
norm of this extension we need to also control the $\uX$ norm.
In the rare (single) case where this is missing, we use the secondary extension.

%-------------------------------------------------------------------------

\begin{proof}[Proof of Proposition \ref{ext_prop}]
%For the seed bounds \eqref{initial_bnds1} and \eqref{initial_bnds2}
%we compare $\phi$ with the solution $\Phi$ to
%\[
% \Box \Phi = 0, \qquad \Phi[0] = \phi[0]
%\]
%For $t_0 > 0$ we can use \eqref{energy_est} to write
%\[
% \| \phi\|_{S[-t_0,t_0]} \lesssim
%\| \Phi\|_{S[-t_0,t_0]}+ \| \phi-\Phi\|_{S[-t_0,t_0]}
%\lesssim \|\phi[0]\|_{\dot H^1 \times L^2} + \|\Box \phi\|_{L^1_t L^2_x[-t_0,t_0]}
%\]
%Since $\phi$ is Schwartz, the last term on the right goes to $0$ as
%$t_0$ goes to $0$, and \eqref{initial_bnds1}. The proof of
%\eqref{initial_bnds2} is similar.
The canonical extension will be defined dyadically for each $\phi_k$. By rescaling we only work with $k=0$. \ret

\step{1}{The canonical extension and estimates
\eqref{can_S}--\eqref{can_N}, \eqref{can_X}--\eqref{can_N_prod}} The
obvious candidate $\Phi^I$ for the extension is obtained by solving
the homogeneous wave equation to the left of $-i_0$ and to the right
of $i_0$, with Cauchy data $P_0\phi[-i_0]$, respectively
$P_0\phi[i_0]$. Denoting the complement of $I$ by $I^- \cup I^+$, we
have:
\begin{equation}
    \Box \Phi^I  \ =  \ 0 \ ,
    \qquad \Phi^I[\pm i_0]  \ =  \ P_0\phi[\pm i_0]  \ , \quad \text{in }
    I^{\pm} \ . \notag
\end{equation}

It is relatively easy to verify that the extension $\Phi^I$ satisfies
all the properties  \eqref{can_S}--\eqref{can_N}, \eqref{can_X}--\eqref{can_N_prod}.
However, there is a core issue with
\eqref{can_Linfty}, as this bound can easily fail because nonconcentration
at time $\pm i_0$, say, does not guarantee nonconcentration at all later times.
To avoid this problem, we truncate $\Phi^I$ outside a compact set and define:
\begin{equation}
    \Phi^{I,\eta}  \ =  \ \chi_I^\eta \Phi^I \ , \notag
\end{equation}
where $\chi^\eta_I$ is a smooth cutoff with
$|\partial_t^k\chi^\eta_I|\lesssim \eta^{k}$, such that
$\chi^\eta_I\equiv 1$ on $I$ and vanishing outside of the extended
interval  $\td{I}=[-i_0-\eta^{-1},i_0+\eta^{-1}]$.  Furthermore,
in $I^\pm$ we have the identity:
\begin{equation}
        \Box \Phi^{I,\eta} \ = \ 2  \partial_t(\chi^\eta_{I})\cdot\partial_t \Phi^I
        +  \partial_t^2 (\chi^\eta_{I})\cdot \Phi^I \ . \notag
\end{equation}
This allows us to estimate:
\begin{equation}
    \lp{ P_0 \Box \Phi^{I,\eta}}{L^1_t( L^2_x)[I^\pm]}  \ \lesssim  \
    \lp{ P_0 \phi[\pm i_0]}{\dot H^1 \times L^2} \ ,
    \label{ipmrhs}
\end{equation}
which in turn leads to:
\begin{equation}
    \lp{ P_0 \Phi^{I,\eta}}{S[I^\pm]} +\lp{P_0 \Box \Phi^{I,\eta}}{N[I^\pm]}
     \ \lesssim  \ \lp{P_0 \phi[\pm i_0]}{ \dot H^1 \times L^2} \ . \notag
\end{equation}

Then the bound \eqref{can_S} follows from \eqref{S_decomp},
while \eqref{can_E} follows from energy estimates for $\Phi^{I,\eta}$
in $I^{\pm}$. The bound \eqref{can_N} is also straightforward,
while for \eqref{can_X} we need to compute:
\begin{equation}
    \lp{ P_0 \Box \Phi^{I,\eta}}{L^2_t L^2_x[I^\pm]}  \ \lesssim\
    \eta^\frac12 \lp{P_0 \phi[\pm i_0]}{\dot H^1 \times L^2}  \ . \notag
\end{equation}
To prove \eqref{can_N_prod} we use Bernstein to estimate:
\begin{align}
    \lp{ P_0 \Box\Phi^{I,\eta}\cdot\psi_j}{L^1_t( L^2_x)[I^\pm]}
    &\lesssim \ \lp{P_0 \Box\Phi^{I,\eta}}{L^1_t( L^\infty_x)[I^\pm]}
    \cdot\lp{\psi_j}{L^\infty_t( L^2_x)[I^\pm]} \notag\\
    &\lesssim \ \lp{ P_0 \Box\Phi^{I,\eta}}{L^1_t( L^2_x)[I^\pm]}\cdot 2^{-j}
    \lp{P_0 \psi_j}{S} \ , \notag
\end{align}
and conclude with \eqref{ipmrhs}.

\ret
%----------------------------------------------------------------
\step{2}{The $L^\infty_t(L^\infty_x)$ estimate \eqref{can_Linfty}}
We now turn our attention to the most interesting part, namely
\eqref{can_Linfty}.
The desired bound follows from a reverse dispersive estimate for the 2D
wave equation:
\begin{equation}
    \lp{ e^{\pm i t|D_x|} P_0 f }{L^\infty_x}
     \ \lesssim  \ \sqrt{1+t}\,  \lp{ P_0 f }{L^\infty_x} \ , \notag
\end{equation}
provided that we can first establish the ``elliptic'' estimate
(setting $P_0\phi=\phi_0$):
\begin{equation}
    \lp{\partial_t\phi_0}{L^\infty_t(L^\infty_x)[I]} \ \lesssim \
    \eta^{-2}\lp{\phi_0}{L^\infty_t(L^\infty_x)[I]}
    + \eta \lp{\phi_0}{\uX[I]} \ , \label{simp_ell}
\end{equation}
provided that $|I|\geqslant \eta^2$. Without loss of generality we may assume we are
in the worst case scenario $|I|=\eta^2$. We begin with the Poincare type inequality:
\begin{equation}
    \lp{\partial_t\phi_0 - (\partial_t\phi_0)^{av}}{L^\infty_t[I]}
    \ \lesssim \ \eta\lp{\partial_t^2\phi_0}{L^2_t[I]} \ , \notag
\end{equation}
where $(\partial_t\phi_0)^{av}=\eta^{-2}\int_I \partial_t\phi_0dt$, so
in particular:
\begin{equation}
     \lp{(\partial_t\phi_0)^{av}}{L^\infty_x[I]}
    \ \lesssim \ \eta^{-2}\lp{\phi_0}{L^\infty_t(L^\infty_x)[I]} \ . \notag
\end{equation}
Therefore, taking the $\sup$ over all space in the second to last
line above we have:
\begin{equation}
    \lp{\partial_t\phi_0}{L^\infty_t(L^\infty_x)[I]} \ \lesssim \
    \eta^{-2}\lp{\phi_0}{L^\infty_t(L^\infty_x)[I]}
    + \eta\lp{\Box\phi_0}{L^2_t(L^\infty_x)[I]}
    + \eta \lp{\Delta \phi_0}{L^2_t(L^\infty_x)[I]} \ . \notag
\end{equation}
The proof of \eqref{simp_ell} is now concluded via Bernstein in
space for the second term on the RHS above, and Cauchy-Schwartz in
time along with the fact that $\Delta P_0$ is bounded on
$L^\infty_x$ to control the third.

\ret
%---------------------------------------------------------------------------------------

\step{3}{The Secondary Extension}
We next turn our attention to the second extension. Again we set
$k=0$. The additional difficulty we face here is that we
 no longer have pointwise bounds for  $\partial_t P_0\phi(\pm i_0)$.
We split the function $\Phi^I$ outside $I$ into two parts:
\begin{equation}
    \Phi^I  \ =  \ \Phi^I_0 + \Phi^I_1 \ , \notag
\end{equation}
corresponding to the two different components of its
Cauchy data at $\pm i_0$:
\begin{equation}
    \Box  \Phi^I_0  \ =  \ 0 \ ,
    \qquad \Phi^I_0[\pm i_0]  \ =  \ (P_0\phi(\pm i_0),0) \qquad
    \text{in  } I^{\pm} \ , \notag
\end{equation}
respectively:
\begin{equation}
    \Box  \Phi^I_1  \ =  \ 0 \ ,
    \qquad \Phi^I_1[\pm i_0]  \ =  \
    (0,\partial_t P_0\phi(\pm i_0)) \quad  \ \text{in  } I^{\pm} \ . \notag
\end{equation}
Then we define the extension $\td{\Phi}^{I,\eta}$ by truncating
the two components on different scales:
\begin{equation}
    \td{\Phi}^{I,\eta}  \ =  \ \chi_I^\eta \Phi^I_0 + \chi_I^{1/\eta} \Phi^I_1 \ . \notag
\end{equation}
For the first component we argue as before. For the second,
we begin with a fixed time  $L^2$ bound:
\begin{equation}
    \lp{  P_0 \Phi^I_1(t)}{L^2_x}  \ \lesssim  \
    |t\mp i_0| \cdot\lp{  \partial_t P_0\phi(\pm i_0)}{L^2_x}
    \qquad \text{in  } I^{\pm}
    \label{daw} \ ,
\end{equation}
which follows at once from   integrating the quantity $\partial_t
P_0 \Phi^I_1$ and energy estimates. This leads to:
\begin{equation}
    \lp{ P_0  \Box (\chi_I^{1/\eta}\Phi^I_1)}{L^1_t( L^2_x)[I^\pm]}
     \ \lesssim\  \lp{\partial_t P_0 \phi(\pm i_0)}{L^2_x} \ , \notag
\end{equation}
which helps us establish  bounds of the type \eqref{can_S}--\eqref{can_N}.
On the other hand, the improved pointwise bound \eqref{can_Linftyb} follows
simply by using Bernstein's inequality in \eqref{daw} to give:
\begin{equation}
    \lp{P_0 \Phi^I_1(t)}{L^\infty_x}  \ \lesssim  \
    |t\mp i_0|\cdot \lp{  \partial_t P_0\phi(\pm i_0)}{L^2_x}
    \qquad \text{in  } I^{\pm} \ . \notag
\end{equation}
The proof of the proposition is thus concluded.
\end{proof}\ret

%-------------------------------------------------------------------------

\subsection{Strichartz and Wolff type bounds}
In this section we prove the estimate \eqref{energy_est} for the
$\underline{S}$ component on line \eqref{phys_str}, as well as a key
$L^2$  bilinear estimate for transverse waves which takes advantage
of the small energy dispersion. The tools we use for these purpose
are the  $V^p_{\pm|D_x|}$ and $U^p_{\pm|D_x|}$ spaces associated to the two
half-wave evolutions.

Precisely, $V^p_{\pm|D_x|}$ is the space of right continuous $L^2_x$ valued
functions with bounded $p$-variation along the half-wave flow:
\begin{equation}
        \lp{ u}{V^p_{\pm|D_x|}}  \ =  \ \lp{e^{\mp it|D_x|}
        u(t)}{V^p(L^2_x)} \ , \notag
\end{equation}
or in expanded form:
\begin{equation}
        \lp{ u}{V^p_{\pm|D_x|}}  \ :=  \ \lp{u}{L^\infty_t( L^2_x)}^p
        + \sup_{t_k \nearrow}
        \sum_{k \in \Z} \lp{u(t_{k+1}) - e^{\pm i(t_{k+1} - t_k)|D_x|}
        u(t_k)}{L^2_x}^p \ , \notag
\end{equation}
where the supremum is taken over all increasing sequences $t_k$. We
note that if $p < \infty$ then $V^p$ functions can have at most
countably many discontinuities as $L^2_x$ valued functions.

On the other hand the slightly smaller space $U^p_{\pm|D_x|}$ is
defined as the atomic space generated by a family $\mathcal A_p$ of atoms
$a$ which have the form:
\begin{equation}
        a(t)  \ =  \ e^{\pm it|D_x|} \sum_k 1_{[t_k,t_{k+1})} u^{(k)} \ ,
        \notag %\label{atom}
\end{equation}
where the sequence $t_k$ is increasing and:
\begin{equation}
        \sum_k \lp{u^{(k)}}{L^2_x}^p  \ \leqslant\  1 \ . \notag
\end{equation}
Precisely, we have:
\begin{equation}
        U^p_{\pm|D_x|} \ = \ \{ u = \sum c_k a_k;\ \sum_k |c_k| < \infty, a_k \in
        \mathcal A_p\} \ . \notag
\end{equation}
The above sum converges uniformly in $L^2_x$; it also converges in
the stronger $V^p_{\pm|D_x|}$ topology. The $U^p_{\pm|D_x|}$ norm is
defined by:
\begin{equation}
        \lp{ u}{U^p_{\pm|D_x|}} \ := \ \inf \{ \sum_k |c_k|;\  u = \sum_k c_k a_k,\ a_k
        \in \mathcal A_p\} \ . \notag
\end{equation}
These spaces are related as follows:
\begin{equation}
    U^p_{\pm|D_x|}  \ \subset  \ V^p_{\pm|D_x|} \ \subset \ U^q_{\pm|D_x|} \ ,
    \qquad 1  \ \leqslant  \ p  \ <  \ q
    \ \leqslant \ \infty  \ . \label{upemb}
\end{equation}
The first inclusion is straightforward. The second is not, and plays
a role similar to the Christ-Kiselev lemma. These spaces were first introduced
in unpublished work of the second author, and have proved their usefulness
as scale invariant substitutes of $X^{s,\frac12}$ type spaces in several problems,
see \cite{KTucp}, \cite{KTnls}, \cite{HHK}, \cite{BIKT}.
\ret

We use these spaces first in the context of the Strichartz
estimates, which for frequency localized homogeneous half-waves can
be expressed in the form:
\begin{equation}
   \lp{e^{\pm it |D_x|}  u_k }{L^q_t( L^r_x)}  \ \lesssim  \ 2^{-(\frac1q + \frac2r -1)k}
    \lp{ u_k}{L^2_x} \ , \qquad \frac{2}q + \frac{1}r \leqslant \frac12 \
    . \notag
\end{equation}
Applying this bound for each segment in $U^q_{\pm|D_x|}$ atoms, one
directly obtains embeddings of these spaces into Strichartz spaces:

\begin{lem}
The following estimates hold:
\begin{equation}
    \lp{ \phi_k}{L^q_t( L^r_x)}  \ \lesssim  \
    2^{-(\frac1q + \frac2r -1)k} \lp{ \phi_k}{U^q_{\pm|D_x|}},
    \qquad
    \frac{2}q + \frac{1}r \leqslant \frac12 \ .
    \label{upst}
\end{equation}
\end{lem}\ret

The second place where these spaces come into play is in the context
of bilinear $L^2_t(L^2_x)$ estimates for transversal waves. The
classical estimate here (see for instance \cite{KS}) has the form:
\begin{equation}
   \lp{e^{\pm it |D_x|}  u_{k_1}^{(1)}
    \cdot e^{\pm it |D_x|}  u_{k_2}^{(2)}}{L^2_t(L^2_x)}
     \ \lesssim  \ 2^{\frac12 \min\{k_1,k_2\}} \theta^{-\frac12}
    \lp{ u_{k_1}^{(1)}}{L^2_x} \lp{u_{k_2}^{(2)}}{L^2_x} \ ,
\end{equation}
provided the $u_{k_i}^{(i)}$ have angular separation in frequency,
namely $|\theta_1-\theta_2|>\theta$ in the $(++)$ or $(--)$ cases,
and $|\pi+(\theta_1-\theta_2)|>\theta$ in the $(+-)$ or $(-+)$
cases.  In subsequent work, Wolff \cite{Wolff} was able to replace
the $L^2_t(L^2_x)$ bound on the left with $L^p_t(L^p_x)$ for $p >
\frac53$. The endpoint $p = \frac53$ was later obtained by Tao
\cite{Taobi}. Our aim here is to first use the Wolff-Tao estimate to
strengthen the classical $L^2_t(L^2_x)$ bound in a way which takes
advantage of the small energy dispersion, and then phrase it in the
set-up of the $V^2_{\pm|D_x|}$ spaces:

\begin{lem}\label{wolff_lem0}
Let $\phi^{(i)}_{k_i} \in V^2_{\pm |D_x|}$ be two test functions
which have angular separation in frequency, namely
$|\theta_1-\theta_2|>\theta$ in the $(++)$ or $(--)$ cases, and also
$|\pi+(\theta_1-\theta_2)|>\theta$ in the $(+-)$ or $(-+)$ cases.
Then for $c < \frac{3}{29}$ we have:
\begin{equation}
       \lp{\phi^{(1)}_{k_1} \phi^{(2)}_{k_2}}{L^2_t(L^2_x)}   \
       \lesssim\
        2^{\frac{1}{2}\max\{k_i\}} \theta^{-1}
        \lp{\phi^{(1)}_{k_1}}{V^2_{\pm |D_x|}} \lp{\phi^{(2)}_{k_2}}{V^2_{\pm |D_x|}}^{1-c}
        \big(2^{-k_2} \lp{\phi^{(2)}_{k_2}}{L^\infty_t(L^\infty_x)}\big)^c
        \label{bii}
\end{equation}
\end{lem}\ret

We remark that we did not make an effort to optimize $c$, the
balance of the frequencies, or the power of $\theta$, as these play
no role in the present paper.

\begin{proof}
Without loss of generality, let us assume we are in the $(++)$ case.
If both $\phi^{(i)}_{k_i}$ were free waves, then Wolff's estimate
(with Tao's endpoint) would yield (see \cite{Taobi} Proposition
17.2):
\begin{equation}
        \lp{\phi^{(1)}_{k_1}\phi^{(2)}_{k_2}}{L^\frac{5}{3}_t(L^\frac{5}{3}_x)} \
        \lesssim \ 2^{\frac15 \max\{ k_i\} }\theta^{-1}
        \lp{\phi^{(1)}_{k_1}(0)}{L^2_x} \lp{
        \phi^{(2)}_{k_2}(0)}{L^2_x} \ . \notag
\end{equation}
Applying this for each intersection of two segments in a product of atoms,
we obtain:
\begin{equation}
         \lp{\phi^{(1)}_{k_1}\phi^{(2)}_{k_2}}{L^\frac{5}{3}_t(L^\frac{5}{3}_x)}
         \ \lesssim \ 2^{\frac15 \max \{k_i\}}\theta^{-1}
            \lp{\phi^{(1)}_{k_1}}{U^{\frac53}_{|D_x|}}
 \lp{\phi^{(2)}_{k_2}}{U^{\frac53}_{|D_x|}}
        \label{u53}
\end{equation}
On the other hand by \eqref{upst} with $(q,r)= (6,6)$
we have:
\begin{equation}
        \lp{\phi^{(1)}_{k_1} \phi^{(2)}_{k_2}}{L_t^3(L^3_x)}  \ \lesssim\
        2^{\frac{k_1+k_2}2} \lp{\phi^{(1)}_{k_1}}{U^{6}_{|D_x|}}
        \lp{\phi_{k_2}^{(2)}}{U^{6}_{|D_x|}}  \ . \label{u3}
\end{equation}
Interpolating \eqref{u53} with \eqref{u3} (it is bilinear
interpolation but it suffices to do it for atoms, so it only involves
$l^p$ and $L^p$ spaces) we obtain:
\begin{equation}
    \lp{\phi^{(1)}_{k_1} \phi^{(2)}_{k_2}}{L_t^p(L^{p}_x)}
    \ \lesssim \ 2^{(2-\frac3p) \max \{k_i\} - \frac{3}{26}|k_1-k_2| }\theta^{-1}
    \lp{\phi^{(1)}_{k_1}}{U^{2}_{|D_x|}} \lp{\phi_{k_2}^{(2)}}{U^{2}_{|D_x|}}
    \ , \ \  p = \frac{13}{7} \
    . \notag
\end{equation}
We want $V^2_{|D_x|}$ instead, so we use the embedding \eqref{upemb}
with $U^{2+}_{|D_x|}$ and
$L_t^{\frac{13}{7}+}(L^{\frac{13}{7}+}_x)$in this last estimate,
which gives the bound:
\begin{equation}
        \lp{\phi^{(1)}_{k_1} \phi^{(2)}_{k_2} }{L_t^{p}(L^{p}_x)}
        \ \lesssim   \ 2^{(2-\frac3p) \max\{ k_i\}  - \frac{3}{26}|k_1-k_2|} \theta^{-1}
         \lp{\phi^{(1)}_{k_1}}{V^{2}_{|D_x|}} \lp{\phi_{k_2}^{(2)}}{V^{2}_{|D_x|}}
         \ , \ \  p >
    \frac{13}{7} \ . \notag
\end{equation}
On the other hand by using an $L^\infty_t(L^\infty_x)$ bound we get:
\begin{equation}
        \lp{\phi^{(1)}_{k_1}\phi^{(2)}_{k_2}}{L^6_t(L^6_x)}
         \ \lesssim \ 2^{\frac{1}{2}k_1}
         \lp{\phi^{(1)}_{k_1}}{V^{2}_{|D_x|}}  \lp{\phi_{k_2}^{(2)}}{L^\infty_t(L^\infty_x)}
            \ . \notag
\end{equation}
Then \eqref{bii} is obtained interpolating the last two lines.
\end{proof}

In this article we work with the $S$ and $N$ spaces. The next lemma
relates them to the $V^2_{\pm|D_x|}$ spaces.

\begin{lem}
Let $\phi_k[0] \in \dot H^1 \times L^2$ and $F_k \in N$. Then the
solution $\phi_k$ to $\Box \phi_k = F_k$ with initial data
$\phi_k[0]$ satisfies:
\begin{equation}
    2^{\frac{1}{2}k}\lp{Q_{\geqslant k-10}\nabla_{t,x}\phi_k}{L^2_t(L^2_x)} + \sum_\pm
    \lp{Q^{\pm}_{<k-10} \nabla_{t,x} \phi_k}{V^2_{\pm|D_x|}}  \! \lesssim  \!
    \lp{\phi_k[0]}{\dot H^1 \times L^2}
    +\lp{F_k}{N}
     . \label{SNV2}
\end{equation}
\end{lem}\ret

\begin{proof}
By rescaling we may assume that $k=0$, and we'll relabel $\phi_k$
and $F_k$ as $\phi,F$ with the implicit understanding that they are
both at unit frequency.

The estimate for $Q_{\geqslant -10}\nabla_{t,x}\phi$ is immediate
from the structure of $S$ and the estimate \eqref{energy_est} (note
that this was shown in \cite{Tao_WM} for all portions of the norm
\eqref{S_norm} except $\underline{S}$).

The linear wave evolution in the energy space $\dot H^1 \times L^2$
is given by the multiplier:
\begin{equation}
    S(t)  \ =  \ \begin{pmatrix} \cos(t|D_x|)  &  |D_x|^{-1} \sin(t|D_x|) \cr
    -|D_x|\sin(t|D_x|)& \cos(t|D_x|) \end{pmatrix}  \ . \notag
\end{equation}
For any increasing sequence $t_j$ we can use the energy component of
\eqref{energy_est} (again established in \cite{Tao_WM}) and
\eqref{N_fung} to estimate:
\begin{equation}
    \sum_j \lp{ \phi[t_{j+1}] - S(t_{j+1}-t_j) \phi[t_j]}{\dot H^1 \times L^2}^2
     \ \lesssim  \ \sum_{j} \lp{ 1_{[t_j,t_{j+1}]} F }{N}^2
     \ \lesssim  \ \lp{F}{N}^2 \ . \notag
\end{equation}
Diagonalizing, one may write the $L^2\times L^2$ evolution as:
\begin{equation}
        \begin{pmatrix}  |D_x|  &  0  \\
        0 &  1 \end{pmatrix} S(t)\phi[t_0]
         \ = \ \mathcal{U}^*  \begin{pmatrix} e^{it|D_x|}  & 0 \cr
        0 & e^{-it|D_x|} \end{pmatrix} \mathcal{U} \begin{pmatrix}  |D_x|
        \phi(t_0)\\
        \partial_t\phi(t_0)\end{pmatrix}  \ , \notag
\end{equation}
where $\mathcal{U}=\frac{1}{\sqrt{2}}\begin{pmatrix}1&-i\\
1&i\end{pmatrix}$. Thus, the LHS of the previous difference formula
may be rotated via $\mathcal{U}$ to yield:
\begin{align}
    \lp{\phi[t] - S(t-s) \phi[s]}{\dot H^1 \times L^2}^2
     \ = \ &\frac{1}{2}\lp{ (\partial_t + i|D_x|)\phi(t) -
     e^{i(t-s)|D_x|} (\partial_t + i|D_x|)\phi(s)}{L^2}^2 \notag\\
    &\hspace{-.4in}+\ \frac{1}{2}\lp{ (\partial_t - i|D_x|)\phi(t) -
    e^{-i(t-s)|D_x|} (\partial_t - i|D_x|)\phi(s)}{L^2}^2 \ .
    \notag
\end{align}
Hence taking the supremum over all increasing sequences $t_k$ we
obtain the pair of bounds:
\begin{equation}
        \lp{  (\partial_t + i|D_x|)\phi}{V^2_{|D_x|}}^2 +
        \lp{  (\partial_t - i|D_x|)\phi}{V^2_{-|D_x|}}^2
         \ \lesssim  \ \lp{\phi[0]}{\dot H^1 \times L^2}^2+ \lp{F}{N}^2 \
        . \notag
\end{equation}
We conclude \eqref{SNV2} by noting that  one has the following
``elliptic'' estimate:
\begin{equation}
        \lp{\nabla_{t,x}Q^\pm_{<-10}P_0\phi}{Y} \ \lesssim \
        \lp{(\partial_t \pm i|D_x|)\phi}{Y} \ , \notag
\end{equation}
for any translation invariant space-time norm $Y$, which is valid
because the convolution kernel of the frequency
localized ratio $\nabla_{t,x} (\partial_t\pm
i|D_x|)^{-1}Q^\pm_{<-10}P_0$ is in $L^1_t(L^1_x)$.
\end{proof}

As a quick application of these ideas, notice that if $(q,r)$ is any
pair of indices in the range of \eqref{phys_str}, we must have $q
\geqslant 4$. Hence from \eqref{upst} and \eqref{SNV2}, and some
Sobolev embeddings interpolated with the  $L^\infty_t(L^2_x)$
estimate for $\nabla_{t,x}\phi$ to control the first member on
the LHS of \eqref{SNV2}, we obtain:

\begin{cor} \label{snstcor}
Let $\phi[0] \in \dot H^1 \times L^2$ and $F \in N$. Then the
solution $\phi$ to $\Box \phi = F$ with initial data $\phi[0]$
satisfies
\begin{equation}
       \lp{\phi}{\underline{S}}  \ \lesssim \
       \lp{ \phi[0]}{\dot H^1 \times L^2}
        +\lp{F}{N} \ .
        \label{snst}
\end{equation}
\end{cor}\ret

This proves \eqref{energy_est}, therefore completing the linear theory
in the $S$ and $N$ spaces, as needed in view of our modification of Tao's \cite{Tao_WM}
definition of the $S$ space, namely by adding the $\underline{S}$ norm to it.
\ret

In a similar manner, we can combine the bounds \eqref{bii} and \eqref{SNV2}
to obtain:

\begin{lem}\label{wolff_lem}
Let $\phi^{(i)}_{k_i}$ be two test functions normalized to that:
\begin{equation}
        \lp{\phi^{(i)}_{k_i}}{S} + \lp{\Box \phi^{(i)}_{k_i}}{N} \ \leqslant \ 1\ ,
        \  \ \ \ j = 1,2  \ ,
        \qquad\qquad \lp{\phi^{(2)}_{k_2}}{L^\infty_t(L^\infty_x)} \ \leqslant\
        \eta \ . \notag
\end{equation}
Assume in addition that the localizations $Q^\pm_{<k_i-10}
\phi^{(i)}_{k_i}$ have the angular separation
$|\theta_1-\theta_2|>\theta$ in the $(++)$ or $(--)$ cases, and
$|\pi+(\theta_1-\theta_2)|>\theta$ in the $(+-)$ or $(-+)$ cases.
Then for $c<\frac{3}{29}$ one has:
\begin{equation}
       \lp{\phi^{(1)}_{k_1} \phi^{(2)}_{k_2}}{L^2_t(L^2_x)} \ \lesssim
       2^{-\frac{3}{2}\min\{k_i\}}\eta^{c} \theta^{-1} \ . \label{bi}
\end{equation}
\end{lem}\ret

\begin{proof}
By an application of Lemma \ref{wolff_lem0}, we need only consider the case
where one factor is at high modulation, i.e. a factor of $Q_{\geqslant k_i-10}\phi_{k_i}^{(i)}$.
In this case, if the other factor has the improved $L^\infty_t(L^\infty_x)$ bound,
estimate \eqref{bi} is immediate on account of the $L^2_t(L^2_x)$
estimate on the LHS of \eqref{SNV2}. On the other hand, if the factor at high modulation
is also the one with improved $L^\infty_t(L^\infty_x)$ control, then using $L^6_t(L^6_x)$ Strichartz
for the first factor we have:
\begin{equation}
    \lp{\phi_{k_1}^{(1)}Q_{\geqslant k_2-10}\phi_{k_2}^{(2)}}{L^\frac{3}{2}_t(L^\frac{3}{2}_x)}
    \ \lesssim \ 2^{-\frac{1}{2}k_1}2^{-\frac{3}{2}k_2} \ . \notag
\end{equation}
Interpolating this with the pointwise bound
$|\phi_{k_1}^{(1)}Q_{\geqslant k_2-10}\phi_{k_2}^{(2)}|\lesssim
\eta$ we again have \eqref{bi}.
\end{proof}

\ret
%-------------------------------------------------------------------------
%%%%%%%%%%%%%%%%%%%%%%%%%%%%%%%%%%%%%%%%%%%%%%%%%%%%%%%%%%%%%%%%%%%%%%%%%%
%%%%%%%%%%%%%%%%%%%%%%%%%%%%%%%%%%%%%%%%%%%%%%%%%%%%%%%%%%%%%%%%%%%%%%%%%%
%-------------------------------------------------------------------------

\section{Bilinear Null Form Estimates}\label{s:matchfreq}
In this section we prove the estimates \eqref{standard_L2_est_bi},
\eqref{ednullw}, and \eqref{bal_l2_N}. The first of these is
essentially standard, being implicitly contained in the calculations
of \cite{Tao_WM}. We provide the proof here for the sake of
completeness:
\ret
%-------------------------------------------------------------------------
\begin{proof}[Proof of estimate \eqref{standard_L2_est_bi}]
We begin with the estimate:
\begin{equation}
        \lp{P_k Q_{<j}F}{L^2_t(L^2_x)} \ \lesssim \
        2^{\frac{j}{2}}\lp{F}{N} \ . \label{standard_L2_pre}
\end{equation}
To see this, notice that if $\Box^{-1}_0$ inverts the wave equation
with zero Cauchy data, we immediately have from \eqref{energy_est}
the inequality:
\begin{equation}
        \lp{P_k Q_{j'}\Box^{-1}_0 F}{X^{1,\frac{1}{2}}_\infty} \ \lesssim \
        \lp{F}{N} \ , \notag
\end{equation}
which implies the fixed frequency estimate:
\begin{equation}
        \lp{P_k Q_{j'}F}{L^2_t(L^2_x)} \ \lesssim \
        2^{\frac{j'}{2}}\lp{F}{N} \ . \notag
\end{equation}
Summing this last line over all $j'<j$ \eqref{standard_L2_pre} is
achieved. We now split the proof of estimate
\eqref{standard_L2_est_bi} into two cases: \ret

\step{1}{$Low \times High$ interaction} In this case we assume that
$k_1 < k_2 - 10$. The case $k_2 < k_1 - 10$ can be handled via a
similar argument. For relatively low modulations we have from
estimates \eqref{standard_est_bi} and \eqref{standard_L2_pre}:
\begin{equation}
        \lp{ Q_{<k_1+10}( \partial^\alpha \phi^{(1)}_{k_1}
        \partial_\alpha \phi^{(2)}_{k_2} )}{L^2_t(L^2_x)} \ \lesssim \
        2^{\frac{k_1}{2}}\lp{\phi^{(1)}_{k_1}}{S}\cdot \lp{\phi^{(2)}_{k_2}}{S}
        \ . \notag
\end{equation}
Therefore, it suffices to look at output modulations larger than
$k_1+10$. In this case we split the modulations of the low frequency
term according to $ \phi^{(1)}_{k_1}= Q_{<k_1}\phi^{(1)}_{k_1} +
Q_{\geqslant k_1}\phi^{(1)}_{k_1}$. For the first term we have that:
\begin{align}
        \lp{\! Q_{\geqslant k_1+10}( Q_{<k_1} \partial^\alpha \phi^{(1)}_{k_1}
        \cdot \partial_\alpha \phi^{(2)}_{k_2} )\!}{L^2_t(L^2_x)}
          \! &\lesssim\!  \lp{ Q_{<k_1} \partial^\alpha \phi^{(1)}_{k_1}
        \cdot Q_{>k_1}\partial_\alpha \phi^{(2)}_{k_2}
        }{L^2_t(L^2_x)} \  \notag\\
         &\lesssim\!
         \lp{\nabla_{t,x}\phi^{(1)}_{k_1}}{L^\infty_t(L^\infty_x)}
        \lp{Q_{>k_1}\nabla_{t,x}\phi^{(2)}_{k_2}}{L^2_t(L^2_x)}
         \notag\\
          &\lesssim\!  2^{\frac{k_1}{2}}\lp{\nabla_{t,x} \phi_{k_1}^{(1)}}{L^\infty_t(L^2_x)}
        \lp{\nabla_{t,x} \phi_{k_2}^{(2)}}{X_\infty^{0,\frac{1}{2}}}
        \ , \notag
\end{align}
which suffices. For the high modulations of the first factor in the
previous decomposition, we estimate:
\begin{equation}
         \lp{\! Q_{\geqslant k_1+10}( Q_{>k_1} \partial^\alpha \phi^{(1)}_{k_1}
        \cdot \partial_\alpha \phi^{(2)}_{k_2} )\!}{L^2_t(L^2_x)}
          \! \lesssim\!
        \lp{Q_{>k_1} \nabla_{t,x} \phi^{(1)}_{k_1}}{L^2_t(L^\infty_x)}
        \lp{\nabla_{t,x} \phi^{(2)}_{k_2}}{L^\infty_t(L^2_x)}
         \ . \notag
\end{equation}
We then conclude using \eqref{L2Linfty_est} for the first
factor.\ret

\step{2}{$High \times High$ interaction} In this case we consider
the frequency interaction $|k_1-k_2|<5 $, and without loss of
generality we may also assume that $k_1\geqslant k_2$. By  using
estimates \eqref{standard_L2_pre} and \eqref{standard_est_bi},   we
may reduce to considering the case of output modulation larger that
$k+ \delta(k_1-k)+10$, where $\delta$ is from the RHS of
\eqref{standard_est_bi} (this ultimately forces a harmless
redefinition of $\delta$ to suit line \eqref{standard_L2_est_bi}).
For this remaining piece, we will show that:
\begin{equation}
            \lp{ P_k Q_{\geqslant k + \delta(k_1-k)+10}
            (  \partial^\alpha \phi^{(1)}_{k_1}
            \cdot \partial_\alpha \phi^{(2)}_{k_2} )}{L^2_t(L^2_x)}
            \ \lesssim \ 2^{\frac{k}{2}}2^{-\frac{1}{2}\delta(k_1-k)}
            \lp{\phi^{(1)}_{k_1}}{S}\cdot\lp{\phi^{(2)}_{k_2}}{S} \ . \notag
\end{equation}
The key observation here is that the output modulation combined with
the output spatial frequency localization guarantees that at least
one term in the product is at modulation greater than $k +
\delta(k_1-k) -20$. Without loss of generality we may assume this is
the first term in the product, and we estimate via Bernstein:
\begin{align}
            &\lp{ P_k Q_{\geqslant k + \delta(k_1-k)+10}
            (  Q_{\geqslant k + \delta(k_1-k)-20} \partial^\alpha \phi^{(1)}_{k_1}
            \cdot \partial_\alpha \phi^{(2)}_{k_2} )}{L^2_t(L^2_x)}
             \notag\\
            \lesssim \ &2^k \sum_{j> k + \delta(k_1-k)-20}
            \lp{Q_j \nabla_{t,x}\phi^{(1)}_{k_1}}{L^2_t(L^2_x)}
            \cdot\lp{\nabla_{t,x}
            \phi^{(2)}_{k_2}}{L^\infty_t(L^2_x)}  \notag\\
            \lesssim \ &2^{\frac{k}{2}}2^{-\frac{1}{2}\delta(k_1-k)}
            \lp{\nabla_{t,x}\phi^{(1)}_{k_1}}{X_\infty^{0,\frac{1}{2}}}\cdot
            \lp{\nabla_{t,x}\phi^{(2)}_{k_2}}{L^\infty_t(L^2_x)} \ . \notag
\end{align}
This concludes our demonstration of \eqref{standard_L2_est_bi}.
\end{proof}\ret

%-------------------------------------------------------------------------

Our next step is to prepare for the proof of Proposition~\ref{p:matchfreq}.
It will first be useful to have a version of these estimates
under simpler assumptions:

\begin{lem}
a) Let $\phi^{(i)}_{k_i}$ be functions localized at frequency $k_i$. Assume that
these functions are normalized as follows:
\begin{align}
    \lp{\phi^{(i)}_{k_i}}{S[I]} + \lp{\Box \phi^{(i)}_{k_i}}{N[I]}  \ &\leqslant \  1 \ ,
    & \lp{\phi^{(1)}_{k_1}}{L^\infty_t(L^\infty_x)[I]} \ &\leqslant \ \eta \ . \label{e_disp_assz}
\end{align}
Then the following bilinear $L^2$ estimate holds:
    \begin{equation}
        \lp{ \partial^\alpha \phi^{(1)}_{k_1} \partial_\alpha \phi^{(2)}_{k_2}}{L^2_t(L^2_x)[I]}
         \ \lesssim  \ 2^{\frac12 \max\{k_1,k_2\}} \eta^\delta \ .
        \label{ednullwz}
    \end{equation}\ret
b) Assume that in addition to \eqref{e_disp_assz}
    we also have the high modulation bounds:
    \begin{align}
         \lp{ \Box \phi^{(1)}_{k_1}}{L^2_t(L^2_x)[I]}  \
         &\leqslant \ 2^{\frac{k_1}2} \eta \ ,
        &\lp{ \Box \phi^{(2)}_{k_2}}{L^2_t(L^2_x)[I]}  \
        &\leqslant \ 2^{\frac{k_2}2} \eta
         \ . \label{high_mod_assz}
    \end{align}
    Then the following estimate holds:
    \begin{equation}
        \lp{ \partial^\alpha \phi^{(1)}_{k_1} \partial_\alpha  \phi^{(2)}_{k_2} }{N[I]}
        \ \lesssim \ 2^{C|k_1-k_2|}\eta^\delta   \ . \label{bal_l2_Nz}
    \end{equation}
\label{lnore}\end{lem}

\begin{proof}
We may assume that the interval length is such that $|I| \geqslant
2^{-\min \{k_i\}}\eta^{2\delta}$, as otherwise the desired bounds follow
from  integrating energy estimates.

We begin by taking extensions of $\phi^{(1)}_{k_1}$ and
$\phi^{(2)}_{k_2}$ according to Proposition \ref{ext_prop} in such a
way that the $L^\infty_t(L^\infty_x)$ bound in \eqref{e_disp_assz}
is preserved; in the case of part (b), we also insure that
\eqref{high_mod_assz} is preserved.  This is achieved using the second
extension in Proposition \ref{ext_prop} in case (a), respectively
the canonical extension in Proposition \ref{ext_prop} in case (b).
Doing this requires balancing
the parameter $\eta$ in Proposition \ref{ext_prop}, and has the
effect of replacing the $\eta$ in both \eqref{e_disp_assz} and
\eqref{high_mod_assz} with a small power of $\eta$ ($\eta^\frac18$
should suffice). This is harmless given the small constant $\delta$
which we seek to obtain in both \eqref{ednullwz} and \eqref{bal_l2_Nz}.

We fix $m$ to be a large spatial frequency separation parameter.
In the course of proving \eqref{ednullwz} and \eqref{bal_l2_Nz},
we will decompose into several frequency ranges. In all cases we
will show a bound of the form:
\begin{equation}
    \hbox{L.H.S.} \ \lesssim \ 2^{Cm}\eta^c + 2^{-\delta m} \ , \notag
\end{equation}
where $c,\delta$ are suitably small constants depending only on the estimates in Propositions
\ref{standard_prop1}, \ref{standard_prop2}, and \ref{wolff_lem}
above, and where $C$ is a suitable large constant. In what follows we call any bound of this type a
``suitable bound''. By choosing $m$ appropriately, and by (globally) redefining
the small parameter $\delta$ one may produce the RHS of estimates
\eqref{ednullwz} and \eqref{bal_l2_Nz} from such bounds.
%This will be later chosen so that:
%\begin{equation}
%   2^{-m}  \ =  \ \eta^{c_1} \ , \notag
%\end{equation}
%for a sufficiently small suitably chosen $c_1$. Thus in all the
%estimates below we seek to obtain $L^2$, respectively $N$ bounds which
%contain factors which are either of the form
%$2^{-c_2 m}$, $c_2 > 0$ or of the form $\eta^{c_3} 2^{Cm}$ for a
%finite collection of constants $c_2$, $c_3$ and $C$ which
%depend on the constants $\delta$
% from Propositions \ref{standard_prop1},
%\ref{standard_prop2}, and \ref{wolff_lem}.
%The parameter $m$ is chosen so as to balance these two competing
%contributions. Note that this forces a
%harmless global redefinition of $\delta$ to accommodate Proposition
%\ref{p:matchfreq}. In what follows we call any bound of this type a
%suitable bound.
\ret

\step{1}{The unbalanced case $|k_1-k_2| \geqslant  \ m$} Here we neglect
the pointwise bound in  \eqref{e_disp_assz} as well as the
high modulation bound in \eqref{high_mod_assz}.
 From estimate  \eqref{standard_L2_est_bi} we immediately have that:
\begin{align*}
         \lp{ \partial^\alpha \phi^{(1)}_{k_1} \partial_\alpha \phi^{(2)}_{k_2}}{L^2_t(L^2_x)}
         \ &\lesssim  \ 2^{\frac12 \max\{k_1,k_2\}} 2^{-\frac{1}{2}m} \ ,
         &|k_1-k_2| \ \geqslant  \ m \ . \label{endullw1}
\end{align*}
which is a suitable $L^2$ bound. Similarly, from \eqref{standard_est_bi}
we obtain a suitable $N$ bound.

Hence, in what follows it suffices to consider the range
$|k_1-k_2|<m$. For the remainder of the proof we let $k_2=0$. We
split into cases depending on the modulations of the factors and the
output.
\ret

%-------------------------------------------------------------------------

\step{2}{The factor $\phi_{k_2}^{(2)}$ at high modulation} Here we
first  prove a suitable $L^2$ bound:
\begin{equation}
          \lp{ \partial^\alpha
          \phi^{(1)}_{k_1} Q_{>-10m}\partial_\alpha \phi^{(2)}_0}{L^2_t(L^2_x)}
         \ \lesssim \  2^{Cm} \eta + 2^{-m}\ .  \label{endullw2}
\end{equation}
For moderate modulations of the first factor, i.e. for $Q_{<10
m}\phi^{(1)}_{k_1}$, we use \eqref{e_disp_assz} to place it in
$L^\infty$:
\begin{equation}
    \lp{ \partial^\alpha Q_{<10 m}\phi^{(1)}_{k_1}}{L^\infty_t(L^\infty_x)}
    \ \lesssim  \ 2^{Cm} \eta \ , \notag
\end{equation}
while the second factor is placed in $L^2$ via to
the general embedding:
\begin{equation}
        \lp{Q_{>j} \phi^{(i)}_{k_i}}{L^2_t(L^2_x)} \ \lesssim  \
        2^{-\frac{1}{2}j}\lp{\phi^{(i)}_{k_i}}{X_\infty^{0,\frac{1}{2}}} \
        . \label{X_infty_largemod}
\end{equation}

For high modulations of the first factor, i.e. for $Q_{>10
m}\phi^{(1)}_{k_1}$, we reverse the roles and bound the first factor in
$L^2$:
\begin{equation}
    \lp{ \partial^\alpha Q_{>10 m}\phi^{(1)}_{k_1}}{L^2_t(L^2_x)}
    \ \lesssim  \ 2^{-5m} \lp{\phi^{(1)}_{k_1}}{S}
     \ . \label{firstvhigh}
\end{equation}
while the second factor has a $\lesssim 1$ bound in $L^\infty$
thanks to \eqref{Linfty_Q_est}, which leads again to a suitable
bound.

In this case it is even easier to obtain the suitable $N$ bound
because we have access to the high modulation assumption \eqref{high_mod_assz}. We prove:
\begin{equation}
          \lp{ \partial^\alpha
          \phi^{(1)}_{k_1} Q_{>-10m}\partial_\alpha \phi^{(2)}_0}{N}
         \ \lesssim \  2^{Cm} \eta.  \label{endullbv}
\end{equation}
This follows from \eqref{standard_est_bi} combined with:
\begin{equation}
    \lp{ Q_{>-10m}\phi^{(2)}_0 }{S}
    \ \lesssim \ 2^{5m}
    \lp{\Box \phi^{(2)}_0 }{L^2_t(L^2_x)}
    \ \lesssim \ 2^{5m} \eta \ , \notag
\end{equation}
where the first inequality follows from \eqref{XS_est}.
\ret

%-------------------------------------------------------------------------

\step{3}{The factor $\phi_{k_1}^{(1)}$ at high modulation} Here we
can also prove a suitable $L^2$ bound, namely:
\begin{equation}
          \lp{ \partial^\alpha Q_{> -10m}
          \phi^{(1)}_{k_1} Q_{\leqslant -10m}\partial_\alpha \phi^{(2)}_0}{L^2_t(L^2_x)}
         \ \lesssim \ 2^{Cm}\eta^\frac{1}{4} + 2^{-m}\ .  \label{endullw3}
\end{equation}
Reusing \eqref{firstvhigh} we can dispense with the very high modulations
in $\phi^{(1)}_{k_1}$ and  replace the first factor
with $\partial^\alpha Q_{[-10m,10m]}
\phi^{(1)}_{k_1}$. This time we cannot directly use the
$L^\infty_t(L^\infty_x)$ estimate for $\phi^{(1)}_{k_1}$. However,
by applying \eqref{X_infty_largemod} and using the $L_t^6(L_x^6)$
Strichartz estimate contained in \eqref{snst} we
have that:
\begin{equation}
         \lp{ \partial^\alpha  Q_{ [-10m, 10m]}
          \phi^{(1)}_{k_1} Q_{\leqslant  -10m}\partial_\alpha \phi^{(2)}_0}{L^\frac{3}{2}_t
          (L^\frac{3}{2}_x)}
         \ \lesssim \ 2^{Cm} \ . \notag
\end{equation}
Next,  using \eqref{e_disp_assz} and
\eqref{Linfty_Q_est} we directly have:
\begin{equation}
         \lp{ \partial^\alpha  Q_{ [-10m,10m]}
          \phi^{(1)}_{k_1} Q_{\leqslant -10m}\partial_\alpha \phi^{(2)}_0}{L^\infty_t
          (L^\infty_x)}
         \ \lesssim \ 2^{Cm}\eta \ . \notag
\end{equation}
Interpolating these last two estimates yields \eqref{endullw3}. It
is important  to notice that in the above estimates one  looses a
polynomial in $2^{m}$ because the multipliers $P_0 Q_{\leqslant
-10m}$ and $P_{k_1}Q_{[-10m, 10m]}$ are not uniformly disposable on
$L^p$. However, a short calculation shows that the resulting
convolution kernels have $L^1_t(L^1_x)$ bounds on the order of
$2^{Cm}$ which is acceptable.

As in the previous step we also have a suitable $N$ bound:
\begin{equation}
          \lp{ \partial^\alpha
          Q_{>-10m} \phi^{(1)}_{k_1}\partial_\alpha Q_{\leqslant -10m} \phi^{(2)}_0}{N}
         \ \lesssim \  2^{Cm} \eta \ .  \label{endullbva}
\end{equation}
\ret

%-------------------------------------------------------------------------
\step{4}{Low frequency output} This is the case when $k_1=k_2+O(1)$,
and we seek to estimate $P_k (\partial^\alpha\phi^{(1)}_{k_1}
\partial_\alpha  \phi^{(2)}_0)$ for $k < -m$. Then we can use \eqref{standard_L2_est_bi}, respectively \eqref{standard_est_bi}
to obtain a $\lesssim 2^{-\delta m}$ suitable bound in $L^2$,
respectively $N$. Here $\delta$ is the previously defined constant from Proposition
\ref{standard_prop1}.
\ret
%-------------------------------------------------------------------------

\step{5A}{Both $\phi_{k_i}^{(i)}$ at low modulation, output at low
modulation $< -2m$ and high frequency $k > -m$} In this case,
to show \eqref{ednullwz} we prove the bound:
\begin{equation}
         \lp{P_k Q_{< -2m}\big( \partial^\alpha Q_{\leqslant  -10m}
         \phi^{(1)}_{k_1} \partial_\alpha Q_{\leqslant -10m}
         \phi^{(2)}_0\big)}{L^2_t(L^2_x)}
         \ \lesssim \ 2^{-\delta m} \ , \label{endullw4}
\end{equation}
where the $\delta$ is the same as in Propositions \ref{standard_prop1} and
\ref{standard_prop2}. This estimate again uses only the $S$ bounds on
$\phi^{(1)}_{k_1}$ and  $\phi^{(2)}_{k_2}$ and the localization
conditions $|k_1|\leqslant m$ and $k_2=0$.
%which therefore play symmetric roles.
%Hence we can assume that $k_1 \leq k_2= 0$, and we need to consider
%two cases, namely (i) $-m < k_1 < O$ and $k = O(1)$, respectively
%(ii) $k_1 = O(1)$ and $-m < k < O(1)$.
To show \eqref{endullw4}, by \eqref{standard_L2_pre} it suffices to prove the
following set of bounds which together also imply \eqref{bal_l2_Nz} in the present case:
\begin{align}
             \lp{P_k Q_{< -2m}\Box \big( Q_{\leqslant  -10m}
         \phi^{(1)}_{k_1} \cdot Q_{\leqslant -10m}
         \phi^{(2)}_0\big)}{X_1^{0,-\frac{1}{2}}}
         \ &\lesssim \ m2^{-\delta m} \ , \label{endullw4_1}\\
          \lp{P_k Q_{< -2m}\big( \Box  Q_{\leqslant  -10m}
         \phi^{(1)}_{k_1}\cdot Q_{\leqslant -10m}
         \phi^{(2)}_0\big)}{N}
         \ &\lesssim \ 2^{-\delta m} \ , \label{endullw4_2}\\
          \lp{P_k Q_{< -2m}\big(  Q_{\leqslant -10m}
         \phi^{(1)}_{k_1} \cdot\Box Q_{\leqslant -10m}
         \phi^{(2)}_0\big)}{N}
         \ &\lesssim \ 2^{-\delta m} \ . \label{endullw4_3}
\end{align}
The first estimate above follows from \eqref{XS_X_est}, while the
second and third  both follow from \eqref{XS_N_est}. Note that while the multiplier
$Q_{\leqslant -2m}$ is not disposable on $N$ (e.g. on the $N\!\!F\!\!A$ atoms),
one may first replace it by $Q_{<0}$, and separately estimate the contribution
of $Q_{[-2m,0]}$ as an $X_1^{0,-\frac{1}{2}}$ atom via \eqref{XN_est}
at an $O(m)$ loss. A similar method using \eqref{XS_est} allows
one to handle the interior  $Q_{\leqslant -10m}$ multipliers, which
are not disposable on $S$, with another  $O(m)$ loss.
% Notice
% that in disposal of the $Q_{< -10m}$ in $S_k$ we only loose a log
% (factor of $m$) when $k=O(m)$, which is easily seen from the estimate:
% \begin{equation}
%         \lp{Q_{<-10 m}\psi_k}{S} \ \lesssim \
%         \lp{Q_{<k}\psi_k }{S} + \sum_{j:\ -10m<  j < k}
%         \lp{Q_j \psi_0 }{X_1^{1,\frac{1}{2}}} \
%         \lesssim \ m \lp{\psi_k}{S} \ . \notag
% \end{equation}
% Similarly, we need to dispose of the outer $Q_{< -2m }$  multiplier
% in both of estimates \eqref{endullw4_2}--\eqref{endullw4_3} whenever
% $k>-2m$. This can be done with only a log loss via:
% \begin{align}
%         \lp{P_k Q_{< -2m } F}{N} \ &\lesssim \
%         \lp{P_k Q_{< k} F}{N} +  \lp{P_k Q_{ -2m\leqslant \cdot < k }
%         F}{N} \ , \notag\\
%         \ &\lesssim \ \lp{P_k  F}{N} +
%         \sum_{j:\ -2m\leqslant j < k} \lp{P_k Q_j F}{X_1^{0,-\frac{1}{2}}}
%         \ \lesssim \ m\lp{P_k  F}{N} \ . \notag
% \end{align}
% In the last inequality we have used \eqref{XN_est} and the fact that
% $k <  m +O(1)$.\\

% The implication \eqref{endullw4_1}$\Rightarrow$\eqref{endullw4} is
% immediate, while the implications
% \eqref{endullw4_1}$\Rightarrow$\eqref{endullw4} and
% \eqref{endullw4_1}$\Rightarrow$\eqref{endullw4} follow from summing
% over the estimate \eqref{standard_L2_pre}.\\
\ret

\step{5B}{Both $\phi_{k_i}^{(i)}$ at low modulation, output at high
frequency and high modulation} In this step, which is the heart of
the matter, we establish the single bound:
\begin{equation}
         \lp{P_{\geqslant -m}Q_{\geqslant -2m}\big( \partial^\alpha Q_{\leqslant -10m}
         \phi^{(1)}_{k_1} \partial_\alpha Q_{\leqslant -10m}
         \phi^{(2)}_0\big)}{L^2_t(L^2_x)}
         \ \lesssim \ 2^{C m}\eta^c \ . \label{endullw6}
\end{equation}
Here $c$ is the same small constant from the RHS of line \eqref{bi}.
To use that estimate, we only need to establish angular separation
of the two factors. This is a standard ``geometry of the cone''
calculation, and one finds that the angle between the two factors
must satisfy $|\theta| \gtrsim 2^{-m}$ in the $(++)$ or $(--)$
cases, and $|\theta-\pi|\gtrsim 2^{-m}$ in the $(+-)$ or $(-+)$
cases (see for example Lemma 11 in Section 13 of \cite{Tao_WM}). By
decomposing the product on the LHS of \eqref{endullw6} into
$O(2^{Cm})$ angular sectors such that each product has these
separation properties, and by repeatedly applying
estimate \eqref{bi} on each interaction we have \eqref{endullw6}.
The proof of the lemma is concluded.

%
% %-------------------------------------------------------------------------
%
% \step{6}{Choosing the frequency separation} To wrap up
% \eqref{ednullw} we combine \eqref{endullw1}, \eqref{endullw2},
% \eqref{endullw3}, \eqref{endullw4}, \eqref{endullw5_1}, and
% \eqref{endullw6} by choosing $m=\delta^2 |\ln(\eta)|$, where
% $\delta$ is from Propositions \ref{standard_prop1},
% \ref{standard_prop2}, and \ref{wolff_lem}. Note that this forces a
% harmless global redefinition of $\delta$ to accommodate Proposition
% \ref{p:matchfreq}.
\end{proof}\ret
%-------------------------------------------------------------------------

%-----------------------------------------------------------------------
\begin{proof}[Proof of Proposition~\ref{p:matchfreq}]
For this we use Lemma~\ref{lnore}.
We begin using the extensions (this will be modified somewhat shortly)
and the same parameter $m$ as  in the proof of Lemma~\ref{lnore}.
We start with several simplifications.
The key observation is that in the proof of Lemma~\ref{lnore} we have used
the bound on $\lp{\Box \phi_{k_i}^{(i)}}{N}$ just once, namely in STEP 5B.
All other cases carry over to the proof  of Proposition~\ref{p:matchfreq}.
Consequently, it suffices to estimate the expression:
\begin{equation}
    P_k Q_{>-2m} R  \ =  \
    P_k Q_{>-2m} (\partial^\alpha \phi^{(1)}_{k_1}
    \partial_\alpha \phi^{(2)}_{k_2}) \ , \notag
\end{equation}
in both $L^2$ and $N$ under the assumptions $k_2 = 0$, $|k_1|
\leqslant m$, and $|k| \leqslant m+2$.

Furthermore, the contribution $P_kQ_{>-2m}\big((1-\chi_I)R\big)$ of $R$ in the
exterior of $I$ is estimated directly by Lemma \ref{lnore} because the
extensions provided by Proposition \ref{ext_prop} enjoy estimate
\eqref{ipmrhs} in the exterior of $I$. Hence, we only need
 consider the expression $P_kQ_{>-2m}(\chi_I R)$.
For this we will establish the pair of suitable bounds:
\begin{equation*}
    \lp{P_k Q_{>-2m} (\chi_I  R)}{L^2_t(L^2_x)}   \lesssim   2^{Cm} \eta^{\delta} + 2^{-5m} ,
    \ \ \ \lp{P_k Q_{>-2m} (\chi_I  R)  }{N}  \lesssim   2^{Cm} \eta^{\delta} + 2^{-4m}
     ,
\end{equation*}
%respectively
%\begin{equation*}
%\|P_k Q_{>-2m} R\|_{N} \lesssim 2^{Cm} \eta^{\delta} + 2^{-3m}
%\end{equation*}
with $\delta$ as in Lemma~\ref{lnore}.  We remark that due to the
frequency and modulation localization of $R$, the second $N$ bound
follows from the first $L^2$ bound  albeit with a readjusted $C$.  Therefore, we
drop the modulation and spatial frequency localization and simply prove that:
\begin{equation}
   \lp{P_k  R}{L^2_t( L^2_x)[I] }  \ \lesssim  \ 2^{Cm} \eta^{\delta} + 2^{-5m}
    \ . \label{l2suit}
\end{equation}
For this we use the renormalization. On the interval $I$, we may decompose
$\phi^{(i)}_{k_i}$ as follows:
\begin{equation}
    \phi_{k_i}^{(i)} \ = \ (U^{(i)}_{,<k_i})^\dagger w^{(i)}_{,k_i} \ , \notag%\label{phi_uw_exp}
\end{equation}
where by using the definition \eqref{big_S_def} we may assume that the
component pieces separately obey the estimates:
\begin{align}
        \lp{P_k w^{(i)}_{k_i}}{S[I]} + \lp{P_k \Box w^{(i)}_{k_i}}{N[I]}\ &\leqslant \  2^{-|k-k_i|}A^{-1} \ , \notag\\
    \lp{ U^{(i)}_{,<k_i}}{S} + \sup_{k > k_i} 2^{C(k_i-k)} \lp{P_k U^{(i)}_{,<k_i}}{S}
    \ &\leqslant \ A \ , \notag %\label{sonu}
\end{align}
for a possibly large constant $A$. By normalization, we may without loss of generality assume that
$A=1$, as any bounds for these two quantities will always appear as a product.
Since $\lp{\phi_{k_1}^{(1)}}{L^\infty_t(L^\infty_x)[I]} \leqslant \eta$, we obtain a similar relation for
$w^{(1)}_{,k_1} $, namely:
\begin{align}
     \lp{P_k w^{(1)}_{,k_1}}{L^\infty_t(L^\infty_x)[I]}
     \ \leqslant \ \eta \ . \notag%\label{e_disp_ass_red}
\end{align}
Furthermore, by using Proposition \ref{ext_prop}, we may extend the $w^{(i)}_{k_i}$ so
that all of the above listed bounds are global, albeit with a fractional modification of $\eta$.
Thus, we may drop the interval $I$, and again work globally.

We decompose the null-form $R$ (first on $I$, then by extension)
into $R =R_1+ R_2+R_3+R_4$ where:
\begin{align}
    R_1  \ &= \ - \partial^\alpha (U^{(1)}_{,<k_1})^\dagger\cdot  w^{(1)}_{,k_1}
 \partial_\alpha (U^{(2)}_{,<0})^\dagger\cdot w^{(2)}_{,0} \ , \notag\\
    R_2 \ &= \  \partial^\alpha (U^{(1)}_{,<k_1})^\dagger\cdot  w^{(1)}_{,k_1}
    \cdot \partial_\alpha \phi^{(2)}_{0} \ , \notag\\
    R_3 \ &= \    \partial^\alpha  \phi^{(1)}_{k_1} \cdot
    \partial_\alpha (U^{(2)}_{,<0})^\dagger\cdot w^{(2)}_{,0} \ , \notag\\
    R_4 \ &= \  (U^{(1)}_{,<k_1})^\dagger\cdot \partial^\alpha w^{(1)}_{,k_1}\cdot
     (U^{(2)}_{,<0})^\dagger\cdot \partial_\alpha w^{(2)}_{,0} \ . \notag
\end{align}
We successively consider each of these terms.
\ret
%-------------------------------------------------------------------------

\step{1}{Estimating the term $R_1$}
Using the $S$ bounds  for $U^{(1)}_{,<k_1}$ and $U^{(2)}_{,<0}$
in the bilinear $L^2$ null form estimate \eqref{standard_L2_est_bi}, after
dyadic summation we obtain:
\begin{equation}
    \lp{ \partial^\alpha (U^{(1)}_{,<k_1})^\dagger \partial_\alpha
    (U^{(2)}_{,<0})^\dagger}{L^2_t(L^2_x)} \ \lesssim \ m \ . \notag
\end{equation}
Note that the RHS loss is the effect of summing over frequencies $k'\leqslant k_1\leqslant m$
on the first factor.
We combine this with the pointwise bound on $w^{(1)}_{,k_1}$ to achieve:
\begin{equation}
    \lp{ R_1}{L^2_t(L^2_x)} \ \lesssim \  \lp{w^{(1)}_{,k_1}}{L^\infty_t(L^\infty_x)}
    \ \lesssim  \ m\eta \ . \notag
\end{equation}
\ret

\step{2}{Estimating the term $R_2$}
This is essentially same as in the previous step. Here we use
the $S$ bounds  for $U^{(1)}_{,<k_1}$ and for $\phi^{(2)}_0$
in conjunction with \eqref{standard_L2_est_bi}, and we again
use the pointwise bound for $w^{(1)}_{,k_1}$.
%\[
%\| \partial^\alpha (U^{(1)}_{,<k_1})^\dagger \partial_\alpha \phi^{(2)}_{0}\|_{L^2}
%\lesssim m
%\]
%Hence $R_2$ satisfies the same bound as $R_1$,
%\[
%\| R_2 \|_{L^2} \lesssim m \| w^{(1)}_{,k_1}\|_{L^\infty} \lesssim m \eta
%\]
\ret

\step{3}{Estimating the term $R_3$} We begin by splitting
$\phi^{(1)}_{k_1}$ into a low and a high modulation part. For the
high modulation part we have from \eqref{X_infty_largemod} the $L^2$
bound:
\begin{equation}
    \lp{ Q_{> 10 m} \partial^\alpha \phi^{(1)}_{k_1} }{L^2_t(L^2_x)}
    \ \lesssim \ 2^{-5m} \ . \notag
\end{equation}
Furthermore, by summing over the energy estimate for $U^{(2)}_{,<0}$
and using the decay of high frequencies we have the pointwise bound:
\begin{equation}
    \lp{\nabla_{t,x}U^{(2)}_{,<0}}{L^\infty_t(L^\infty_x)} \ \lesssim \
    \sum_k 2^k\lp{P_k \nabla_{t,x}U^{(2)}_{,<0}}{L^\infty_t(L^2_x)}\
    \lesssim \ 1 \ . \label{Linfty_dt_U}
\end{equation}
Combining these two estimates with the pointwise bound for $ w^{(2)}_{,0}$ we can estimate
the corresponding part of $R_3$, call it $R_{31}$, by:
\begin{equation}
    \lp{  R_{31}}{L^2_t(L^2_x)} \ \lesssim\  2^{-5m} \ . \notag
\end{equation}

It remains to consider the contribution of the low modulation
part $Q_{< 10 m} \phi^{(1)}_{k_1}$ in $R_3$, which we will label by
$R_{32}$. Using the $S$ bounds for
$\phi^{(1)}_{k_1}$ and $U^{(2)}_{,<0}$ along with
\eqref{standard_L2_est_bi}, after dyadic summation we obtain
the usual $L^2$ estimate:
\begin{equation}
    \lp{ \partial^\alpha  Q_{< 10 m} \phi^{(1)}_{k_1} \cdot
    \partial_\alpha (U^{(2)}_{,<0})^\dagger}{L^2_t(L^2_x)}
    \ \lesssim \ 1 \ . \notag
\end{equation}
On the other hand, from the Strichartz control \eqref{phys_str} and the
boundedness of the gauge we have:
\begin{equation}
    \lp{ w^{(2)}_{,0}}{L^6_t(L^6_x)} \ \lesssim \
    \lp{\phi^{(2)}_{0}}{L^6_t(L^6_x)} \ \lesssim \ 1 \ . \notag
\end{equation}
therefore we obtain a low index space-time $L^p$ bound for $R_{32}$,
namely:
\begin{equation}
    \lp{R_{32}}{L^\frac{3}{2}_t(L^\frac{3}{2}_x)}
    \ \lesssim  \ 1 \ . \notag
\end{equation}
On the other hand, from the pointwise bound \eqref{e_disp_ass} for
$\phi^{(1)}_{k_1}$ we obtain:
\begin{equation}
    \lp{ \partial^\alpha  Q_{< 10 m} \phi^{(1)}_{k_1}}{L^\infty_t(L^\infty_x)}
    \ \lesssim \ 2^{Cm} \eta \ . \notag
\end{equation}
Combining this with \eqref{Linfty_dt_U} and the pointwise bound for
$ w^{(2)}_{,0}$ we have:
\begin{equation}
    \lp{ R_{32}}{L^\infty_t(L^\infty_x)} \ \lesssim \ 2^{Cm} \eta \ . \notag
\end{equation}
Interpolating the last two lines we obtain:
\begin{equation}
    \lp{ R_{32}}{L^2_t(L^2_x)}  \ \lesssim \ 2^{Cm} \eta^\frac14 \ . \notag
\end{equation}
\ret

%-------------------------------------------------------------------------

\step{4}{Estimating the term $R_4$}
We start by dividing the main part of the product into all spatial frequencies:
\begin{equation}
    \partial^\alpha w^{(1)}_{,k_1}\cdot  \partial_\alpha w^{(2)}_{,0}
     \ =  \ \sum_{j_i} \partial^\alpha P_{j_1} w^{(1)}_{,k_1}\cdot  \partial_\alpha
    P_{j_2} w^{(2)}_{,0}  \ . \notag
\end{equation}
Using the bound \eqref{ednullwz} if $j_1,j_2 < 10 m$, and
\eqref{standard_L2_est_bi} otherwise in conjunction with the $2^{-|j_i-k_i|}$ frequency
separation gains for $P_{j_i} w^{(i)}_{,k_i}$ we have:
\begin{equation}
    \lp{ \partial^\alpha w^{(1)}_{,k_1}\cdot  \partial_\alpha w^{(2)}_{,0}}{L^2_t(L^2_x)}
     \ \lesssim  \ 2^{Cm} \eta^\delta + 2^{-5m} \ . \notag
\end{equation}
This estimate is directly transferred to $R_4$ due to the pointwise
bounds on the gauge factors.
\end{proof}\ret

%-------------------------------------------------------------------------
%%%%%%%%%%%%%%%%%%%%%%%%%%%%%%%%%%%%%%%%%%%%%%%%%%%%%%%%%%%%%%%%%%%%%%%%%%
%-------------------------------------------------------------------------

\section{Proof of the Trilinear Estimates}\label{s:der_tri}

In this section we will prove estimates
\eqref{core_L2}--\eqref{ext_N}. In all cases the desired bounds
follow easily from a combination of the standard estimates
\eqref{standard_L2_est_bi}--\eqref{standard_est_tri}  for widely
separated frequencies, and the improved matched frequency estimates
\eqref{ednullw} and \eqref{bal_l2_N}.
\ret

\begin{proof}[Proof of estimate \eqref{core_L2}]
The proof will be accomplished in a series of steps whose goal
is to reduce things to the matched frequency bilinear estimate \eqref{ednullw}.
\ret

\step{1}{Disposal of the $\phi^{(1)}$ Factor}
As a first step we will show the general estimate:
\begin{equation}
    \lp{\phi\cdot F}{L^2_t(\dot{H}^{-\frac{1}{2}})_c[I]} \ \lesssim \ \lp{\phi}{S[I]}
    \cdot\lp{F}{L^2_t(\dot{H}^{-\frac{1}{2}})_c[I]} \ , \label{SL2_est}
\end{equation}
where $\{c_k\}$ is any $(\delta_0,\delta_0)$-admissible frequency
envelope. To prove this, we split into the three main frequency
interactions. \ret

In the $Low\times High$ case we immediately have:
\begin{equation}
    \lp{P_{k}(\phi_{<k-10}\cdot F)}{L^2_t(\dot{H}^{-\frac{1}{2}}_x)[I]} \ \lesssim \
    \lp{\phi}{S[I]}\cdot 2^{-\frac{k}{2}}\lp{P_{[k-5,k+5]}F }{L^2_t(L^2_x)[I]} \ , \notag
\end{equation}
which is sufficient.\ret

In the $High\times Low $ case, we freeze the dyadic frequency of $F$ and we have a similar estimate:
\begin{equation}
    \lp{P_{k}(\phi\cdot F_{k'})}{L^2_t(\dot{H}^{-\frac{1}{2}}_x)[I]} \ \lesssim \
    2^{\frac{k'-k}{2}}c_{k'}\ \lp{\phi}{S[I]}\cdot \lp{F }{L^2_t(\dot{H}^{-\frac{1}{2}}_x)_c[I]} \ , \notag
\end{equation}
for any $k'\leqslant k-10$. Summing this over all such $k'\leqslant
k-10$ and using \eqref{sum_rule1} we have \eqref{SL2_est} in this
case.\ret

In the $High\times High$ case we freeze the frequency of the inputs and output to estimate:
\begin{equation}
    \lp{P_{k}(\phi_{k_1}\cdot F_{k_2})}{L^2_t(L^2_x)} \ \lesssim \ 2^{k-k_1}\lp{\phi_{k_1}}{S[I]}
    \cdot 2^{\frac{k_2}{2}}c_{k_2} \lp{F}{L^2_t(\dot{H}^{-\frac{1}{2}})_c[I]} \ ,
    \label{HH_XS_est}
\end{equation}
which follows easily from Bernstein's inequality \eqref{Bernstein}
and the bound \eqref{L2Linfty_Q_est}. Multiplying this last line by
$2^{-\frac{1}{2}k}$, and then summing over all $k_1$ and $k_2$ such
that $|k_1-k_2|\leqslant 20$ and $k_1\geqslant k-10$,  and then
using \eqref{sum_rule2}, we arrive at the estimate \eqref{SL2_est}
for this case.

\ret \step{2}{The Bilinear Estimate} In light of estimate
\eqref{SL2_est} above, it suffices to show that:
\begin{equation}
        \lp{\partial^\alpha\phi^{(2)} \partial_\alpha\phi^{(3)} }{L^2_t(\dot{H}^{-\frac{1}{2}})_c[I]}
        \lesssim \ \eta^{\delta_1} \ , \label{red_core_L2}
\end{equation}
assuming the conditions of estimate \eqref{core_L2}. This will be
done in two steps.

\ret \step{2A}{Reduction to Matched Frequencies} Our first step is
to peel off all frequency interactions that cannot be treated by
estimate \eqref{ednullw}. In all of these interactions, we will
exploit the fact that there is a wide separation in the frequency.
This is measured by choosing a large integer $m_0=m_0(\eta)$ such
that:
\begin{equation}
        2^{-\frac{1}{2}\delta m_0} \ = \ \eta^{\delta_1} \ ,
        \label{L2_red_delta_cond}
\end{equation}
where we remind the reader that $\delta$ is the small dyadic savings
from the standard $L^2$ bilinear estimate on line
\eqref{standard_L2_est_bi}, and because of the definition of
$\delta_1$ we have:
\begin{equation}
        m_0 \ \lesssim \ \sqrt{\delta_1}|\ln(\eta)| \ . \label{m_bound}
\end{equation}
Our goal in this step is to show the following fixed frequency
estimate:
\begin{equation}
        \sum_{\substack{k_i\\ \max\{|k_i-k|\}\geqslant m_0}}\!\!\!
        2^{-\frac{1}{2}k}\lp{P_k\big( \partial^\alpha \phi^{(2)}_{k_2}
        \partial_\alpha\phi^{(3)}_{k_3}
        \big) }{L^2_t(L^2_x)[I]}  \lesssim
       2^{-\frac{1}{2}\delta m_0} c_k \lp{\phi^{(2)}}{S[I]}\lp{\phi^{(3)}}{S_c[I]}
       , \label{core_L2_unmatched}
\end{equation}
which in light of \eqref{L2_red_delta_cond} suffices to establish
\eqref{red_core_L2} for all frequency interactions except for the
case $k=k_1 + O(m_0)= k_2 + O(m_0)$. By an application of estimate
\eqref{standard_L2_est_bi},  the two sum rules
\eqref{sum_rule1}--\eqref{sum_rule2}, and the definition
\eqref{freq_env_defn} we immediately have:
\begin{align}
        \hbox{(L.H.S.)}\eqref{core_L2_unmatched}  &\lesssim
         \sum_{\substack{k_i\\ \max\{|k_i-k|\}\geqslant m_0}}  \!\!\!\!\!
         2^{-\frac{1}{2}k}
         2^{\frac{1}{2}\min\{k_i\}}2^{-(\frac{1}{2}+\delta)(\max\{k_i\}-k)}c_{k_3}
         \lesssim  2^{-(\delta-\delta_0) m_0}c_k  , \notag
\end{align}
which  by  using \eqref{L2_red_delta_cond} and the definition of the
$\delta_i$ suffices to establish \eqref{core_L2_unmatched}.

\ret \step{2B}{The Matched Frequency Case} We have now reduced
estimate \eqref{red_core_L2} to showing the  matched frequency
bound:
\begin{equation}
         \sum_{\substack{k_i\\ \max\{|k_i-k|\}< m_0}}
        2^{-\frac{1}{2}k}\lp{P_k\big( \partial^\alpha \phi^{(2)}_{k_2}
        \partial_\alpha\phi^{(3)}_{k_3}
        \big) }{L^2_t(L^2_x)[I]} \ \lesssim \ \eta^{\delta_1}c_k  \ . \notag
\end{equation}
Due to the fact that there are only $O(m_0)\leqslant| \ln(\eta)|$
terms in this sum, it suffices to show:
\begin{align}
        \lp{P_k\big( \partial^\alpha \phi^{(2)}_{k_2}
        \partial_\alpha\phi^{(3)}_{k_3}
        \big) }{L^2_t(L^2_x)[I]} \ &\lesssim \ \eta^{2\delta_1}2^{\frac{1}{2}k}c_k  \ ,
     &\max\{|k_i-k|\}  &\lesssim \sqrt{\delta_1}|\ln(\eta)| \ .
     \notag
\end{align}
But this last estimate follows immediately from \eqref{ednullw} and
the definition of the $\delta_i$.
\end{proof}\ret

%---------------------------------------------------------

\begin{proof}[Proof of estimate \eqref{ext_L2}]
This estimate was essentially established in the previous proof. We
split the estimate into a sum of two pieces:
\begin{align}
        (\hbox{L.H.S.})\eqref{ext_L2} \ \leqslant & \
        \lp{P_k\big[P_{<k+10}\phi^{(1)}\partial^\alpha\phi^{(2)}
        \partial_\alpha\phi^{(3)}\big]}{L_t^2(L^2_x)[I]}\notag
       \\ &  \ \ \ \ + \lp{P_k \big[P_{\geqslant k+10}\phi^{(1)}\partial^\alpha\phi^{(2)}
        \partial_\alpha\phi^{(3)}\big]}{L_t^2(L^2_x)[I]} \ . \notag
\end{align}
For the first term we simply use \eqref{core_L2}. For the second
term, we use the following version of \eqref{HH_XS_est} above:
\begin{equation}
    \lp{P_{k}(\phi_{k_1}\cdot F_{k_2})}{L^2_t(\dot{H}^{-\frac{1}{2}})} \ \lesssim \ 2^{\frac{k-k_1}{2}}c_{k_1}
    \lp{\phi}{S_c[I]}\cdot \lp{F}{L^2_t(\dot{H}^{-\frac{1}{2}})[I]} \ ,
    \notag
\end{equation}
for $|k_1-k_2|\leqslant 5$, along with \eqref{red_core_L2}. This
suffices via the sum rule \eqref{sum_rule2}.
\end{proof}\ret

%----------------------------------------------------------------------

\begin{rem}\label{large_core_L2_rem}
It is possible to prove the frequency envelope estimate \eqref{core_N}
with $\eta=1$ in the case where there is no energy dispersion. As the
previous step shows, one may first reduce to a bilinear estimate.
Then the desired bound follows from summation over \eqref{standard_L2_est_bi}
using the sum rules \eqref{sum_rule1}--\eqref{sum_rule2}. The details
are standard and left to the reader.
\end{rem}

\ret
%---------------------------------------------------------

\begin{proof}[Proof of estimate \eqref{core_N}]
The proof will be accomplished in a series of steps whose goal is to
reduce things to the bilinear estimate \eqref{bal_l2_N}.\ret

\step{0}{A Preliminary Reduction} The first order of business is to
reduce estimate \eqref{core_N} to the case where we replace the
condition on line \eqref{trilin_delta_conds} with a maximal case:
\begin{equation}
    m \ = \ \max\{ \sqrt{\delta_1}|\ln(\eta)|,10\} \ . \label{tri_N_m_cond}
\end{equation}
We claim that a proof of \eqref{core_N} with this choice of $m$
implies \eqref{core_N} for any other choice of $m$ where we turn
\eqref{tri_N_m_cond} into an $\leqslant$ inequality. The only caveat
is that we must replace the multiplier $P_k$ in the definition of
\eqref{basic_trilin} by a version $\td{P}_k$ with a slightly
fattened support, so that one obtains the quasi-idempotence identity
$\td{P}_k P_{[k-5,k+5]} = P_{[k-5,k+5]}$.  To see this, simply
notice that one has the reshuffling identity:
\begin{equation}
    T^{m_0}_k(\phi^{(1)}\!\!,\phi^{(2)}\!\!,\phi^{(3)}) \! =\! T^{m}_k(\phi^{(1)}\!\!,\phi^{(2)}\!\!,\phi^{(3)})
    - \td{T}^{m}_{1;k}(\phi^{(1)}\!\!,\phi^{(2)}\!\!,\phi^{(3)})
    -  \td{T}^{m}_{2;k}(\phi^{(1)}\!\!,\phi^{(2)}\!\!,\phi^{(3)})  \ , \notag
\end{equation}
for any $10\leqslant m_0\leqslant m$, where the
$\td{T}^{m}_{i;k}(\phi^{(1)},\phi^{(2)},\phi^{(3)})$ are the
trilinear forms obtained from applying the definition of $T^m_k$,
with $\td{P}_k$ instead of $P_k$, to the second and third terms
(resp.) on the RHS of \eqref{basic_trilin} in the definition of
$T^{m_0}_k$.

\ret \step{1}{Removal of the Commutator} We are now trying to prove
\eqref{core_N} under the condition \eqref{tri_N_m_cond}.  Our next
step is to use \eqref{commutator} to write \eqref{basic_trilin} in the
form:
\begin{equation}
\begin{split}
    \!\! T^{m}_k(\phi^{(1)},\phi^{(2)},\phi^{(3)})   \ &= \
        P_k \big(\phi^{(1)}\partial^\alpha\phi^{(2)} \partial_\alpha\phi^{(3)}\big)
        - P_k \big(\phi^{(1)}_{<k-m}\partial^\alpha\phi^{(2)}_{<k-m} \partial_\alpha\phi^{(3)}\big)
    \\ \    -  & P_k  \big(\phi^{(1)}_{<k-m}\partial^\alpha\phi^{(2)} \partial_\alpha\phi^{(3)}_{<k-m}\big)
    -  (C_1 +C_2+C_3+C_4)  \ ,
\end{split}\label{big_Tm_red}
\end{equation}
where the commutator terms $C_1$, $C_2$, $C_3$ and $C_4$ have the form
(here the $L_i$ refers to the effect of the commutator, and not
the $L$ in the definition of \eqref{basic_trilin} which has been
dropped):
\begin{align*}
    C_1  \ &= \ 2^{-k} L_1(\nabla_x \phi^{(1)}_{<k-m},\partial^\alpha\phi^{(2)}_{<k-m} ,
    \partial_\alpha\td{P}_k\phi^{(3)}) \ , \\
    C_2 \ &= \  2^{-k}L_2(\phi^{(1)}_{<k-m},\nabla_x\partial^\alpha\phi^{(2)}_{<k-m} ,
    \partial_\alpha\td{P}_k\phi^{(3)})\ ,  \\
    C_3 \ &= \  2^{-k}L_3(\nabla_x\phi^{(1)}_{<k-m},\partial^\alpha\td{P}_k\phi^{(2)} ,
    \partial_\alpha\phi^{(3)}_{<k-m})\ , \\
    C_4 \ &= \  2^{-k}L_4(\phi^{(1)}_{<k-m},\partial^\alpha\td{P}_k\phi^{(2)},\nabla_x
    \partial_\alpha\phi^{(3)}_{<k-m}) \ .
\end{align*}

Here $\td{P}_k$ is the same as in the previous step, and we remind the reader that the
$L_i$ are disposable. The goal of this step is to prove the estimates:
\begin{equation}
    \lp{C_i}{N[I]} \ \lesssim \ \eta^{\delta_1} c_k \ , \label{comm_Li_bound}
\end{equation}
which suffice to establish \eqref {core_N} for all but the first
three terms on the RHS of the equation for $T^{m}_k$ above. It
suffices to work with the case of $i=3,4$; the cases $i=1,2$ are
similar but simpler because the frequency envelope is on the high
term.\ret

For the trilinear form $C_3$ we decompose into all possible frequencies
and use \eqref{standard_est_tri}, which gives:
\begin{equation}
    \lp{C_3}{N[I]}
   \ \lesssim \ 2^{-k}\sum_{k_1,k_3< k-m} 2^{k_1}2^{-\delta(k_1-k_3)_+}  c_{k_3}
     \ \lesssim \  2^{-m} c_{k-m}  \ \lesssim  \ 2^{-(1-\delta_0)m} c_k \ , \notag
\end{equation}
which suffices to show \eqref{comm_Li_bound} in light of the definition \eqref{tri_N_m_cond}
for $m$.

To prove the bound \eqref{comm_Li_bound} for $C_4$ we only split $\phi^{(3)}$ into
separate frequencies, and we use \eqref{N_prod_est1} and \eqref{standard_est_bi} to bound:
\begin{equation}
     \lp{C_4}{N[I]}
     \ \lesssim \  2^{-k}\sum_{k_3< k-m} 2^{k_3}  c_{k_3} \
     \ \lesssim \  2^{-m}c_{k-m}  \ \lesssim  \ 2^{(1-\delta_0)m} c_k \ . \notag
\end{equation}

%----------------------------------------------------------------------------------------------

\ret
\step{2}{Reduction to Matched Frequencies}
We are now trying to bound the sum of the first three terms on the RHS of \eqref{big_Tm_red} above.
Here we write:
\begin{equation}
    (\hbox{First three terms on R.H.S.})\eqref{big_Tm_red} \ = \ A_1 + A_2+ A_3 + B_1 + B_2  \ , \notag
\end{equation}
where $A_1$, $A_2$ and $A_3$ account for the unmatched frequency interactions:
\begin{align*}
        A_1 \ &= \ \sum_{k_1 \geqslant k-m} \ \sum_{\max\{k_2,k_3\} \geqslant k+m}
    P_k \big(\phi^{(1)}_{ k_1}
        \partial^\alpha\phi^{(2)}_{ k_2} \partial_\alpha\phi^{(3)}_{ k_3}\big) \ , \\
    A_2 \ &= \ \sum_{k_1 \geqslant k-m} \ \sum_{\min\{k_2,k_3\} \leqslant k-2m}
    P_k \big(\phi^{(1)}_{ k_1}
        \partial^\alpha\phi^{(2)}_{ k_2} \partial_\alpha\phi^{(3)}_{ k_3}\big) \ , \\
    A_3 \ &= \ \sum_{\max\{k_2,k_3\} \geqslant k+m}
    P_k \big(\phi^{(1)}_{< k-m}
        \partial^\alpha\phi^{(2)}_{ k_2} \partial_\alpha\phi^{(3)}_{k_3}\big) \ ,
%   A_3 \ &=  \ \sum_{k_1 \geq k-m} \ \sum_{k_2,k_3 \leq k+m}
%   P_k \big(\phi^{(1)}_{ k_1}
%   \partial^\alpha\phi^{(2)}_{ k_2} \partial_\alpha\phi^{(3)}_{ k_3}\big) \ ,
\end{align*}
while $B_1$ and $B_2$ account for the matched frequency interactions:
\begin{align*}
        B_1 \ &= \    \sum_{k-2m < k_2,k_3 < k+m}
    P_k \big(\phi^{(1)}_{\geqslant k-m}
        \partial^\alpha\phi^{(2)}_{k_2 }
        \partial_\alpha\phi^{(3)}_{ k_3}\big)
        \ , \notag\\
        B_2 \ &= \ \sum_{k-m < k_2,k_3 < k+m} P_k \big(\phi^{(1)}_{<k-m}\partial^\alpha\phi^{(2)}_{k_2}
        \partial_\alpha\phi^{(3)}_{k_3}\big)
        \ . \notag
\end{align*}
The goal of this step is to prove the set of estimates:
\begin{equation}
    \lp{A_{i}}{N[I]} \ \lesssim \ 2^{-\frac{1}{2}\delta m} c_k\ , \label{Aij_est}
\end{equation}
which is sufficient to establish \eqref{core_N} for these terms because of the
definition \eqref{tri_N_m_cond}.

To prove \eqref{Aij_est} for the term $A_{1}$  we use \eqref{standard_est_tri}.
The two highest frequencies can only differ by $O(1)$, therefore we
get three distinct contributions if the highest pairs are $\{12\}$,  $\{13\}$,  or $\{23\}$
respectively:
\begin{align*}
    \lp{ A_1}{N[I]}  \ &\lesssim \ \sum_{k_2 \geqslant k+m} \sum_{k_3 \leqslant k_2}\!\!
    2^{\delta(k_3-k_2)} 2^{\delta(k-k_2)} c_{k_3}
    + \!\!
    \sum_{k_3 \geqslant k+m}  \sum_{k_2 \leqslant k_3} \!\!\!    2^{\delta(k_2-k_3)} 2^{\delta(k-k_3)} c_{k_3}
    \\ & + \sum_{k_3 \geqslant k+m}  \sum_{k_1 =k-m}^{k_3} 2^{\delta(k-k_3)} c_{k_3}
     \ \lesssim  \ m 2^{-\frac{3}{4}\delta m} c_{k+m} \ \lesssim \ 2^{-\frac{1}{2}\delta m} c_k \ .
\end{align*}

In the case of the term $A_{2}$ we must have either the condition
$\max\{k_2,k_3\}>k-10$, or the conditions $\max\{k_2,k_3\}\leqslant
k-10$ and $k_1 >k-10$. This gives two distinct contributions using
estimate  \eqref{standard_est_tri}, which after summing out the $k_1$
index  may be (resp) written as:
\begin{equation}
     \lp{ A_2}{N[I]} \ \lesssim \ S_1 + S_2 \ , \notag
\end{equation}
with:
\begin{align}
    S_1 \ &= \ \sum_{\min\{k_2,k_3\}<k-2m}\ \sum_{\max\{k_2,k_3\}> k-10}
     2^{\delta(k-\max\{k_2,k_3\})}2^{\delta(\min\{k_2,k_3\}-k+m)}c_{k_3}
     \ , \notag\\
    S_2 \ &= \ \sum_{\min\{k_2,k_3\}<k-2m}\ \sum_{\max\{k_2,k_3\}\leqslant k-10}
    2^{\delta(\min\{k_2,k_3\}-k)}c_{k_3} \ . \notag
\end{align}
For the sum $S_1$ we split into cases depending on which index is minimal, and then sum out
$k_2$ which yields:
\begin{equation}
    S_1 \! \lesssim \!\!\!\!\!\! \sum_{k_3<k-2m} 2^{\delta(k_3-k+m)}c_{k_3} + 2^{-\delta m}\!\!\!\!\sum_{k_3>k-10}
    2^{\delta(k-k_3)}c_{k_3} \! \lesssim \!
    2^{-\delta m}(c_{k-2m} + c_k) \! \lesssim \! 2^{-\frac{1}{2}\delta m}c_k
    \ . \notag
\end{equation}
For the sum $S_2$ we again split into cases depending on which index is minimal:
\begin{equation}
    S_2 \lesssim \ \sum_{k_3<k-2m}
    (k-k_3)2^{\delta(k_3-k)}c_{k_3} + \sum_{k_2<k-2m}\ \sum_{k_3=k_2}^{k-10}
    2^{\delta(k_2-k)}c_{k_3} \ . \notag
\end{equation}
For the first sum on the RHS above we get $2^{-\delta m}c_{k-2m}$ which is acceptable.
For the second sum we further split the range into $k_3<k-2m$ and $k-2m\leqslant k_3 < k-10$.
In the first case we again get $2^{-\delta m}c_{k-2m}$, while in the second
we are left with $m 2^{-2\delta m}\sup_{k-2m\leqslant k_3<k-10}c_{k_3}$, which again suffices.
%\begin{equation}
%   2^{-2\delta m}\sum_{k_3=k-2m}^{k-10} c_{k_3} \ \lesssim \
%\end{equation}

Finally, in the term $A_3$ we must have $|k_2 - k_3|<10$ and only the
frequency $2^k$ part of the null form $\partial^\alpha\phi^{(2)}_{ k_2}
\partial_\alpha\phi^{(3)}_{ k_3}$ will contribute. Then we use
\eqref{standard_est_bi} for the null form, combined with \eqref{N_prod_est1}:
\begin{equation*}
    \lp{A_3}{N[I]}   \ \lesssim\  \sum_{k_3 > k+m} 2^{\delta(k-k_3)} c_{k_3}
     \ \lesssim  \ 2^{-\delta m} c_{k+m}  \ \lesssim \
     2^{-\frac{1}{2}\delta m} c_k \ .
\end{equation*}

%----------------------------------------------------------------------------------

\ret
\step{3}{The Matched Frequency Estimate}
After the last step, it remains to bound the remaining  two terms $B_i$. In both  cases,
by an application of either \eqref{N_prod_est1} or \eqref{N_prod_est2},
 we only need to show the more general matched frequency estimate:
\begin{equation}
    \lp{\partial^\alpha\phi^{(2)}_{[k-O(m),k+O(m)]}
    \partial_\alpha\phi^{(3)}_{ [k-O(m),k+O(m)]}}{N[I]}
    \ \lesssim \ \eta^{\delta_1} c_k
        %+\eta^{-C\delta_1} \underline{c}_k^\delta c_k^{1-\delta}
    \ , \label{core_N_RED}
\end{equation}
under the conditions of Proposition \ref{big_tri_prop}. Using the bound on
$m$ \eqref{tri_N_m_cond} and the definition \eqref{freq_env_defn},
it suffices to establish the fixed frequency estimate:
\begin{equation}
    \lp{\partial^\alpha\phi^{(2)}_{k_2 }
    \partial_\alpha\phi^{(3)}_{ k_3}}{N[I]}
    \ \lesssim \ 2^{Cm} \eta^\delta c_{k_3}
    %+ \underline{c}_{k_3}^\delta c_{k_3}^{1-\delta} )
    \ , \notag
\end{equation}
where we are restricting $|k_2-k_3|\lesssim m$. This follows
immediately from \eqref{bal_l2_N}.
\end{proof}\ret

%---------------------------------------------------------

\begin{rem}\label{large_core_N_rem}
We remark here that one may prove estimate \eqref{large_core_N} by a
quick application of the above work. To see this, notice
the above proof up to \textbf{Step 3} does not use
Proposition \ref{p:matchfreq}. Thus, we are left with showing estimate
\eqref{core_N_RED} in this case, and by inspection of
\textbf{Step 2} we may assume the gap between $k_2$ and $k_3$ is no larger than $3m$.
By directly applying estimate \eqref{standard_est_bi} we have
\eqref{large_core_N} in this case.
%First, assume that $m=10$. Then
%by following the process of the above proof through \textbf{Step 3},
%with line \eqref{tri_N_m_cond} replaces with $m=10$, we have the
%desired bound once we estimate the contribution on line
%\eqref{core_N_RED} with $m=10$. This follows immediately from
%\eqref{standard_est_bi}.\ret

%To prove estimate \eqref{large_core_N} for larger $m$, we may reduce
%to the case of $m=10$ by estimating differences of the form:
%\begin{align}
%        T \ &= \
%        L(\phi^{(1)}_{<k-10},\partial^\alpha\phi^{(2)}_{<k-10}, \partial_\alpha\phi^{(3)}_k)
 %       -L(\phi^{(1)}_{<k-m},\partial^\alpha\phi^{(2)}_{<k-m}, \partial_\alpha\phi^{(3)}_k)
 %        \ , \notag\\
 %        &= \ L(\phi^{(1)}_{<k-10},\partial^\alpha\phi^{(2)}_{[k-m,k-10]},
%         \partial_\alpha\phi^{(3)}_k) +
%         L(\phi^{(1)}_{[k-m,k-10]},
%         \partial^\alpha\phi^{(2)}_{<k-m},
%         \partial_\alpha\phi^{(3)}_k) \ , \notag\\
%         &= \ T_1 + T_2 \ , \notag
%\end{align}
%with a similar set of expressions where the roles of $\phi^{(2)}$ and $\phi^{(3)}$ are switched.
%Estimate \eqref{large_core_N} for the first term $T_1$ follows
%directly from \eqref{N_prod_est1} and then summing over
%\eqref{standard_est_bi} at a loss of $O(m)$. The bound for $T_2$ is
%even better, and there is no $O(m)$ loss by an application of
%\eqref{standard_est_tri}. The details are left to the reader.
%The case of these decompositions with  $\phi^{(3)}$ at low frequency
%is similar, except that this time one looses $O(2^{\delta_0 m})$
%because the frequency envelope $\{c_k\}$ is on the low frequency term.
\end{rem}\ret

%---------------------------------------------------------

\begin{proof}[Proof of estimate \eqref{ext_N}]
The proof of this estimate follows from some simple manipulations of
the bounds used to produce \eqref{core_N}. A quick review of the previous
proof shows that all bounds were achieved with RHS
$\lesssim \eta^{\frac{1}{2}\sqrt{\delta_1}\delta}$.
Thus, by a direct application of those bounds and using the $(\delta_0,\delta_0)$
variance condition on $\{c_k\}$ we have:
\begin{equation}
         \lp{T^m_k(P_{<k+m}\phi^{(1)},\phi^{(2)},\phi^{(3)})}{N[I]} \ \lesssim \
        \eta^{\delta_1}\big(c_k + \lp{P_{<k}\phi^{(1)}}{S[I]}\big)  \ ,
        \notag
\end{equation}
where $m$ is from line \eqref{tri_N_m_cond}.\ret

To bound the contribution with a $P_{\geqslant
k+m}\phi^{(1)}$ factor we directly apply \eqref{standard_est_tri}
which yields the sum:
\begin{align}
    \lp{T^m_k(P_{>k+m}\phi^{(1)},\phi^{(2)},\phi^{(3)})}{N[I]}  &\lesssim \!\!
    \sum_{k_i: \ k_1>k+m} 2^{-\delta(\max\{k_i\}-k)}2^{-\delta(k_1-\min\{k_2,k_3\})_+}c_{k_1}\notag\\
    &\lesssim \sum_{k_1: \  k_1>k+m} 2^{-\delta(k_1-k)}c_{k_1} \ \lesssim \ 2^{-(\delta-\delta_0) m}c_k \ , \notag
\end{align}
which suffices.
\end{proof}\ret

\begin{rem}\label{large_ext_N_rem}
To prove \eqref{large_ext_N} we follow a similar procedure as in the previous proof,
except this time applied to estimate \eqref{large_core_N} instead of \eqref{core_N}.
Here it suffices to split cases according to $P_{< k+10}\phi^{(1)}$ or
$P_{\geqslant k+10}\phi^{(1)}$ contributions. The details are left to the reader.
\end{rem}

\ret
%-------------------------------------------------------------------------
%%%%%%%%%%%%%%%%%%%%%%%%%%%%%%%%%%%%%%%%%%%%%%%%%%%%%%%%%%%%%%%%%%%%%%%%%%
%-------------------------------------------------------------------------

\section{The Gauge Transformation}
\label{s:renormalize}

In this section we prove Proposition~\ref{propphiu}.  The proof is
divided into several portions which deal with different aspects of the
problem. \ret

%For convenience, we employ a locally defined small constant
%$0 < \delta_2\leqslant \frac{1}{2}\delta$, where $\delta$ is
%as in Propositions \ref{standard_prop1} and \ref{standard_prop2}.
%We will show the bounds \eqref{env_est1}--\eqref{env_est0}
%with this $\delta_2$, and then one may globally redefine $\delta$
%in Propositions  \ref{standard_prop1} and \ref{standard_prop2} to uniformize notation.

\subsection{ Bounds for B}
Here we transfer the bounds from $\phi$ to $B$. Precisely, we
have:

\begin{lem}
Let $\phi$ be a wave map as in Proposition \ref{propphiu}. Then the
matrix $B$ defined on line \eqref{B_def} has an antisymmetric
extension off the interval $I$, which satisfies the following global bounds:
\begin{align}
    \lp{\nabla_{t,x}B_k}{L^\infty_t(L^2_x)} \ &\lesssim_E \ \td{c}_k \ , \label{bken}\\
    \lp{B_k}{S \cap \uX}  \ &\lsf  \ c_k \ , \label{bksux}\\
    \lp{\Box B_{k'}\cdot \psi_k }{N} \ &\lsf  \ 2^{\delta (k'-k)}c_{k'}\lp{\psi_k}{S} \ ,
    \ \ \ \ k' < k-10 \ .  \label{boxbkpsi}
\end{align}
\end{lem}\ret

\begin{proof}
By definition we have that:
\begin{equation}
        B_k \ = \ \mathcal{S}(\phi)_{<k-10}\phi_k \ , \notag
\end{equation}
where $\mathcal{S}$ is the antisymmetric part of the original second
fundamental form. The bound \eqref{bken} on $I$ follows from the
same bound for $\phi$ combined with Leibnitz's rule. Furthermore,
an $S[I]$ norm bound as on line \eqref{bksux} follows from the
algebra property \eqref{S_prod2} combined with the
Moser estimate \eqref{basic_moser_c}.\ret

For the estimates involving $\Box B_k$ we remark that the function
$\mathcal{S}(\phi)$ solves a similar wave equation to $\phi$ on the
interval $I$, which we write schematically as:
\begin{equation}
        \Box \mathcal{S}(\phi)  \ = \ \mathcal{F}(\phi) \partial^\alpha \phi
        \partial_\alpha \phi \ . \notag
\end{equation}
By the version of estimate \eqref{core_L2} in Remark
\ref{trilinear_rem}, we have the pair of bounds:
\begin{equation}
       \lp{\Box P_k \mathcal{S}(\phi)}{L^2_t(L^2_x)[I]}
        + \lp{\Box P_k \phi}{L^2_t(L^2_x)[I]} \lesssim_F 2^{\frac{k}2} c_k
        \label{boxsl2}
\end{equation}
By Leibnitz's rule we have:
\begin{equation}
        \Box B_k = \Box \mathcal{S}(\phi)_{<k-10} \phi_{k}
        + 2  \partial^\alpha \mathcal{S}(\phi)_{<k-10} \partial_\alpha \phi_{k}
        + \mathcal{S}(\phi)_{<k-10}\Box \phi_{k} \ .
        \label{Box_B_exp}
\end{equation}
Hence using \eqref{boxsl2} for the first and last term, and again
using Remark \ref{trilinear_rem} for the null form in the middle
term, we obtain an $\uX[I]$ bound as on line \eqref{bksux}.\ret

It remains to prove the estimate \eqref{boxbkpsi} localized to $I$.
We again use the expression \eqref{Box_B_exp} for $\Box B_k$.
We multiply the RHS of this line by a function $\psi_j$ of frequency
$j > k+10$. The contribution of the middle term can be estimated by
\eqref{standard_est_tri}:
\begin{equation}
    \lp{   \partial^\alpha \mathcal{S}(\phi)_{<k-10} \partial_\alpha \phi_{k}
    \psi_j}{N[I]} \  \lesssim_F  \
        2^{-\delta(j-k)}c_{k} \lp{ \psi_j}{S[I]} \ . \notag
\end{equation}
In both other cases, by the $S\cdot N$ algebra property,  it suffices to prove the estimate:
\begin{equation}
        \lp{ P_{k}\big( \phi^{(1)} \partial^\alpha \phi^{(2)} \partial_\alpha \phi^{(3)}\big)
        \psi_j}{N[I]} \  \lesssim_F  \
        2^{-\delta(j-k)}c_{k} \lp{ \psi_j}{S[I]} \ , \notag
\end{equation}
for any set of test functions $\phi^{(i)}$ with $S[I]$ norm of size $F$, and
frequency envelopes $\{c_k\}$. To show this, we let $T^{10}_k$ be
the trilinear form defined on line \eqref{basic_trilin}, built up
out of the $\phi^{(i)}$  in the above estimate. Then by a
combination of \eqref{large_core_N} and estimate \eqref{N_prod_est2}
we have:
\begin{equation}
        \lp{ T^{10}_k\cdot \psi_j}{N[I]} \  \lesssim_F  \
        2^{-\delta(j-k)}c_{k} \lp{ \psi_j}{S[I]} \ . \notag
\end{equation}
It remains to show a bound of the form:
\begin{align}
        \lp{  \psi_j\cdot \phi^{(1)}_{<k-10} \partial^\alpha \phi^{(2)}_{<k-10}
        \partial_\alpha \phi^{(3)}_k }{N[I]} \  \lesssim_F  \
        2^{-\delta(j-k)}c_{k} \lp{ \psi_j}{S[I]} \ , \notag
\end{align}
which encapsulates the remainder from $T^{10}_k$ (there are two such
remainders, but they are essentially symmetric). This bound follows
immediately from by applying the algebra estimate \eqref{S_prod2} to
the first two terms, and then summing the
resulting trilinear via estimate \eqref{standard_est_tri}.\ret

To conclude our proof of the estimates
\eqref{bken}--\eqref{boxbkpsi} we simply need to extend off the
interval $I$ in a simultaneous way. To do this we use the canonical
extension defined in  Proposition \ref{ext_prop}.
\end{proof}\ret

%------------------------------------------------------------------------

\subsection{The gauge construction}
Here we construct the gauge transformation $U$ and obtain estimates
on $U$ in $S$ and $\uX$. For comparison purposes we note that in the
small data results of \cite{Tao_WM}, \cite{Tataru_WM1} the function $U$ is
constructed iteratively by setting:

\begin{equation}
    U_{<k}  \ =  \ \sum_{k'<k} U_{k'} \ , \qquad
    U_k \ = \ U_{< k-C}  B_k \ , \notag
\end{equation}
with $k \in \Z$. This insures that $U_k$ are localized at frequency
$2^k$, while  the  smallness of $\phi$ in $S$ is used to prove that
$\Ut U - I$ is small.
\ret

Such a construction is no longer satisfactory here, as $\phi$ can be
large in $S$ and thus $U$ may fail to be almost orthogonal. Instead
we switch to a continuous version of the above construction where we
seek $U$ and its ``frequency localized'' version $U_{,<k}$ in the
integrated form:
\begin{equation}
        U  \ =  \ \int_{-\infty}^\infty U_{,k} dk \ , \qquad
        U_{,<j} = \int_{-\infty}^j U_{,k} dk \ , \notag
\end{equation}
where each $U_{,k}$ is defined by:
\begin{equation}
        U_{,k}  \ =  \ U_{,< k} B_k \ . \notag
\end{equation}
In other words, $U_{,<k}$ solves the Cauchy problem:
\begin{equation}
        \frac{d}{dk} U_{,<k}  \ =  \ U_{,<k} B_k \ , \qquad
        U_{,<-\infty}  \ =  \ I_N \ .
        \label{uode}
\end{equation}
Owing to the antisymmetry of the $B_k$, solutions to this ODE enjoy
the conservation law $U_{,<k}U^\dagger_{,<k}=I_N$, so they are
automatically \emph{exactly} orthogonal. However, the price one pays
is that the exact frequency localization of each $U_{,k}$ is lost.
In spite of this, we will prove that $U_{,k}$ is approximately
localized at frequency $2^k$ modulo rapidly decreasing tails:
\begin{equation}
        \lp{P_{j} U_{,k}}{S}  \ \lesssim_F  \ 2^{-\delta|k-j|} c_k \ . \label{pjuks}
\end{equation}
Note however that arbitrarily high frequencies are immediately
introduced, and their evolution is not easy to track. In particular
a bootstrap argument for the above $S$ norm bound would seem to fail
due to the lack of smallness of the $B_k$'s. We proceed with the
proof in several steps aimed at building up to the full $S$ norm
estimate by using the conservation law of \eqref{uode} in a crucial
way:
\ret

%------------------------------------------------------------------------

\step{1}{$L^\infty_t(L^\infty_x)$ and $L^\infty_t( L^2_x)$ bounds
for $U_{,k}$} We will work exclusively with the energy frequency
envelope $\{\td{c}_k\}$ for $B$ in this step. Without loss of
generality we may assume that this is bounded by the $S$ norm
frequency envelope $\{c_k\}$. We start with the pointwise and energy
bounds:
\begin{equation}
        \lp{B_k}{L^\infty_t(L^\infty_x)} \lse \td{c}_k  , \
        \lp{ \nabla_{t,x} B_k}{L^\infty_t(L^\infty_x)} \lse 2^k \td{c}_k ,
        \
        \lp{ \nabla_{t,x} B_k}{L^\infty_t( L^2_x)}  \lse   \td{c}_k , \label{philil2}
\end{equation}
all derived from \eqref{bken}. We claim that this implies the
following energy type for $U_{,k}$ itself:
\begin{equation}
        \lp{P_{k'} \nabla_{t,x} U_{,k} }{L^\infty_t(L^2_x)}  \ \lse  \
        2^{-|k'-k|}2^{C(k-k')_+} \td{c}_k \ . \label{ulil2}
\end{equation}
To show this, notice that by construction of $U_{,<k}$ we
immediately obtain:
\begin{equation}
        \lp{ U_{,k}}{L^\infty_t(L^\infty_x)}  \ \lse  \ \td{c}_k \ , \qquad
        \lp{U_{,k} }{L^\infty_t( L^2_x)} \ \lse  \ 2^{-k} \td{c}_k \ .
        \notag
\end{equation}
We estimate $\nabla_{t,x} U_{,k}$ by differentiating \eqref{uode}:
\begin{equation}
        \frac{d}{dk} \nabla_{t,x} U_{,<k}  \ =  \
        \nabla_{t,x} U_{,<k}\cdot B_k +  U_{,<k}\nabla_{t,x} B_k \ ,
        \qquad \nabla_{t,x} U_{,<-\infty}  \ = \  0 \ . \label{dUk}
\end{equation}
In view of the second estimate on line \eqref{philil2}, we have good
bounds for the second term on the RHS of the above expression, and
we wish to transfer these to $\nabla_{t,x} U_{,<k}$. In order to do
this, we employ the following device that will be used repeatedly in
the sequel:

\begin{lem}[Unitary Variation of Parameters]\label{var_param_lem}
Let $V_{,<k}$ be given by the ODE:
\begin{equation}
             \frac{d}{dk} V_{,<k}  \ =  \
             V_{,<k} B_k +  W_{,k} \ \qquad
             V_{,<-\infty} \ = \ W_{,-\infty} \ = \ 0 \ , \label{V_form}
\end{equation}
where $B_k$ is antisymmetric and the forcing term $W_{,k}$ is
arbitrary. Then in any mixed Lebesgue space $L^q_t(L^r_x)$ space we
have the following bound:
\begin{equation}
        \lp{ V_{,<k}}{L^q_t(L^r_x)} \ \lesssim \ \int_{-\infty}^k
        \lp{W_{,k'}}{L^q_t(L^r_x)}dk' \ . \label{unitary_var_param}
\end{equation}
\end{lem}\ret

\begin{proof}
We write the formula for $V_{,<k}$ via variation of parameters as
follows:
\begin{equation}
        V_{,<k} \ = \ \int_{-\infty}^k \mathcal{P}(k,k')W_{,k'}dk' \
        , \notag
\end{equation}
where $\mathcal{P}$ is the propagator of the unitary problem:
\begin{equation}
        \frac{d}{dk} \mathcal{P}(k,k') \ = \ B_k\mathcal{P}(k,k')  \
        , \qquad
        \mathcal{P}(k',k') \ = \ I_N \ . \notag
\end{equation}
In particular,
$\lp{\mathcal{P}(k,k')}{L^\infty_t(L^\infty_x)}\lesssim 1$. The proof is
concluded via an application of Minkowski's integral inequality.
\end{proof}\ret

We now use estimate \eqref{unitary_var_param} to integrate
\eqref{dUk}, which yields:
\begin{equation}
        \lp{ \nabla_{t,x} U_{,<k}}{L^\infty_t(L^\infty_x)}  \ \lse  \ 2^{k}
        \td{c}_k \ , \notag
\end{equation}
through a direct application of the sum rule \eqref{sum_rule1}. From
the differentiated equation for $U_{,k}$ this shows that:
\begin{equation}
        \lp{\nabla_{t,x} U_{,k}}{L^\infty_t( L^2_x)}  \ \lse  \  \td{c}_k
        \ . \notag
\end{equation}
Repeating the process for all possible spatial derivatives of
$\nabla_{t,x} U_{,<k}$ we get the inductive bounds:
\begin{equation}
        \lp{ \nabla_x^J \nabla_{t,x} U_{,<k}}{L^\infty_t(L^\infty_x)}
         \lse   2^{(|J|+1) k} \td{c}_k \ , \quad
        \lp{  \nabla_x^J \nabla_{t,x}  U_{,k}}{L^\infty_t(L^2_x)}
         \lse   2^{|J| k} \td{c}_k \ . \label{highebd}
\end{equation}
The second relation shows in particular that:
\begin{equation}
        \lp{P_{k'} \nabla_{t,x} U_{,k} }{L^\infty_t( L^2_x)}
        \ \lse  \ 2^{C(k-k')} \td{c}_k \ , \notag
\end{equation}
for any positive constant $C$, and therefore by integration that:
\begin{equation}
        \lp{ P_{k'} \nabla_{t,x} U_{,<k}}{L^\infty_t( L^2_x)}
        \ \lse \ 2^{C(k-k')} \td{c}_k \ , \notag
\end{equation}
which suffices for \eqref{ulil2} if $k' \geqslant k$.

It remains to bound the low frequencies in $U_{,k}$, and so we
write:
\begin{equation}
        P_{k'} U_{,k}  \ =  \
        P_{k'} ( P_{[k-10, k+10]}U_{<k}\cdot B_k ) \ , \qquad
        k'  \ <  \ k-20 \ . \notag
\end{equation}
Using Bernstein's inequality \eqref{Bernstein} we obtain:
\begin{equation}
        \lp{P_{k'} \nabla_{t,x} U_{,k}}{L^\infty_t( L^2_x)}
        \ \lesssim\  2^{k'} \lp{ \nabla_{t,x} ( P_{[k-10,k+10]}
        U_{,<k}B_k)}{L^\infty_t( L^1_x)} \ \lse  \ \td{c}_k 2^{k'-k} \ .
        \notag
\end{equation}
Hence \eqref{ulil2} is proved.
\ret

%-------------------------------------------------------------------------------

\step{2}{Strichartz bounds for $U_{,k}$} This section largely mimics
the previous one, so we will be more terse here. By \eqref{bksux} we
have the Strichartz bounds:
\begin{equation}
        \lp{ B_k}{D\underline{S}_k} \ \lsf \ 2^{-k}c_k \ , \qquad
        \lp{ \nabla_{t,x} B_k}{D\underline{S}_k} \ \lsf \  c_k  \ ,
        \notag
\end{equation}
where we recall that $D\underline{S}_k$ is the space of Strichartz
admissible $L^q_t(L^r_x)$ norms from line \eqref{dphys_str} with
appropriate dyadic weight (note that this norm does \emph{not}
include frequency localization, which will be notationally useful
here).

Using the bounds for $B_k$ with equation \eqref{dUk} or its
derivatives, we directly have:
\begin{equation}
        \lp{\nabla_x^J \nabla_{t,x}U_{,k}}{D\underline{S}_k}
        \ \lesssim_F \ 2^{|J|k}c_k \ ,
        \qquad |J|\geqslant 0 \ . \notag
\end{equation}
By using this last set of estimates for high frequencies, and
\eqref{ulil2} and Bernstein's inequality for low frequencies, we
have:
\begin{equation}
        \lp{P_{k'}\nabla_{t,x}U_{,k} }{D\underline{S}_k} \
        \lesssim_F \ 2^{-|k'-k|} 2^{C(k'-k)_+}c_k \ . \notag
\end{equation}
In particular, one has the inequality:
\begin{equation}
        \lp{P_{k'} \nabla_{t,x} U_{,k} }{L^4_t(L^\infty_x)}
        \ \lsf \
        2^{\frac{3}{4}k}2^{-|k'-k|}c_k \ ,  \label{ul4li}
\end{equation}
which will be useful later in this section. Finally, by
interpolating this last bound with \eqref{ulil2} and recalling the
definition from line \eqref{phys_str}, we have the following
$\underline{S}$ norm portion of estimate \eqref{pjuks}:
\begin{equation}
        \lp{P_{k'}U_{,k}}{\underline{S}} \ \lesssim_F \
        2^{-\frac{1}{4}|k'-k|}c_k \ . \notag
\end{equation}\ret

%By directly using the $L^\infty_t(L^\infty_x)$ bound for $U_{,<k}$,
%or by integrating the derivative equation \eqref{dUk} via Lemma
%\ref{var_param_lem}, we have:
%\begin{equation}
%        \lp{U_{,k}}{\underline{S}} \ \lsf \ c_k \ , \qquad
%        \lp{\nabla_x^J \nabla_{t,x} U_{,k}}{\underline{S}}
%        \ \lsf \ 2^{(|J|+1)k} c_k \ . \notag
%\end{equation}
%As above we can bound higher frequencies directly, and lower
%frequencies via Bernstein's inequality, to obtain the full estimate:
%\begin{equation}
 %      \lp{P_{k'} \nabla_{t,x} U_{,k} }{\underline{S}}  \ \lsf  \
 %       2^{k-|k'-k|} 2^{C(k'-k)_+}c_k  \ . \notag
%\end{equation}
%for any $C\geqslant 0$. Because it will be particularly used in the
%sequel, we state the following two specific instances of the last
%line above: \ret

\step{3}{High modulation bounds for $U_{,k}$} Here we will show
that:
\begin{equation}
        \lp{P_{j} \Box U_{,k} }{L^2_t(L^2_x)} \ \lsf \
        2^{\frac{k}2} 2^{-(\frac12+\delta) |k-j|}c_k \ .
        \label{uhighmod}
\end{equation}
Differentiating the equation for $U_{,<k}$ we obtain the evolution
equation for $ \Box U_{,<k}$:
\begin{equation}
        \frac{d}{dk} \Box U_{,<k}  \ =   \ \Box U_{,<k} B_k + U_{,<k} \Box B_{,k} +
        2\partial^\alpha  U_{,<k} \partial_\alpha B_k \ .
        \label{boxueq}
\end{equation}
Our first goal will be to use Lemma \ref{var_param_lem}
to show that:
\begin{equation}
        \lp{\Box U_{,<k}}{L^2_t(L^2_x)} \ \lsf \
        2^{\frac{k}2} c_k \ . \label{boxu<k}
\end{equation}
By estimate \eqref{unitary_var_param} and the $\uX$ control for
$B_k$ from line \eqref{bksux}, it suffices to have the null-form
bound:
\begin{equation}
        \lp{ \partial^\alpha U_{,<k} \partial_\alpha B_k} {L^2_t(L^2_x)}  \
        \lsf\ 2^{\frac{k}2} c_k  \ . \label{boxaunull}
\end{equation}
Expanding the term on the LHS of this last line we have:
\begin{equation}
\begin{split}
        \partial^\alpha  U_{,<k} \partial_\alpha B_k = &\
        \int_{-\infty}^k \partial^\alpha  (U_{,<k'} B_{k'} )  \partial_\alpha
        B_k dk' \ , \\
        = &\  \int_{-\infty}^k U_{,<k'} \partial^\alpha B_{k'} \partial_\alpha  B_k dk'
        + \int_{-\infty}^k   \partial^\alpha U_{,<k'} B_{k'} \partial_\alpha
        B_k dk' \ .
\end{split}\label{expa}
\end{equation}
Estimate \eqref{boxaunull} for the first term on the RHS of this
last line follows by summing over the bound
\eqref{standard_L2_est_bi}. For the second term on the RHS of the
last line above we may take a product of two $L^4_t(L^\infty_x)$
estimates for the terms at frequency $k'$ and $,<k'$, and one energy
type bound for $B_k$. This again yields \eqref{boxaunull}. A similar
argument allows us to prove the analog of estimates
\eqref{boxaunull} and \eqref{boxu<k} for higher spatial derivatives:
\begin{equation}
\begin{split}
        \lp{ \nabla_x^J (\partial^\alpha U_{,<k} \partial_\alpha B_k)}{L^2_t(L^2_x)}
         \lsf &\
        2^{(\frac{1}{2} +|J|) k} c_k \ ,  \\
        \lp{ \nabla_x^J \Box U_{,<k}}{L^2_t(L^2_x)}   \lsf &\
        2^{(\frac{1}{2} +|J|) k} c_k\ . \label{boxaudiff}
\end{split}
\end{equation}\ret

Turning our attention to $U_k$ we have the identity:
\begin{equation}
        \Box U_{,k}  \ =   \ U_{,<k}  \Box B_k+ \Box U_{,<k} B_k  +
        2 \partial^\alpha U_{,<k} \partial_\alpha B_k \ .
        \label{boxukexp}
\end{equation}
By estimates \eqref{boxu<k}, \eqref{boxaunull}, and the analogous
bound for $B_k$ from line \eqref{bksux} we directly have:
\begin{equation}
        \lp{ \Box U_{,k}}{L^2_t(L^2_x)}  \ \lsf  \
        2^{\frac{k}2} c_k \ , \notag
\end{equation}
while the estimates on line \eqref{boxaudiff} combined with the
energy and $L^\infty_t(L^\infty_x)$ bounds for derivatives of
$U_{,<k}$ proved in the first step allow us to prove:
\begin{equation}
        \lp{ \nabla_x^J  \Box U_{,k}}{L^2_t(L^2_x)}
        \ \lsf \  2^{(\frac{1}{2}+|J|)k} c_k \ . \notag
\end{equation}
This suffices to give \eqref{uhighmod} for all but the low
frequencies.
\ret

It remains to obtain improved low frequency bounds, i.e. prove
\eqref{uhighmod} in the case when $j < k-10$.  The first two terms
in \eqref{boxukexp} are easy to estimate, combining the
$L^2_t(L^2_x)$ bound for one factor with the $L^\infty_t( L^2_x)$
energy type bound for the other, while using Bernstein's inequality
at low frequency.

The third term on the RHS of \eqref{boxukexp} has already been
estimated before using \eqref{expa}, but now we need to be more
careful to gain from small $j$.  The first term on the RHS of
\eqref{expa}, call it $T_1$, can at low frequency be split into
three contributions, $P_j T_1 = T_{11} + T_{12} + T_{13}$ where
\begin{align}
    T_{11}   \ &=  \
    P_j \int_{k-4}^k P_{<j+4}U_{,<k'}\cdot
    P_{<j+8} (\partial^\alpha B_{k'} \partial_\alpha  B_k) dk' \ , \notag\\
    T_{12} \ &=  \
    P_j \int_{k-4}^k \int_{j+4}^{\infty} P_l U_{,<k'}
    \cdot
    P_{[l-4,l+4]} (\partial^\alpha B_{k'} \partial_\alpha  B_k) dl dk' \ , \notag\\
    T_{13} \ &= \  P_j \int_{-\infty}^{k-4} \! P_{[k-10,k+10]} U_{,<k'}\cdot
        \partial^\alpha B_{k'} \partial_\alpha  B_k dk' \ . \notag
\end{align}
We explain the estimates for each of these terms. In the case of
$T_{11}$ we bound $P_{<j+4}U_{,<k'}$ in $L^\infty_t(L^\infty_x)$ and
then apply \eqref{standard_L2_est_bi} for the remaining null form.
In the case of $T_{12}$ we use \eqref{ulil2} to bound $P_l U_{,<k'}$
in $L^\infty_t( L^2_x)$, \eqref{standard_L2_est_bi} for the
remaining null form, and then conclude with Bernstein's inequality.
Finally, the bound $T_{13}$ is obtained in the same way as in the
case of $T_{12}$.

Finally, we need to prove the low frequency part of the estimate
\eqref{uhighmod} for the second term on line \eqref{expa} above,
which we denote by $T_2$. This cannot be done directly, because
there is no extra room in the application of Strichartz estimates to
use Bernstein's inequality. Therefore we reexpand as follows:
\begin{equation*}
\begin{split}
        T_2 \ = \ T_{21} + T_{22} \ = &\
        \int_{-\infty}^k \int_{-\infty}^{k_1} U_{,<k_2} \partial^\alpha
        B_{k_2}   B_{k_1} \partial_\alpha B_k dk_2 dk_1\\
   &  + \int_{-\infty}^k \int_{-\infty}^{k_1} \partial^\alpha  U_{,<k_2}
        B_{k_2} B_{k_1} \partial_\alpha  B_k dk_2 dk_1 \ . \notag
\end{split}
\end{equation*}
The first term $T_{21}$ on the RHS above has a structure very
similar to the whole of $T_1$ above. The only new development is
that extra factor of $B_{k_1}$, but it is harmless due to the fact
that its frequency is always greater than the differentiated term
$\partial^\alpha B_{k_2}$. Therefore, one can use the same methods
as in the previous paragraph to bound this term (one could as well
use the procedure we are about to describe for bounding the second
term $T_{22}$). To handle $T_{22}$ above, we split it further as:
\begin{equation*}
\begin{split}
        P_j T_{22} \ = \  T_{221} +  T_{222} \ = &\
        P_j \int_{-\infty}^{k-8} \int_{-\infty}^{k_1} P_{>k-20}\partial^\alpha U_{,<k_2}
        \cdot B_{k_2} B_{k_1} \partial_\alpha  B_k dk_2 dk_1\\
       & + P_j \int_{k-8}^k \int_{-\infty}^{k_1} \partial^\alpha  U_{,<k_2}
        B_{k_2} B_{k_1} \partial_\alpha  B_k dk_2 dk_1 \ . \notag
\end{split}
\end{equation*}
For the first term above we put the two (i.e. first and fourth) high
frequency terms in $L^\infty_t(L^2_x)$, while the middle two terms
are both estimated with $L^4_t(L^\infty_x)$; then we use Bernstein's inequality. One is forced to loose
in the low frequencies this way, but this is made up for by the
arbitrary gain in the difference $(k-k_2)$ coming from estimate
\eqref{ulil2}:
\begin{equation}
        \lp{T_{221}}{L^2_t(L^2_x)}  \lesssim_F \
         c_k \int_{-\infty}^{k-8}\! \! \int_{-\infty}^{k_1}
        2^j\cdot2^{C(k_2-k)}\cdot2^{-\frac{1}{4}k_2}\cdot 2^{-\frac{1}{4}k_1}
        dk_2 dk_1 \
        \lesssim_F \ 2^{\frac{1}{2}j} 2^{\frac{1}{2}(j-k)}c_k . \notag
\end{equation}
To bound the term $T_{222}$ we put both the $k_2$ indexed terms in
$L^4_t(L^\infty_x)$, and the other two factors in
$L^\infty_t(L^2_x)$ while using Bernstein's inequality at low
frequency. This gives the inequality:
\begin{equation}
         \lp{T_{222}}{L^2_t(L^2_x)} \ \lesssim_F \
         c_k \int_{k-8}^k\int_{-\infty}^{k_1} 2^j\cdot2^{\frac{1}{2}k_2}\cdot
        2^{-k_1} dk_2 dk_1 \
        \lesssim_F \  2^{\frac{1}{2}j} 2^{\frac{1}{2}(j-k)} c_k \ . \notag
\end{equation}
This completes our demonstration of the estimate \eqref{uhighmod}.

\ret
%----------------------------------------------------------------------------------

\step{4}{High frequency bounds for $U_{,k}$.}
Here we show that the high frequencies in $U_{,k}$ can be estimated
in a much more favorable way:
\begin{equation}
    \lp{  \nabla_{t,x}^J P_{k'}  U_{,k}}{L^1_t(L^1_x)} \ \lsf\
    2^{(-3+|J|)k}2^{ -C(k'-k)} c_k \ , \qquad k' \ > \ k+10 \ ,
    \label{highpk}
\end{equation}
where $|J|\leqslant 2$. For this we expand with $D = \{k_5 < k_4 <
k_3 < k_2 < k_1 < k\}$:
\begin{equation}
    P_{k'} U_{,k}  \ =  \ P_{k'}
    \int_{D}
    U_{,<k_5} B_{k_5} B_{k_4} B_{k_3} B_{k_2} B_{k_1} B_k\ dk_5
    dk_4 dk_3 dk_2 dk_1 \ . \notag
\end{equation}
Due to the frequency localizations we can replace $U_{,<k_5}$ by
$P_{[k'-10,k'-10]}U_{,<k_5}$, for which we may use the
$L^\infty_t(L^\infty_x)$ bound derived from \eqref{ulil2}. Hence by
the Strichartz estimates alone for the $B_{k_i}$'s we obtain:
\begin{equation}
    \lp{ P_{k'}  U_{,k}}{L^1_t(L^1_x)} \lsf
    c_k \int_{D}
    2^{C(k_5-k') -\frac14(k_5+k_4+k_3+k_2)-k_1-k} c_k dk_5
    dk_4 dk_3 dk_2 dk_1 \ . \notag
\end{equation}
The bound \eqref{highpk} with $J=0$ follows after integration. The
cases $|J|=1$ is treated similarly. A minor variation is needed in
the case $|J|=2$ when two time derivatives occur. There one writes
$\partial_t^2 = \Box + \Delta_x$, using either \eqref{bksux} or
\eqref{boxaudiff} for the factor containing the d'Alembertian.

\ret
%------------------------------------------------------------------------------

\step{5}{Full $S$ norm bounds for $U_{,k}$.} Here we prove that:
\begin{equation}
    \lp{P_{k'} U_{,k} }{S} \ \lsf  \ 2^{-\delta|k-k'| } c_k \ .
    \label{full_U_S}
\end{equation}
In view of the previous step it suffices to consider the case $ k' <
k+10$.

Here we encounter the main difficulty compared to \cite{Tao_WM},
\cite{Tataru_WM1}.  The inductive bound used there grows exponentially in
$k$ due to lack of smallness, so it is useless. Bootstrapping fails
for a similar reason.  Instead we consider iterated expansions. There
are two bounds we need to prove, namely for $\lp{P_{k'}Q_j U_{,k}}{X^{1,\frac12}_\infty}$
and $\lp{ P_{k'} Q_j U_{,k} }{S[k';j]}$.  Due to the high modulation bound
\eqref{uhighmod} and the high frequency bound \eqref{highpk}
 it suffices to consider the case $j < k'-20 < k - 10$.  The key
technical step asserts that in either case we can bound the
contribution of $U_{,< j-20}$ using only pointwise and high modulation
bounds:

\begin{lem}\label{lvvlow}
Let $j < k-10$. Then the following estimate holds for test functions
$u$ and $\phi$:
\begin{equation}
    \lp{ Q_j (  u_{<j-10} \phi_k) }{X^{1,\frac12}_\infty} +
     \lp{ Q_{<j} ( u_{<j-10}  \phi_k)  }{S[k;j]}
     \ \lesssim  \
     \lp{u }{L^\infty_t(L^\infty_x)\cap\uX} \cdot \lp{ \phi_k}{S}
     \ . \label{u<jphik}
\end{equation}
\end{lem}

\begin{proof}
We write
\begin{equation}
    u_{<j-10} \phi_k  \ =  \ Q_{> j-10} u_{<j-10}
    \cdot \phi_k + Q_{< j-10} u_{<j-10}\cdot  \phi_k \ . \notag
\end{equation}
For the first term we obtain an $L^2_t(L^2_x)$ bound, which by \eqref{XS_est}
suffices for both norms on the left in \eqref{u<jphik}:
\begin{align}
    \lp{ Q_{> j-10} u_{<j-10}\cdot \phi_k }{L^2_t(L^2_x)} \ &\lesssim \
    \lp{ Q_{> j-10}
    u_{<j-10}}{L^2_t( L^\infty_x) }\lp{ \phi_k }{L^\infty_t( L^2_x)}   \notag\\
    &\lesssim \ 2^{j} \lp{ Q_{> j-10}
    u_{<j-10}}{L^2_t(L^2_x)}  \lp{\phi_k }{L^\infty_t( L^2_x)}   \notag\\
    &\lesssim \ 2^{-\frac{j}2-k} \lp{\Box u}{L^2_t(\dot{H}^{-\frac12})}\lp{\nabla_{t,x}
    \phi_k}{L^\infty_t(L^2_x)} \ . \notag
\end{align}
For the second term we consider separately the two cases.
On one hand:
\begin{equation}
    Q_j (Q_{< j-10} u_{<j-10}\cdot \phi_k)  \ =  \ Q_j (Q_{< j-10} u_{<j-10}\cdot
    Q_{>j-10} \phi_k) \ . \notag
\end{equation}
therefore we directly have:
\begin{align}
    \lp{ Q_j (Q_{< j-10} u_{<j-10}\cdot  \phi_k)}{X^{1,\frac12}_\infty}
     &\lesssim
    2^{\frac{j}2+k} \lp{ Q_{< j-10}
    u_{<j-10}}{L^\infty_t(L^\infty_x)} \lp{  Q_{>j-10} \phi_k}{L^2_t(L^2_x)}
      \notag\\
    &\lesssim \ \lp{ u}{L^\infty_t(L^\infty_x)}\lp{ \phi_k}{X^{1,\frac{1}{2}}_\infty} \ . \notag
\end{align}
On the other hand, by a direct application of estimate \eqref{S_prod4}
we have:
\begin{equation}
    \lp{  Q_{<j}( Q_{< j-10} u_{<j-10}\cdot  \phi_k)}{S[k;j]}
    \ \lesssim \
    \lp{ u}{L^\infty_t(L^\infty_x)} \lp{\phi_k}{S} \ . \notag
\end{equation}
The proof of the lemma is concluded.
\end{proof}\ret

We now  return to the main proof, and consider the two bounds we
need in order to bound $U_{,k}$ in $S$, namely:
\begin{equation}
    \lp{ P_{k'}Q_j U_{,k}}{X^{1,\frac12}_\infty} +
    \lp{P_{k'}Q_{<j} U_{,k}}{S[k';j]}  \ \lesssim_F \ 2^{-\frac{1}{4}|k-k'| }c_k
    \qquad j < k'-20 < k- 10 \ . \notag
\end{equation}
For each fixed modulation index $j$, we expand $U_{,k}$ in the form:
\begin{equation}
    U_{,k}  = U_{,<j-20} B_k \! +\! \int_{j-20}^k\!\!\!\! U_{,<j-20} B_{k_1} B_k dk_1
    \!+\! \int_{j-20}^k\!\! \int_{j-20}^{k_1}\!\!\!\!  U_{,<k_2} B_{k_2} B_{k_1} B_k
    dk_2dk_1 \ . \label{uexp}
\end{equation}\ret

%--------------------------------------------------------------------------

\step{5A}{Contribution of the first term in \eqref{uexp}} We write:
\begin{equation}
        U_{,<j-20} B_k = P_{<j-10} U_{,<j-20}\cdot B_k
        +  P_{>j-10} U_{,<j-20}\cdot B_k \ . \notag
\end{equation}
The first component has output at frequency $k$, and its
contribution is accounted for due to Lemma~\ref{lvvlow}. The second
can have both high and low frequency output, so we need to split it
further.

For the high frequency output we estimate:
\begin{equation}
    \lp{ P_{>j-10} U_{,<j-20}\cdot B_k}{L^2_t(L^2_x)}
    \lesssim  \!
    \lp{ P_{>j-10} U_{,<j-20}}{L^2_t(L^\infty_x)}
    \lp{ B_k}{L^\infty_t(L^2_x)}   \lsf
    2^{-\frac{j}{2}-k} c_k \ , \notag
\end{equation}
where the $L^2_t( L^\infty_x)$ norm is estimated  by interpolating
the (summed version of the) energy bound \eqref{ulil2} with the
$L^1_t(L^1_x)$ high frequency bound \eqref{highpk} for $U_{,<j-20}$,
and by using Bernstein's inequality.

In the case of low frequency output $k'<k-20$, the first factor is
further restricted to high frequencies so we may bound:
\begin{equation}
\begin{split}
    \lp{P_{k'}( P_{>k-10} U_{,<j-20}\cdot B_k)}{L^2_t(L^2_x)} \
    & \lesssim   \lp{P_{>k-10} U_{,<j-20}}{L^2_t(L^\infty_x)}
    \lp{ B_k}{L^\infty_t(L^2_x)} \\ & \lsf
    2^{-C(k-j)}2^{-\frac{3}{2}k} c_k \ , \label{large_sep_bound}
\end{split}
\end{equation}\ret
where we have followed the same procedure as in the previous
estimate. The restriction $j<k'$ then suffices to produce
\eqref{full_U_S} for this term.
\ret

%----------------------------------------------------------------------

\step{5B}{Contribution of the second term in \eqref{uexp}} We need
to split this into several subcases:
\ret

%----------------------------------------------------------------------

\step{5B.1}{Contribution of high frequencies in $ U_{<j-20}$} This
term may have both low and high frequency output. In the case of
high frequencies we estimate directly in $L^2_t(L^2_x)$ using
Strichartz estimates as follows:
\begin{align}
    \lp{\! P_{>j-10} U_{,<j-20}\cdot  B_{k_1} B_k\!}{L^2_{t,x}}
     \! &\lesssim
    \lp{ P_{>j-10} U_{,<j-20}}{L^4_t(L^\infty_x)}
    \lp{\! B_{k_1}\!}{L^4_t(L^\infty_x)}
    \lp{\! B_k\!}{L^\infty_t( L^2_x)}  \notag\\
    &\lsf \ 2^{-\frac{j+k_1}{4}-k} c_k \ , \notag
\end{align}
where the $j-20< k_1< k$ integration is now straightforward and
yields a RHS expression of the form $\lsf
2^{-k}2^{-\frac{1}{2}j}c_k$ which suffices.
\ret

In the case of low frequency output where $k'<k-20$, we further
split the integrand as follows:
\begin{align*}
        P_{k'} (P_{>j-10} U_{,<j-20}\cdot  B_{k_1} B_k) \ = & \
        P_{k'} (P_{>j-10} U_{,<j-20}\cdot P_{>k-10} B_{k_1}\cdot B_k)\\
       & + P_{k'} (P_{>j-10} U_{,<j-20}\cdot  P_{<k-10}B_{k_1}\cdot B_k)
        \ . \notag
\end{align*}
The first term is estimated as above with a gain of
$2^{-\frac{1}{4}(k-j)}$ due to the restriction on $k_1$ (which in
particular restricts the range of integration for this term). This
suffices to show \eqref{full_U_S} for this term. To handle the
second term, we use the fact that the first factor is now forced to
be at large frequency, which gives an $L^2_t(L^2_x)$ bound as on
line \eqref{large_sep_bound} above. Notice that the additional
integration in $j-20<k_1<k$ may be absorbed via the factor of
$2^{-C(k-j)}$.
\ret

%----------------------------------------------------------------------

\step{5B.2}{Contribution of low frequencies but high modulations in
$U_{<j-20}$} In this case the only possible low frequency
contribution comes when $|k_1-k|<10$. Therefore we may proceed as
above using  the high modulation bound \eqref{boxu<k} for the first
factor as follows:
\begin{align}
    \ \ &\lp{P_{k'}( Q_{> j-10} P_{<j-10} U_{,<j-20}\cdot B_{k_1}
    B_k)}{L^2_t(L^2_x)} \  \notag\\
    \lesssim \ \ \
    &\lp{ Q_{> j-10} P_{<j-10} U_{,<j-20}}{L^4_t(L^\infty_x)}
    \lp{ B_{k_1}}{L^4_t(L^\infty_x)}
    \lp{ B_k}{L^\infty_t( L^2_x)}  \  \notag\\
    \lsf  \ &2^{-\frac{j+k_1}{4}-k} c_k \ , \notag
\end{align}
and the integral in $k_1$ is the same as above depending on whether
$|k'-k|<10$ or $k'<k-10$. In either case one gains a RHS factor of
$\lsf 2^{-\frac{1}{4}(k-k')}2^{-k}2^{-\frac{1}{2}j}c_k$.
\ret

%----------------------------------------------------------------------

\step{5B.3}{Contribution of low frequencies and low modulations in
$U_{,<j-20}$} Here we deal with the expression $Q_{< j-10}
P_{<j-10} U_{,<j-20}\cdot B_{k_1} B_k$. We consider two subcases:
\ret

%----------------------------------------------------------------------

\step{5B.3.a}{Contribution of the range $k_1 > k-10$} Under this
restriction, we may group the product $B_{k_1}B_k$ as a single term,
which we further decompose into all frequencies $k'<k+10$. For each
such localized term we have from the algebra bound \eqref{S_prod3}
the estimate:
\begin{equation}
    \lp{ P_{k'}(  B_{k_1} B_k) }{S}  \ \lesssim  \ 2^{-|k'-k|}
    \lp{B_{k_1}}{S}\,
    \lp{B_k}{S} \ \lsf  \ 2^{-|k'-k|} c_k \ . \notag
\end{equation}
Therefore, in the range $j<k'$ the resulting term may be estimated
in essentially the same way the first term on the RHS of
\eqref{uexp} was estimates in \textbf{Step 5A} above with the
additional simplification that the low frequency gains are already
implicit in the $P_{k'}(B_{k_1}B_k)$ localization.
\ret
%----------------------------------------------------------------------

\step{5B.3.b}{Contribution of the range $j-20 < k_1 < k-10$} In this
case the output is automatically at frequency $2^k$. Notice that if
we argue as in the previous case then we run into trouble with the
$k_1$ integration. Instead, we observe that one has access to the
additional localization:
\begin{equation}
        Q_j \big( Q_{< j-10} P_{<j-10} U_{,<j-20}\cdot B_{k_1} B_k\big) \
        = \
        Q_j \big( Q_{< j-10} P_{<j-10} U_{,<j-20}\cdot Q_{< j+4}(B_{k_1} B_k)\big)
        \ , \notag
\end{equation}
and according to estimate \eqref{XS_X_est} we may bound the entire
contribution of the second factor as:
\begin{equation}
    \lp{ Q_{< j+4}(B_{k_1} B_k)}{S}  \ \lesssim  \ 2^{-\delta |j-k_1|}
    \lp{ B_{k_1}}{S} \lp{B_k}{S} \ , \qquad j < k_1 \leqslant k \ . \notag
\end{equation}
This provides  the needed additional gain that enables us to
integrate with respect to $k_1$.
\ret

%--------------------------------------------------------------------------------

\step{5C}{Contribution of the last term in \eqref{uexp}} As in the
previous step, we need to split into two further subcases depending
on the range of integration:
\ret

%-------------------------------------------------------------------------

\step{5C.1}{Contribution of the range $k_1 > k-10$} A direct application of Strichartz bounds  gives the
estimate:
\begin{equation}
        \lp{  B_{k_2} B_{k_1} B_k }{L^2_t(L^2_x)}
          \lesssim  \ \lp{ B_{k_2}} {L^4_t( L^\infty_x)}\lp{
        B_{k_1}}{L^4_t( L^\infty_x)} \lp{ B_k}{L^\infty_t( L^2_x)}
          \lsf  \ 2^{-k} 2^{-\frac{k_1+k_2}{4}} c_k \ . \notag
\end{equation}
The integration with respect to $k_1,k_2$ over the region $j-20 <
k_2 < k_1< k$ with the additional restriction that $k_1=k+O(1)$is
straightforward and yields the RHS term $\lsf
2^{-\frac{1}{4}(k-j)}2^{-k}2^{-\frac{1}{2}j}c_k$ which suffices to
produce \eqref{full_U_S} for this term in light of the restriction
$j<k'$.
\ret

%-----------------------------------------------------------------------

\step{5C.2}{Contribution of the range $j-20 < k_1 < k-10$} In this
case with high frequency output we may proceed as in the previous
step. Notice that integration over the full range $j-20 < k_2 < k_1<
k$ with no additional work still yields a RHS of the form $\lsf
2^{-k}2^{-\frac{1}{2}j}c_k$.

The contribution of this range with low frequency output forces the
first term in the product to have localization in the range
$k+O(1)$. One may again proceed as in the last case of \textbf{Step
5A} above to produce an $L^2_t(L^2_x)$ estimate via
\eqref{large_sep_bound}. Notice that the integration in both $k_1$
and $k_2$ is safely absorbed by the factor $2^{-C(k-j)}$. This
concludes our demonstration of the estimate \eqref{full_U_S}.\ret

%------------------------------------------------------------------------------------

\step{6}{Proof of the bound \eqref{env_est3}}
By the algebra estimates \eqref{N_prod_est1} and \eqref{N_prod_est2} it suffices to do this for
$|k'-k|>20$. We rescale to $k=0$. There are two cases:\ret

%-----------------------------------------------------------------------------------
\step{6.A}{Low frequencies; $k'<-20$}
Here we may further localize the transformation
matrix to $P_{[-10,10]}U_{,<-20}$. Therefore, we have access to
\eqref{env_est1h}.
For the lower modulations in $G_0$ we estimate via Bernstein:
\begin{multline}
    \lp{ P_{k'} ( P_{[-10,10]}U_{,<-20}\cdot Q_{<20 }G_0) }{L^1_t(L^2_x)} \\
    \ \lesssim \ 2^{k'}
    \lp{ P_{[-10,10]}U_{,<-20} }{L^2_t(L^2_x)} \lp{Q_{<20 }G_0}{L^2_t(L^2_x)}
    \ . \notag%\label{L12_Bern_line}
\end{multline}
This suffices by estimate \eqref{standard_L2_pre} in Section \ref{s:matchfreq} above.

For the high modulation contribution, we split  $Q_{>20}G_0=G^{(1)}+G^{(2)}$,  a sum (resp) of
an $L^1_t(L^2_x)$ atom  and an $X^{0,-\frac{1}{2}}_1$ atom. For  $G^{(1)}$ the bound
\eqref{env_est3} follows by taking $P_{[-10,10]}U_{,<-20}$ in $L^\infty_t(L^2_x)$ and using Bernstein.

For the  $X^{0,-\frac{1}{2}}_1$ atom $G^{(2)}$, we may assume we are working with a single modulation
$Q_j G^{(2)}$ where $j>20$. For modulations $Q_{<j-10}P_{[-10,10]}U_{,<-20}$, estimate
\eqref{env_est3} follows by again putting the first factor in $L^\infty_t(L^2_x)$ and using Bernstein
to estimate the product as a $X^{0,-\frac{1}{2}}_1$ atom with a $2^{k'}$ gain.

For high modulations of the first factor, we estimate:
\begin{multline}
    \lp{  P_{k'} (Q_{>j-10} P_{[-10,10]}U_{,<-20}\cdot Q_{j }G_0^{(2)}) }{L^1_t(L^2_x)} \\
    \ \lesssim \ 2^{k'}
    2^{-j}
    \lp{ \partial_t Q_{>j-10}P_{[-10,10]}U_{,<-20} }{L^2_t(L^2_x)} \lp{Q_{j }G_0^{(2)}}{L^2_t(L^2_x)}
    \ , \notag
\end{multline}
which is sufficient to place the second factor in $X^{0,-\frac{1}{2}}_1$.

\ret
%-----------------------------------------------------------------------------------
\step{6.B}{High frequencies; $k'>20$}
Here we may further localize the transformation
matrix to $P_{[k'-5,k'+5]}U_{,<-20}$. Therefore, we again have access to
\eqref{env_est1h}. In this case we may proceed exactly as above, using
at each step the same estimates, which in every case suffice due to the exponential
decay in \eqref{env_est1h} for large frequencies.\ret

%--------------------------------------------------------------------------------------

\step{7}{Proof of the bound \eqref{env_est2}} Here we establish the
estimate:
\begin{equation}
    \lp{ P_{k} (\Box U_{,k_1}\cdot\psi_{k_2})}{N}  \ \lsf \
    2^{-|k-k_2|} 2^{-\delta(k_2-k_1)}c_{k_1} \lp{\psi}{S} \ ,
    \qquad k_1  \ <  \ k_2-10 \ . \label{boxu}
\end{equation}
We use the expansion \eqref{boxueq} for $\Box U_{,<k_1}$ to write
the expression we are estimating via a linear combination of the
following three terms:
\begin{align}
        G_1  =  \Box U_{,<k_1}\cdot B_{k_1}\psi_{k_2} \ , \quad
        G_2  =  \partial^\alpha U_{,<k_1} \partial_\alpha B_{k_1}\cdot
        \psi_{k_2} \ , \quad
        G_3  =  U_{,<k_1} \Box B_{k_1}\cdot \psi_{k_2} \ .
        \notag
\end{align}
We prove the bound \eqref{boxu} separately for each of these in
reverse order. Without loss of generality we will assume that
$\lp{\psi_{k_2}}{S}=1$.
\ret
%------------------------------------------------------------------------------------

\step{7A}{Estimating the term $G_3$} We directly have from \eqref{boxbkpsi} and \eqref{env_est3}
the product estimate:
\begin{equation}
        \lp{P_{k}( U_{,<k_1}\cdot \Box B_{k_1}\cdot \psi_{k_2})}{N} \ \lsf \
        2^{-|k-k_2|}2^{-\delta(k_2-k_1)}c_{k_1} \ . \notag
\end{equation}

%------------------------------------------------------------------------------------

\step{7B}{Estimating the term $G_2$} In this case the bound
\eqref{boxu} follows by applying the trilinear null-form estimate
\eqref{standard_est_tri} along with the bound \eqref{full_U_S} shown
for the first factor in the previous section. One can again split
into medium, high and low output frequency cases as in the previous
step. The details are left to the reader. Notice that the gains from
frequencies higher that $k_2$ in the first factor are essential for
maintaining the separation $2^{-\delta(k_2-k_1)}$.\ret

%------------------------------------------------------------------------------------

\step{7C}{Estimating the term $G_1$} We break this term into two
further contributions:
\begin{equation}
        G_1 \ = \ G_{11} + G_{12} \ = \
        P_{<k_2-10} \Box U_{,<k_1}\cdot B_{k_1}\psi_{k_2}
        +   P_{>k_2-10} \Box U_{,<k_1}\cdot B_{k_1}\psi_{k_2} \ .
        \notag
\end{equation}
The first term has output localized to frequency $k_2$, and we
estimate it directly via Strichartz estimates and \eqref{boxu<k}:
\begin{equation}
       \lp{G_{11}}{L^1_t(L^2_x)} \ \lesssim \
        \lp{\Box U_{,<k_1}}{L^2_t(L^2_x)}\lp{B_{k_1}}{L^4_t(L^\infty_x)}
        \lp{\psi_{k_2}}{L^4_t(L^\infty_x)} \ \lsf \
        2^{-\frac{1}{4}(k_2-k_1)}c_{k_1} \ . \notag
\end{equation}
The second term $G_{12}$ can have both high and low frequency
outputs. When the output is in the range $k<k_2+10$ we use
\eqref{highpk} and Bernstein's inequality  to bound it as follows:
\begin{align}
        \lp{P_k G_{12}}{L^1_t(L^2_x)} \ &\lesssim \
        2^{k}\lp{P_{>k_2-10} \Box U_{,<k_1}}{L^1_t(L^1_x)}
        \lp{B_{k_1}}{L^\infty_t(L^\infty_x)}\lp{\psi_{k_2}}{L^\infty_t(L^\infty_x)}
        \  \notag\\
        &\lesssim_F \ 2^{k}2^{-k_1}2^{-C(k_2-k_1)}c_{k_1} \ , \notag
\end{align}
which suffices to show \eqref{boxu} in this case. When $G_{12}$ has
output in the high range $k>k_2+10$ we have further high frequency
localization of the first factor and we may estimate via the same
procedure:
\begin{align}
        \lp{P_k G_{12}}{L^1_t(L^2_x)} \ &\lesssim \
        2^{k}\lp{P_{[k-5,k+5]} \Box U_{,<k_1}}{L^1_t(L^1_x)}
        \lp{B_{k_1}}{L^\infty_t(L^\infty_x)}\lp{\psi_{k_2}}{L^\infty_t(L^\infty_x)}
        \ , \notag\\
        &\lesssim_F \ 2^{k}2^{-k_1}2^{-C(k-k_1)}c_{k_1} \ , \notag
\end{align}
which is again sufficient to show \eqref{boxu} in this case. This
concludes our demonstration of Proposition \ref{propphiu}.

\ret
%-------------------------------------------------------------------------
%%%%%%%%%%%%%%%%%%%%%%%%%%%%%%%%%%%%%%%%%%%%%%%%%%%%%%%%%%%%%%%%%%%%%%%%%%
%-------------------------------------------------------------------------

\section{The Linear Paradifferential Flow}\label{s:lin_sect}

We now proceed with the proof of Proposition \ref{p:para}. The main
difficulty here is that we do not necessarily have smallness of the
constant from line \eqref{A_tdphi_est}, which would otherwise make
estimate \eqref{linearized_est} consequence of Propositions
\ref{standard_prop1} and \ref{propphiu}. Instead of proceeding
directly, we shall follow a more measured approach of building up
our estimate piece by piece. Since this is a lengthy argument, we
begin with a brief outline.

The first step of the proof is to take advantage of the
antisymmetry of $A_\alpha$, which makes our paradifferential equation
almost conservative. Precisely, the only nontrivial contributions
to  energy estimates arise from terms where one derivative falls on the
coefficients. But such terms are small due to the
large frequency gap $m$. Consequently, we are able to prove a favorable
estimate:
\begin{equation}
        \lp{\psi_k}{\underline{E}[I]} \ \lesssim_F \
        \lp{\psi_k[0]}{\dot H^1 \times L^2}  +
        2^{\delta m}\lp{G_k}{N[I]}
        + 2^{-\delta m}\lp{\psi_k}{S[I]} \ , \label{para_ap_est}
\end{equation}
for the energy \eqref{uE_norm} on both time slices and characteristic surfaces.

We still need an estimate on the $S$ norm of $\psi_k$, for which we renormalize
the equation \eqref{reduced_lin_eq} using an orthogonal
gauge transformation $U_{,<k-m}$ obtained by Proposition~\ref{propphiu}.
The function $w_{,k} = U_{,<k-m} \psi_k$ solves a perturbed wave
equation of the form:
\begin{equation}
    \Box w_{,k}  \ =  \ \mathcal{R}^{pert}_{,<k-m} \psi_k + U_{,<k-m} G_k \ .
    \label{renormw}\end{equation}
In the analysis of the small data problem in \cite{Tao_WM}, \cite{Tataru_WM1}
one uses a perturbative bound of the form:
\begin{equation}
    \lp{ \mathcal{R}^{pert}_{,<k-m} \psi_k}{N[I]}  \
    \lesssim_F  \ \lp{\psi_k}{S[I]} \ , \notag
\end{equation}
where the implicit constant is at least quadratic in $F$, for $F$ small.
This is no longer sufficient here. Instead, we observe that we can
rebalance the above estimate and use only the energy norm of $\psi_k$
to estimate the bulk of the LHS above. Thus, we prove that
for $ 0 \leqslant m_0 \ll m$ we have:
\begin{equation}
    \lp{ P_j (\mathcal{R}^{pert}_{,<k-m} \psi_k)}{N[I]}  \ \lesssim_F  \ 2^{-|j-k|}
    \big( 2^{-\delta_1 m_0}\lp{\psi_k}{S[I]} + 2^{2m_0}\lp{\psi_k}{\underline{E}[I]}\big) \ .
    \label{qpert}
\end{equation}
By the linear solvability bound \eqref{energy_est} we have:
\begin{equation}
    \lp{ w_{,k}}{S[I]}  \ \lesssim  \
    \lp{w_{,k}[0]}{\dot H^1 \times L^2} +\lp{\Box w_{,k}}{N[I]} \ . \notag
\end{equation}
Since both $U_{,<m-k}$ and $(U_{,<m-k})^{-1}=U_{,<m-k}^{\dagger}$
are in $S$ with norm $\lesssim_F 1$, by the $S$ algebra property and estimate
\eqref{env_est3} we have the gauge removal bounds:
\begin{equation}
    \lp{ \psi_k}{S[I]}  \ \lesssim_F  \ \lp{w_{,k}}{S[I]} \ ,
    \quad \lp{P_{k'}( U_{,<k-m} G_k)}{N[I]}  \ \lesssim_F  \ 2^{-|k-k'|}\lp{G_k}{N[I]}
    \ . \label{basic_para_ests}
\end{equation}
On the other hand using the energy component of
\eqref{env_est1} we obtain:
\begin{equation}
    \lp{ P_j w_{,k}[0]}{\dot H^1 \times L^2}  \ \lesssim_F\
    2^{-|j-k|}\lp{ \psi_{k}[0]}{\dot H^1 \times L^2} \ .
    \label{dataforw}
\end{equation}
Summing up the estimates on the last four lines we obtain the
$S$ bound for $\psi_k$:
\begin{equation}
    \lp{ \psi_k}{S[I]} \lesssim_F \lp{\psi_k[0]}{\dot H^1 \times L^2}
    + \lp{ G_k}{N[I]} +  2^{-\delta_1 m_0}\lp{\psi_k}{S[I]}
    + 2^{2m_0}\lp{\psi_k}{\underline{E}[I]}
    \ .
    \label{psik-s}
\end{equation}
Now all we have to do is combine this with \eqref{para_ap_est},
carefully balancing the constants.
Assuming that $m_0 = m_0(F)$ for a large enough  $m_0(F)\sim\ln(F)$,
the third term on the right can be absorbed on the left to obtain:
\begin{equation}
    \lp{\psi_k}{S[I]}  \ \lesssim_F  \ \lp{\psi_k[0]}{\dot H^1 \times L^2}
    + \lp{ G_k}{N[I]}  + \lp{\psi_k}{\underline{E}[I]} \ . \notag
    %\label{psik-sf}
\end{equation}
Substituting \eqref{para_ap_est} for the third term on the RHS of this last line
we arrive at:
\begin{equation}
    \lp{\psi_k}{S[I]} \ \lesssim_F \
    \lp{ \psi_{k}[0]}{\dot H^1 \times L^2}
    + 2^{\delta m}\lp{G_k}{S[I]} +  2^{-\delta m}\lp{\psi_k}{S[I]} \ , \notag
\end{equation}
so now assuming $m > m(F)$ for a larger $m(F)\sim\ln(F)$,
the last term on the RHS is again absorbed on the left:
\begin{equation}
    \lp{\psi_k}{S[I]} \ \lesssim_F \
    \lp{ \psi_{k}[0]}{\dot H^1 \times L^2}  + \lp{G_k}{S[I]} \ . \notag
    %\label{psik-ef}
\end{equation}
To conclude the proof of \eqref{linearized_est} we need to improve
the $S$ bound above to a $\W$ bound.
Returning to $w_{,k}$, we have the estimate:
\begin{equation}
    \lp{P_{j} \Box w_{,k}}{N[I]} +
    \lp{ P_j w_{,k}[0]}{\dot H^1 \times L^2}
     \ \lesssim_F  \ 2^{-|j-k|} \big(\lp{\psi_{k}[0]}{\dot H^1 \times L^2}
     + \lp{G_k}{S[I]}\big) \ . \notag
\end{equation}
This follows from \eqref{qpert}, \eqref{dataforw}, and
the second member on line \eqref{basic_para_ests}.\ret

It remains to prove the two main estimates above, namely
\eqref{para_ap_est} and \eqref{qpert}. In  the proof we shall make use
of three auxiliary Lemmas whose proofs we postpone until the end of
this section. The first one is used to estimate perturbative expressions
which are small due to the large frequency gap $m$.

\begin{lem}[Some auxiliary estimates]\label{para_aux_lem}
Let $A_\alpha$ be the connection one-form defined on line
\eqref{red_con} above with estimates \eqref{A_tdphi_est}. Then the
following bounds hold:
\begin{align}
        \lp{A_\alpha}{L^\infty_t(L^\infty_x)} \ &\lesssim_F \ 2^{k-m} \ ,
        \label{A_bound}\\
        \lp{A_\alpha \psi_k}{DS[I]} \ &\lesssim_F \
        2^{-m}\, \lp{\psi_k}{S[I]} \ . \label{A0_bound}
\end{align}
Also, for three test functions $\phi^{(i)}$ normalized with
$S\cap\underline{E}[I]$ size one, the following list of multilinear
estimates holds:
\begin{align}
        \lp{\phi^{(1)} \partial^\alpha \phi^{(2)}_{<k-m}
        \partial_\alpha \phi^{(3)}_{<k-m}\cdot\psi_k}{N[I]}
        \ &\lesssim \ 2^{-\delta m}\,
        \lp{\psi_k}{S[I]} \ , \label{I_est_1}\\
    \lp{P_{<k-m}\left(\phi^{(1)}\partial^\alpha\phi^{(2)}\partial_\alpha
        \phi^{(3)}\right)\psi_k }{N[I]} \ &\lesssim \ 2^{-\delta
        m}\,
        \lp{\psi_k}{S[I]} \ , \label{I_est_2}\\
        \lp{ \nabla_{t,x}\phi^{(1)}_{<k-m}
        \partial^\alpha \phi^{(2)}_{<k-m} \partial_\alpha\psi_k}{N[I]}
        &\lesssim \ 2^{-\delta m}2^k\,
        \lp{\psi_k}{S[I]} \ , \label{I_est_3}\\
        \lp{ \nabla_{t,x}\phi^{(1)}_{<k-m}
        \partial^\alpha \phi^{(2)}_{<k-m} \partial_\alpha \phi^{(3)}_{<k-m}
        \cdot \psi_k}{N[I]}
        &\lesssim \ 2^{-\delta m}2^k\,
        \lp{\psi_k}{S[I]} \ . \label{I_est_4}
%         \lp{ \phi^{(1)}_{<k-m}
%        \partial_t \partial^\alpha \phi^{(2)}_{<k-m} \partial_\alpha\psi_k}{N[I]}
%        &\lesssim \ 2^{-\delta m}2^k\,
%        \lp{\psi_k}{S[I]} \ , \label{I_est_4}
\end{align}
\end{lem}\ret

In proving energy estimates we need to restrict integration
to half-spaces. This is where the next lemma comes handy:

\begin{lem}[Half-space duality estimate]\label{half_dual_lem}
Let $\psi_k\in S$ and $H_k\in N$ be  frequency localized
functions. Then for any time-slab $I$, any unit vector  $\omega$,
and any spatial point $x_0\in\mathbb{R}^2$ the following truncated
duality estimate holds uniformly:
\begin{equation}
        \Big| \dint_{ I\cap\{t>\omega\cdot(x-x_0)\}}\
        H_k\cdot  \psi_k \ dxdt
        \Big|
        \lesssim \ \lp{H_k}{N[I]}\cdot
       \lp{\psi_k}{DS[I]}
         \ . \label{half_dual_est}
\end{equation}
\end{lem}\ret

Finally, for the bulk of the estimate \eqref{qpert} we need the
following lemma, which improves upon the trilinear bound \eqref{standard_est_tri} in the case of balanced low frequencies
$|k_1 -k_2|\ll m$, $k_1,k_2 < k_3-m$:

\begin{lem}[An improved trilinear estimate]\label{para_trilin_lem}
There exists a universal constant $C> 0$ such that for any integer
$m\geqslant 0$ and  $S[I]$ unit normalized test functions
$\phi^{(i)}_{k_i}$ with $k_i\leqslant k-m$, and $\psi_k$  any
additional test function defined on $I$,  one has the following
imbalanced trilinear estimate:
\begin{equation}
    \lp{\phi^{(1)}_{k_1}\partial_\alpha \phi^{(2)}_{k_2} \partial^\alpha \psi_k }{N[I]} \ \lesssim \
    2^{C|k_1-k_2|}
    \left(2^{-\delta^2 m}\lp{\psi_k}{S[I]} + 2^{m}\lp{\psi_k}{\underline{E}[I]}\right)
    \  , \label{imbal_trilin}
\end{equation}
%Here $\delta$ is the same as in Propositions \ref{standard_prop1}
%and \ref{standard_prop2}.
\end{lem}\ret

\noindent Assuming these estimates, we give a proof of
\eqref{para_ap_est} and \eqref{qpert} in a series of steps. To close the argument properly, we will  employ our chain of small constants
\eqref{delta_list} (although there use here is independent of their
use in other sections).\ret

%---------------------------------------------------------------------

\step{1}{A-priori control of the  energy norm of $\psi_k$: proof of \eqref{para_ap_est}}
We begin by writing the equation \eqref{reduced_lin_eq} for $\psi_k$ on the
interval $I$ in a covariant form:
\begin{equation}
        \Box_A \psi_k \ =  \ \mathcal{R}^{en}_{<k-m} \psi_k + G_k
        \label{ext_psi_eq}
\end{equation}
where $\Box_A = (\partial + A)^\alpha(\partial + A)_\alpha$ is the
gauge covariant wave equation with the connection $A_\alpha$ is
given by the formula on line \eqref{red_con} and the function $\mathcal{R}^{en}_{<k-m}$
has the form:
\begin{equation}
    \mathcal{R}^{en}_{<k-m}   \! = \!
    \partial^\alpha
        \mathcal{S}(\phi)_{<k-m}\partial_\alpha\phi_{<k-m}
       \! + \mathcal{S}(\phi)_{<k-m}
    P_{<k-m}\left(\mathcal{S}(\phi)\partial^\alpha\phi\partial_\alpha
        \phi\right) + A^\alpha A_\alpha \ .
\end{equation}
Note that in the RHS of
this last line, the matrix $\mathcal{S}(\phi)$ is either the pure
second fundamental form $\mathcal{S}(\phi)_{ab}^c$, or its
antisymmetric version as it appears in the formula for $A_\alpha$.
The distinction will not be important for us here. Also, notice that
we have used the Wave-Map equation for $\phi$ on the interval $I$,
which we may do by the assumptions of Proposition \ref{p:para}.\ret

To obtain the energy estimates we proceed via a simple
integration-by-parts argument. First, we form the gauge-covariant
energy momentum density:
\begin{equation}
        Q_{\alpha\beta}[\psi_k] \ = \ (\Anabla_\alpha \psi_k)^\dagger
        \Anabla_\beta \psi_k - \frac{1}{2}g_{\alpha\beta}
        (\Anabla^\gamma \psi_k)^\dagger \Anabla_\gamma \psi_k \ . \notag
\end{equation}
Here we are writing $\Anabla_\alpha = \partial_\alpha + {A}_\alpha$.
A quick calculation shows that (notice that this identity crucially
uses the antisymmetry of $A_\alpha$, which is the main source of the
cancelation that makes \eqref{para_ap_est} possible):
\begin{equation}
    \nabla^\alpha Q_{\alpha\beta}[\psi_k] \ = \
    (\Box_{A}\psi_k)^\dagger \Anabla_\beta \psi_k +
    (F_{\gamma\beta} \psi_k)^\dagger \Anabla^\gamma \psi_k \ ,
    \label{div_iden}
\end{equation}
where $F_{\alpha\beta} \ = \ \partial_\alpha A_\beta -
\partial_\beta A_\alpha + [A_\alpha,A_\beta]$ is the curvature of
$A_\alpha$. Next, we form the linear momentum one-form
$P_\alpha[\psi_k]=Q_{\alpha 0}[\psi_k]$.  Integrating $\nabla^\alpha
P_\alpha[\psi_k]=\nabla^\alpha Q_{\alpha 0}[\psi_k]$ over all possible half spaces
of the form $[0,t_0]\cap\{t>\omega\cdot(x-x_0)\}$ we have the bound:
\begin{multline}
    \!\!    \lp{\! \Anabla_{t,x}\psi_k}{L^\infty_t(L^2_x)[I]}^2 + \sup_\omega
        \lp{\ \, \! \slash\!\!\!\!\!\!\Anabla_{t,x}\psi_k}
        {L^\infty_{t_\omega}(L^2_{x_\omega})[I]}^2
        \lesssim
        \lp{\! \Anabla_{x,t} \psi_k(0)}{L^2_x}^2 + I_1 + I_2
         \label{raw_energy_est}
\end{multline}
where:
\begin{align}
        I_1 \ &= \  \sup_{I,\omega,x_0}
        \Big| \dint_{ I\cap\{t>\omega\cdot(x-x_0)\}}(\Box_{A}\psi_k)^\dagger  \Anabla_0 \psi_k \ dxdt
        \Big| \ , \notag\\
        I_2 \ &= \  \sup_{I,\omega,x_0}
        \Big| \dint_{ I\cap\{t>\omega\cdot(x-x_0)\}} (F_{0\gamma} \psi_k)^\dagger \Anabla^\gamma \psi_k
        \ dxdt
        \Big| \ . \notag
\end{align}
Our task is to estimate $I_1$ and $I_2$ and to show that we can
replace covariant differentiation by regular differentiation in
\eqref{raw_energy_est}. For the right hand side of \eqref{raw_energy_est}
we claim that both:
\begin{align}
        \lp{\Anabla_{x,t} \psi_k(0)}{L^2_x}^2 \ &\lesssim_F \
        \lp{\nabla_{x,t}\psi_k(0)
        }{L^2_x}^2\
        , \label{data_bound}\\
        I_1 + I_2 \ &\lesssim_F \ 2^{-2\delta m}\lp{\psi_k}{S[I]}^2 +
        2^{2\delta m}\lp{G_k}{S[I]}^2 \ . \label{I_12_bound}
\end{align}

The proof of \eqref{data_bound} is an immediate consequence of expanding the
covariant derivative
$\Anabla$ and using the triangle inequality, followed by the $L^\infty_t(L^\infty_x)$
bound for  $A_\alpha$ in \eqref{A_bound}.\ret

To  obtain \eqref{I_12_bound}, we use the half-space duality
estimate \eqref{half_dual_est} and Young's inequality for the term
involving $G_k$. For the other terms, we again use half-space
duality, and then conclude with an application of the estimates
\eqref{s-ds}, \eqref{A0_bound}--\eqref{I_est_4}. It suffices to
establish the bounds:
\begin{align}
    \lp{\mathcal{R}^{en}_{<k-m} \psi}{N[I]} &\lesssim_F
    2^{-\delta m} \lp{\psi_k}{S[I]} \ ,
    &\lp{F_{0\gamma}\cdot \Anabla^\gamma \psi_k }{N[I]} &\lesssim_F
    2^{-\delta m} 2^k \lp{\psi_k}{S[I]}  . \notag
\end{align}
The first estimate above follows from applying
\eqref{I_est_1}--\eqref{I_est_2} to each of the terms in
$\mathcal{R}^{en}_{<k-m}$. The second estimate follows from the
bounds \eqref{I_est_3}--\eqref{I_est_4} applied to the definition of
the curvature. Notice that these two multilinear estimates suffice
because there are never any terms in $I_2$ with a single factor
containing more than one derivative thanks to the skew symmetry of
the curvature. \ret

The bound \eqref{para_ap_est} will now follow once we can rid
ourselves of the gauge covariant derivatives $\Anabla_{t,x}$ on the
LHS of \eqref{raw_energy_est} in favor of the usual derivatives
$\nabla_{t,x}$. This can be done with a successive application of the two
estimates:
\begin{align}
        \lp{A\psi_k}{L^\infty_t(L^2_x)[I]} \ &\lesssim_F \ 2^{-m}
        \lp{\nabla_{t,x}\psi_k}{L^\infty_t(L^2_x)[I]}
        \ , \notag\\
        \sup_\omega
        \lp{\slash\!\!\!\!A\psi_k}{L^\infty_{t_\omega}(L^2_{x_\omega})[I]}
        \ &\lesssim_F \ \|\psi_k\|_{L^\infty_t(L^\infty_x)[I]}
    \ \lesssim  \ \lp{\nabla_{t,x}\psi_k}{L^\infty_t(L^2_x)[I]}
        \ . \notag
\end{align}
The first of these follows immediately from the bound
\eqref{A_bound}, while the second uses the
characteristic energy estimates we are assuming for $\phi$.
We remark that using the first bound above requires $m$ to be
large enough, i.e. $2^m \gg_F 1$.

% Adding
% these two together and using \eqref{raw_energy_est} in conjunction
% with \eqref{data_bound}--\eqref{I_12_bound} we have the bound:
% \begin{multline}
%          \lp{\nabla_{t,x}\psi_k}{L^\infty_t(L^2_x)[I]} + M^{-1}\sup_\omega
%         \lp{\ \, \snabla_{t,x}\psi_k}
%         {L^\infty_{t_\omega}(L^2_{x_\omega})[I]} \\
%         \lesssim_F \ 2^{-m}\lp{\nabla_{t,x}\psi_k}{L^\infty_t(L^2_x)[I]}
%         +
%         M^{-1}
%         \lp{\nabla_{t,x}\psi_k}{L^\infty_t(L^2_x)[I]}\\
%         + \lp{(\nabla_x
%         f_k,\dot{f}_k)}{L^2_x} + 2^{-\delta_2 m}\lp{\psi_k}{S[I]} +
%         \lp{G_k}{S[I]}\ . \notag
% \end{multline}
% When $m$ is sufficiently large, we may take $M\gg
% Q\big(\lp{\phi}{S\cap \underline{E}[I]}\big)$, where $Q$ is the
% implicit polynomial on the RHS of the above inequality, to yield
% \eqref{para_ap_est}.\ret

\ret
%-----------------------------------------------------------------------------

\step{2}{The $S$ bound for $\psi$: Proof of \eqref{qpert}}
 The first thing we need to do is to rewrite the equation
\eqref{reduced_lin_eq} in a gauged formulation (we have no further use for \eqref{ext_psi_eq}). As usual, we write:
\begin{equation}
    A^\alpha_{<k-m} \ = \ \partial^\alpha B_{<k-m} + D^\alpha_{<k-m} \ , \notag
\end{equation}
where the RHS is given by the integrated terms:
\begin{align}
    B_{<k-m} \ &= \ \int_{k'<k-m}\big(S^a_{cb}(\phi) - S^b_{ca}(\phi)\big)_{\leqslant k'-10}
    \phi^c_{k'}dk' \ , \notag\\
      D^\alpha_{<k-m} \ &= \
       \int_{k'<k-m}\big(S^a_{cb}(\phi) - S^b_{ca}(\phi)\big)_{k'-10< \cdot < k-m}
    \partial^\alpha\phi^c_{k'}dk' \notag\\
      &\ \ \  - \ \int_{k'<k-m}
      \partial^\alpha \big(S^a_{cb}(\phi) - S^b_{ca}(\phi)\big)_{\leqslant k'-10}
    \phi^c_{k'}dk' \ . \notag
\end{align}
The connection $\partial^\alpha B_{<k-m}$ is of the form in
Proposition \ref{propphiu}, and we define the $SO(N)$ matrix
$U=U_{,<k-m}$ accordingly. We also set $C_\alpha = \nabla_\alpha B -
U^\dagger \nabla_\alpha U  $, which is given by the second term on
the RHS of formula \eqref{renorm_form}. Finally, we  denote by
$w_{,k} = U_{,<k-m}\psi_k$. Then $w_k$ obeys the gauged
equation \eqref{renormw} with:
\begin{equation}
    \mathcal{R}^{pert}_{,<k-m} \psi_k   \ = \ -2U_{,<k-m}(C^\alpha + D^\alpha)
    \partial_\alpha \psi_k   + \Box U_{,<k-m}\cdot \psi_k
    \ . \label{w_k_eq}
\end{equation}
The second term on the right is easy to estimate using \eqref{env_est2},
which yields:
\begin{equation}
    \lp{P_{k'}( \Box U_{,<k-m}\cdot\psi_k) }{N[I]} \ \lesssim_F
        \ 2^{-|k-k'|}2^{-\delta_1 m}\lp{\psi_k}{S[I]} \ . \label{box_U_psi_est}
\end{equation}
It remains to estimate the first term in $ \mathcal{R}^{pert}_{,<k-m} \psi_k $,
for which we will show the bound:
\begin{equation}
        \lp{\! P_{k'}\big(
        U_{,<k-m}(C^\alpha + D^\alpha) \partial_\alpha \psi_k\big)\!}{N[I]} \! \lesssim_F \!
       2^{-|k-k'|}\big(2^{-\delta_1 m_0}\lp{\!\psi_k\!}{S[I]} +
       2^{2m_0}\lp{\!\psi_k\!}{\underline{E}[I]}\big)  \ , \label{CD_phi_est}
\end{equation}
for $0 < m_0\ll m$. We will prove this in a further series of steps. \ret

%------------------------------------------------------------------------------

\step{2A}{Removal of the gauge and high frequency connection} Here we write $ C^\alpha_{low} $
for the second term on the RHS of line \eqref{renorm_form} with each gauge factor replaced
by $P_{<k-10}U_{,<k'}$. Thus  $C^\alpha_{low}=P_{<k-10}C^\alpha_{low}$.
Notice that the connection $D^\alpha$ also has frequency $<k-m+10$. Therefore,
from estimate \eqref{env_est3} we have:
\begin{equation}
        \lp{P_{k'}\big( U_{,<k-m}(C^\alpha_{low} + D^\alpha)\partial_\alpha \psi_k\big)}{N[I]}
        \ \lesssim_F \ 2^{-|k-k'|}  \lp{(C^\alpha_{low} + D^\alpha)\partial_\alpha \psi_k}{N[I]}
    \ . \notag
\end{equation}
Furthermore, we claim the remainder estimate:
\begin{equation}
    \lp{P_{k'}\big( U_{,<k-m}(C^\alpha - C^\alpha_{low} )\partial_\alpha \psi_k\big)}{N[I]}
        \ \lesssim_F \ 2^{-|k-k'|}  2^{-Cm}\lp{\psi_k}{S[I]}
    \ . \notag
\end{equation}
Setting $R=U_{,<k-m}(C^\alpha - C^\alpha_{low} )$, this follows at once
from Bernstein's inequality and the improved bounds:
\begin{equation}
    \lp{R}{L^1_t(L^1_x)[I]} \ \lesssim_F \ 2^{-2k}2^{-Cm}
    \ , \qquad
    \lp{P_{j}R}{L^1_t(L^1_x)[I]} \ \lesssim_F \ 2^{-2k}2^{-C(j-k)}2^{-Cm}
    \ . \notag
\end{equation}
These estimates are a consequence of the improved estimate \eqref{env_est1h}, and the fact that at
least one of the gauge factors in the $C^\alpha - C^\alpha_{low}$ integral is localized to
$P_{>k-10}U_{,<k'}$.

\ret

%-----------------------------------------------------------------------------------

\step{2B}{Estimation of the main term} The purpose of this step is
to prove the remaining estimate:
\begin{equation}
        \lp{ (C^\alpha_{low} + D^\alpha)\partial_\alpha \psi_k}{N[I]}
        \ \lesssim_F \
       2^{-\delta_1 m_0}\lp{\psi_k}{S[I]} +
       2^{2m_0}\lp{\psi_k}{\underline{E}[I]} \ . \label{red_CD_phi_est}
\end{equation}
 We'll do this separately
for each of the two terms on the left.\ret

%-----------------------------------------------------------------------------------

\step{2B.1}{Estimation of $D^\alpha$ term} The plan is to use Lemma
\ref{para_trilin_lem}. To do this we need to separate the connection
$D^\alpha$ into two pieces, one with essentially matched frequencies
and one with wide frequency separation. We write
$D^\alpha=D^\alpha_{(\delta)} + \td{D}^\alpha$ where:
\begin{align}
        D^\alpha_{(\delta)} \ &= \  \int_{k'<k-m}\big(S^a_{cb}(\phi)
        - S^b_{ca}(\phi)\big)_{[k'-10, k'+c\delta^2 m_0]}
        \partial^\alpha\phi^c_{k'}dk' \notag\\
        &\ \ - \ \int_{k'<k-m}
        \partial^\alpha \big(S^a_{cb}(\phi) - S^b_{ca}(\phi)\big)_{[k'-c\delta^2 m_0,k'-10]}
        \phi^c_{k'}dk' \ . \notag
\end{align}
Here $c\ll 1$ is an additional small constant. By a direct
application of estimate \eqref{imbal_trilin} we have:
\begin{equation}
        \lp{D_{(\delta)}^\alpha \partial_\alpha \psi_k}{N[I]} \ \lsf \
        2^{Cc\delta^2 m_0}
        \left(2^{-\delta^2 m_0}\lp{\psi_k}{S[I]} +
        2^{m_0}\lp{\psi_k}{\underline{E}[I]}\right) \ . \label{delta_D_est}
\end{equation}
For $c$ small enough in relation to $C$ we have
\eqref{red_CD_phi_est} for this term. The remainder term is  in
the range where the standard trilinear estimate
\eqref{standard_est_tri} gives additional savings. A quick
computation shows that for this term we in fact have:
\begin{equation}
        \lp{\td{D}^\alpha \partial_\alpha \psi_k}{N[I]} \ \lsf \
        2^{-c\delta^3 m_0} \lp{\psi_k}{S[I]} \ . \label{tdD_est}
\end{equation}
The details of the dyadic summation are left to the reader.\ret

%--------------------------------------------------------------------------------

\step{2B.2}{Estimation of $C^\alpha_{low}$ term} We follow the same
strategy as in the previous argument. We split
$C^\alpha_{low}=C^\alpha_{(\delta)} + \td{C}^\alpha$ where:
\begin{equation}
        C^\alpha_{(\delta)} \ = \
        \int_{-\infty}^{k-m} \big[B_{k'},P_{<k-10} U^\dagger_{,< k'}
        \nabla^\alpha P_{[k'-c\delta^2m_0,k'+c\delta^2m_0]}
        P_{<k-10} U_{,<k'}\big] \, dk' \ . \notag
\end{equation}
The factors $P_{<k-10} U_{,<k'}$ are bounded on $N$ via estimate
\eqref{N_prod_est1}, and can therefore be neglected. Again, by
summing over the bound \eqref{imbal_trilin} with the help of
\eqref{env_est1} we have the analog of \eqref{delta_D_est} (but this
time with a factor of $2^{2Cc\delta^2 m_0}$ instead) for the
contraction $ C^\alpha_{(\delta)}\partial_\alpha \psi_k$. Similarly,
we have the analog of \eqref{tdD_est} for the contraction
$\td{C}^\alpha \partial_\alpha \psi_k$, which also uses estimate
\eqref{env_est1h}.\ret

%--------------------------------------------------------------------------------

%--------------------------------------------------------------------------------

\begin{rem}\label{direct_para_rem_proof}
The above process can also be used to show that if one already has
$\lp{\psi_k}{S[I]}$ norm control, then one may conclude
normalization bounds $\lp{\psi_k}{\W[I]}$ under the much less
restrictive assumption that $m\geqslant 20$. In this case, one
simply skips all of \textbf{Step 1} above, and carry out
 \textbf{Step 2} without introducing at all the terms $C^\alpha_{(\delta)}$
and $D^\alpha_{(\delta)}$.
\end{rem}\ret

%-----------------------------------------------------------------------------

%To end this section, we give the proofs of Lemmas \ref{para_aux_lem}--\ref{para_trilin_lem}.

\begin{proof}[Proof of Lemma \ref{para_aux_lem}]
The estimate \eqref{A_bound} follows from the energy bounds for
$\phi$ combined with Bernstein's inequality. On the other hand
\eqref{A0_bound} is a consequence of \eqref{S_DS_prod} and
\eqref{s-ds}. \ret

Estimate \eqref{I_est_1} follows from an application of
\eqref{N_prod_est1}--\eqref{N_prod_est2}, and then summation over
the trilinear bound \eqref{standard_est_tri}. The relevant detail is
that one has the dyadic sum:
\begin{equation}
        \sum_{k_2,k_3: \ k_i< k-m} 2^{-\delta(k-\min\{k_2,k_3\})} \
        \lesssim \ 2^{-\delta m} \ . \notag
\end{equation}
Estimate \eqref{I_est_2} is a more elaborate use of such summations,
but it is standard and left to the reader.\ret

Consider now \eqref{I_est_3}. For  modulations
at most comparable to the frequencies in the first factor we
can replace the time derivative with a frequency factor
and prove the estimate
\eqref{I_est_3}  by summing over
\eqref{standard_est_tri}. The relevant detail is that one has the
dyadic sum:
\begin{equation}
        \sum_{k_1,k_2: \ k_i< k-m} 2^{k_2}2^{-\delta(k_2-k_1)_+} \
        \lesssim \ 2^{- m}2^k \ . \notag
\end{equation}
It remains to bound the expression when the first factor is at high
modulation. In this case we take a product of the two bounds:
\begin{align}
        \lp{ \nabla_{t,x}Q_{|\xi|\ll |\tau|}\phi^{(1)}_{<k-m}}{L^2(L^\infty)}
        \ &\lesssim \ 2^{\frac{1}{2}(k-m)}    \ , \notag\\
        \lp{\partial^\alpha \phi^{(2)}_{<k-m}\partial_\alpha \psi_k}{L^2_x(L^2_x)}
         \ &\lesssim \  2^{\frac{1}{2}(k-m)}\lp{\psi_k}{S[I]}
         \ , \notag
\end{align}
the first of which follows from summation over \eqref{L2Linfty_est}
and the second of which follows from summation over
\eqref{standard_L2_est_bi}. The estimate \eqref{I_est_4} follows
from similar reasoning and is left to the reader.
\end{proof}\ret

\begin{proof}[Proof of Lemma \ref{half_dual_lem}]
The bound we seek is scale invariant, so without loss of generality
we may assume that $k=0$, and we may rotate and center the estimate
so that $\omega=(1,0)$ and $x_0=(0,0)$. In light of
\eqref{duality_est} we see that the main point of
\eqref{half_dual_est} is to be able to drop half space cutoffs of
the form $\chi_{t<0}$ and $\chi_{t<x^1}$. The required boundedness
of   cutoffs with discontinuities across space-like hypersurfaces
was already shown in \eqref{dscut}. Therefore, we seek an analog of
\eqref{dscut} in the null case. Due to the frequency localization of
both factors on the LHS of \eqref{half_dual_est}, it suffices to
prove the  following product estimate:
\begin{align}
    \lp{P_0 \big(\chi_{t<x^1}\cdot
     \psi_0) }{DS} \ &\lesssim \ \lp{\psi_0}{DS} \ , \label{char_dual}\
\end{align}

To save notation we will write $\chi=\chi_{t<x^1}$. Our point of
view will be to observe that  $\chi$ is a singular solution to the
wave equation, so one can hope that \eqref{char_dual} is in some
sense a version of the standard product estimate \eqref{S_prod2}.
While this is true, the demonstration requires a bit of care because
the $P\!W$ norm of $P_{k}(\chi)$ does not gain the usual  weight
from $L^1$ summation over angles, even though its Fourier support is
well localized in the angular variable. In fact, a quick calculation
shows that:
\begin{equation}
    \widetilde{\chi}(\tau,\xi) \ = \ \begin{cases}
                    \frac{c_+}{\xi_1}\delta(\tau+|\xi|)\delta(\xi_2) \ ,
                    &\xi_1 \ > \ 0 \ ;\\
                    \frac{c_-}{\xi_1}\delta(\tau-|\xi|)\delta(\xi_2) \ ,
                    &\xi_1 \ < \ 0 \ .
                           \end{cases}\notag
\end{equation}
Here $c_\pm$ are appropriate constants depending on ones in the
definition of the Fourier transform.  The above formulas show that
the $(+)$ wave portion of $\chi$ is a measure concentrated on the
ray $(1,-1,0)$, and opposite for the $(-)$ wave portion. We have the
frequency localized $P\!W$ type bound:
\begin{equation}
    \lp{Q^\pm P_k \chi }{L^2_{t_{(1,0)}}(L^\infty_{x_{(1,0)}})}
    \ \lesssim \ 2^{-\frac{1}{2}k} \ . \label{PW_chi}
\end{equation}
Finally, note that due to the frequency localization in
\eqref{char_dual}, we may replace the cutoff
with $Q_{<10}P_{<10}(\chi)$. Also, if $\phi_0$ is at high modulation
$ \> 10$ then $P_0 \big(\chi_{t<x^1}\cdot \psi_0)$ is at comparable
modulation, therefore  \eqref{char_dual} is immediate due to the $L^\infty$ estimate for $\chi$. We now proceed to prove \eqref{char_dual} in a series of steps:\ret

\step{1}{Controlling the Strichartz norms} Due to the boundedness
of $\chi$,  we easily have:
\begin{equation}
        \lp{P_0(P_{<10}\chi\cdot \psi_0)}{L^q_t(L^r_x)}
        \ \lesssim \ \lp{\psi_0}{DS} \ . \notag
\end{equation}\ret

\step{2}{Controlling the $X^{s,b}$ norm} Our first order of
business is to bound the $X_\infty^{0,\frac{1}{2}}$ part of the norm
\eqref{dtS_norm}. Freezing the outer modulation, our goal is to show
that:
\begin{equation}
    \lp{ Q_jP_0\big(P_{<10}\chi\cdot \psi_0 )}{L^2_t(L^2_x)} \ \lesssim \
    2^{-\frac{j}{2}}\lp{\psi_0}{DS} \ . \label{X_est}
\end{equation}
We now split into subcases.\ret

\step{2.A}{Output far from cone} In this step we consider the
contribution of  output modulations $j > 20$. In this case, we may
further localize the product to $Q_jP_0(P_{<10}\chi\cdot Q_{j +
O(1)}\psi_0)$. Estimate \eqref{X_est} follows immediately from
$L^\infty$ control of $\chi$.\ret

\step{2.B}{$\chi$ at low frequency ($\leq j-10$)} In this case $\phi_0$ must be at
modulation $2^j$ therefore we consider the
contribution of the expression  $Q_jP_0(P_{<j-10}\chi\cdot Q_{j +
O(1)}\psi_0)$. Then \eqref{X_est} is immediate from
$L^\infty$ control of $\chi$.\ret

\step{2.C}{$\chi$ at medium frequency, $\psi$ at larger modulation}
In this case we consider the contribution of the term $Q_jP_0(P_{[j-10,10]}
\chi\cdot Q_{\geqslant j - 20}\psi_0)$.
Again, only the boundedness of $\chi$ is used.
\ret

\step{2.D}{$\chi$ at medium frequency, $\psi$ at low modulation} The
contribution of $Q_jP_0(P_{[j-10, 10]}\chi\cdot Q_{< j
- 20}\psi_0)$ is considered here. This is the main term. Without loss
of generality, we may assume that we are in a $(++)$ interaction,
which we decompose into all possible angular sectors of cap size
$|\kappa|\sim 2^{\frac{1}{2}j-10}$, respectively
$|\kappa'|\sim 2^{\frac{1}{2}j}$:
\begin{equation}
    Q_jP_0(Q^+ P_{[j-10,10]}\chi\cdot Q^+_{< j - 20}\psi_0)
     =\hspace{-.15in}
    \sum_{j-10 \leqslant k < 10}\sum_{\kappa,\kappa'}
    Q_j P_{0,\kappa'}\big(Q^+ P_k(\chi)\cdot P_{0,\kappa}Q^+_{<j-20}\psi_0\big) \ . \notag
\end{equation}
The main difficulty here is that we cannot really sum over $k$,
because $\chi$ is only in an $\ell^\infty$ type Besov space.
However, using Lemma 11 of \cite{Tao_WM} we see that the above sum
is both essentially \emph{diagonal} in $\kappa,\kappa'$, and
essentially \emph{frequency disjoint} in its contribution of angles
for each fixed $k$. Precisely, two sectors $\kappa,\kappa'$ and a
frequency $k$ can provide nonzero output if and only if:
\begin{equation}
    \hbox{dist}(\kappa,\kappa')\sim 2^{\frac{j}{2}-10} \ , \qquad
    \hbox{dist}(\kappa,(-1,0))\sim 2^\frac{j-k}{2}\ . \notag
\end{equation}
In particular the sector $\kappa'$ centered at $(1,0)$ does not
yield any output. Taking this into account we may bound:
\begin{align}
    &\lp{Q_jP_0(Q^+P_{j-10 \leqslant \cdot < 10}\chi\cdot
    Q^+_{< j -20}\psi_0)}{L^2_t(L^2_x)}^2 \  \notag\\
    \lesssim \
    &\sum_{k=j-10}^{10} \sum_{\substack{
    \hbox{dist}(\kappa,(-1,0))\sim 2^\frac{j-k}{2}\\
    \hbox{dist}(\kappa,\kappa')\sim 2^\frac{j}{2}-10 }}
    \lp{P_{0,\kappa'}\big(Q^+ P_k(\chi)\cdot P_{0,\kappa}
    Q^+_{<j-20}\psi_0\big)}{L^2_t(L^2_x)}^2 \  \notag\\
    \lesssim \
    &\sum_{k=j-10}^{10}\sum_{\hbox{dist}(\kappa,(-1,0))\sim 2^\frac{j-k}{2}}
     \!\!\!    \lp{Q^+ P_k(\chi)}{L^2_{t_{(1,0)}}(L^\infty_{x_{(1,0)}})}^2
    \lp{P_{0,\kappa}Q^+_{<j-20}\psi_0}{L^\infty_{t_{(1,0)}}(L^2_{x_{(1,0)}})}^2  \notag\\
    \lesssim \ &\sum_{k=j-10}^{10} \sum_{\hbox{dist}(\kappa,(-1,0))\sim 2^\frac{j-k}{2}}
    2^{-k}\sup_{\hbox{dist}(\omega,\kappa)\sim 2^\frac{j-k}{2}}
    \lp{P_{0,\kappa}Q^+_{<j-20}\psi_0}{L^\infty_{t_\omega}(L^2_{x_\omega}) }^2  \notag\\
    \lesssim \ &2^{-j}\sum_\kappa
    \lp{P_{0,\kappa}Q^+_{<j-20}\psi_0}{S[0,\kappa]}^2 \ . \notag
\end{align}
From the definition of the $S$ norm \eqref{S_norm}, this suffices to
prove \eqref{X_est}.\ret

\step{3}{Controlling the square sum of $S[0,\kappa]$ norms} Again
freezing $j < -10$ we need to demonstrate that:
\begin{equation}
        \sup_\pm \sum_\kappa \lp{Q^\pm_{<j} P_{0,\pm\kappa}(P_{<10}\chi\cdot \psi_0)}
        {S[0,\kappa]}^2 \ \lesssim \ \lp{\psi_0}{S}^2 \ , \label{char_ang}
\end{equation}
where angular sector size is $|\kappa|\sim 2^{\frac{1}{2}j}$. The
subcases repeat \textbf{Case 2} above with little difference, and
are mostly left to
the reader:\ret

\step{3.A}{$\chi$ at low frequency} This is the contribution of
the expression
\[
Q_{<j}P_0(P_{<j-10}\chi\cdot Q_{<j + O(1)}\psi_0) \ .
\]
In this case \eqref{char_ang} is
immediate from the $L^\infty$ control of $\chi$.
\ret

\step{3.B}{$\chi$ at medium frequency, $\psi$ at larger modulation}
As before, this is the term $Q_{<j}P_0(P_{j-10,10]}\chi\cdot Q_{\geqslant j - 20}\psi_0)$ for which we have a stronger
$L^2$ bound:
\begin{equation}
    \lp{Q_{<j}P_0(P_{[j-10,
    10]}\chi\cdot Q_{\geqslant j - 20}\psi_0)}{L^2_t L^2_x}  \ \lesssim \
    2^{-\frac{1}{2}j} \lp{\psi_0}{S} \ . \notag
\end{equation}
\ret

\step{3.C}{$\chi$ at medium frequency, $\psi$ at low modulation}
Here we consider the contribution of $Q_{<j}P_0(P_{[j-10,10]}\chi\cdot Q_{< j - 20}\psi_0)$. This is again the main term.
Without loss of generality  we may assume that we are in a $(++)$
interaction in terms of output and $\psi_0$ modulation (in
particular, from the estimate in \textbf{step 2} above we may
dispense with the case $(+)$ output and $(-)$ input from $\psi_0$),
and we again use Lemma 11 of \cite{Tao_WM} to decompose into a diagonal
sum over caps of size $|\kappa|\sim 2^{\frac{1}{2}j-C}$, respectively
$|\kappa'|\sim 2^{\frac{1}{2}j}$:
\begin{equation}
    Q_{<j}^+P_0(Q^+ P_{[j-C,C]}\chi\cdot Q^+_{< j - C}\psi_0)
 =
    \sum_{\hbox{dist}(\kappa,\kappa')\sim 2^{\frac{j}{2}}}\!\!\!\!
    Q_{<j}^+ P_{0,\kappa'}\big(P_{[j-C,C]}(\chi)
    \cdot P_{0,\kappa}Q^+_{<j-C}\psi_0\big)  \notag
\end{equation}
Notice that we do no need to frequency localize the factor $P_{[j-10,10]
}\chi$ to obtain this diagonally, which is a good
thing because the rougher bounds on the output modulation and that
of $\psi_0$ do not win us disjoint angular contributions in the
$k$-sum of $P_k\chi$. Plugging the above decomposition into the LHS
side of estimate, \eqref{char_ang} the RHS bound follows at once
from $L^\infty$ control of $\chi$.
\end{proof}\ret

\begin{proof}[Proof of Lemma \ref{para_trilin_lem}]
We  begin by  extending $\psi_k$ via the universal extension in
Proposition \ref{ext_prop} in
such a way that we simultaneously maintain the $\underline{E}$ and
$S$ norm control. The functions $\phi_{k_i}^{(i)}$ are similarly extended.
Thus, it suffices to prove the bound on all of
space-time.\ret

The constant $C$ will be fixed in the proof in just a moment. Let
$m\geqslant 0$ be any fixed integer. Without loss of generality, we
may assume that $k=0$. Furthermore, we may also assume that
$|k_1-k_2| <  \delta m$, for otherwise the estimate follows
immediately from an application of the standard trilinear bound
\eqref{standard_est_tri}, and taking $C>1$ on the RHS of
\eqref{imbal_trilin}. The proof will be accomplished in a series of
steps:\ret

\step{1}{Reduction to a bilinear estimate} In this step we consider
the contribution of  $\phi^{(1)}_{k_1}Q_{< k_2
-\delta m }\big(\partial_\alpha \phi^{(2)}_{k_2} \partial^\alpha
\psi_0\big)$. By an application of the estimates \eqref{XN_est},
\eqref{standard_est_bi}, and \eqref{XS_N_est} we easily have that
for $j<k_2-\delta m$ (which also implies $j<k_1$) and $\delta\ll 1$
sufficiently small:
\begin{equation}
    \lp{\phi^{(1)}_{k_1}Q_{ j }
    \left(\partial_\alpha \phi^{(2)}_{k_2} \partial^\alpha \psi_0\right)}{N} \
    \lesssim \ 2^{-\delta (k_1-j)} \lp{\phi^{(1)}_{k_1}}{S}\, \lp{\phi^{(2)}_{k_2}}{S}\cdot
    \lp{\psi_0}{S} \ . \notag
\end{equation}
Summing over all $j<k_2-\delta m$ we directly have
\eqref{imbal_trilin} for this component. It remains to estimate
the contribution of $\phi^{(1)}_{k_1}Q_{> k_2
-\delta m }\big(\partial_\alpha \phi^{(2)}_{k_2} \partial^\alpha
\psi_0\big)$.  We peel off the factor
$\phi^{(1)}_{k_1}$ from the trilinear estimate via the bound
\eqref{N_prod_est1}. It remains to prove the bilinear bound
\begin{equation}
 \| Q_{> k_2
-\delta m }\big(\partial_\alpha \phi^{(2)}_{k_2} \partial^\alpha
\psi_0\big)\|_{N} \lesssim_F
    \left(2^{-\delta^2 m}\lp{\psi_k}{S[I]} + 2^{m}\lp{\psi_k}{\underline{E}[I]}\right)
\end{equation}
\ret

\step{2}{$\psi_0$ is far from the cone} In this step we consider the
contribution of $Q_{\geqslant k_2 -\delta m } \big(\partial_\alpha
\phi^{(2)}_{k_2} Q_{\geqslant 0}\partial^\alpha \psi_0\big)$.  We will prove that the remaining
null-form is an $X^{0,-\frac{1}{2}}_{1}$ atom. In the present case,
we freeze the output modulation $j $ and then estimate:
\begin{align}
    \lp{Q_{j } \left(\partial_\alpha
    \phi^{(2)}_{k_2} Q_{\geqslant 0}\partial^\alpha \psi_0\right)}{L^2_t(L^2_x)}
    \ &\lesssim \ \lp{\nabla_{t,x}\phi^{(2)}_{k_2} }{L^\infty_t(L^\infty_x)}
    \cdot\lp{Q_{\geqslant 0}\nabla_{t,x}\psi_0}{L^2_t(L^2_x)} \ , \notag\\
    &\lesssim \ 2^{k_2}\lp{\phi^{(2)}_{k_2}}{S}\cdot \lp{\psi_0}{S} \ .
    \notag
\end{align}
Multiplying both sides of this bound by $2^{-\frac{j}{2}}$ and then
summing  over all dyadic $j\geqslant k_2 -\delta m$ we arrive at:
\begin{equation}
        \lp{Q_{\geqslant k_2 -\delta m } \left(\partial_\alpha
        \phi^{(2)}_{k_2} Q_{\geqslant 0}\partial^\alpha \psi_0\right)}{N} \
        \lesssim \ 2^{\frac{k_2}{2}+ \delta m}\lp{\phi^{(2)}_{k_2}}{S}\cdot \lp{\psi_0}{S}
        \ , \notag
\end{equation}
which suffices due to the condition $k_2<-m$.\ret

\step{3}{$\phi^{(2)}_{k_2}$ is far from the cone} In this step we
consider the contribution of $Q_{\geqslant k_2 -\delta m }
\big(Q_{>k_2-8\delta m}\partial_\alpha \phi^{(2)}_{k_2}
Q_{<0}\partial^\alpha \psi_0\big)$. In this case, we again freeze
the output modulation $j$ and proceed to bound:
\begin{align}
        &\lp{Q_{j}
        \big(Q_{>k_2-8\delta m}\partial_\alpha \phi^{(2)}_{k_2}
        Q_{<0}\partial^\alpha \psi_0\big)}{L^2_t(L^2_x)}
       \notag\\
        \lesssim \ &\lp{Q_{>k_2-8\delta m}
        \nabla_{t,x}\phi^{(2)}_{k_2} }{L^2_t(L^\infty_x)}
        \cdot\lp{\nabla_{t,x}\psi_0}{L^\infty_t(L^2_x)} \ .  \label{step3_frozenj}
\end{align}
By summing over all $j > k_2-8\delta m$ in  estimate
\eqref{L2Linfty_est} we have that:
\begin{equation}
        \lp{Q_{>k_2-8 \delta m}
        \nabla_{t,x}\phi^{(2)}_{k_2} }{L^2_t(L^\infty_x)}
        \ \lesssim \ 2^{\frac{1}{2}k_2 +4 \delta m} \lp{\phi^{(2)}_{k_2}}{S} \ . \notag
\end{equation}
Substituting this into the RHS of \eqref{step3_frozenj}, multiplying
the result by $2^{-\frac{1}{2}j}$, and then summing over all $ j
\geqslant k_2 - \delta m $ we have the estimate:
\begin{equation}
        \lp{Q_{\geqslant k_2 - \delta m }
        \big(Q_{>k_2-8 \delta m}\partial_\alpha \phi^{(2)}_{k_2}
        Q_{<0}\partial^\alpha \psi_0\big)}{N} \ \lesssim \
        2^{5 \delta m}\lp{\phi^{(2)}_{k_2}}{S}\cdot\lp{\psi_0}{\underline{E}} \ . \notag
\end{equation}\ret

\step{4}{The core contribution} In this step we consider the
contribution of the expression $Q_{\geqslant k_2 -\delta m }
\big(Q_{<k_2-8\delta m}\partial_\alpha \phi^{(2)}_{k_2}
Q_{<0}\partial^\alpha \psi_0\big)$. This is the main case, and
requires a decomposition into angular sectors of cap size
$|\kappa|\sim 2^{-4\delta m}$. Without loss of generality we may
assume we are in the $(++)$ configuration. The other cases $(--)$,
$(+-)$, and $(-+)$ are the same with only minor modifications and
are therefore left to the reader. We break the entire contribution
into a $Q_{\geqslant k_2 - \delta  m}$ localized sum of two
principle terms $T_1$ and $T_2$, where:
\begin{align}
        T_1 \ &= \ \sum_{\kappa\in K_l} Q^+_{ <k_2-8 \delta m}P_{k_2,\kappa}
        \partial_\alpha \phi^{(2)}_{k_2}\cdot(I-P_{0,2\kappa})Q_{<0}^+\partial^\alpha
        \psi_0 \ , \notag\\
        T_2 \ &= \ \sum_{\kappa\in K_l} Q^+_{ <k_2-8 \delta m}P_{k_2,\kappa}
        \partial_\alpha \phi^{(2)}_{k_2}\cdot P_{0,2\kappa}Q_{<0}^+\partial^\alpha
        \psi_0 \ . \notag
\end{align}\ret

To help state the estimates, we introduce the following weaker
version of the $N\!F\!A^*$ portion of the $S[k,\kappa]$ norm from
line \eqref{str_norm}:
\begin{equation}
        \lp{\psi}{\td{S}_k} \ := \ \sup_{\substack{l > 10 \\
        \kappa\in K_l}} \lp{P_k \psi}{\td{S}[k,\kappa]}
        \ , \notag
\end{equation}
where:
\begin{equation}
        \lp{\psi}{\td{S}[k,\kappa]} \ := \
        \sup_\pm \sup_{ \omega\in\frac{1}{2}\kappa}
        2^k |\kappa|\cdot
        \lp{Q_{<k}^\pm(I-P_{k,\pm\kappa})\psi}{L^\infty_{t_\omega}(L^2_{x_\omega})}
        \ . \notag
\end{equation}
Notice that we do not use the more eccentric $Q_j$ multipliers for
$j<k-10$ in this definition, and there is no square-summing over
angles. The reason this notation is useful is that we have the
relation: $\lp{\psi_k}{\td{S}} \lesssim \lp{\psi_k}{\underline{E}}$.
This is shown through an application of the estimate:
\begin{equation}
          \lp{Q_{<k}^\pm(I-P_{k,\pm\kappa})\psi_k}{L^\infty_{t_\omega}(L^2_{x_\omega})}
         \ \lesssim \ 2^{-k} |\kappa|^{-1}\cdot
         \lp{\snabla_{t,x} \psi_k}{L^\infty_{t_\omega}(L^2_{x_\omega})}
         \ . \label{ang_est}
\end{equation}
Such an inequality may be proved by decomposing the multiplier
$Q_{<k}^\pm(I-P_{k,\pm\kappa})$ into a dyadic sum of angular sectors
of increasing size and spread from $\pm\kappa$. Without loss of
generality, we may assume we are in the ``+'' case, and we decompose
$Q_{<k}^+(I-P_{k,\kappa})=\sum_{1 > 2^j > |\kappa|}Q_{<k}^+
P_{k,\kappa_j}$ where each sector size is $|\kappa_j|\sim 2^{j}$
with distance $\hbox{dist}(\kappa,\kappa_j)\sim 2^j$. For each of
these sectors we use the uniform multiplier bounds:
\begin{equation}
    \lp{Q_{<k}^+ P_{k,\kappa_j}\psi_k}{L^\infty_{t_\omega}(L^2_{x_\omega})} \ \lesssim \
    2^{-k} |\kappa_j|^{-1}\cdot\lp{\snabla_{t,x}\psi_k}{L^\infty_{t_\omega}(L^2_{x_\omega})} \ , \notag
\end{equation}
which is an easy consequence of the fact that the  kernels
associated to the operators:
\begin{equation}
    \mathcal{L} \ = \ 2^{k} |\kappa_j|\snabla_{t,x}^{-1} Q_{<k}^+
    P_{k,\kappa_j} \ , \notag
\end{equation}
 are uniformly in $L^1_t(L^1_x)$. The  estimate
\eqref{ang_est} now follows from simply summing over this last bound
overall all dyadic $1<|\kappa_j|^{-1}<|\kappa|^{-1}$.\ret

Returning to the main thread, we first bound the term $T_1$ above.  In
this case, we are going to loose a large constant because the sum is
not well localized in the second factor an therefore we cannot use
orthogonality with respect to $\kappa$. Furthermore, we will not
bother to gain anything from the null-structure, because the frequency
localization of this term eliminates parallel interactions. To
compensate for the large number of non-orthogonal sectors, we may use the
$\td{S}$ norm for the second factor. Using the product estimate \eqref{PW_product} we may
bound:
\begin{align}
        \lp{Q_j T_1}{L^2_t(L^2_x)} \
        \lesssim  & \
        \sum_{\kappa\in K_l} |\kappa|^\frac{1}{2}2^{-\frac{1}{2}k_2}
        \lp{Q^+_{ <k_2-8 \delta m}P_{k,\kappa}
        \partial_\alpha \phi^{(2)}_{k_2}}{S[k_2,\kappa]} \notag
        \\ & \ \ \ \ \ \ \cdot\sup_{\omega\in \kappa}\lp{(I-P_{0,2\kappa}) Q_{<0}^+\partial^\alpha
        \psi_0}{L^\infty_{t_\omega}(L^2_{x_\omega})} \notag\\
        \lesssim & \  2^{\frac{1}{2}k_2}\lp{\phi^{(2)}_{k_2}}{S}\cdot
        2^{4\delta m}\lp{\psi_0}{\td{S}} \ . \notag
\end{align}
Multiplying both sides of this last estimate by the factor
$2^{-\frac{1}{2}j}$ and then summing over all $j > k_2- \delta m$ we
have:
\begin{equation}
        \lp{Q_{> k_2- \delta m} T_1}{N} \ \lesssim \
        2^{5 \delta
        m}\lp{\phi_{k_2}^{(2)}}{S}\cdot\lp{\psi_0}{\td{S}} \ , \notag
\end{equation}
which is sufficient.\ret

Our final task here is to bound the term $Q_{> k_2- \delta m} T_2$
in the space $X_1^{0,\frac{1}{2}}$. Notice that because of the
angular and $(++)$ localization, as well as the fact that $j> k_2-
\delta m$, for each $Q_j T_2$ we may freely insert the multiplier
$Q_{> j-10}$ in front of the second factor, because the complement
vanishes (see Lemma 11 of \cite{Tao_WM}). In this case the resulting
sum  is both diagonal and orthogonal in $\kappa$, so freezing $Q_j
T_2$ we have with the aid of Bernstein's inequality
\eqref{Bernstein} the estimate:
\begin{align}
        \lp{Q_j T_2}{L^2_t(L^2_x)}^2  &\lesssim
        \sum_{\kappa\in K_l} \lp{Q^+_{ <k_2-8 \delta m}P_{k,\kappa}
        \partial_\alpha \phi^{(2)}_{k_2}}{L^\infty_t(L^\infty_x)}^2
        \cdot\lp{P_{0,2\kappa}Q_{[j-10,0]}^+\partial^\alpha
        \psi_0}{L^2_t(L^2_x)}^2 \notag\\
        &\lesssim  \left(
        2^{k_2-2 \delta m} \lp{\phi_{k_2}^{(2)}}{S}\cdot
        2^{-\frac{1}{2}j}\lp{\psi_0}{S}
        \right)^2 \ . \notag
\end{align}
Multiplying the root of this inequality by the factor
$2^{-\frac{1}{2}j}$ and then summing over all $j > k_2 - \delta m$
we finally have:
\begin{equation}
        \lp{Q_{> k_2-\delta m} T_2}{N} \ \lesssim \
        2^{- \delta  m} \lp{\phi_{k_2}^2}{S}\cdot\lp{\psi_0}{S} \ . \notag
\end{equation}
This concludes our proof of estimate \eqref{imbal_trilin}.
\end{proof}

\ret
%-------------------------------------------------------------------------
%%%%%%%%%%%%%%%%%%%%%%%%%%%%%%%%%%%%%%%%%%%%%%%%%%%%%%%%%%%%%%%%%%%%%%%%%%
%-------------------------------------------------------------------------

\section{Structure of Finite $S$ Norm Wave-Maps and Energy Dispersion}\label{s:wm_struct}

In this section we prove Proposition~\ref{p:wm_struct}. There is
almost nothing to do for \eqref{ext__phi_bnds}. The $\uX$ bound
follows from the reduced version of \eqref{core_L2} in Remark \ref{trilinear_rem},
while the $\underline{E}$ bound follows
from energy estimates on null surfaces.\ret

%---------------------------------------------------------------------------

\subsection{Renormalization}
Here we establish the  renormalization bound
\eqref{renorm_phi_bnds}. Our starting point is the construction of
the renormalization matrix $U$ in Proposition~\ref{propphiu}. The
frequency localized wave-map equation for $\phi$ is given by:
\begin{equation}
        \Box \phi_k  \ = \
        -P_k \big( \mathcal{S}(\phi)\partial^\alpha
        \phi\partial_\alpha \phi\big) \ . \label{vanilla_wm}
\end{equation}
For each index $m$ the RHS of this expression can be written in
terms of the trilinear from $T^{m}_k$ from line \eqref{basic_trilin}
as follows:
\begin{equation}
        P_k \big( \mathcal{S}(\phi)\partial^\alpha
        \phi\partial_\alpha \phi\big) \ = \
    2 \mathcal{S}(\phi)_{<k-m}\partial^\alpha_{<k-m}
        \phi\partial_\alpha \phi_k
     + T_{1;k}^{m} \big( \mathcal{S}(\phi),\partial^\alpha
        \phi,\partial_\alpha \phi\big) \ . \notag
\end{equation}
Using the identity \eqref{sff_iden} we have:
\begin{equation}
     P_k\big(\mathcal{S}^\dagger(\phi)\partial^\alpha\phi\partial_\alpha\phi_{<k-m+2}\big)
     \ = \ 0 \ . \notag
\end{equation}
Thus, we may further write:
\begin{equation}
        \mathcal{S}(\phi)_{<k-m}\partial^\alpha_{<k-m}
        \phi\partial_\alpha \phi_k \ = \
        \big( \mathcal{S}(\phi)_{<k-m} -\mathcal{S}(\phi)_{<k-m}^\dagger \big)\partial^\alpha
        \phi_{<k-m}\partial_\alpha \phi_k
    +T^m_{2;k} \ , \notag
\end{equation}
where $T^m_{2;k}$ is obtained by applying the decomposition
\eqref{basic_trilin} to the previous line.
Therefore, we have written
the original frequency localized wave-map equation in the form:
\begin{equation}
        \Box \phi_k \ = \ -2A^\alpha_{<k-m}\partial_\alpha \phi_k + \sum_i T_{i;k}^{m} \ ,
        \label{WM_Tm_form}
\end{equation}
where the $T_{i;k}^{m}$ are trilinear forms as on line
\eqref{basic_trilin} with $O(m)$ gap indices. By an
application of estimates \eqref{easy_linearized_est} and
\eqref{large_core_N} with $m=20$ we have the bound:
\begin{equation}
        \llp{\phi_k}{\W[I]} \ \lesssim_F \  c_k \ , \notag
\end{equation}
where $\{c_k\}$ is some $S[I]$ frequency envelope for $\phi_k$. This
proves \eqref{renorm_phi_bnds}.\ret

%-------------------------------------------------------------------------------

\subsection{Partial fungibility of the $S$ norm}
Here we prove that there is always a decomposition of intervals
$I=\cup_{i,l}^{K(F)} I_{il}$ where $K(F)$ is some polynomial in the
$S[I]$ norm of $\phi$, and where \eqref{fung_bnd} holds in each
subinterval. Our starting point is the series of frequency localized
equations \eqref{WM_Tm_form}.
For a fixed $\phi_k$ we  use \eqref{WM_Tm_form} with $m=20$.
As in the previous section, we can find a renormalization
$w_{,k} = U_{,<k-20}\phi_k$
on all of $I$ such that:
\begin{equation}
        \lp{P_{k'}\Box w_{,k} }{N[I]} \ \lsf \ 2^{-|k-k'|}c_k \ . \label{box_wk_line}
\end{equation}
Let $\eta \ll 1$ to be chosen later. By the fungibility property
\eqref{N_fung} (and continuity) there exists a polynomial $K_1$ in
$F\eta^{-1}$ such that $I=\cup_i^{K_1}I_i$ such that:
\begin{equation}
        \lp{P_{k'}\Box w_{,k}}{N[I_i]} \ \leqslant \ 2^{-\frac12|k-k'|}\eta c^i_k \ , \notag
\end{equation}
where $\{c_k^i\}$ are now some unit normalized frequency envelope
which may depend on the interval $I_i$. We label each time
interval as $I_i=[t_i,t_{i+1}]$, and on each of these time slabs we
write $ w_{,k}= w_{,k}^{free} + w_{,k}^{source}$ where
$w_{,k}^{free}$ is a free wave with data $w_{,k}[t_i]$. By the
previous line and the energy estimate \eqref{energy_est} we have on
$I_i$ the bound:
\begin{equation}
        \lp{P_{k'} w_{,k}^{source}}{S[I_i]} \ \lesssim \
        2^{-\frac{1}{2}|k-k'|}\eta c_k^i \ . \label{w_source_est}
\end{equation}
Consequently, for the corresponding part $ U^\dagger_{,<k-20}w_{,k}^{source}$
of $\phi_k$ we obtain:
\begin{equation}
    \lp{ U^\dagger_{,<k-20}w_{,k}^{source} }{S[I_i]}  \ \lesssim_F  \ \eta c_k^i \ . \notag
\end{equation}
By choosing $\eta$ as the reciprocal of an appropriate polynomial in $F$, we have:
%from \eqref{S_prod2}--\eqref{S_prod3} and \eqref{env_est1} the estimate:
\begin{equation}
        \lp{U^\dagger_{,<k-20}w_{,k}^{source}}{S[I_i]} \ \leqslant\ c_k^i \ . \label{w_source_subest}
\end{equation}
It remains to bound the free wave contribution  $U^\dagger_{,<k-20}w_{,k}^{free}$
on each of the intervals $I_i$, or on some further subdivision thereof.
\ret

Unfortunately we do not directly know that $U^\dagger_{,<k-20}$ is
manageable on $I_i$. However, we do have from estimate
\eqref{env_est0} and the energy bound for $\phi$ that:
\begin{equation}
        \lp{P_{k'}w_{,k}^{free}[t_i]}{\dot{H}\times L^2} \ \lesssim_E \
2^{-|k-k'|}c_k \ , \notag
\end{equation}
uniformly with respect to $i$ where we may choose the unit frequency
envelope $\{c_k\}$ to be the same as  on line \eqref{box_wk_line} above.
In particular, we have the uniform control:
\begin{equation}
 \lp{P_{k'}w_{,k}^{free}}{S[I_{i}]} \ \lesssim_E \ 2^{-|k-k'|}c_k \ . \label{w_free_est}
\end{equation}

Now we turn our attention to the $U_{,<k-20}$'s. Given a large parameter
$m$ to be chosen later, we consider the sections $P_{[j-m,j+m]}
U_{,<k-20}$ of $U_{,<k-20}$. Recall that from Remark \ref{direct_para_rem}
each $U_{,<k-20}$ is built up out of the same connection
\eqref{B_def}, and therefore the bounds \eqref{env_est1} for each
$U_{,<k-20}$ may be taken in terms of the \emph{same} frequency
envelope. Hence, except for  a polynomial in $mF$ number of indices $j$
we already have:
\begin{equation}
    \sup_{k}  \lp{P_{[j-m,j+m]}U_{,<k-20}}{S[I]} \ \leqslant\ 1 \ . \label{good j}
\end{equation}
Such indices $j$ are called ``good $j$'s''; the remainder (of which
we have at most a polynomial in $mF$) are called  ``bad $j$'s''.
We also introduce the corresponding parts of  $U^\dagger_{,<k-20}w_{,k}^{free}$:
\begin{equation}
        \phi_k(j) \ = \ P_{[j-m,j+m]}U_{,<k-20}^\dagger \cdot w_{,k}^{free} \
        . \notag
\end{equation}
% In particular,
%by the bounds \eqref{env_est1} and the algebra estimates
%\eqref{S_prod2}--\eqref{S_prod3} and \eqref{w_free_est},

The goal of the argument is now to choose a polynomial in $mF$
collection of subintervals $I_{il}$, partitioning the $I_i$, such that on
each there is the uniform control over all $k$ and $j$:
\begin{equation}
        \lp{P_{k'} \phi_k(j)}{S[I_{il}]} \ \lesssim_E \ 2^{-\frac{1}{2}\delta|k-k'|}c_k^{jl} \ .
        \label{subinterval_est}
\end{equation}
for some additional set of unit normalized frequency envelopes $\{c_{k}^{jl}\}$.  For good $j$'s
this is straightforward in view of \eqref{w_free_est} and \eqref{good
j}. Since there are $\lesssim mF$ bad $j$'s, it suffices to consider
a fixed such bad $j$. The equation for each fixed $\phi_k(j)$ is:
\begin{equation}
        \Box\phi_k(j) \ = \ P_{[j-m,j+m]}\Box U_{,<k-20}^\dagger \cdot w_{,k}^{free}
        + 2  \partial^\alpha P_{[j-m,j+m]} U_{,<k-20}^\dagger \partial_\alpha  w_{,k}^{free}
        \ . \notag
\end{equation}
Therefore, by a direct application of the estimates
\eqref{w_free_est}, \eqref{env_est1}, \eqref{standard_est_bi}, and
\eqref{env_est1h}--\eqref{env_est2} we have on all of $I_i$ the bound:
\begin{equation}
        \lp{\Box P_{k'} \phi_k(j)}{N[I_i]} \ \lesssim_{mF} \
        2^{-\delta |k-k'|}c_k \ , \notag
\end{equation}
and from the energy norm control giving \eqref{w_free_est} and
estimate \eqref{env_est0} we also have the uniform energy control:
\begin{equation}
        \lp{P_{k'}\phi_k(j)[t]}{\dot{H}^1\times L^2}
        \ \lesssim_E \ 2^{-|k-k'|}c_k \ . \notag
\end{equation}
Thus, by again using the property \eqref{N_fung} we obtain the desired
partition $\{I_{il}\}$ of $I$, with estimate \eqref{subinterval_est} uniformly,
at a cost of at most $\lesssim_{mF} 1$ subdivisions.\ret

To conclude the proof we need to estimate
$U^\dagger_{,<k-20}w_{,k}^{free}$ on each subinterval $J=I_{il}$,
which is now fixed with the property that \eqref{subinterval_est}
holds.  We split $U^\dagger_{,<k-20}$ into:
\begin{equation}
    U^\dagger_{,<k-20} \ =  \ P_{<k-m} U^\dagger_{,<k-20}
    + P_{[k-m,k+m]}  U^\dagger_{,<k-20}
    +  P_{>k+m} U^\dagger_{,<k-20} \ . \notag
\end{equation}
For the high frequency part we use \eqref{env_est1h} in conjunction with the
product bounds \eqref{S_prod2}--\eqref{S_prod3} to obtain:
\begin{equation}
    \lp{ P_{k'}( P_{>k+m} U^\dagger_{,<k-20} w_{,k}^{free})}{S[I]}  \
    \lesssim_F \ 2^{-Cm}2^{-|k'-k|} c_k \ , \notag
\end{equation}
which suffices provided $m$ is large enough, $m \sim \ln{F}$.  For the
medium frequency part we can use directly \eqref{subinterval_est} with
$j=k$.  Thus, we are reduced to providing good $S[J]$ norm bounds for
the quantities $P_{<k-m}U^\dagger_{,<k-20}\cdot w_{,k}^{free}$ which
are localized at frequency $2^k$. We do this in a series of steps
depending on what component of the $S[J]$ norm is being
considering:\ret

%-----------------------------------------------------------------------------------
\step{1}{Energy and Strichartz norm control} For any of the Strichartz norms
we immediately have from Leibnitz rule, estimates \eqref{w_free_est} and \eqref{env_est0} the
bound:
\begin{equation}
    \lp{\nabla_{t,x}\big( P_{<k-m}U^\dagger_{,<k-10}\cdot
    w_{,k}^{free}\big)}{D\underline{S}[J]} \  \lesssim_E c_k
     \ , \notag
\end{equation}
which is sufficient.\ret

%-----------------------------------------------------------------------------------
\step{2}{$X_\infty^{0,\frac{1}{2}}$ norm control} Fix a modulation
$Q_{j}$. Without loss of generality we will assume that $j<k$, as the
complimentary region is easier to treat using the high modulation
bounds in \eqref{env_est1} and \eqref{env_est1h}.  We decompose as
follows:
\begin{equation}
    Q_j \big( P_{<k-m}U^\dagger_{,<k-20}\cdot
    w_{,k}^{free}) \ =  \
    \begin{cases}
        Q_j \phi_k(j) + Q_j R_{,k}  \ , & j < k-2m;\\
    Q_j \phi_k(k-2m) +   Q_j R_{,k} \ , & j > k-2m.
    \end{cases}
    \label{Qj_wdecomp}
\end{equation}
where:
\begin{equation}
    R_{,k} \ =  \
    \begin{cases}
        P_{<j-m}U^\dagger_{,<k-20}\cdot w_{,k}^{free} +
        P_{[j+m,k-m]}U^\dagger_{,<k-20}\cdot
        w_{,k}^{free}  \ , &j < k-2m;\\
    P_{<k-3m}U^\dagger_{,<k-20}\cdot
    w_{,k}^{free} \ , & j > k-2m.
    \end{cases}\notag%\\
     %&= \ R_1 + R_2 + R_3 \ . \notag
\end{equation}
By estimate \eqref{subinterval_est} we already control the first terms
on the RHS  of \eqref{Qj_wdecomp}, so we only need to
bound the contribution of $Q_jP_k R_{,k}$. This is given by the following
analog of Lemma \ref{lvvlow}:

\begin{lem}\label{lvvlow_imp}
  Let $j < k-10$ and $m>10$ an integer. Then the following estimates
  hold for test functions $u=u_{<k-10}$ and $\phi_k$:
\begin{align}
    \lp{\! Q_j (  u_{<j-m} \phi_k)\! }{X^{1,\frac12}_\infty}\!\! +\!\!
     \lp{\! Q_{<j} ( u_{<j-m}  \phi_k)\!  }{S[k;j]}
     \! \lesssim \! \big(
     \lp{\! u\! }{L^\infty_t(L^\infty_x)}\!\!  +\!\!  2^{-\delta m}
     \lp{\! u\! }{\uX}\big)  \lp{\! \phi_k\! }{S}  ,
     \notag\\
     \lp{\! Q_j ( P_{>j+m} u_{<k-10} \phi_k)\! }{X^{1,\frac12}_\infty}\! +\!
     \lp{\! Q_{<j} ( P_{>j+m} u_{<k-10} \phi_k)\!  }{S[k;j]}
        \lesssim   2^{-\delta m}
     \lp{\! u\! }{S}  \lp{\! \phi_k \! }{S}
      . \notag\end{align}
\end{lem}\ret

\begin{proof}%[Proof of estimates \eqref{u<jphik_imp}--\eqref{u<jphik_imph}]
The proof of the first bound is immediate from
\eqref{S_prod4} and  the product bounds \eqref{S_prod2}--\eqref{S_prod3}
in conjunction with the following easy estimate for very high modulations:
\begin{equation}
    \lp{Q_{>j-\frac{1}{2}m}u_{<j-m}}{S} \
    \lesssim \ 2^{-\frac{1}{4}m}\lp{u}{\uX} \ . \notag
\end{equation}
The second estimate  is just a summed version of \eqref{XS_X_est}
which also incorporates \eqref{XS_est}.
\end{proof}\ret

Using a combination of the estimates  in this last Lemma, and
\eqref{w_free_est}, we have:
\begin{equation}
    \lp{Q_jP_k R_{,k}}{X_\infty^{0,\frac{1}{2}}} \ \lesssim_E \
    \big(1 + 2^{-\delta m}Q_{2}(F)\big) c_k  \ , \notag
\end{equation}
which suffices.\ret

%-----------------------------------------------------------------------------------
\step{3}{$S[k;j]$ norm  control} This is immediate from the
decomposition \eqref{Qj_wdecomp}, the estimate
\eqref{subinterval_est}, and Lemma \ref{lvvlow_imp}.\ret

%-----------------------------------------------------------------------------------

\subsection{The role of the energy dispersion}
By applying estimate \eqref{renorm_phi_bnds}
and then using \eqref{core_L2} on equation \eqref{vanilla_wm}
we have \eqref{small_mod_est}.\ret

Suppose now that $\{c_k\}$ is a frequency envelope for the
initial data of $\phi$ in $\dot{H}\times L^2$. Then by the
seed bounds \eqref{seed_S} we have the full control:
\begin{equation}
    \lp{\phi}{S_c[J]} \ \leqslant \ K_1(F) \ , \label{phi_boot_better}
\end{equation}
on some sufficiently small subinterval $J\subseteq I$. Here $K_1$
is a universal polynomial that will be chosen in a moment. The goal now
is to bootstrap this control and show that if:
\begin{equation}
    \lp{\phi}{S_c[J]} \ \leqslant \ 2K_1(F) \ , \label{phi_boot_hyp}
\end{equation}
then we have \eqref{phi_boot_better}. By Proposition \ref{propboot}
we may continue and finally close this last estimate on all of $I$.

By applying estimates \eqref{renorm_phi_bnds} and \eqref{small_mod_est}
to \eqref{phi_boot_hyp}, we have:
\begin{equation}
    \lp{\phi}{\W_c[J]} \ \leqslant \ K_2(F)K_1(F) \ ,
    \qquad \lp{\phi}{\uX_c[J]} \ \leqslant \ \epsilon^{\delta_1}K_2(F)K_1(F)\ ,
    \notag
\end{equation}
for a universal polynomial $K_2$.

Next, choose the gap $m\lesssim \ln(F)$ in equation \eqref{WM_Tm_form}
in a way that is consistent
with the assumptions of Proposition \ref{p:para}, and apply estimate \eqref{linearized_est}
to \eqref{WM_Tm_form}, while using \eqref{core_N} via the last two bounds.
This gives:
\begin{equation}
    \lp{\phi_k}{S[I]} \ \leqslant \ K_3(F)\big(1 + \epsilon^{\delta_1^2}K_2(F)K_1(F)\big)c_k
     \ . \notag
\end{equation}
The proof is concluded by choosing $K_1=2K_3$ and assuming $\epsilon$ is sufficiently small.

\ret
%-------------------------------------------------------------------------
%%%%%%%%%%%%%%%%%%%%%%%%%%%%%%%%%%%%%%%%%%%%%%%%%%%%%%%%%%%%%%%%%%%%%%%%%%
%-------------------------------------------------------------------------

\section{Initial Data Truncation}
\label{s:cut}

Here we prove that for each initial data set with small energy
dispersion we can continuously regularize it. In a sufficiently
small tubular neighborhood $V(\mathcal{M})$ of the surface
$\mathcal{M}\subset \R^N$ we introduce a projection operator:
\begin{equation}
        \Pi: V(\mathcal{M})  \ \to  \ \mathcal{M} \ . \notag
\end{equation}
This also induces a projection operator on the tangent bundle:
\begin{equation}
    \Pi: TV(\mathcal{M})  \ \to  \ T\mathcal{M} \ , \notag
\end{equation}
which is a product of $\Pi$ in $\mathbb{R}^N$ and Euclidean linear
orthogonal projection onto each fiber in the second factor. Given an
initial data set:
\begin{equation}
        \phi[0] \ = \ (\phi_0,\phi_1): \R^{2}
        \ \to  \ T\mathcal{M} \ , \notag
\end{equation}
which belongs to $ \dot H^{1} \times L^2$, we regularize it as
follows:
\begin{equation}
        \phi_{,<k}[0]  \ =  \ \Pi(P_{<k}\phi[0]) \ . \notag
\end{equation}
The following result asserts that if $\phi[0]$ has small energy
dispersion then its regularizations are well defined, and stay close
to the corresponding Littlewood-Paley projections:\ret

\begin{prop} \label{p:cuta}
For each $E > 0$ there exists $\epsilon_0 > 0$ so that for each
initial data set $\phi[0]$ for \eqref{main_eq} with energy $E$ and
energy dispersion $\epsilon \leqslant \epsilon_0$ and $k,k_* \in \Z$
we have:
\begin{equation}
        \| P_{k}  (P_{< k_*} \phi[0] -
        \phi_{,<k_*}[0]) \|_{\dot H^1
        \times L^2}   \ \lesssim_E
        \ \min\{ \epsilon |\ln \epsilon|^2, 2^{- |k-k_*|}\}
        \ . \label{cutphi_red}
\end{equation}
\end{prop}\ret

\begin{proof}
By rescaling we assume that $k_*=0$. We begin with two simple Moser
type estimates which  we will repeatedly use in the sequel.
Precisely, for each smooth and bounded function $G$ with bounded
derivatives we have:
\begin{equation}
        \lp{ \nabla_{x}^J G(P_{<k} \phi_0)}{L^\infty_x} \ \lesssim_E \ 2^{|J| k}
         \ , \label{moser1}
\end{equation}
and:
\begin{equation}
    \lp{ \nabla_{x}^J G(P_{<k} \phi_0)}{L^2_x}
    \ \lesssim_E \ 2^{(|J|-1)k} \ ,
    \qquad |J|  \ \geqslant  \ 1 \ ,
    \label{moser2}
\end{equation}
which are easily proved using the chain rule and Bernstein's
inequality.\ret

We first show that if $\epsilon$ is small enough then the
projection $\Pi  P_{< 0} \phi[0]$ is well defined:

\begin{lem}
Under the assumptions of Proposition \ref{p:cuta} we have:
\begin{equation}
        \hbox{dist}\;( P_{< 0} \phi_0,\mathcal{M})  \ \lesssim_E \
        \epsilon  |\log\epsilon| \ . \notag
\end{equation}
\end{lem}
\begin{proof}
By translation invariance, it suffices to show that:
\begin{equation}
        I  \ = \ \int_{|x| \leqslant 1} |P_{< 0} \phi_0(0) - \phi_0(x)| dx
        \ \lesssim_E \
        \epsilon |\log\epsilon| \ . \label{close}
\end{equation}
We use a positive parameter $m$ and a Littlewood-Paley decomposition
to estimate $I$ as follows:
\begin{equation}
        I  \ \lesssim \
        \lp{ \nabla_x P_{<-m} \phi_0} {L^\infty_x} + \lp{ P_{[-m,m]} \phi_0}{L^\infty_x}
        + \lp{ P_{ > m} \phi_0}{L^2_x} \ . \notag
\end{equation}
Using Sobolev embeddings for the first term, energy dispersion for
the second, and the $\dot{H}^1$ norm for the third we obtain:
\begin{equation}
        I  \ \lesssim\  2^{-m}E + m \epsilon + 2^{-m} E \ . \notag
\end{equation}
Then \eqref{close} is obtained by choosing $2^m = \epsilon^{-1}$.
\end{proof}\ret

To continue the proof of the proposition, we remark that $\Pi$ can
be expressed as:
\begin{equation} \label{piform}
        \Pi(\phi^{(1)},\phi^{(2)}) \ = \
        \big(G(\phi^{(1)}),H(\phi^{(1)}) \phi^{(2)}\big) \ ,
\end{equation}
where $G$ is some smooth extension of $\Pi$ to all of
$\mathbb{R}^N$, and $H$ is some extension of the fiber projection
composed with $G$. Note that both $G$ and $H$ may be chosen as
bounded functions with bounded derivatives. We separately estimate
the high frequencies, middle frequencies and low frequencies of the
difference $P_{< k_*} \phi[0] - \phi_{,<k_*}[0]$.\ret

\step{1}{High frequency bounds, the contribution of $k > 0$} For the
high frequencies we do not use at all the fact that $\phi[0]$ takes
values in $T\mathcal{M}$. Instead, we use \eqref{moser1} to directly
estimate:
\begin{equation}
        \lp{ P_k G(P_{<0} \phi_0)}{L^2_x}
        \ \lesssim_E \ 2^{-(C+1)k} \ , \notag
\end{equation}
where $C$ is a large integer. Similarly we have:
\begin{equation}
        \lp{ P_k \big( H(P_{<0} \phi_0) P_{<0} \phi_1\big)}{L^2_x}
        \ \lesssim_E \ 2^{-Ck} \ . \notag
\end{equation}
Thus we  obtain:
\begin{equation}
        \lp{ P_{k}  (P_{<0} \phi[0] -
        \phi_{,<0}[0]) }{\dot{H}^1
        \times L^2} \  \lesssim_E \ 2^{-Ck} \ .
        \label{mllh}
\end{equation}\ret

%-----------------------------------------------------------------------------

\step{2}{Low frequencies bounds, the contribution of $k < 0$} Here
we take advantage of the identity $\Pi \phi[0] = \phi[0]$. Then we
can write:
\begin{equation}
        P_k \big( P_{< 0} \phi[0] -  \phi_{,<0}[0]\big)  \ = \
        P_k \big(\Pi \phi[0] - \Pi ( P_{< 0} \phi[0])\big)
        \ := \ \psi[0] \ . \notag
\end{equation}
To estimate the last difference we use an integral expansion as
follows:
\begin{align}
        \psi_0  \ &=  \ P_k\int_0^{\infty}  \frac{d}{dk_1}
        G({P_{<k_1}}\phi_0) dk_1 \ , \notag\\
        &=  \ P_k \int_0^{\infty} G'(P_{<k_1}\phi_0) P_{k_1} \phi_0 dk_1
        \ , \notag \\
        &= \ P_k \int_0^{\infty} P_{> k_1-10}
        G'(P_{<k_1}\phi_0)\cdot
        P_{k_1} \phi_0 dk_1 \ . \notag
\end{align}
Next, we use Bernstein's inequality and \eqref{moser2} to estimate:
\begin{align}
    \lp{ P_k\big( P_{> k_1-10}  G'(P_{<k_1}\phi_0)\cdot P_{k_1}\phi_0\big)}{L^2_x}
    \ &\lesssim \ 2^k \lp{P_{> k_1-10} G'(P_{<k_1}\phi_0)\cdot  P_{k_1}\phi_0}{L^1_x}
    \  \notag\\
    &\lesssim \ 2^k \lp{P_{> k_1-10} G'(P_{<k_1}\phi_0)}{L^2_x}
    \lp{P_{k_1}  \phi_0}{L^2_x}
    \  \notag\\
    &\lesssim_E \ 2^{k-2k_1}  \ . \label{moserhhl}
\end{align}
Hence after integration with respect to $k_1\geqslant 0$ we obtain:
\begin{equation}
   \lp{ \psi_0}{L^2_x} \ \lesssim_E \ 2^{k} \ . \notag
\end{equation}
A similar computation shows that:
\begin{equation}
        \psi_1 \ = \ P_k \int_0^\infty
        H'(P_{<k_1}\phi_0)P_{k_1}\phi_0\cdot
        P_{<k_1} \phi_1 dk_1 + P_k\int_0^\infty H(P_{<k_1}\phi_0)
        P_{k_1} \phi_1 dk_1 \ . \notag
\end{equation}
Then proceeding as above, we may estimate both integrands on the RHS
in terms of $\lesssim_E 2^{k-k_1}$, which upon integration over all
$k_1\geqslant 0$ yields a similar bound:
\begin{equation}
        \lp{ \psi_1}{L^2_x} \ \lesssim_E \ 2^{k} \ . \notag
\end{equation}
Thus we have proved that:
\begin{equation}
        \lp{ P_{k}  \big(P_{< 0} \phi[0] -
        \phi_{,<0}[0]\big) }{\dot H^1
        \times L^2}   \ \lesssim_E \ 2^{k}
        \ . \label{mhhl}
\end{equation}\ret

%-----------------------------------------------------------------------------

\step{3}{Intermediate frequency bounds, the contribution of $-m < k<
m$} Here $m$ is some fixed large integer. The goal of the argument
here is to show the estimate:
\begin{equation}
        \lp{P_k \big(P_{< 0} \phi[0] -
        \phi_{,<0}[0]\big)}{\dot{H}^1\times L^2} \ \lesssim_E \
        m^2\epsilon + 2^{-m} \ , \qquad |k|\ \leqslant \ m \ .  \label{mhhh}
\end{equation}
This is used with $m$ chosen so that $2^{-m} \approx \epsilon$.  Due
to the identity $\Pi \phi[0] = \Pi\phi[0]$ we can rewrite \eqref{mhhh}
in the form
\begin{equation}
\lp{P_k \big(P_{< 0} \Pi \phi[0] -
        \Pi P_{<0} \phi[0]\big)}{\dot{H}^1\times L^2} \ \lesssim_E \
        m^2\epsilon + 2^{-m}.
 \label{mhhha}\end{equation}
This is a direct consequence of the following paradifferential relation:

\begin{lem} \label{pipara}
Let $\Pi$ be as in \eqref{piform}, and  $D\Pi$ be its differential. Then for each $\psi[0] \in \dot H^1\times L^2$
with energy $E$ and energy dispersion $\epsilon$ and each $k \in \R$ we have
\begin{equation}
\| P_k \Pi \psi[0] - D\Pi( P_{< k-m} \psi[0]) P_k \psi[0]\|_{\dot H^1
      \times L^2}  \lesssim_E  m^2 \epsilon + 2^{-m}, \qquad m > 4
\label{mhhhaux}
\end{equation}
where $P_k$ can be substituted by any multiplier whose symbol has
similar size, localization and regularity.
\end{lem}
We remark that in \eqref{mhhhaux} there is no geometry left.
That is to say, $\psi[0]$ in \eqref{mhhhaux} need not satisfy the identity
$\Pi \psi[0] = \psi[0]$.

It is easy to see that \eqref{mhhhaux} implies \eqref{mhhha}. Indeed,
if $k \geq 2$ then the first term $P_k P_{< 0} \Pi \phi[0]$ in
\eqref{mhhha} does not contribute, while for the second we use
\eqref{mhhhaux} with $\psi[0] = P_{<0} \phi[0]$.  On the other hand if
$k < 2$ then for the first term $P_k P_{< 0} \Pi \phi[0]$ in
\eqref{mhhha} we use \eqref{mhhhaux} with $P_k$ replaced by $P_k P_{<
  0}$ and $\psi[0] = \phi[0]$, while for the second term we again use
\eqref{mhhhaux} with $\psi[0] = P_{<0} \phi[0]$. It remains to prove
the lemma.
\begin{proof}[Proof of Lemma~\ref{pipara}]
We write
\[
 P_k \Pi \psi[0] - D\Pi( P_{< k-m} \psi[0]) P_k \psi[0] = (w_0,w_1),
\]
with
\[
w_0 =  P_k G(\psi_0) - (\nabla G)(P_{< k-m} \psi_0) P_k \psi_0,
\]
\[
w_1 =  P_k ( H(\psi_0)\psi_1)  - ( H(P_{< k-m} \psi_0) P_k \psi_1+
(\nabla H)(P_{< k-m} \psi_0) P_k \psi_0 P_{< k-m} \psi_1   ).
\]
Then we need to prove that
\[
\|w_0\|_{\dot H^1} + \|w_1\|_{L^2} \lesssim_E   m^2 \epsilon + 2^{-m}.
\]
We observe that the expression for $\nabla w_0$ coincides
with the expression for $w_1$ with $H = \nabla G$ and $\psi_1 = \nabla \psi_0$.
Hence it suffices to prove the bound for $w_1$. Furthermore,
the last term in $w_1$ is directly estimated as
\[
\| (\nabla H)(P_{< k-m} \psi_0) P_k \psi_0 P_{< k-m} \psi_1 \|_{L^2}
\lesssim \|  P_k \psi_0\|_{L^\infty} \|\psi_1\|_{L^2} \lesssim_E \epsilon.
\]
It remains to show that
\begin{equation}
\|   P_k ( H(\psi_0)\psi_1)  -  H(P_{< k-m} \psi_0) P_k \psi_1\|_{L^2}
\lesssim_E   m^2 \epsilon + 2^{-m}.
\label{paralast}
\end{equation}
We use the integral representation
\[
\begin{split}
P_k ( H(\psi_0)\psi_1) = &\ P_k \int_{-\infty}^\infty \frac{d}{dk_1}
( H(P_{<k_1}\psi_0)P_{<k_1} \psi_1) dk_1 \\ = &\ P_k \int_{-\infty}^\infty
 H(P_{<k_1}\psi_0)P_{k_1} \psi_1   +
\nabla H(P_{<k_1}\psi_0)  P_{k_1}\psi_0  P_{< k_1} \psi_1  dk_1.
\end{split}
\]
The integrals from $-\infty$ to $k-m$, respectively from $k+m$ to $\infty$
can be bounded $\lesssim_E 2^{-m}$ as in STEP 1, respectively STEP 2 above.
For the integral from $k-m$ to $k+m$ we consider the two terms in the integrand
separately. The second term is estimated directly,
\[
\| \nabla H(P_{<k_1}\psi_0)  P_{k_1}\psi_0  P_{< k_1} \psi_1\|_{L^2}
\lesssim \|  P_{k_1}\psi_0\|_{L^\infty} \| P_{< k_1} \psi_1\|_{L^2} \lesssim_E \epsilon.
\]
Thus so far we have
\[
P_k ( H(\psi_0)\psi_1) =  P_k \int_{k-m}^{k+m}
 H(P_{<k_1}\psi_0)P_{k_1} \psi_1   dk_1 + O_{L^2}(m \epsilon + 2^{-m}).
\]
The remaining integrand is further expanded,
\[
 H(P_{<k_1}\psi_0)P_{k_1} \psi_1 = H(P_{<k-m}\psi_0)P_{k_1} \psi_1+
\int_{k-m}^{k_1} \nabla H(P_{<k_2}\psi_0)  P_{k_2}\psi_0  P_{ k_1} \psi_1 dk_2.
\]
The second term can be estimated as above by $\lesssim_E m \epsilon$. We arrive at
\[
P_k ( H(\psi_0)\psi_1) = P_k ( H(P_{<k-m}\psi_0) P_{[k-m,k+m]} \psi_1 ) +
O_{L^2}(m^2 \epsilon + 2^{-m})
\]
This implies \eqref{paralast} via a commutator bound, see \eqref{commutator}:
\[
\begin{split}
\|[P_k,H(P_{<k-m}\psi_0)]   P_{[k-m,k+m]} \psi_1\|_{L^2}
 \lesssim  &\ 2^{-k} \| \nabla H(P_{<k-m}\psi_0)\|_{L^\infty}
\| P_{[k-m,k+m]} \psi_1\|_{L^2}
\\
\lesssim_E &\  2^{-m}
\end{split}
\]
The proof of the lemma is complete.
\end{proof}

 This concludes our demonstration of Proposition \ref{p:cuta}.
\end{proof}

\ret
%-------------------------------------------------------------------------
%%%%%%%%%%%%%%%%%%%%%%%%%%%%%%%%%%%%%%%%%%%%%%%%%%%%%%%%%%%%%%%%%%%%%%%%%%
%%%%%%%%%%%%%%%%%%%%%%%%%%%%%%%%%%%%%%%%%%%%%%%%%%%%%%%%%%%%%%%%%%%%%%%%%%
%-------------------------------------------------------------------------
\bibliography{wm}

\bibliographystyle{plain}

%-------------------------------------------------------------------------
%%%%%%%%%%%%%%%%%%%%%%%%%%%%%%%%%%%%%%%%%%%%%%%%%%%%%%%%%%%%%%%%%%%%%%%%%%
%%%%%%%%%%%%%%%%%%%%%%%%%%%%%%%%%%%%%%%%%%%%%%%%%%%%%%%%%%%%%%%%%%%%%%%%%%
%-------------------------------------------------------------------------

\end{document}